\documentclass[11pt,twoside,leqno,openany]{amsart}
\usepackage{amssymb,amsbsy,amsmath,amsfonts,amscd,times}
\usepackage{graphics,color,xy,xypic,footmisc,fancyhdr,multicol}
\usepackage{fancybox,graphicx,mathrsfs}
\usepackage[T1]{fontenc} 
\newcommand{\explain}[1]{\text{\scriptsize\sf [#1]}}
\let\mathcal\mathscr
\newcommand{\C}{\mathbb{C}}\newcommand{\K}{\mathbb{K}}
\newcommand{\N}{\mathbb{N}}
\newcommand{\Q}{\mathbb{Q}}\newcommand{\R}{\mathbb{R}}
\newcommand{\Z}{\mathbb{Z}}

\newcommand{\zero}[1]{\underline{#1}_{\circ}}

\newtheorem{The}{Theorem}[section]

\newtheorem{Theorem}{Theorem}[section]
\newtheorem{Proposition}[The]{Proposition}
\newtheorem{Lemma}[The]{Lemma}
\newtheorem{Corollary}[The]{Corollary}

\theoremstyle{definition}
\newtheorem{Definition}[The]{Definition}

\newcommand{\HEAD}[2]{%
\pagestyle{fancy}
\fancyhead[RO]{\scriptsize\sf\thepage}
\fancyhead[CO]{{\scriptsize\sf \thesection.\,\,\,#1}}
\fancyhead[LE]{\scriptsize\sf\thepage}
\fancyhead[CE]{{\scriptsize\sf #2}}
\fancyfoot{}}

\begin{document}

\title[]{
Effective Cartan-Tanaka connections 
\\
on $\mathcal{ C}^6$ strongly pseudoconvex hypersurfaces 
$M^3 \subset \C^2$}

\address{Department of Mathematical Sciences,
Isfahan University of Technology, Isfahan, IRAN}
\email{m.aghasi@cc.iut.ac.ir}

\address{D\'epartment de Math\'ematiques d'Orsay,
B\^atiment 425, Faculté des Sciences, Université Paris XI - Orsay, 
F-91405 Orsay Cedex, FRANCE}
\email{merker@dma.ens.fr}

\address{Department of Mathematical Sciences,
Isfahan University of Technology, Isfahan, IRAN}
\email{sabzevari@math.iut.ac.ir}

\date{\number\year-\number\month-\number\day}

\maketitle

\begin{center}
{\sc Mansour Aghasi, Jo\"el Merker, and
Masoud Sabzevari}
\end{center}



\begin{center}
\begin{minipage}[t]{11.75cm}
\baselineskip =0.35cm {\scriptsize

\centerline{\bf Table of contents}

{\bf \ref{statement-results}. Synthetic Statement of the Results
\dotfill~\pageref{statement-results}.}

{\bf \ref{general-infinitesimal-CR}. 
Infinitesimal CR Automorphisms of Real Analytic Generic
Submanifolds of $\C^{n+d}$
\dotfill~\pageref{general-infinitesimal-CR}.}

{\bf \ref{Infinitesimal-CR-section}. 
CR Symmetries of the Heisenberg sphere
$\mathbb H^3\subseteq\mathbb C^2$ 
\dotfill~\pageref{Infinitesimal-CR-section}.}

{\bf \ref{Tanaka-Prolongation}. Tanaka prolongation 
\dotfill~\pageref{Tanaka-Prolongation}.}

{\bf \ref{Cohomology-section}. 
Second Cohomology of the Heisenberg Lie Algebra
\dotfill~\pageref{Cohomology-section}.}

{\bf \ref{initial-frame}. Initial Frame on a Strongly Pseudoconvex 
$M^3 \subset \C^2$ 
\dotfill~\pageref{initial-frame}.}

{\bf \ref{free-Lie-algebras}. 
Free Lie Algebras of Rank Two 
and Relations Between Brackets of Length $\leqslant 6$
\dotfill~\pageref{free-Lie-algebras}.}

{\bf \ref{Cartan-connections-coordinates}.
Cartan Connections in Terms of Coordinates and Bases
\dotfill~\pageref{Cartan-connections-coordinates}.}

{\bf \ref{Cartan-construction}. Effective Construction of 
a Normal, Regular Cartan-Tanaka connection 
\dotfill~\pageref{Cartan-construction}.}

{\bf \ref{second-cohomology-graded}. 
General Formulas for the Second Cohomology of Graded Lie Algebras
\dotfill~\pageref{second-cohomology-graded}.}

}\end{minipage}
\end{center}

\section{Synthetic Statement of the Results}
\label{statement-results}

\HEAD{Synthetic Statement of the Results}{
Mansour Aghasi, Joël Merker, and Masoud Sabzevari}

\begin{Theorem}
Let $M^3 \subset \C^2$ be an arbitrary local Levi nondegenerate
$\mathcal{ C}^6$-smooth real $3$-dimensional hypersurface of $\C^2$
which is represented in coordinates $(z,w) = (x + iy, \, u+ iv)$ as a
graph:
\[
v
=
\varphi(x,y,u)
=
x^2+y^2+{\rm O}(3),
\]
and whose complex tangent bundle $T^c M = {\rm Re}\, T^{ 1, 0} M$
is generated by the two explicit intrinsic vector fields:
\[
H_1 
:= 
{\textstyle{\frac{\partial}{\partial x}}} 
+ 
\big(
{\textstyle{\frac{\varphi_y-\varphi_x\,\varphi_u}{
1+\varphi_u^2}}} 
\big)
{\textstyle{\frac{\partial}{\partial u}}} 
\ \ \ \ \ \ \ \ \ \ 
\text{\rm and} 
\ \ \ \ \ \ \ \ \ \
H_2 
:= 
{\textstyle{\frac{\partial}{\partial y}}} 
+ 
\big(
{\textstyle{\frac{-\varphi_x-\varphi_y\,\varphi_u}{
1+\varphi_u^2}}} 
\big)
{\textstyle{\frac{\partial}{\partial u}}},
\]
satisfying $H_1 \vert_0 = \frac{ \partial}{ \partial x} \big\vert_0$
and $H_2 \vert_0 = \frac{ \partial}{ \partial y} \big\vert_0${\em ;} 
introduce the three $\mathcal{ C}^5$-smooth functions:
\[
\Delta
:=
1+\varphi_u^2,
\ \ \ \ \ \ \ \ \ \
\Lambda_1
:=
\varphi_y-\varphi_x\,\varphi_u,
\ \ \ \ \ \ \ \ \ \
\Lambda_2
:=
-\,\varphi_x-\varphi_y\,\varphi_u,
\]
so that:
\[
H_1 
= 
{\textstyle{\frac{\partial}{\partial x}}} 
+
{\textstyle{\frac{\Lambda_1}{\Delta}}}\,
{\textstyle{\frac{\partial}{\partial u}}} 
\ \ \ \ \ \ \ \ \ \ 
\text{\rm and} 
\ \ \ \ \ \ \ \ \ \
H_2 
= 
{\textstyle{\frac{\partial}{\partial y}}} 
+ 
{\textstyle{\frac{\Lambda_2}{\Delta}}}\,
{\textstyle{\frac{\partial}{\partial u}}};
\]
consider the third, Levi form-type Lie-bracket vector field:
\[
\aligned 
T 
:= 
&\,
{\textstyle{\frac{1}{4}}}\, [H_1,H_2]
\\
= 
&\,
\Big({\textstyle{\frac{1}{4}}}\,
{\textstyle{\frac{1}{(1+\varphi_u^2)^2}}}
\big\{
-
\varphi_{xx}-\varphi_{yy}-2\,\varphi_y\,\varphi_{xu}
- 
\varphi_x^2\,\varphi_{uu} 
+
2\,\varphi_x\,\varphi_{yu} 
- 
\varphi_y^2\,\varphi_{uu} 
+
\\
& \ \ \ \ \ \ \ \ \ \ \ \ \ \ \ \ \ \ \ \ \ \ \ \ +
2\,\varphi_y\,\varphi_u\,\varphi_{yu} +
2\,\varphi_x\,\varphi_u\,\varphi_{xu} - \varphi_u^2\,\varphi_{xx} -
\varphi_u^2\,\varphi_{yy}
\big\} \Big)\, \frac{\partial}{\partial u},
\endaligned
\]
satisfying $T \big\vert_0 = - \frac{ \partial}{ \partial u}
\big\vert_0$ which produces, jointly with $H_1$ and $H_2$ of which it
is locally linearly independent, a frame for $TM$ in a neighborhood of
the origin; introduce the $\mathcal{ C}^4$-smooth function:
\[
\aligned
\Upsilon
&
:=
-
\varphi_{xx}-\varphi_{yy}-2\,\varphi_y\,\varphi_{xu}
- 
\varphi_x^2\,\varphi_{uu} 
+
2\,\varphi_x\,\varphi_{yu} 
- 
\varphi_y^2\,\varphi_{uu} 
+
\\
& \ \ \ \ \ \ \ \ \ \ \ \ \ \ \ \ \ \ \ \ \ \ \ \ +
2\,\varphi_y\,\varphi_u\,\varphi_{yu} +
2\,\varphi_x\,\varphi_u\,\varphi_{xu} - \varphi_u^2\,\varphi_{xx} -
\varphi_u^2\,\varphi_{yy},
\endaligned
\]
which specifies the numerator of the Levi form-type Lie-bracket:
\[
T 
=
{\textstyle{\frac{1}{4}}}\, [H_1,H_2]
=
\frac{1}{4}\,
\frac{\Upsilon}{\Delta^2}\,
\frac{\partial}{\partial u};
\]
allow the two notational coincidences: $x_1 \equiv x$, $x_2 \equiv y$;
introduce the two length-three brackets:
\[
\big[H_i,\,T\big]
=
{\textstyle{\frac{1}{4}}}\,
\big[H_i,
[H_1,H_2]\big]
=:
\Phi_i\,T
\ \ \ \ \ \ \ \ \ \ \ \ \ 
{\scriptstyle{(i\,=\,1,\,2)}},
\]
which are both multiples of $T$ by means of two functions:
\[
\Phi_i
:=
\frac{A_i}{\Delta^2\,\Upsilon}
\ \ \ \ \ \ \ \ \ \ \ \ \ 
{\scriptstyle{(i\,=\,1,\,2)}}
\]
whose numerators are explicitly given by:
\[
A_i
:=
\Delta^2\,\Upsilon_{x_i}
+
\Delta
\big(
-2\,\Delta_{x_i}\,\Upsilon
+
\Lambda_i\,\Upsilon_u
-
\Upsilon\,\Lambda_{i,u}
\big)
-
\Lambda_i\,\Upsilon\,\Delta_u
\ \ \ \ \ \ \ \ \ \ \ \ \ 
{\scriptstyle{(i\,=\,1,\,2)}};
\]
introduce furthermore the following 
$4 + 8 + 16$ iterated brackets for $i, k_1, k_2, k_3 = 1, 2${\em :}
\[
\small
\aligned
\big[H_{k_1},\big[H_i,\,T\big]\big]
&
=
{\textstyle{\frac{1}{4}}}\,
\big[H_{k_1},\big[H_i,\,
[H_1,H_2]\big]\big]
\\
\big[H_{k_2},\big[H_{k_1},\big[H_i,\,T\big]\big]\big]
&
=
{\textstyle{\frac{1}{4}}}\,
\big[H_{k_2},\big[H_{k_1},\big[H_i,\,
[H_1,H_2]\big]\big]\big]
\\
\big[H_{k_3},\big[H_{k_2},\big[H_{k_1},\big[H_i,\,T\big]\big]\big]\big]
&
=
{\textstyle{\frac{1}{4}}}\,
\big[H_{k_3},\big[H_{k_2},\big[H_{k_1},\big[H_i,\,
[H_1,H_2]\big]\big]\big]\big],
\endaligned
\]
up to length $6$ that are all multiples of $T$:
\[
\footnotesize
\aligned
\big[H_{k_1},\big[H_i,\,T\big]\big]
&
=
\big(
H_{k_1}(\Phi_i)
+
\Phi_i\,\Phi_{k_1}
\big)\,T,
\\
\big[H_{k_2},\big[H_{k_1},\big[H_i,\,T\big]\big]\big]
&
=
\big(
H_{k_2}(H_{k_1}(\Phi_i))
+
\Phi_{k_1}\,H_{k_2}(\Phi_i)
+
\Phi_i\,H_{k_2}(\Phi_{k_1})
+
\\
&
\ \ \ \ \ \ \
+
\Phi_{k_2}\,H_{k_1}(\Phi_i)
+
\Phi_i\,\Phi_{k_1}\Phi_{k_2}
\big)\,T,
\\
\big[H_{k_3},\big[H_{k_2},\big[H_{k_1},\big[H_i,\,T\big]\big]\big]\big]
&
=
\big(
H_{k_3}(H_{k_2}(H_{k_1}(\Phi_i)))
+
\Phi_{k_1}H_{k_3}(H_{k_2}(\Phi_i))
+
\Phi_i\,H_{k_3}(H_{k_2}(\Phi_{k_1}))
+
\\
&
\ \ \ \ \ \ \
+
\Phi_{k_2}\,H_{k_3}(H_{k_1}(\Phi_i))
+
\Phi_{k_3}H_{k_2}(H_{k_1}(\Phi_i))
+
\\
&
\ \ \ \ \ \ \
+
H_{k_3}(\Phi_{k_1})\,H_{k_2}(\Phi_i)
+
H_{k_3}(\Phi_i)\,H_{k_2}(\Phi_{k_1})
+
H_{k_3}(\Phi_{k_2})\,H_{k_1}(\Phi_i)
+
\\
&
\ \ \ \ \ \ \
+
\Phi_{k_1}\Phi_{k_2}H_{k_3}(\Phi_i)
+
\Phi_i\,\Phi_{k_2}H_{k_3}(\Phi_{k_1})
+
\Phi_i\,\Phi_{k_1}H_{k_3}(\Phi_{k_2})
+
\\
&
\ \ \ \ \ \ \
+
\Phi_{k_1}\Phi_{k_3}H_{k_2}(\Phi_i)
+
\Phi_i\,\Phi_{k_3}H_{k_2}(\Phi_{k_1})
+
\Phi_{k_2}\Phi_{k_3}H_{k_1}(\Phi_i)
+
\\
&
\ \ \ \ \ \ \
+
\Phi_i\,\Phi_{k_1}\Phi_{k_2}\Phi_{k_3}
\big)\,T
\endaligned
\]
and in the expressions of which the 
$H_k$-iterated derivatives of the functions $\Phi_i$ up
to order $3$: 
\[
\aligned
H_{k_1}(\Phi_i)
&
=
\frac{A_{i,k_1}}{\Delta^4\,\Upsilon^2}
\ \ \ \ \ \ \ \ \ \ \ \ \ 
{\scriptstyle{(i,\,k_1\,=\,1,\,2)}},
\\
H_{k_2}(H_{k_1}(\Phi_i))
&
=
\frac{A_{i,k_1,k_2}}{\Delta^6\,\Upsilon^3}
\ \ \ \ \ \ \ \ \ \ \ \ \ 
{\scriptstyle{(i,\,k_1,\,k_2\,=\,1,\,2)}},
\\
H_{k_3}(H_{k_2}(H_{k_1}(\Phi_i)))
&
=
\frac{A_{i,k_1,k_2,k_3}}{\Delta^8\,\Upsilon^4}
\ \ \ \ \ \ \ \ \ \ \ \ \ 
{\scriptstyle{(i,\,k_1,\,k_2,\,k_3\,=\,1,\,2)}},
\endaligned
\]
have numerators $A_{ i,k_1}$, $A_{ i, k_1, k_2}$, $A_{ i, k_1, k_2,
k_3}$ that are certain differential polynomials whose {\em completely
explicit} expressions in terms of the jets
$J_{ x, y, u}^4 \varphi$, $J_{ x, y, u}^5 \varphi$, $J_{ x, y, u}^6
\varphi$ of the graphing function $\varphi ( x, y, u)$ of orders
$4$, $5$, $6$ (respectively) are provided through the induction
formulas:
\[
\footnotesize
\aligned
A_{i,k_1}
&
:=
\Delta^2
\big(
\Upsilon\,A_{i,x_{k_1}}
-
\Upsilon_{x_{k_1}}\,A_i
\big)
+
\Delta
\big(
-2\,\Delta_{x_{k_1}}\,\Upsilon\,A_i
+
\Upsilon\,\Lambda_{k_1}\,A_{i,u}
-
\Upsilon_u\,\Lambda_{k_1}\,A_i
\big)
-
\\
&
\ \ \ \ \ 
-
2\,\Delta_u\,\Upsilon\,\Lambda_{k_1}\,A_i 
\ \ \ \ \ \ \ \ \ \ \ \ \ \ \ \ \ \ \ \ \ \ \ \ \ \ \ \ \
\ \ \ \ \ \ \ \ \ \ \ \ \ \ \ \ \ \ \ \ \ \ \ \ \ \ \ \ \
\ \ \ \ \ \ \ \ \ \ \ \ \ \ \ \ \ \ \ \ \ \ \ \ \ \ \ 
{\scriptstyle{(i,\,k_1\,=\,1,\,2)}},
\\
A_{i,k_1,k_2}
&
:=
\Delta^2
\big(
\Upsilon\,A_{i,k_1,x_{k_2}}
-
2\,\Upsilon_{x_{k_2}}\,A_{i,k_1}
\big)
+
\Delta
\big(
-3\,\Delta_{x_{k_2}}\,\Upsilon\,A_{i,k_1}
+
\Upsilon\,\Lambda_{k_2}\,A_{i,k_1,u}
-
\\
&
\ \ \ \ \
-
2\,\Upsilon_u\,\Lambda_{k_2}\,A_{i,k_1}
\big)
-
3\,\Delta_u\,\Upsilon\,\Lambda_{k_2}\,A_{i,k_1}
\ \ \ \ \ \ \ \ \ \ \ \ \ \ \ \ \ \ \ \ \ \ \ \ \ \ \ \ \
\ \ \ \ \ \ \ \ \ \ \ \ \ \ \
{\scriptstyle{(i,\,k_1,\,k_2\,=\,1,\,2)}},
\\
A_{i,k_1,k_2,k_3}
&
=
\Delta^2
\big(
\Upsilon\,A_{i,k_1,k_2,x_{k_3}}
-
\Upsilon_{x_{k_3}}\,A_{i,k_1,k_2}
\big)
+
\Delta
\big(
-6\,\Delta_{x_{k_3}}\,\Upsilon\,A_{i,k_1,k_2}
+
\Upsilon\,\Lambda_{k_3}\,A_{i,k_1,k_2,u}
-
\\
&
\ \ \ \ \
-\,3\,\Upsilon_u\,\Lambda_{k_3}\,A_{i,k_1,k_2}
\big)
-
6\,\Delta_u\,\Upsilon\,\Lambda_{k_3}\,A_{i,k_1,k_2}
\ \ \ \ \ \ \ \ \ \ \ \ \ \ \ \ \ \ \ \ \ \ \ \ \ \ \ \ \ \
{\scriptstyle{(i,\,k_1,\,k_2,\,k_3\,=\,1,\,2)}}.
\endaligned
\]
Then (first, preliminary effective assertion): these iterated
derivatives identically satisfy: 
\[ 
\small
\aligned
H_2(\Phi_1) 
\equiv 
H_1(\Phi_2),
\endaligned
\]
together with the following four third-order relations:
\[
\small
\aligned
0
&
\equiv
-\,H_1(H_2(H_1(\Phi_2)))
+
2\,H_2(H_1(H_1(\Phi_2)))
-
H_2(H_2(H_1(\Phi_1)))
-
\\
&
\ \ \ \ \ 
-\,\Phi_2\,H_1(H_2(\Phi_1))
+
\Phi_2\,H_2(H_1(\Phi_1)),
\endaligned
\]
\[
\small
\aligned
0
&
\equiv
-\,H_2(H_1(H_1(\Phi_2)))
+
2\,H_1(H_2(H_1(\Phi_2)))
-
H_1(H_1(H_2(\Phi_2)))
-
\\
&
\ \ \ \ \ 
-\,\Phi_1\,H_2(H_1(\Phi_2))
+
\Phi_1\,H_1(H_2(\Phi_2)),
\endaligned
\]
\[
\small
\aligned
0
&
\equiv
-\,H_1(H_1(H_1(\Phi_2)))
+
2\,H_1(H_2(H_1(\Phi_1)))
-
H_2(H_1(H_1(\Phi_1)))
+
\\
&
\ \ \ \ \ 
+\Phi_1\,H_1(H_1(\Phi_2))
-
\Phi_1\,H_2(H_1(\Phi_1)),
\endaligned
\]
\[
\small
\aligned
0
&
\equiv
H_2(H_2(H_1(\Phi_2)))
-
2\,H_2(H_1(H_2(\Phi_2)))
+
H_1(H_2(H_2(\Phi_2)))
-
\\
&
\ \ \ \ \
-\,\Phi_2\,H_2(H_1(\Phi_2))
+
\Phi_2\,H_1(H_2(\Phi_2)).
\endaligned
\]
Moreover (second, well known effective assertion), the model
Heisenberg sphere $\mathbb{ H}^3 \subset \C^2$ whose graphing function
simply has no ${\rm O} (3)$ remainder:
\[
v
=
x^2+y^2
\]
possesses an eight-dimensional graded Lie algebra:
\[
\mathfrak{ hol}(\mathbb{H}^3) 
= 
\mathfrak{h}_{-2}
\oplus 
\mathfrak{h}_{-1} 
\oplus
\mathfrak{h}_0
\oplus 
\mathfrak{h}_1
\oplus
\mathfrak{h}_2
\]
of
(local or global) holomorphic vector fields ${\sf X}$ whose real parts
$\frac{ 1}{ 2}\big( {\sf X} + \overline{ \sf X} \big)$ are tangent
to it, having components:
\[
\aligned \mathfrak{h}_{-2} & = \mathbb R\,{\sf T}, 
\ \ \ \ \ 
\mathfrak{h}_{-1} 
=
\mathbb R\,{\sf
H}_1 \oplus \mathbb R\,{\sf H}_2
\\
\mathfrak{h}_0 & = \mathbb R\,{\sf D} \oplus \mathbb R\,{\sf R},
\\
\mathfrak{h}_1 & = \mathbb R\,{\sf I}_1 \oplus \mathbb R\,{\sf I}_2, \ \ \ \ \
\mathfrak{h}_2 
=
\mathbb R\,{\sf J}.
\endaligned
\]
where:
\[
\aligned 
& 
\mathfrak{h}_{-2}\colon \left\{ {\sf T} 
:=
\partial_w
\right. 
\ \ \ \ \ \ \ \ \ \ \ 
\mathfrak{h}_{-1}\colon \left\{ \aligned {\sf H}_1 &
:=
\partial_z
+ 2iz\,\partial_w
\\
{\sf H}_2 
& 
:= 
i\,\partial_z + 2z\,\partial_w
\endaligned
\right.
\\
& \mathfrak{h}_0\colon 
\left\{ 
\aligned 
{\sf D} & 
:= 
z\,\partial_z + 2w\,\partial_w
\\
{\sf R} & := iz\,\partial_z
\\
\endaligned\right.
\\
& \mathfrak{h}_1\colon \left\{ \aligned {\sf I}_1 & := (w+2iz^2)\,\partial_z +
2izw\,\partial_w
\\
{\sf I}_2 & := (iw+2z^2)\,\partial_z + 2zw\,\partial_w
\endaligned\right.
\ \ \ \ \ \ \ \ \ \ \ \mathfrak{h}_2\colon \left\{ {\sf J} := zw\,\partial_z +
w^2\,\partial_w,
\right.
\endaligned
\]
and these eight holomorphic fields
enjoy the following commutator table with
real (in fact, integer) structure constants:
\medskip
\begin{center}
\begin{tabular} [t] { l | l l l l l l l l }
& ${\sf T}$ & ${\sf H}_1$ & ${\sf H}_2$ & ${\sf D}$ & ${\sf R}$ & ${\sf I}_1$ &
${\sf I}_2$ & ${\sf
J}$
\\
\hline ${\sf T}$ & $0$ & $0$ & $0$ & $2\,{\sf T}$ & $0$ & ${\sf H}_1$ & ${\sf H}_2$
& ${\sf D}$
\\
${\sf H}_1$ & $*$ & $0$ & $4\,{\sf T}$ & ${\sf H}_1$ & ${\sf H}_2$ & $6\,{\sf R}$ &
$2\,{\sf D}$ &
${\sf I}_1$
\\
${\sf H}_2$ & $*$ & $*$ & $0$ & ${\sf H}_2$ & $-{\sf H}_1$ & $-2\,{\sf D}$ &
$6\,{\sf R}$ & ${\sf
I}_2$
\\
${\sf D}$ & $*$ & $*$ & $*$ & $0$ & $0$ & ${\sf I}_1$ & ${\sf I}_2$ & $2\,{\sf J}$
\\
${\sf R}$ & $*$ & $*$ & $*$ & $*$ & $0$ & $-{\sf I}_2$ & ${\sf I}_1$ & $0$
\\
${\sf I}_1$ & $*$ & $*$ & $*$ & $*$ & $*$ & $0$ & $4\,{\sf J}$ & $0$
\\
${\sf I}_2$ & $*$ & $*$ & $*$ & $*$ & $*$ & $*$ & $0$ & $0$
\\
${\sf J}$ & $*$ & $*$ & $*$ & $*$ & $*$ & $*$ & $*$ & $0$.
\end{tabular}
\end{center}

\noindent
Lastly (third, main effective assertion), to any $\mathcal{ C}^6$
strongly pseudoconvex $M^3 \subset \C^2$ is uniquely associated an
{\em effective} local Cartan connection: 
\[
\omega
\colon\,\,\,\,
T\mathcal{P}
\longrightarrow
\mathfrak{g}
\]
on the local principal bundle:
\[
\mathcal{P}
:=
M^3
\times
H^5
\]
which is the Cartesian product of $M$ with the unique (connected and
simply connected) local $5$-dimensional Lie group $H$ equipped with
some $5$ real coordinates:
\[
(a,b,c,d,e)\in\R^5,
\] 
that is associated to the isotropy Lie subalgebra:
\[
\mathfrak{hol}(\mathbb{H}^3,0)
=
\R\,{\sf D}
\oplus
\R\,{\sf R}
\oplus
\R\,{\sf I}_1
\oplus
\R\,{\sf I}_2
\oplus
\R\,{\sf J}
\]
of the origin $0 \in \mathbb{ H}^3${\em ;} this Cartan connection 
$\omega \colon T\mathcal{ P} \longrightarrow \mathfrak{ g}$
is valued in the eight-dimensional {\em abstract} real Lie algebra:
\[
\mathfrak{g}
:=
\R\,{\sf t}
\oplus
\R\,{\sf h}_1
\oplus
\R\,{\sf h}_2
\oplus
\R\,{\sf d}
\oplus
\R\,{\sf r}
\oplus
\R\,{\sf i}_1
\oplus
\R\,{\sf i}_2
\oplus
\R\,{\sf j}
\]
spanned by some eight abstract vectors enjoying the same commutator table:
\medskip
\begin{center}
\begin{tabular} [t] { l | l l l l l l l l }
& ${\sf t}$ & ${\sf h}_1$ & ${\sf h}_2$ & ${\sf d}$ & ${\sf r}$ &
${\sf i}_1$ & ${\sf i}_2$ & ${\sf j}$
\\
\hline ${\sf t}$ & $0$ & $0$ & $0$ & $2\,{\sf t}$ 
& $0$ & ${\sf h}_1$ & ${\sf h}_2$
& ${\sf d}$
\\
${\sf h}_1$ & $*$ & $0$ & $4\,{\sf t}$ 
& ${\sf h}_1$ & ${\sf h}_2$ & $6\,{\sf r}$ &
$2\,{\sf d}$ &
${\sf i}_1$
\\
${\sf h}_2$ & $*$ & $*$ & $0$ & ${\sf h}_2$ 
& $-{\sf h}_1$ & $-2\,{\sf d}$ &
$6\,{\sf r}$ & ${\sf
i}_2$
\\
${\sf d}$ & $*$ & $*$ & $*$ & $0$ 
& $0$ & ${\sf i}_1$ & ${\sf i}_2$ & $2\,{\sf j}$
\\
${\sf r}$ & $*$ & $*$ & $*$ & $*$ 
& $0$ & $-{\sf i}_2$ & ${\sf i}_1$ & $0$
\\
${\sf i}_1$ & $*$ & $*$ & $*$ & $*$ 
& $*$ & $0$ & $4\,{\sf j}$ & $0$
\\
${\sf i}_2$ & $*$ & $*$ & $*$ & $*$ & $*$ & $*$ & $0$ & $0$
\\
${\sf j}$ & $*$ & $*$ & $*$ & $*$ & $*$ & $*$ & $*$ & $0${\rm ;}
\end{tabular}
\end{center}
this Cartan connection $\omega \colon T \mathcal{ P} \longrightarrow
\mathfrak{ g}$ is normal and regular in the sense of Tanaka (\cite{
Cap, EMS}, \pageref{c3-c4-anticipate} below), and if one denotes the
Lie algebra of $H$ by:
\[
\mathfrak{h}
:=
\R\,{\sf d}
\oplus
\R\,{\sf r}
\oplus
\R\,{\sf i}_1
\oplus
\R\,{\sf i}_2
\oplus
\R\,{\sf j},
\]
with corresponding five left-invariant vector fields on $H$ of the form:
\[
\aligned
D
&
:=-
a\,{\textstyle{\frac{\partial}{\partial a}}}
-
b\,{\textstyle{\frac{\partial}{\partial b}}}
-
c\,{\textstyle{\frac{\partial}{\partial c}}}
-
d\,{\textstyle{\frac{\partial}{\partial d}}}
-
2e\,{\textstyle{\frac{\partial}{\partial e}}}
\\
R
&
:=-b\,{\textstyle{\frac{\partial}{\partial a}}}
+
a\,{\textstyle{\frac{\partial}{\partial
b}}}
+
d\,{\textstyle{\frac{\partial}{\partial c}}}
-
c\,{\textstyle{\frac{\partial}{\partial d}}}
\\
I_1
&
:={\textstyle{\frac{\partial}{\partial a}}}
-
b\,{\textstyle{\frac{\partial}{\partial e}}}
\\
I_2
&
:={\textstyle{\frac{\partial}{\partial b}}}
+
a\,{\textstyle{\frac{\partial}{\partial e}}}
\\
J
&
:={\textstyle{\frac{1}{2}}}\,
{\textstyle{\frac{\partial}{\partial e}}}
\endaligned
\]
near the origin $(a_0, b_0, c_0, d_0, e_0) := (0, 0, 1, 1, 0)$, 
then the curvature function:
\[
\kappa
\in
\mathcal{C}^0
\big(
\mathcal{P},\,
\Lambda^2(\mathfrak{g}^*/\mathfrak{h}^*)
\otimes
\mathfrak{g}
\big)
\]
of the Cartan connection $\omega \colon \mathcal{ P}
\longrightarrow \mathfrak{ g}$, a function $\kappa(p)$ of the
eight real variables:
\[
\mathcal{P}
\ni
p
:=
(x,y,u,a,b,c,d,e)
\in
M^3\times H
\]
has an algebraic expression which reduces to:
\[
\aligned
\kappa(p)
&
=
\kappa^{h_1t}_{i_1}(p)\,
{\sf h}_1^\ast\wedge{\sf t}^\ast
\otimes
{\sf i}_1
+
\kappa^{h_1t}_{i_2}(p)\,
{\sf h}_1^\ast\wedge{\sf t}^\ast
\otimes 
{\sf i}_2
+
\kappa^{h_2t}_{i_1}(p)\,
{\sf h}_2^\ast\wedge{\sf t}^\ast
\otimes 
{\sf i}_1
+
\\
&
\ \ \ \ \
+
\kappa^{h_2t}_{i_2}(p)\,
{\sf h}_2^\ast\wedge{\sf t}^\ast
\otimes 
{\sf i}_2
+
\kappa^{h_1t}_j(p)\,
{\sf h}_1^\ast\wedge{\sf t}^\ast
\otimes
{\sf j}
+
\kappa^{h_2t}_j(p)\,
{\sf h}_2^\ast\wedge{\sf t}^\ast
\otimes
{\sf j},
\endaligned
\]
where the two main curvature coefficients are {\em explicitly} given
by:
\[
\footnotesize
\aligned
\kappa_{i_1}^{h_1t}(p)
&
=
-\,\mathbf{\Delta_1}\,c^4
-
2\,\mathbf{\Delta_4}\,c^3d
-
2\,\mathbf{\Delta_4}\,cd^3
+
\mathbf{\Delta_1}\,d^4,
\\
\kappa_{i_2}^{h_1t}(p)
&
=
-\,\mathbf{\Delta_4}\,c^4
+
2\,\mathbf{\Delta_1}\,c^3d
+
2\,\mathbf{\Delta_1}\,cd^3
+
\mathbf{\Delta_4}\,d^4,
\endaligned
\]
in which the two functions $\mathbf{ \Delta_1}$ and $\mathbf{ \Delta_4}$
of only the three variables $(x, y, u)$ are {\em explicitly} given by
the symmetric expressions:
\[
\footnotesize\aligned
\mathbf{\Delta_1}
&
=
{\textstyle{\frac{1}{384}}}
\Big[
H_1(H_1(H_1(\Phi_1)))
-
H_2(H_2(H_2(\Phi_2)))
+
11\,H_1(H_2(H_1(\Phi_2)))
-
11\,H_2(H_1(H_2(\Phi_1)))
+
\\
&
\ \ \ \ \ \ \ \ \ \ \ \ \
+6\,\Phi_2\,H_2(H_1(\Phi_1))
-
6\,\Phi_1\,H_1(H_2(\Phi_2))
-
3\,\Phi_2\,H_1(H_1(\Phi_2))
+
3\,\Phi_1\,H_2(H_2(\Phi_1))
-
\\
&
\ \ \ \ \ \ \ \ \ \ \ \ \
-\,3\,\Phi_1\,H_1(H_1(\Phi_1))
+
3\,\Phi_2\,H_2(H_2(\Phi_2))
-
2\,\Phi_1\,H_1(\Phi_1)
+
2\,\Phi_2\,H_2(\Phi_2)
-
\\
&
\ \ \ \ \ \ \ \ \ \ \ \ \
-\,2\,(\Phi_2)^2\,H_1(\Phi_1)
+
2\,(\Phi_1)^2\,H_2(\Phi_2)
-
2\,(\Phi_2)^2\,H_2(\Phi_2)
+
2\,(\Phi_1)^2\,H_1(\Phi_1)
\Big],
\\
\mathbf{\Delta_4}
&
=
{\textstyle{\frac{1}{384}}}
\Big[
-\,3\,H_2(H_1(H_2(\Phi_2)))
-
3\,H_1(H_2(H_1(\Phi_1)))
+
5\,H_1(H_2(H_2(\Phi_2)))
+
5\,H_2(H_1(H_1(\Phi_1)))
+
\\
&
\ \ \ \ \ \ \ \ \ \ \ \ \
+4\,\Phi_1\,H_1(H_1(\Phi_2))
+
4\,\Phi_2\,H_2(H_1(\Phi_2))
-
3\,\Phi_2\,H_1(H_1(\Phi_1))
-
3\,\Phi_1\,H_2(H_2(\Phi_2))
-
\\
&
\ \ \ \ \ \ \ \ \ \ \ \ \
-\,7\,\Phi_2\,H_1(H_2(\Phi_2))
-
7\,\Phi_1\,H_2(H_1(\Phi_1))
-
2\,H_1(\Phi_1)\,H_1(\Phi_2)
-
2\,H_2(\Phi_2)\,H_2(\Phi_1)
+
\\
&
\ \ \ \ \ \ \ \ \ \ \ \ \
+4\,\Phi_1\Phi_2\,H_1(\Phi_1)
+
4\,\Phi_1\Phi_2\,H_2(\Phi_2)
\Big],
\endaligned
\]
and where the remaining four secondary curvature coefficients are
given by:
\[
\aligned
\kappa_{i_1}^{h_2t}
&
=
\kappa_{i_2}^{h_1t},
\\
\kappa_{i_2}^{h_2t}
&
=
-\,\kappa_{i_1}^{h_1t},
\\
\kappa^{h_1t}_j
&
=
\widehat{H}_1\big(\kappa^{h_2t}_{i_2}\big)
-
\widehat{H}_2\big(\kappa^{h_1t}_{i_2}\big),
\\
\kappa^{h_2t}_j
&
=
-\widehat{H}_1\big(\kappa^{h_2t}_{i_1}\big)
+
\widehat{H}_2\big(\kappa^{h_1t}_{i_1}\big),
\endaligned
\]
if one denotes the eight constant vector fields on
$\mathcal{ P}$ associated to the Cartan connection by:
\[
\aligned
\widehat{T}
&
:=
\omega^{-1}({\sf t}),
\ \ \ \ \ \
\widehat{H}_1
:=
\omega^{-1}({\sf h}_1),
\ \ \ \ \ 
\widehat{H}_2
:=
\omega^{-1}({\sf h}_2),
\\
\widehat{D}
&
:=
\omega^{-1}({\sf d}),
\ \ \ \ \
\widehat{R}
:=
\omega^{-1}({\sf r}),
\ \ \ \ \ \
\widehat{I}_1
:=
\omega^{-1}({\sf i}_1),
\ \ \ \ \ 
\widehat{I}_2
:=
\omega^{-1}({\sf i}_2),
\ \ \ \ \
\widehat{J}
&
:=
\omega^{-1}({\sf j});
\endaligned
\]
furthermore and for completeness, 
the $22$ coefficients $\alpha_{ {}_\bullet 
{}_\bullet}$ of these eight constant fields with respect to 
the frame $\{ T, H_1, H_2, D, R, I_1, I_2, J\}$ on $\mathcal{ P}$:
\[
\footnotesize
\aligned
\widehat{T}
&
=
\alpha_{tt}\,T
+\ \
\alpha_{th_1}\,H_1
+\ \
\alpha_{th_2}\,H_2
+\ \
\alpha_{td}\,D
+\
\alpha_{tr}\,R
+\ \
\alpha_{ti_1}\,I_1
+\ \ \
\alpha_{ti_2}\,I_2
+\ \
\alpha_{tj}\,J
\\
\widehat{H}_1
&
=
\ \ \ \ \ \ \ \ \ \ \ \ \
\alpha_{h_1h_1}\,H_1
+
\alpha_{h_1h_2}\,H_2
+
\alpha_{h_1d}\,D
+
\alpha_{h_1r}\,R
+
\alpha_{h_1i_1}\,I_1
+
\alpha_{h_1i_2}\,I_2
+
\alpha_{h_1j}\,J
\\
\widehat{H}_2
&
=
\ \ \ \ \ \ \ \ \ \ \ \ \
\alpha_{h_2h_1}\,H_1
+
\alpha_{h_2h_2}\,H_2
+
\alpha_{h_2d}\,D
+
\alpha_{h_2r}\,R
+
\alpha_{h_2i_1}\,I_1
+
\alpha_{h_2i_2}\,I_2
+
\alpha_{h_2j}\,J
\\
\widehat{D}
&
=
\ \ \ \ \ \ \ \ \ \ \ \ \ \ \ \ \ \ \ \ \ \ \ \ \ \ \ \ \ \ \ \
\ \ \ \ \ \ \ \ \ \ \ \ \ \ \ \ \ \ \ \ \ \ \ \ \ \ \ \ \ \ 
D
\\
\widehat{R}
&
=
\ \ \ \ \ \ \ \ \ \ \ \ \ \ \ \ \ \ \ \ \ \ \ \ \ \ \ \ \ \ \ \
\ \ \ \ \ \ \ \ \ \ \ \ \ \ \ \ \ \ \ \ \ \ \ \ \ \ \ \ \ \ \ \
\ \ \ \ \ \ \ \ \ \ \ \ \ \ \
R
\\
\widehat{I}_1
&
=
\ \ \ \ \ \ \ \ \ \ \ \ \ \ \ \ \ \ \ \ \ \ \ \ \ \ \ \ \ \ \ \
\ \ \ \ \ \ \ \ \ \ \ \ \ \ \ \ \ \ \ \ \ \ \ \ \ \ \ \ \ \ \ \
\ \ \ \ \ \ \ \ \ \ \ \ \ \ \ \ \ \ \ \ \ \ \ \ \ \ \ \ \ \ \ \ \
I_1
\\
\widehat{I}_2
&
=
\ \ \ \ \ \ \ \ \ \ \ \ \ \ \ \ \ \ \ \ \ \ \ \ \ \ \ \ \ \ \ \
\ \ \ \ \ \ \ \ \ \ \ \ \ \ \ \ \ \ \ \ \ \ \ \ \ \ \ \ \ \ \ \
\ \ \ \ \ \ \ \ \ \ \ \ \ \ \ \ \ \ \ \ \ \ \ \ \ \ \ \ \ \ \ \
\ \ \ \ \ \ \ \ \ \ \ \ \ \ \ \ \ \ \
I_2
\\
\widehat{J}
&
=
\ \ \ \ \ \ \ \ \ \ \ \ \ \ \ \ \ \ \ \ \ \ \ \ \ \ \ \ \ \ \ \
\ \ \ \ \ \ \ \ \ \ \ \ \ \ \ \ \ \ \ \ \ \ \ \ \ \ \ \ \ \ \ \
\ \ \ \ \ \ \ \ \ \ \ \ \ \ \ \ \ \ \ \ \ \ \ \ \ \ \ \ \ \ \ \
\ \ \ \ \ \ \ \ \ \ \ \ \ \ \ \ \ \ \ \ \ \ \ \ \ \ \ \ \ \ \ \ \ \
J
\endaligned
\]
are given, in terms of the five variables $(a, b, c, d, e)$ of the
structure group $H$ of the principal bundle $\mathcal{ P}$ and in
terms of the $H_k$-derivatives (up to order $3$) of the fundamental
coefficient functions $\Phi_i$, explicitly by:
\[
\footnotesize\aligned \alpha_{tt} 
&
= 
c^2+d^2, \ \ \ \ \ \alpha_{th_1}
= 
bd-ac, \ \ \
\ \
\alpha_{th_2}
= 
-ad-bc, 
\ \ \ \ \ \ 
\alpha_{h_1h_1} 
& 
= 
c, 
\ \ \ \ \ \ \ \ \ \ \ \ \ \ \ \ \ \ \ \
\ \ \ \ \ \ \ \ \ \ \ \ \ \ \ \ \ \ \ \ \ \ \ \ \
\\
\alpha_{h_1h_2} &= d, \ \ \ \ \ \alpha_{h_2h_1} = -d, \ \ \ \ \ \ \ \ \ \
\alpha_{h_2h_2}=c,
\\
\alpha_{h_1d}&=-2b+{\textstyle{\frac{1}{2}}}\Phi_1c+{\textstyle{\frac{1}{2}}}\Phi_2d,
\ \ \ \
\alpha_{h_2d}=2a+{\textstyle{\frac{1}{2}}}\Phi_2c-{\textstyle{\frac{1}{2}}}\Phi_1d,
\\
\alpha_{h_1r}&=
-6a-{\textstyle{\frac{1}{2}}}\Phi_2c+{\textstyle{\frac{1}{2}}}\Phi_1d, \ \ \ \
\alpha_{h_2r}=-6b+{\textstyle{\frac{1}{2}}}\Phi_1c+{\textstyle{\frac{1}{2}}}\Phi_2d,
\endaligned
\]
\[
\footnotesize\aligned
\alpha_{td}&={\textstyle{\frac{1}{2}}}(bd-ac)\Phi_1-{\textstyle{\frac{1}{2}}}\Phi_2(bc+ad)\Phi_2-2e,
\\
\alpha_{tr}&={\textstyle{\frac{1}{32}}}(H_1(\Phi_1)H_2(\Phi_2))c^2+
 {\textstyle{\frac{1}{32}}}(H_1(\Phi_1)+H_2(\Phi_2))d^2-{\textstyle{\frac{1}{2}}}(ad+bc)\Phi_1+{\textstyle{\frac{1}{2}}}(ac-bd)\Phi_2
 +3a^2+3b^2,
\\
\alpha_{h_1i_1}&={\textstyle{\frac{1}{2}}}(bd+ac)\Phi_1-{\textstyle{\frac{1}{2}}}(bc-ad)\Phi_2-4ab-2e,
\\
\alpha_{h_1i_2}&={\textstyle{\frac{1}{32}}}(H_1(\Phi_1)+H_2(\Phi_2))c^2+{\textstyle{\frac{1}{32}}}(H_1(\Phi_1)+
H_2(\Phi_2))d^2+
{\textstyle{\frac{1}{2}}}(bc-ad)\Phi_1+{\textstyle{\frac{1}{2}}}(ac+bd)\Phi_2+3a^2-b^2,
\\
\alpha_{h_2i_1}&=-{\textstyle{\frac{1}{32}}}H_1((\Phi_1)+H_2(\Phi_2))c^2-{\textstyle{\frac{1}{32}}}(H_1(\Phi_1)+
H_2(\Phi_2))d^2+{\textstyle{\frac{1}{2}}}(bc-ad)\Phi_1+{\textstyle{\frac{1}{2}}}(ac+bd)\Phi_2+
a^2-3b^2,
\\
\alpha_{h_2i_2}&=-{\textstyle{\frac{1}{2}}}(ac+bd)\Phi_1-{\textstyle{\frac{1}{2}}}(ad-bc)\Phi_2+4ab-2e,
\endaligned
\]
\[
\footnotesize\aligned
\ \
\alpha_{ti_1}&={\textstyle{\frac{1}{192}}}\big[-H_2(H_2(\Phi_2))+
H_1(\Phi_1)\Phi_2-{\textstyle{5}}H_2(H_1(\Phi_1))+H_2(\Phi_2)\Phi_2+{\textstyle{4}}H_1(H_1(\Phi_2))\big]d^3+
\\
&+ {\textstyle{\frac{1}{192}}}\big[{\textstyle{4}}H_2(H_1(\Phi_2))+
H_2(\Phi_2)\Phi_1-H_1(H_1(\Phi_1))-{\textstyle{5}}H_1(H_2(\Phi_2))+H_1(\Phi_1)\Phi_1\big]c^3+
\\
&+ {\textstyle{\frac{1}{192}}}\big[{\textstyle{4}}H_2(H_1(\Phi_2))+
H_2(\Phi_2)\Phi_1-H_1(H_1(\Phi_1))-{\textstyle{5}}H_1(H_2(\Phi_2))+H_1(\Phi_1)\Phi_1\big]cd^2+
\\
&+ {\textstyle{\frac{1}{16}}}\big[H_2(\Phi_2)+
H_1(\Phi_1)\big]bc^2+{\textstyle{\frac{1}{192}}}\big[-H_2(H_2(\Phi_2))+H_1(\Phi_1)\Phi_2-{\textstyle{5}}H_2(H_1(\Phi_1))+
\\
&+ H_2(\Phi_2)\Phi_2+{\textstyle{4}}H_1(H_1(\Phi_2))\big]c^2d+
{\textstyle{\frac{1}{16}}}\big[H_2(\Phi_2)+H_1\Phi_1\big]bd^2+
\\
&+
{\textstyle{\frac{1}{2}}}\big[-\Phi_1a^2c+{\textstyle{4}}b^3-\Phi_1b^2c+{\textstyle{4}}ba^2-\Phi_2b^2d-\Phi_2a^2d\big],
\endaligned
\]
\[\footnotesize\aligned
\alpha_{ti_2}&={\textstyle{\frac{1}{192}}}\big[-H_2(H_2(\Phi_2))+H_1(\Phi_1)\Phi_2-{\textstyle{5}}H_2(H_1(\Phi_1))+H_2(\Phi_2)\Phi_2+
{\textstyle{4}}H_1(H_1(\Phi_2))\big]c^3-
\\
&-{\textstyle{\frac{1}{16}}}\big[H_1(\Phi_1)+H_2(\Phi_2)\big]ac^2-{\textstyle{\frac{1}{16}}}\big[H_1(\Phi_1)+H_2(\Phi_2)\big]ad^2+
{\textstyle{\frac{1}{192}}}\big[-H_2(H_2(\Phi_2))+H_1(\Phi_1)\Phi_2-
\\
&-{\textstyle{5}}H_2(H_1(\Phi_1))+H_2(\Phi_2)\Phi_2+{\textstyle{4}}H_1(H_1(\Phi_2))\big]cd^2+
{\textstyle{\frac{1}{192}}}\big[-{\textstyle{4}}H_2(H_1(\Phi_2))-H_1(\Phi_1)\Phi_1-
\\
&-H_2(\Phi_2)\Phi_1+{\textstyle{5}}H_1(H_2(\Phi_2))+
H_1(H_1(\Phi_1))\big]c^2d+{\textstyle{\frac{1}{192}}}\big[-{\textstyle{4}}H_2(H_1(\Phi_2))-H_1(\Phi_1)\Phi_1-
H_2(\Phi_2)\Phi_1+
\\
&+{\textstyle{5}}H_1(H_2(\Phi_2))+H_1(H_1(\Phi_1))\big]d^3-
{\textstyle{\frac{1}{2}}}\big[\Phi_2a^2c+\Phi_2b^2c-\Phi_1b^2d-
\Phi_1a^2d-{\textstyle{4}}ab^2+{\textstyle{4}}a^3\big],
\endaligned
\]
\[
\footnotesize\aligned
\alpha_{h_1j}&={\textstyle{\frac{1}{96}}}\big[-H_2(H_2(\Phi_2)+H_2(\Phi_2)\Phi_2+H_1(\Phi_1)\Phi_2+
{\textstyle{7}}H_2(H_1(\Phi_1))-{\textstyle{8}}H_1(H_1(\Phi_2))\big]c^3-
\\
&+{\textstyle{\frac{1}{96}}}\big[-H_1(\Phi_1)\Phi_1+
{\textstyle{8}}H_2(H_1(\Phi_2))-{\textstyle{7}}H_1(H_2(\Phi_2))-H_2(\Phi_2)\Phi_1+
H_1(H_1(\Phi_1))\big]c^2d+
\\
&+{\textstyle{\frac{1}{96}}}\big[-H_2(H_2(\Phi_2))+H_2(\Phi_2)\Phi_2+H_1(\Phi_1)\Phi_2+
{\textstyle{\frac{7}{16}}}H_2(H_1(\Phi_1))- {\textstyle{8}}H_1(H_1(\Phi_2))\big]cd^2+
\\
&+{\textstyle{\frac{1}{96}}}\big[-H_1(\Phi_1)\Phi_1+{\textstyle{8}}H_2(H_1(\Phi_2))-{\textstyle{7}}H_1(H_2(\Phi_2))-H_2(\Phi_2)\Phi_1+
H_1(H_1(\Phi_1))\big]d^3-
\\
&-{\textstyle{\frac{1}{8}}}\big[H_2(\Phi_2)+H_1(\Phi_1)\big]ac^2-{\textstyle{\frac{1}{8}}}\big[\frac{1}{8}H_2(\Phi_2)+H_1(\Phi_1)\big]ad^2
-\Phi_2a^2c-\Phi_2b^2c+2\Phi_1ce-
\\
&-8be+2\Phi_2de+\Phi_1b^2d+ \Phi_1a^2d-4ab^2-4a^3,
\endaligned
\]
\[
\footnotesize\aligned
\alpha_{h_2j}&={\textstyle{\frac{1}{96}}}\big[-H_1(\Phi_1)\Phi_1+{\textstyle{8}}H_2(H_1(\Phi_2))-{\textstyle{7}}H_1(H_2(\Phi_2))-H_2(\Phi_2)\Phi_1+
H_1(H_1(\Phi_1))\big]c^3+
\\
&-{\textstyle{\frac{1}{8}}}\big[H_2(\Phi_2)+
H_1(\Phi_1)\big]bd^2+{\textstyle{\frac{1}{96}}}\big[-H_2(\Phi_2)\Phi_2-
H_1(\Phi_1)\Phi_2+H_2(H_2(\Phi_2))+8H_1(H_1(\Phi_2))-
\\
&-7H_2(H_1(\Phi_1))\big]d^3-
{\textstyle{\frac{1}{8}}}\big[H_2(\Phi_2)+H_1(\Phi_1)\big]bc^2+{\textstyle{\frac{1}{96}}}\big[-H_1(\Phi_1)\Phi_1+8H_2(H_1(\Phi_2))-
\\
&-
7H_1(H_2(\Phi_2))-H_2(\Phi_2)\Phi_1+H_1(H_1(\Phi_1))\big]cd^2+{\textstyle{\frac{1}{96}}}\big[-H_2(\Phi_2)\Phi_2-
H_1(\Phi_1)\Phi_2+H_2(H_2(\Phi_2))+
\\
&+8H_1(H_1(\Phi_2))-7H_2(H_1(\Phi_1))\big]c^2d+
\Phi_1a^2c-2\Phi_1de+\Phi_2b^2d+\Phi_2a^2d+8ae{\textstyle{4}}b^3+\Phi_1b^2c-
\\
&-4a^2b+2\Phi_2ce,
\endaligned
\]
\[
\footnotesize
\aligned
& 
\alpha_{tj}
=
3a^4+3b^4-4e^2-\Phi_1a^2bc+ca\Phi_2b^2-\Phi_1ab^2d-
\Phi_2a^2bd-2\Phi_2bce-2\Phi_1ace-2\Phi_2ade+2\Phi_1bde-
\\
&-\Phi_1a^3d+\Phi_2a^3c-\Phi_1b^3c-\Phi_2b^3d+6a^2b^2+\big[{\textstyle{\frac{3}{16}}}H_1(\Phi_1)+{\textstyle{\frac{3}{16}}}H_2(\Phi_2)\big]b^2d^2+
\\
&+\big[-{\textstyle{\frac{11}{1536}}}H_2(\Phi_2)H_1(\Phi_1)-
{\textstyle{\frac{1}{192}}}H_1(H_1(\Phi_1))\Phi_1-{\textstyle{\frac{11}{3072}}}H_2({\Phi_2}^2)+{\textstyle{\frac{1}{384}}}{\Phi_2}^2H_2(\Phi_2)-
{\textstyle{\frac{11}{3072}}}H_1({\Phi_1^2})+
\\
&+{\textstyle{\frac{1}{384}}}\Phi_1^2H_1(\Phi_1)+{\textstyle{\frac{1}{48}}}H_1(H_2(H_1(\Phi_2)))+
{\textstyle{\frac{1}{384}}}H_2(H_2(H_2(\Phi_2)))+{\textstyle{\frac{1}{384}}}H_1(H_1(H_1(\Phi_1)))+
{\textstyle{\frac{1}{384}}}{\Phi_2}^2H_1(\Phi_1)-
\endaligned
\]
\[
\footnotesize
\aligned
&-{\textstyle{\frac{1}{192}}}H_2(H_2(\Phi_2))\Phi_2+
{\textstyle{\frac{1}{48}}}H_2(H_1(H_1(\Phi_2)))+{\textstyle{\frac{1}{64}}}H_2(H_1(\Phi_1))\Phi_2-{\textstyle{\frac{1}{48}}}\Phi_1H_2(H_1(\Phi_2))+
{\textstyle{\frac{1}{384}}}{\Phi_1}^2H_2(\Phi_2)-
\\
&-{\textstyle{\frac{7}{384}}}H_2(H_2(H_1(\Phi_1)))+{\textstyle{\frac{1}{64}}}H_1(H_2(\Phi_2))\Phi_1-
{\textstyle{\frac{7}{384}}}H_1(H_1(H_2(\Phi_2)))-{\textstyle{\frac{1}{48}}}\Phi_2H_1(H_1(\Phi_2))\big]d^4+
\\
&+\big[-{\textstyle{\frac{11}{768}}}H_2(\Phi_2)H_1(\Phi_1)-{\textstyle{\frac{7}{192}}}H_2(H_2(H_1(\Phi_1)))+
{\textstyle{\frac{1}{192}}}H_2(H_2(H_2(\Phi_2)))+
{\textstyle{\frac{1}{192}}}H_1(H_1(H_1(\Phi_1)))+
\endaligned
\]
\[
\footnotesize
\aligned
&+{\textstyle{\frac{1}{24}}}H_1(H_2(H_1(\Phi_2)))-
{\textstyle{\frac{1}{96}}}H_2(H_2(\Phi_2))\Phi_2+{\textstyle{\frac{1}{32}}}H_1(H_2(\Phi_2))\Phi_1+
{\textstyle{\frac{1}{192}}}\Phi_2^2H_1(\Phi_1)-{\textstyle{\frac{7}{192}}}H_1(H_1(H_2(\Phi_2)))+
\\
&+ {\textstyle{\frac{1}{192}}}\Phi_2^2H_2(\Phi_2)-
{\textstyle{\frac{11}{1536}}}H_1(\Phi_1^2)-{\textstyle{\frac{1}{24}}}\Phi_2H_1(H_1(\Phi_2))-
{\textstyle{\frac{11}{1536}}}H_2(\Phi_2^2)+{\textstyle{\frac{1}{32}}}H_2(H_1(\Phi_1))\Phi_2-{\textstyle{\frac{1}{96}}}H_1(H_1(\Phi_1))\Phi_1+
\\
&+
{\textstyle{\frac{1}{192}}}\Phi_1^2H_2(\Phi_2)+{\textstyle{\frac{1}{192}}}\Phi_1^2H_1(\Phi_1)-
{\textstyle{\frac{1}{24}}}\Phi_1H_2(H_1(\Phi_2))+{\textstyle{\frac{1}{24}}}H_2(H_1(H_1(\Phi_2)))\big]c^2d^2+
\big[-{\textstyle{\frac{1}{32}}}H_1(H_1(\Phi_1))+
\endaligned
\]
\[
\footnotesize
\aligned
&+{\textstyle{\frac{1}{32}}}H_2(\Phi_2)\Phi_1-{\textstyle{\frac{1}{32}}}H_1(H_2(\Phi_2))+
{\textstyle{\frac{1}{32}}}H_1(\Phi_1)\Phi_1\big]bcd^2+
\big[{\textstyle{\frac{1}{32}}}H_2(H_1(\Phi_1))+{\textstyle{\frac{1}{32}}}H_2(H_2(\Phi_2))-
\\
&-{\textstyle{\frac{1}{32}}}H_2(\Phi_2)\Phi_2-
{\textstyle{\frac{1}{32}}}H_1(\Phi_1)\Phi_2\big]acd^2+
\big[-{\textstyle{\frac{1}{32}}}H_1(H_1(\Phi_1))+{\textstyle{\frac{1}{32}}}H_2(\Phi_2)\Phi_1-{\textstyle{\frac{1}{32}}}H_1(H_2(\Phi_2))+
\\
& +{\textstyle{\frac{1}{32}}}H_1(\Phi_1)\Phi_1\big]ad^3+
\big[{\textstyle{\frac{1}{32}}}H_2(H_1(\Phi_1))+{\textstyle{\frac{1}{32}}}H_2(H_2(\Phi_2))-{\textstyle{\frac{1}{32}}}H_2(\Phi_2)\Phi_2-
{\textstyle{\frac{1}{32}}}H_1(\Phi_1)\Phi_2\big]ac^3+
\endaligned
\]
\[
\footnotesize
\aligned
&+ {\textstyle{\frac{3}{16}}}\big[H_1(\Phi_1)+H_2(\Phi_2)\big]a^2d^2+
{\textstyle{\frac{1}{32}}}\big[H_2(\Phi_2)\Phi_2-H_2(H_1(\Phi_1))-H_2(H_2(\Phi_2))+H_1(\Phi_1)\Phi_2\big]bd^3+
\\
&+\big[-{\textstyle{\frac{1}{32}}}H_1(H_1(\Phi_1))+
{\textstyle{\frac{1}{32}}}H_2(\Phi_2)\Phi_1-
{\textstyle{\frac{1}{32}}}H_1(H_2(\Phi_2))+{\textstyle{\frac{1}{32}}}H_1(\Phi_1)\Phi_1\big]bc^3+
\\
&+{\textstyle{\frac{3}{16}}}\big[H_1(\Phi_1)+H_2(\Phi_2)\big]a^2c^2+
{\textstyle{\frac{3}{16}}}\big[H_1(\Phi_1)+H_2(\Phi_2)\big]b^2c^2+{\textstyle{\frac{1}{32}}}\big[H_2(\Phi_2)\Phi_2-
H_2(H_1(\Phi_1))-
\endaligned
\]
\[
\footnotesize
\aligned
&-
H_2(H_2(\Phi_2))+H_1(\Phi_1)\Phi_2\big]dbc^2+{\textstyle{\frac{1}{32}}}\big[-H_1(H_1(\Phi_1))+
H_2(\Phi_2)\Phi_1- H_1(H_2(\Phi_2))+H_1(\Phi_1)\Phi_1\big]ac^2d+
\\
&+\big[-{\textstyle{\frac{11}{1536}}}H_2(\Phi_2)H_1(\Phi_1)-
{\textstyle{\frac{1}{192}}}H_1(H_1(\Phi_1))\Phi_1-{\textstyle{\frac{11}{3072}}}H_2(\Phi_2^2)+
{\textstyle{\frac{1}{384}}}\Phi_2^2H_2(\Phi_2)-{\textstyle{\frac{11}{3072}}}H_1(\Phi_1^2)+
\\
&+
{\textstyle{\frac{1}{384}}}\Phi_1^2H_1(\Phi_1)+{\textstyle{\frac{1}{48}}}H_1(H_2(H_1(\Phi_2)))+{\textstyle{\frac{1}{384}}}H_2(H_2(H_2(\Phi_2)))+
{\textstyle{\frac{1}{384}}}H_1(H_1(H_1(\Phi_1)))+{\textstyle{\frac{1}{384}}}\Phi_2^2H_1(\Phi_1)-
\endaligned
\]
\[
\footnotesize
\aligned
&-
{\textstyle{\frac{1}{192}}}H_2(H_2(\Phi_2))\Phi_2+{\textstyle{\frac{1}{48}}}H_2(H_1(H_1(\Phi_2)))+{\textstyle{\frac{1}{64}}}H_2(H_1(\Phi_1))\Phi_2-
{\textstyle{\frac{1}{48}}}\Phi_1H_2(H_1(\Phi_2))+{\textstyle{\frac{1}{384}}}\Phi_1^2H_2(\Phi_2)-
\\
&-{\textstyle{\frac{7}{384}}}H_2(H_2(H_1(\Phi_1)))+{\textstyle{\frac{1}{64}}}H_1(H_2(\Phi_2))\Phi_1-{\textstyle{\frac{7}{384}}}H_1(H_1(H_2(\Phi_2)))-
{\textstyle{\frac{1}{48}}}\Phi_2H_1(H_1(\Phi_2))\big]c^4.
\endaligned
\]
\end{Theorem}

For a conceptional, motivational and historical introduction to the
domain that it would be quite useless to reproduce here, the reader is
referred to the excellent expository article (\cite{ EMS}) by Vladimir
Ezhov, Ben McLaughlin and Gerd Schmalz which appeared recently in the
{\em Notices of the American Mathematical Society}. By performing the
above choice $\{ H_1, H_2, T\}$ of an initial frame for $TM$ which is
explicit in terms of the graphing function $\varphi( x, y, u)$, we
deviate from the initial normalization made in~\cite{ EMS} (with a more
geometric-minded approach), our computational objective being to
provide a Cartan-Tanaka connection all elements of which are
completely effective in terms of 
$\varphi (x, y, u)$\,\,---\,\,assuming only
$\mathcal{ C}^6$-smoothness of $M$.

It took about fifteen years (between 1810 and 1826) to Gauss to derive
what he considered to be a completely convincing proof that the
(Gaussian) curvature $\kappa = \kappa (u,v)$ of a surface equipped
with a metric:
\[
ds^2 
=
E(u,v)\,du^2 
+ 
2F(u,v)\,dudv 
+ 
G(u,v)\,dv^2
\]
is a completely {\em intrinsic} invariant through infinitesimal
isometries {\em because} it expresses ({\em Theorema Egregium}, \cite{
Gauss}) for {\em any} surface as the following explicit rational
differential in the second-order jet of the three elements $E, F, G$, 
as Gauss showed:
\[
\footnotesize
\aligned
\kappa 
=
& \ 
\frac{1}{4\,(EG-F^2)^2}
\bigg\{\,
E 
\bigg[
\frac{\partial E}{\partial v} 
\cdot
\frac{\partial G}{\partial v}
-
2\,\frac{\partial F}{\partial u} 
\cdot 
\frac{\partial G}{\partial v} 
+
\frac{\partial G}{\partial u}
\cdot
\frac{\partial G}{\partial u}
\bigg] 
+
\\
&
+
F\bigg[
\frac{\partial E}{\partial u}
\cdot\frac{\partial G}{\partial v}
-
\frac{\partial E}{\partial v} 
\cdot
\frac{\partial G}{\partial u} 
-
2\,\frac{\partial E}{\partial v} 
\cdot\frac{\partial F}{\partial v}
+
4\, 
\frac{\partial F}{\partial u} 
\cdot
\frac{\partial F}{\partial v} 
- 
2\,\frac{\partial F}{\partial u} 
\cdot\frac{\partial G}{\partial u}
\bigg] 
+
\\
& 
\ \ \ \ \ \ \ \ \ \ \ \ \ \ \ \ \ \ \ \ \ \ \ \ \ 
+
G\bigg[
\frac{\partial E}{\partial u} 
\cdot
\frac{\partial G}{\partial u} 
- 
2\,\frac{\partial E}{\partial u} 
\cdot 
\frac{\partial F}{\partial v} 
+
\frac{\partial E}{\partial v}
\cdot
\frac{\partial E}{\partial v}
\bigg]
+
\\
&
+2\,
\big(EG-F^2\big) 
\bigg[
-\frac{\partial^2 E}{\partial v^2} 
+
2\,\frac{\partial^2 F}{\partial u \partial v } 
- 
\frac{\partial^2 G}{\partial u^2}
\bigg]
\,\bigg\}.
\endaligned
\]

However, for what is commonly considered to constitute the {\em
simplest} instance of Cauchy-Riemann geometry, namely for the case of
(embedded) Levi nondegenerate real hypersurfaces $M^3 \subset \C^2$,
what would correspond to the above {\em formula egregia} concerning CR
curvature seems not to have yet ever been achieved, perhaps due to the
fact that in modern differential geometry\,\,---\,\,and in Cartan's
theory of the problem of equivalence as well\,\,---, 
the {\em causality} of
`intrinsic-ness' never relies upon some elimination computations
(Gauss' proof), but it
is set {\em ab initio} in theories.
Nonetheless, a folklore yet unresolved question seems to 
remain: can one
characterize the vanishing of curvature explicitly in terms of
$\varphi ( x, y, u)$?

\begin{Corollary}
The local biholomorphic equivalence to the Heisenberg sphere: $v' =
{x'}^2 + {y'}^2$ of an arbitrary real analytic Levi
nondegenerate hypersurface $M^3 \subset \C^2$ represented as a graph
of the form:
\[
v
=
\varphi(x,y,u)
\]
with $\varphi_{ xx} (0) + \varphi_{ yy} (0) \neq 0$ is explicitly
characterized by the identical vanishing:
\[
0
\equiv
\mathbf{\Delta_1}
\equiv
\mathbf{\Delta_4},
\]
of the two main functions of $(x, y, u)$ appearing in the curvature 
function of
the above effective Cartan-Tanaka connection.
\qed
\end{Corollary}

Of course, thanks to our extensive theorem stated in length,
expansions of these two principal functions $\mathbf{
\Delta_1}$ and $\mathbf{ \Delta}_4$ can straightforwardly be
achieved on a computer by just applying the induction formulas to
which the numerators $A_i$, $A_{ i, k_1}$, $A_{ i, k_1, k_2}$, $A_{ i,
k_1, k_2, k_3}$ of the above fundamental functions $\Phi_i$, $H_{ k_1}
( \Phi_i)$, $H_{ k_2} ( H_{ k_1} ( \Phi_i))$, $H_{ k_3} ( H_{ k_2} (
H_{ k_1} ( \Phi_i)))$ are subjected. As a mathematically satisfactory
fact, the numerators of 
both expressions of $\mathbf{ \Delta_1}$ and $\mathbf{
\Delta}_4$ then become a completely explicit differential polynomial
in the sixth-order jet $J_{ x, y, u}^6 \varphi$ of the graphing
function (the same can be done of course for the functions $\alpha_{
{}_\bullet {}_\bullet}$ too). But their prohibitive
lengths\,\,---\,\,nearly one thousand pages long on a computer, not
copied in this {\LaTeX} file\,\,---\,\,presumably explain why no
reference in the literature ({\em cf.} {\em e.g.}~\cite{Slovak,
Cartan, Chern-Moser, EMS, Isaev, Jacobowitz, Le, Merker5}) succeeded
to be fully effective on the topic, whence unexpectedly and a bit
paradoxically also, {\sc die Gaussche Strenge} (the Gaussian
requirement) for total computational effectiveness in mathematics
happens to be unsatisfiable at human scale 
even in the case of the simplest $M^3 \subset \C^2$.

\subsection*{Acknowledgments}
One year ago, Gerd Schmalz kindly provided us with a pdf copy of the
accepted version of~\cite{ EMS}, and this was of great help during the
(painful) preparation of the present paper. Another article~\cite{
BES} of Valerii Beloshapka, Vladimir Ezhov and Gerd Schmalz was also
used to enter the theory in the right way. The
authors would also like to thank Gerd Schmalz and Ben McLaughlin for
very helpful explanations through e-mail exchanges.


\section{Infinitesimal CR Automorphisms 
\\
of Real Analytic Generic Submanifolds of $\C^{n+d}$}
\label{general-infinitesimal-CR}

\HEAD{Infinitesimal CR Automorphisms 
of Real Analytic generic Submanifolds of $\C^{n+d}$}{
Mansour Aghasi, Joël Merker, and Masoud Sabzevari}

\subsection{Real and complex local equations for generic submanifolds}
Consider a local generic CR submanifold $M\subset \C^{ n+d}$
of positive codimension $d \geqslant 1$ and of positive
CR dimension $n \geqslant 1$ and let $p$ be a point of
$M$. In any system of local holomorphic coordinates $(z, w) = (z_1,
\ldots , z_n, w_1, \ldots ,w_d)$ decomposed in real and imaginary
parts as $(z, w) = (x + iy, u+iv)$ which vanish at $p$ and for
which $T_p M = \{ {\rm Re}\, w_j = 0, j = 1, \dots, d\}$, the generic
submanifold $M \subset \C^{ n+d}$ is locally represented by $d$ real
equations:
\begin{equation}
\label{real-equations-M}
u_j
=
\varphi_j(x,y,v)
\ \ \ \ \ \ \ \ \ \ \ \ \
{\scriptstyle{(j\,=\,1\,\cdots\,d)}}
\end{equation}
as a graph over the $d$-codimensional plane $T_pM$ with of course the
property that the first order jet of each graphing function $\varphi_j$
is zero at the origin:
\[
0
=
\varphi_j(0)
=
\partial_{x_k}\varphi_j(0)
=
\partial_{y_k}\varphi_j(0)
=
\partial_{v_{j'}}\varphi_j(0)
\ \ \ \ \ \ \ \ \ \ \ \ \
{\scriptstyle{(k\,=\,1\,\cdots\,n\,;\,\,\,
j,\,\,j'\,=\,1\,\cdots\,d)}}.
\]
We shall assume in this
section that $M$ is {\em real analytic}, so that the
functions $\varphi_1, \dots, \varphi_d$ are all expandable in Taylor
series converging in a certain neighborhood of the origin
in $\C^n \times \C^n \times \C^d$.

In fact, the adequate invariants of (CR mappings between)
CR manifolds can be viewed mostly when $M$ is represented by 
$d$ so-called {\sl complex defining equations}. Such equations may be
obtained by rewriting the above real equations just as:
\[
{\textstyle{\frac{w_j+\overline{w}_j}{2}}}
=
\varphi_j
\big(
{\textstyle{\frac{z+\overline{z}}{2}}},\,
{\textstyle{\frac{z-\overline{z}}{2{\scriptstyle{\sqrt{-1}}}}}},\,
{\textstyle{\frac{w-\overline{w}}{2{\scriptstyle{\sqrt{-1}}}}}}
\big)
\ \ \ \ \ \ \ \ \ \ \ \ \
{\scriptstyle{(j\,=\,1\,\cdots\,d)}},
\]
and then by solving the so written equations with respect to the
variables $w_j$ by means of the {\em analytic} implicit function
theorem; in this way, one obtains a collection of $d$ equations of the
shape\footnote{\,
Recall that $\C
\{ {\sf x}_1, \dots, {\sf x}_\nu \}^d$ denotes the ring of power
series $\sum_{ \alpha_1, \dots, \alpha_n \in \N}\, C_{ \alpha_1,
\dots, \alpha_n} \cdot {\sf x}_1^{ \alpha_1} \cdots {\sf x}_n^{
\alpha_n}$ with $\C^d$-valued complex coefficients $C_{ \alpha_1,
\dots, \alpha_n}$ which converge in some neighborhood of the origin.
} 
(written in vectorial notation):
\[
w
=
\Theta\big(z,\,\overline{z},\,\overline{w}\big)
=
\sum_{\alpha\in\N^n,\,\beta\in\N^n,\,\gamma\,\in\,\N^d\atop
\vert\alpha\vert+\vert\beta\vert+\vert\gamma\vert\geqslant 1}\,
\Theta_{\alpha,\beta,\gamma}\,
z^\alpha\,\overline{z}^\beta\,\overline{w}^\gamma
\in
\C\big\{z,\,\overline{z},\,\overline{w}\big\}^d,
\]
whose right-hand side converges of course near the origin $(0, 0, 0)
\in \C^n \times \C^n \times \C^d$ and whose (vector) coefficients
$\Theta_{ \alpha, \beta, \gamma} \in \C^d$ are {\em complex}. Since
$d\varphi ( 0) = 0$, one has $\Theta = - \overline{ w} + {\sf
order}\,2\,{\sf terms}$. 

The paradox that any such $d$ {\em complex} equations provide in fact
$2d$ real defining equations for the {\em real} generic submanifold $M
\subset \C^{ n+d}$ which is $d$-codimensional, and also in addition
the fact that one could as well have chosen to solve the above
equations with respect to the $\overline{ w}_j$, instead of the $w_j$,
these two apparent ``contradictions'' are corrected by means of a
fundamental, elementary statement that transfers to $\Theta$ (in a
natural way) the condition of reality enjoyed by the initial defining
$\R^d$-valued map $\varphi$:
\[
\small
\aligned
\overline{\varphi(x,y,v)}
=
\sum_{\vert\alpha\vert+\vert\beta\vert+\vert\gamma\vert\geqslant 1}\,
\overline{\varphi_{\alpha,\beta,\gamma}}\,
\overline{x}^\alpha\overline{y}^\beta\overline{v}^\gamma
=
\sum_{\vert\alpha\vert+\vert\beta\vert+\vert\gamma\vert\geqslant 1}\,
\varphi_{\alpha,\beta,\gamma}\,x^\alpha y^\beta v^\gamma
=
\varphi(x,y,v).
\endaligned
\] 
In the sequel, we shall work exclusively with the complex
graphing functions $\Theta_j$, so we recall a basic
result\footnote{\,
According to a general, common convention, given a power series $\Phi
(Z) = \sum_{ \delta \in \N^N}\, \Phi_\delta\, Z^\delta$, $Z \in \C^N$,
$\Phi_\delta \in \C$, one defines the series $\overline{ \Phi} (Z) :=
\sum_{ \delta \in \N^N }\, \overline{ \Phi }_\delta\, Z^\delta$ by
conjugating only its complex coefficients. Then the complex
conjugation operator distributes oneself simultaneously on functions
and on variables: $\overline{ \Phi(Z)} \equiv \overline{ \Phi}
(\overline{ Z})$, a trivial property which is nonetheless frequently
used in the formal CR reflection principle (\cite{ Merker2,
Merker3, MerkerPorten}). }.
The complex analytic $\C^d$-valued map $\Theta = \Theta
( z, \overline{ z}, \overline{ w} )$ with $\Theta = - \overline{ w} +
{\rm O}(2)$ together with its complex conjugate:
\[
\overline{\Theta}
=
\overline{\Theta}\big(\overline{z},z,w)
=
\sum_{\alpha\in\N^n,\,\beta\in\N^n,\,\gamma\in\N^d}\,
\overline{\Theta}_{\alpha,\beta,\gamma}\,
\overline{z}^\alpha\,z^\beta\,w^\gamma
\in
\C\big\{\overline{z},\,z,\,w\big\}^d
\]
satisfy the two (equivalent by conjugation) collections of
$d$ functional equations:
\begin{equation}
\label{reality-Theta}
\aligned
\overline{w}_j
\equiv
& \
\overline{\Theta}_j
\big(\overline{z},z,\Theta(z,\overline{z},\overline{w})\big)
\ \ \ \ \ \ \ \ \ \ \ \ \ {\scriptstyle{(j\,=\,1\,\cdots\,d)}},
\\
w_j
\equiv
& \
\Theta_j\big(z,\overline{z},
\overline{\Theta}(\overline{z},z,w)\big)
\ \ \ \ \ \ \ \ \ \ \ \ \ {\scriptstyle{(j\,=\,1\,\cdots\,d)}};
\endaligned
\end{equation}
conversely, given a local holomorphic $\C^d$-valued map $\Theta ( z,
\overline{ z}, \overline{ w} ) \in \C \{ z, \overline{ z}, \overline{
w} \}^d$, $\Theta = - \, \overline{ w} + {\rm O} ( 2)$ which, in
conjunction with its complex conjugate $\overline{ \Theta} (
\overline{ z}, z, w)$, satisfies this pair of equivalent identities,
then the two zero-sets:
\[
\big\{
0=-\,w+\Theta\big(z,\,\overline{z},\,\overline{w}\big)
\big\}
\ \ \ \ \ \ \ \ \ \
\text{\rm and}
\ \ \ \ \ \ \ \ \ \
\big\{
0=-\,\overline{w}+
\overline{\Theta}\big(\overline{z},\,z,\,w\big)
\big\}
\]
coincide and define a local generic $d$-codimensional real analytic
submanifold passing through the origin in $\C^{ n + d}$. In fact more
precisely, one may show (\cite{ Merker2, Merker3}) that there is an
invertible $d \times d$ matrix $a ( z, w, \overline{ z}, \overline{
w})$ of analytic functions defined near the origin such that one has:
\[
w-\Theta(z,\overline{z},\overline{w})
\equiv
a(z,w,\overline{z},\overline{w})\,
\big[
\overline{w}-\overline{\Theta}(\overline{z},z,w)
\big],
\]
identically in $\C \{ z, w, \overline{ z}, \overline{ w} \}$,
whence the coincidence of the two zero-sets immediately follows.

\subsection{Extrinsic complexification} As is known in local analytic
CR geometry,
it is natural to introduce new independent complex variables
$(\underline{ z}, \underline{ w}) \in \C^n \times 
\C^d$\,\,---\,\,underlining
should {\em not} be confused here with conjugating\,\,---\,\,and 
to define the so-called
{\sl extrinsic complexification} $M^{e_c}$ of $M$ as being
the complex analytic $d$-codimensional submanifold of $\C^{n+d}\times
\C^{n+d}$ equipped with the $2n + 2d$ coordinates 
$(z, w, \underline{ z}, \underline{ w})$ which is
defined by the $d$ equations:
\[
w_j
=
\Theta_j(z,\underline{z},\underline{w})
\ \ \ \ \ \ \ \ \ \ \ \ \ {\scriptstyle{(j\,=\,1\,\cdots\,d)}}.
\] 
Notice that the replacement $( \overline{ z}, \overline{ w})$ by $(
\underline{ z}, \underline{ w})$ in the Taylor series of $\Theta$ is
really meaningful:
\[
\Theta_j(z,\underline{z},\underline{w})
:=
\sum_{\alpha\in\N^n,\,\beta\in\N^n,\,\gamma\in\N^d}\,
\Theta_{j,\alpha,\beta,\gamma}\,
z^\alpha\,\underline{z}^\beta\,\underline{w}^\gamma,
\]
thanks to the fact that the series converges locally. Equivalently,
$M^{ e_c}$ is defined by the $d$ equations $\underline{ w}_j =
\overline{ \Theta}_j ( \underline{ z}, z, w)$. Then $M$ is recovered
from $M^{ e_c}$ by just replacing these independent variables $(
\underline{ z}, \underline{ w})$ by the original conjugates $(
\overline{ z}, \overline{ w})$. The following standard uniqueness
principle is useful.

\begin{Lemma}
\label{uniqueness-lemma}
Consider a complex-valued converging power series:
\[
\Phi
= 
\Phi(z,w,\overline{z},\overline{ w})
=
\sum_{\alpha\in\N^n,\,\beta\in\N^d,\,\gamma\in\N^n,\,\delta\in\N^d}\,
\Phi_{\alpha,\beta,\gamma,\delta}\,
z^\alpha\,w^\beta\,\overline{z}^\gamma\,\overline{w}^\delta
\]
in $\C \{ z, w, \overline{ z}, \overline{ w} \}$ having complex
coefficients $\Phi_{ \alpha, \beta, \gamma, \delta} \in \C$. Then 
the following four properties are equivalent:

\begin{itemize}

\smallskip\item[{\bf (i)}]
$\Phi$ takes only the value zero when the point $(z, w)$ varies (without
restriction) on $M \subset \C^n$;

\smallskip\item[{\bf (ii)}]
the extrinsic complexification of $\Phi$:
\[
\Phi^{e_c}
= 
\Phi^{e_c}(z,w,\underline{z},\underline{ w})
:=
\sum_{\alpha\in\N^n,\,\beta\in\N^d,\,\gamma\in\N^n,\,\delta\in\N^d}\,
\Phi_{\alpha,\beta,\gamma,\delta}\,
z^\alpha\,w^\beta\,\underline{z}^\gamma\,\underline{w}^\delta
\]
takes only the value zero when the point $(z, w, \underline{
z}, \underline{ w})$ varies (without restriction) on 
the complexification $M^{ e_c}
\subset \C^{ 2n+2d}$;

\smallskip\item[{\bf (iii)}]
after replacing $\underline{ w}$ by $\overline{ \Theta} ( \underline{
z}, z, w)$ in the extrinsic complexification $\Phi^{ e_c}$ of $\Phi$,
the result is an identically zero series in $\C\{ \underline{ z}, z, w
\}^d$, namely:
\[
0
\equiv
\sum_{\alpha\in\N^n,\,\beta\in\N^d,\,\gamma\in\N^n,\,\delta\in\N^d}\,
\Phi_{\alpha,\beta,\gamma,\delta}\,
z^\alpha\,w^\beta\,\underline{z}^\gamma\,
[\overline{\Theta}(\underline{z},z,w)]^\delta;
\] 

\item[{\bf (iv)}]
after replacing $w$ by $\Theta ( z, \underline{ z}, \underline{ w})$
in the extrinsic complexification $\Phi^{ e_c}$, the result
is an identically zero power series in $\C\{ z, \underline{ z}, 
\underline{ w} \}^d$, namely:
\[
0
\equiv
\sum_{\alpha\in\N^n,\,\beta\in\N^d,\,\gamma\in\N^n,\,\delta\in\N^d}\,
\Phi_{\alpha,\beta,\gamma,\delta}\,
z^\alpha\,
[\Theta(z,\underline{z},\underline{w})]^\beta\,
\underline{z}^\gamma\,\underline{w}^\delta.
\qed
\]

\end{itemize}\smallskip

\end{Lemma}

Let $Z \in \C^N$. A converging power series $\Phi ( Z) \in \C \{ Z\}$
will be called a {\sl holomorphic} function. A converging power series
$\Pi ( \overline{ Z}) \in \C \{ \overline{ Z} \}$ will be called an
{\sl antiholomorphic} function. But in general, in local analytic
Cauchy-Riemann geometry, some variables $Z$ and $\overline{ Z}$ are
mixed or considered together. Because any converging power series
$\Psi ( Z, \overline{ Z}) \in \C\{ Z, \overline{ Z}\}$ may also be
considered as the series:
\[
\Psi^\sim({\rm Re}\,Z,\,{\rm Im}\,Z) 
:= 
\Psi({\rm Re}\,Z+i\,{\rm Im}\,Z,\,\, 
{\rm Re}\,Z-i\,{\rm Im}\,Z), 
\]
belonging to $\C\{ {\rm Re}\, Z, \, {\rm Im}\, Z\}$ (in terms of the
basic {\em real} $2N$ variables ${\rm Re}\, Z$ and ${\rm Im}\, Z$),
such a series $\Psi$ will be called a {\sl real analytic} function,
not only when it has purely real values, namely when $\Psi^\sim \in \R
\{ {\rm Re}\, Z, \, {\rm Im}\, Z\}$, but also when it has complex
values, namely when $\Psi^\sim \in \C \{ {\rm Re}\, Z, \, {\rm Im}\,
Z\}$. Thus, the terminology ``{\sl real analytic}'' is used for $(Z,
\overline{ Z})$-dependence.

\subsection{Holomorphic and antiholomorphic tangent vector fields}
In such coordinates $(z, w) \in \C^n \times \C^d$, we claim that the
bundle $T^{ 1,0}M$ and its conjugate $T^{ 0, 1} M = \overline{ T^{ 1,
0}M}$ are generated, respectively, by the two collections of mutually
independent $(1, 0)$ and $(0, 1)$ vector fields:
\[
\aligned
\mathcal{L}_k
&
=
\frac{\partial}{\partial z_k}
+
\sum_{j=1}^d\,
\frac{\partial\Theta_j}{\partial z_k}(z,\overline{z},\overline{w})\,
\frac{\partial}{\partial w_j}
\ \ \ \ \ \ \ \ \ \ \ \ \ {\scriptstyle{(k\,=\,1\,\cdots\,n)}},
\\
\overline{\mathcal{L}}_k
&
=
\frac{\partial}{\partial\overline{z}_k}
+
\sum_{j=1}^d\,
\frac{\partial\overline{\Theta}_j}{\partial\overline{z}_k}
(\overline{z},z,w)\,
\frac{\partial}{\partial\overline{w}_j}
\ \ \ \ \ \ \ \ \ \ \ \ \ {\scriptstyle{(k\,=\,1\,\cdots\,n)}}
\endaligned
\]
having real analytic coefficients. Indeed, one checks immediately
that $0 \equiv \mathcal{ L}_k \big( w_j - \Theta_j(z, \overline{ z},
\overline{ w}) \big)$ and that $0 \equiv \overline{ \mathcal{ L}}_k
\big( \overline{ w_j} - \overline{ \Theta} ( \overline{ z}, z, w)
\big)$, and since the two complex vector bundles $T^{ 1, 0} M$ and
$T^{ 0, 1} M$ are known (\cite{ Boggess, BER, MerkerPorten}) 
to be of rank $n = {\rm CRdim}\, M$ (which
truly means that $M$ is generic), the claim is clear. Of course,
these two collections of vector fields have extrinsic
complexifications:
\[
\aligned
\mathcal{L}_k^{e_c}
&
=
\frac{\partial}{\partial z_k}
+
\sum_{j=1}^d\,
\frac{\partial\Theta_j}{\partial z_k}
(z,\underline{z},\underline{w})\,
\frac{\partial}{\partial w_j}
\ \ \ \ \ \ \ \ \ \ \ \ \ {\scriptstyle{(k\,=\,1\,\cdots\,n)}},
\\
\underline{\mathcal{L}}_k^{e_c}
&
=
\frac{\partial}{\partial\underline{z}_k}
+
\sum_{j=1}^d\,
\frac{\partial\overline{\Theta}_j}{\partial\underline{z}_k}
(\underline{z},z,w)\,
\frac{\partial}{\partial\underline{w}_j}
\ \ \ \ \ \ \ \ \ \ \ \ \ {\scriptstyle{(k\,=\,1\,\cdots\,n)}}.
\endaligned
\]

\subsection{Intrinsic generators of $T^cM$}
It is also useful to write the $(1, 0)$ and $(0, 1)$ vector fields
tangent to $M$ in terms of the {\em real} graphed defining equations
$u_j = \varphi_j ( x, y, v)$ of $M$. For any $k = 1, \dots, n$, a $(1,
0)$ vector field of the general form:
\[
\mathcal{L}_k
=
\frac{\partial}{\partial z_k} 
+ 
\sum_{l=1}^d\,{\tt A}_{k,l}\,
\frac{\partial}{\partial w_l}
\]
is tangent to the $d$ real equations of $M$:
\[
0
=
-\,u_j+\varphi_j(x,y,u)
\ \ \ \ \ \ \ \ \ \ \ \ \ {\scriptstyle{(j\,=\,1\,\cdots\,d)}}
\]
if and only if its $d$ complex coefficients ${\tt A}_{k, l}$ satisfy, 
on restriction to $M$, the following $n\, d$ scalar equations:
\[
0
=
-\,
\frac{1}{2}\,{\tt A}_{k,j}
-
\frac{i}{2}\,
\sum_{l=1}^d\,{\tt A}_{k,l}\,\varphi_{j,v_l}
+
\varphi_{j,z_k}
\ \ \ \ \ \ \ \ \ \ \ \ \ 
{\scriptstyle{(k\,=\,1\,\cdots\,n\,;\,\,\,j\,=\,1\,\cdots\,d)}}.
\]
Fixing $k$, if one introduces the column matrix ${\tt A}_k :=
\big( {\tt A}_{k,1}, \dots, {\tt A}_{k, d} \big)^{\tt t}$, the column
matrix $\varphi_{ z_k} := \big( \varphi_{ 1, z_k}, \dots,
\varphi_{ d, z_k} \big)^{\tt t}$ and the $d \times d$ matrix $\varphi_v :=
\big( \varphi_{ j, v_l}\big)_{ 1 \leqslant j \leqslant d}^{ 1
\leqslant l \leqslant d}$ in which $j$ is the index of lines, the
corresponding $d$ equations, when rewritten as:
\[
2\,\varphi_{z_k}
=
\Big(
\big(\delta_{j,l}+i\,\varphi_{j,v_l}\big)_{
1\leqslant j\leqslant d}^{1\leqslant l\leqslant d}
\Big)
\cdot{\tt A}_k
\ \ \ \ \ \ \ \ \ \ \ \ \ {\scriptstyle{(j\,=\,1\,\cdots\,d)}}
\]
constitute a linear system of $d$ equations in the $d$ unknowns
${\tt A}_{ k, 1}, \dots, {\tt A}_{ k, d}$ which may 
be solved by means of a matrix inversion:
\[
{\tt A}_k
=
2\,\big(I+i\,\varphi_v\big)^{-1}\cdot\varphi_{z_k}.
\]
Then the decomposition in real and imaginary parts: 
\[
{\tt A}_k
=
{\tt A}_k'
+
i\,{\tt A}_k''
\ \ \ \ \ \ \ \ \ \ \ \ \ {\scriptstyle{(k\,=\,1\,\cdots\,n)}}
\]
of these coefficients writes:
\[
\aligned
{\tt A}_k'
&
=
\big(I+i\,\varphi_v\big)^{-1}\cdot\varphi_{z_k}
+
\big(I-i\,\varphi_v\big)^{-1}\cdot\varphi_{\overline{z}_k},
\\
{\tt A}_k''
&
=
-\,i\,\big(I+i\,\varphi_v\big)^{-1}\cdot\varphi_{z_k}
+
i\,\big(I-i\,\varphi_v\big)^{-1}\cdot\varphi_{\overline{z}_k}.
\endaligned
\]
In this way, transposing the column matrix of
basic $\frac{ \partial}{ \partial w_l}$ derivations, we
obtain precisely the right
number $n = {\rm CRdim}\, M$ linearly independent
generators of $T^{ 1, 0} M$:
\[
\aligned
\mathcal{L}_k
&
=
\frac{1}{2}\,\frac{\partial}{\partial x_k}
-
\frac{i}{2}\,\frac{\partial}{\partial y_k}
+
\bigg(
\frac{\partial}{\partial w}
\bigg)^{\tt t}
\cdot
{\tt A}_k
\\
&
=
\frac{1}{2}\,\frac{\partial}{\partial x_k}
-
\frac{i}{2}\,\frac{\partial}{\partial y_k}
+
\bigg(
\frac{1}{2}\,
\frac{\partial}{\partial u}
-
\frac{i}{2}\,
\frac{\partial}{\partial v}
\bigg)^{\tt t}
\cdot
\big(
{\tt A}_k'+i\,{\tt A}_k''
\big).
\endaligned
\]
However, such $n$ generators $\mathcal{ L}_k$ of $T^{ 1, 0} M$ are
still {\em extrinsic}, namely they involve the coordinates $u_j$, and
if we want to pull-back them to the generic submanifold $M$ that needs
only its intrinsic coordinates $(x, y, v)$:
\[
\aligned
H_k^1
&
:=
2\,{\rm Re}\,\big(\mathcal{L}_k\big\vert_M\big),
\\
H_k^2
&
:=
-\,2\,{\rm Im}\,\big(\mathcal{L}_k\big\vert_M\big),
\endaligned
\]
we just have to drop the $\frac{ \partial}{ \partial u}$ above, and we
receive in this way $2n$ generators for the intrinsic complex tangent
bundle $T^c M = {\rm Re} \, T^{ 1, 0} M$:

\[
\aligned
H_k^1
&
=
\frac{\partial}{\partial x_k}
+
\bigg(
\frac{\partial}{\partial v}
\bigg)^{\tt t}
\cdot
{\tt A}_k''
\ \ \ \ \ \ \ \ \ \ \ \ \ {\scriptstyle{(k\,=\,1\,\cdots\,n)}},
\\
H_k^2
&
=
\frac{\partial}{\partial y_k}
+
\bigg(
\frac{\partial}{\partial v}
\bigg)^{\tt t}
\cdot
{\tt A}_k'
\ \ \ \ \ \ \ \ \ \ \ \ \ {\scriptstyle{(k\,=\,1\,\cdots\,n)}}.
\endaligned
\]

\subsection{Infinitesimal CR automorphisms}
According to~\cite{ Stanton, Beloshapka, BES}, a {\sl (local)
infinitesimal CR-automorphism of $M$}, when understood extrinsically,
is a holomorphic vector field:
\[
{\sf X}
=
\sum_{k=1}^n\,Z^k(z,w)
\frac{\partial}{\partial z_k}
+
\sum_{j=1}^d\,W^j(z,w)
\frac{\partial}{\partial w_j}
\]
whose real part ${\rm Re}\, X = \frac{ 1}{ 2} ( {\sf X} +
\overline{\sf X} )$ is tangent to $M$. (One should mind that, contrary
to the above $(1, 0)$ generators $\mathcal{ L}_k$ of $T^{ 1, 0} M$,
such an ${\sf X}$ is supposed to have purely holomorphic coefficients,
whereas the $\frac{ \partial \Theta_j}{ \partial z_k} ( z, \overline{
z}, \overline{ w})$ are\,\,---\,\,most 
of the time\,\,---\,\,neither purely holomorphic,
nor purely antiholomorphic, but only real analytic.) Determining all
such ${\sf X}$'s is the same as knowing the {\em CR symmetries} of
$M$, a question which lies at the heart of the problem of classifying
all local analytic CR manifolds up to biholomorphisms.

By integration, the {\em real} flow:
\[
(t,z,w)
\longmapsto
\exp(t\,{\sf X})(z,w)
\ \ \ \ \ \ \ \ \ \ \ \ \
{\scriptstyle{(t\,\in\,\R\,\,\,
\text{\rm small})}}
\]
constitutes a local one-parameter group of local biholomorphisms of
$\C^n$; because ${\sf X}$ is tangent to $M$, this flow leaves $M$
invariant (locally), that is to say: through this flow, points of $M$ are
transferred to points of $M$. We note {\em passim} that this
real flow coincides with restricting the consideration of the 
complex (holomorphic) flow:
\[
(\tau,z,w)
\longmapsto
\exp(\tau\,{\sf X})(z,w)
\ \ \ \ \ \ \ \ \ \ \ \ \
{\scriptstyle{(\tau\,\in\,\C\,\,\,
\text{\rm small})}}
\]
to real time parameters $\tau := t \in \R$, a fact
we will not use. 
Conversely, one may show:

\begin{Lemma} 
{\em (\cite{ Stanton, BER})}
If $M \subset \C^{n+d}$ is a generic submanifold and if $(z, w)
\longmapsto \phi_t ( z, w)$ is a local real one-parameter group of {\em
holomorphic} self-transformations of $\C^{n+d}$ which stabilizes $M$
locally, then the vector field:
\[
{\textstyle{\frac{d}{dt}}}\big\vert_0
\big(
\phi_t(z,w)
\big)
\]
has holomorphic coefficients and its real part is tangent to $M$. 
\qed
\end{Lemma}

Since holomorphy of coefficients and tangency to a submanifold is
preserved under taking Lie brackets, the collection $\mathfrak{ hol}
(M)$ of all such ${\sf X}$ is obviously a Lie algebra. Also, when
$\mathfrak{ hol} (M)$ is finite-dimensional (which occurs except in
degenerate situations, {\em see} {\em e.g.}~\cite{ GaussierMerker}),
the corresponding finite-dimensional local Lie group is real, whence
$\mathfrak{ hol} (M)$ constitutes a {\em real} Lie algebra. So
according to one of Lie's fundamental theorems (\cite{ Merker7},
Chap.~9), if ${\sf X}_1, \dots, {\sf X}_r$ denote any basis of
$\mathfrak{ hol}( M)$ as a vector space, there are {\em real}
structure constants $c_{ jk}^s \in \R$ such that:
\begin{equation}
\label{real-structure-constants}
\big[{\sf X}_j,\,{\sf X}_k\big]
=
\sum_{s=1}^r\,c_{jk}^s\,{\sf X}_s.
\end{equation}
For an explicitly given $M \subset \C^{ n+d}$, determining a basis of
the Lie algebra $\mathfrak{ hol} ( M)$ is a natural problem for which
systematic computational procedures exists, as we will establish in a
while. The groundbreaking works of Sophus Lie and his followers
(Friedrich Engel, Georg Scheffers, Gerhard Kowalewski, Ugo Amaldi and
others) showed that the most fundamental question in concern here is to
draw up lists of possible Lie algebras $\mathfrak{ hol} ( M)$ which
would classify all possible $M$'s according to their CR symmetries.

Alternatively, if one prefers to view the CR manifold $M$ in a
purely intrinsic way, one may consider the local group ${\rm Aut}_{
CR} ( M)$ of automorphisms of the CR structure, namely of local
$\mathcal{ C}^\infty$ diffeomorphisms $g \colon M \to M$ (close to the
identity mapping) which satisfy:
\[
dg_p\big(T_p^cM\big) 
= 
T_{g(p)}^cM
\ \ \ \ \
\text{\rm and}
\ \ \ \ \ 
dg_p\big(J(v_p)\big)
= 
J_{g(p)}\big(dg_p(v_p)\big)
\] 
at any point $p \in M$ and for any complex-tangent
vector $v_p \in T_p^c M$. In other
words, $g$ belongs to ${\rm Aut}_{ CR} ( M)$ if and only if it is a
(local) {\sl CR-diffeomorphism} of $M$, namely a diffeomorphism which
respects the (intrinsic) CR structure of $M$. As did Lie most of the
time in his original theory (\cite{ Merker7, EngelLie}), we shall
consider only a neighborhood of the identity mapping, hence all our
groups will be {\sl local Lie groups}; the reader is again referred
to~\cite{ Merker7, Olver1, GaussierMerker}) 
for fundamentals about local Lie groups in
general, especially concerning the fact that it is essentially useless
to point out open sets and domains in which mappings and
transformations are defined.

Accordingly, let: 
\[
\mathfrak{aut}_{CR}(M)
\] 
denote the collection of all (real) vector fields ${\sf Y}$ on $M$ the
flow of which $(t, p) \mapsto \exp ( t {\sf Y}) ( p)$ becomes a local
CR diffeomorphism of $M$. When ${\rm Aut}_{ CR} (M)$ is a
finite-dimensional Lie group, $\mathfrak{ aut}_{ CR} ( M)$ is just its
Lie algebra.

\begin{Lemma}
{\rm (\cite{Merker7}, Chap.~8; \cite{BER})} A local real
analytic vector field ${\sf Y}$ on $M$ belongs to $\mathfrak{ aut}_{
CR} ( M)$, if and only if for every local section $L$ of the complex
tangent bundle $T^cM$, the Lie bracket $[ {\sf Y}, \, L]$ is again a
section of $T^cM$.
\qed
\end{Lemma}

In all cases which are of interest, namely when $M$ is nondegenerate
in some sense (the interested
reader is referred to~\cite{ BER, GaussierMerker, 
Merker2, Merker3, MerkerPorten}, for we prefer
not to dwell on that topic here), such $\mathcal{
C}^\infty$ flows $(t, p) \mapsto \exp ( t {\sf Y}) ( p)$ happen to be
in fact {\em real analytic}, whence, according to a classical theorem,
they extend as local {\em biholomorphic} maps from a neighborhood of
$M$ in $\C^{n+d}$. If follows that any such
intrinsic ${\sf Y}$
happens in fact to
be a restriction ${\sf Y} = {\sf X} \big\vert_M$ to $M$ of
some extrinsic ${\sf X} \in \mathfrak{ hol} (M)$. 
In all these circumstances which cover a broad universe of
yet unstudied CR structures, one has the fundamental relation:
\[
\boxed{
\mathfrak{aut}_{CR}(M)
=
{\rm Re}
\big(
\mathfrak{hol}(M)
\big)}\,,
\]
where both sides are finite-dimensional, spanned by vector fields
whose coefficients are expandable in converging power series. Thus,
one may work exclusively with the {\em holomorphic} vector fields
generating $\mathfrak{ hol} ( M)$, as we will do from now on. And in
any case, there will be no confusion to call an {\em infinitesimal CR
automorphism} either the holomorphic vector field ${\sf X} \in
\mathfrak{ hol} ( M)$ or its real part $\frac{ 1}{ 2} ( {\sf X} +
\overline{\sf X}) \in \mathfrak{ aut}_{CR} ( M)$.

Since holomorphic vector fields obviously commute with antiholomorphic
vector fields, we deduce from~\thetag{ \ref{real-structure-constants}}
that when $\mathfrak{ hol} ( M) = \R {\sf X}_1 \oplus \cdots \oplus \R
{\sf X}_r$ is $r$-dimensional, the real parts of the ${\sf X}_j$ which
generate $\mathfrak{ aut}_{ CR} ( M)$ simply have the same (real)
structure constants:
\begin{equation}
\big[{\sf X}_j+\overline{\sf X}_j,\,
{\sf X}_k+\overline{\sf X}_k\big]
=
\big[{\sf X}_j,\,{\sf X}_k\big]
+
\big[\overline{\sf X}_j,\,\overline{\sf X}_k\big]
=
\sum_{s=1}^r\,c_{jk}^s\,
\big({\sf X}_s+\overline{\sf X}_s\big).
\end{equation}

To conclude these generalities, at any fixed point $p \in M$, one may
also consider the Lie subalgebras $\mathfrak{ hol}( M, p)$ of
$\mathfrak{ hol} ( M)$ and $\mathfrak{ aut}_{ CR}( M, p)$ of
$\mathfrak{ aut}_{ CR} ( M)$ consisting of those vector fields whose
values vanish at $p$. Then $\mathfrak{ hol}( M, p)$ and $\mathfrak{
aut}_{ CR}( M, p)$ are the Lie algebra of the subgroups ${\rm Hol} (
M, p)$ of ${\rm Hol} ( M)$ and ${\rm Aut}_{ CR} ( M, p)$ of ${\rm
Aut}_{ CR} ( M)$ consisting of only the maps that fix the point
$p$. Of course, one has $\mathfrak{ aut}_{ CR} ( M, p) = {\rm Re}
\big( {\sf hol} ( M, p) \big)$.

\subsection{Effective tangency equations}
In order to compute $\mathfrak{ hol} ( M)$ for an explicitly given
generic submanifold $M \subset \C^{n+d}$, it is most convenient, as
already pointed out, to work with {\sl complex defining equations} of
the specific shape (\cite{ MerkerPorten, Merker2}):
\[
\overline{w}_j+w_j
=
\overline{\Xi}_j(\overline{z},z,w)
\ \ \ \ \ \ \ \ \ \ \ \ \
{\scriptstyle{(j\,=\,1\,\cdots\,d)}},
\]
so that, compared to the notations introduced a moment ago, one should
consider that the following 
{\sl notational coincidence} holds:
\[
\overline{\Theta}_j
(\overline{z},z,w)
\equiv
-w_j
+
\overline{\Xi}_j(\overline{z},z,w)
\ \ \ \ \ \ \ \ \ \ \ \ \
{\scriptstyle{(j\,=\,1\,\cdots\,d)}}.
\]
Concretely and precisely, the condition that a general holomorphic
vector field ${\sf X} = \sum_{ k=1}^n\, Z^k ( z, w)\, \frac{
\partial}{ \partial z_k} + \sum_{ l=1}^d\, W^l( z, w) \, \frac{
\partial}{ \partial w_l}$ belongs to $\mathfrak{ hol} ( M)$, namely
that $\overline{\sf X} + {\sf X}$ is tangent to $M$, means that each
one of the following $d$ differentiated equation:
\[
\footnotesize
\aligned
0
=
(\overline{\sf X}+{\sf X})
&
\big[
\overline{w}_j+w_j
-
\overline{\Xi}_j(\overline{z},z,w)
\big]
=
\\
&
=
\overline{\sf X}\,\big[\overline{w}_j+w_j
-
\overline{\Xi}_j(\overline{z},z,w)
\big]
+
{\sf X}\,
\big[\overline{w}_j+w_j
-
\overline{\Xi}_j(\overline{z},z,w)
\big]
\\
&
=
\overline{W}^j(\overline{z},\overline{w})
-
\sum_{k=1}^n\,\overline{Z}^k(\overline{z},\overline{w})\,
\frac{\partial\overline{\Xi}_j}{\partial
\overline{z}_k}
(\overline{z},z,w)
+
\\
&
\ \ \ \ \
+
W^j(z,w)
-
\sum_{k=1}^n\,
Z^k(z,w)\,
\frac{\partial\overline{\Xi}_j}{\partial z_k}(\overline{z},z,w)
-
\sum_{l=1}^d\,W^l(z,w)\,
\frac{\partial\overline{\Xi}_j}{\partial w_l}(\overline{z},z,w)
\\
&
\ \ \ \ \ \ \ \ \ \ \ \ \ \ \ \ \ \ \ \ \ \ \ \ \ \ 
\ \ \ \ \ \ \ \ \ \ \ \ \ 
{\scriptstyle{(j\,=\,1\,\cdots\,d)}}
\endaligned
\]
should vanish for every $(z, w) \in M$. According to
Lemma~\ref{uniqueness-lemma}, this condition holds true if and only if,
after extrinsic complexification and replacement of 
$\underline{ w}$ by
$- w + 
\overline{ \Xi} ( \underline{ z}, z, w)$, the $d$ power
series obtained in $\C \{ \underline{ z}, z, w \}$
vanish identically, namely if and only if:
\[
\footnotesize
\aligned
0
\equiv
\bigg[
&
\overline{W}^j(\underline{z},\underline{w})
-
\sum_{k=1}^n\,\overline{Z}^k(\underline{z},\underline{w})\,
\frac{\partial\overline{\Xi}_j}{\partial
\underline{z}_k}
(\underline{z},z,w)
+
\\
&
+
W^j(z,w)
-
\sum_{k=1}^n\,
Z^k(z,w)\,
\frac{\partial\overline{\Xi}_j}{\partial z_k}(\underline{z},z,w)
-
\sum_{l=1}^d\,W^l(z,w)\,
\frac{\partial\overline{\Xi}_j}{\partial w_l}(\underline{z},z,w)
\bigg]_{\underline{w}=-w+\overline{\Xi}(\underline{z},z,w)}
\\
&
\ \ \ \ \ \ \ \ \ \ \ \ \ \ \ \ \ \ \ \ \ \ \ \ \ \ 
\ \ \ \ \ \ \ \ \ \ \ \ \ \ \ \ \ \
{\scriptstyle{(j\,=\,1\,\cdots\,d)}},
\endaligned
\]
or else in greater details, when one really performs the said
substitution:
\begin{equation}
\label{initial-tangency}
\footnotesize
\aligned
0
&
\equiv
\overline{W}^j\big(\underline{z},\,-w+\overline{\Xi}(\underline{z},z,w)\big)
-
\sum_{k=1}^n\,\overline{Z}^k\big(
\underline{z},\,-w+\overline{\Xi}(\underline{z},z,w)\big)\,
\frac{\partial\overline{\Xi}_j}{\partial\underline{z}_k}
(\underline{z},z,w)
+
\\
&
\ \ \ \ \
+
W^j(z,w)
-
\sum_{k=1}^n\,Z^k(z,w)\,
\frac{\partial\overline{\Xi}_j}{\partial z_k}
(\underline{z},z,w)
-
\sum_{l=1}^d\,
W^l(z,w)\,
\frac{\partial\overline{\Xi}_j}{\partial w}
(\underline{z},z,w)
\\
&
\ \ \ \ \ \ \ \ \ \ \ \ \ \ \ \ \ \ \ \ \ \ \ \ \ \ 
\ \ \ \ \ \ \ \ \ \ \ \ \ \ \ \ \ \
{\scriptstyle{(j\,=\,1\,\cdots\,d)}}.
\endaligned
\end{equation}
Interestingly enough, this condition may also be interpreted as saying that
the complexified sum of vector fields:
\[
\footnotesize
\aligned
\big(\overline{\sf X}\big)^{e_c}
+
{\sf X}^{e_c}
:=
&\,
\sum_{k=1}^n\,\overline{Z}_k(\underline{z},\underline{w})\,
\frac{\partial}{\partial\underline{z}_k}
+
\sum_{l=1}^d\,\overline{W}^l(\underline{z},\underline{w})\,
\frac{\partial}{\partial\underline{w}_l}
+
\sum_{k=1}^n\,Z^k(z,w)\,
\frac{\partial}{\partial z_k}
+
\sum_{l=1}^d\,W^l(z,w)\,
\frac{\partial}{\partial w_l}
\\
=:
&\,
\underline{\sf X}
+
{\sf X}
\endaligned
\]
is tangent to $M^{e_c}$, {\em cf.}~\cite{Merker4} for similar
considerations in the broader context of completely integrable
analytic systems of partial differential equations. But we must now
analyze further what this condition really means.

To this aim, we may at first introduce the expansions of the
coefficients of such a sought ${\sf X}$ with respect to the powers of $z$:
\[
\aligned
Z^k(z,w)
=
\sum_{\alpha\in\N^n}\,z^\alpha\,Z^{k,\alpha}(w)
\ \ \ \ \
\text{\rm and}
\ \ \ \ \
W^l(z,w)
=
\sum_{\alpha\in\N^n}\,z^\alpha\,W^{l,\alpha}(w),
\endaligned
\]
where the $Z^{ k, \alpha}( w)$ and the $W^{ l, \alpha} ( w)$ are local
holomorphic functions. We will show that the identical vanishing of
the $d$ equations~\thetag{ \ref{initial-tangency}} in $\C \{
\underline{ z}, z, w\}$ is equivalent to a certain (in general
complicated) linear system of partial differential equations involving
the
$\frac{ \partial^{ \vert \gamma \vert} Z^{k, \alpha}}{
\partial w^\gamma}(w)$, 
the
$\frac{ \partial^{ \vert \gamma' \vert} Z^{ k', \alpha'}}{
\partial w^{ \gamma'}}(w)$, 
the
$\frac{ \partial^{ \vert \gamma'' \vert} W^{ l, \alpha'''}}{
\partial w^{ \gamma''}}(w)$
and the
$\frac{ \partial^{ \vert \gamma''' \vert} W^{ l', \alpha'''}}{
\partial w^{ \gamma'''}}(w)$.

Applying these expansions with respect to the
powers of $z$, we get:
\[
\footnotesize
\aligned
0
&
\equiv
\sum_{\alpha\in\N^n}\,\underline{z}^\alpha\,
\overline{W}^{j,\alpha}
\big(-w+\overline{\Xi}\big)
-
\sum_{k=1}^n\,\sum_{\alpha\in\N^n}\,
\underline{z}^\alpha\,
\overline{Z}^{k,\alpha}
\big(-w+\overline{\Xi}\big)\,
\frac{\partial\overline{\Xi}_j}{\partial\underline{z}_k}
(\underline{z},z,w)
+
\\
&
\ \ \ \ \
+
\sum_{\beta\in\N^n}\,
z^\beta\,W^{j,\beta}(w)
-
\sum_{k=1}^n\,\sum_{\beta\in\N^n}\,
z^\beta\,Z^{k,\beta}(w)\,
\frac{\partial\overline{\Xi}_j}{\partial z_k}
(\underline{z},z,w)
-
\sum_{l=1}^d\,\sum_{\beta\in\N^n}\,
z^\beta\,W^{l,\beta}(w)\,
\frac{\partial\overline{\Xi}_j}{\partial w_l}
(\underline{z},z,w)
\\
&
\ \ \ \ \ \ \ \ \ \ \ \ \ \ \ \ \ \ \ \ \ \ \ \ \ \ 
\ \ \ \ \ \ \ \ \ \ \ \ \ \ \ \ \ \ \ \ \ \ \ \ \ \
\ \ \ \ \ \ \ \ 
{\scriptstyle{(j\,=\,1\,\cdots\,d)}}.
\endaligned
\]
Since in these equations, $w$ is the argument both of all the 
$Z^{ k, \beta}$ and of all the $W^{ l,
\beta}$ appearing in the second line,
one should arrange that the same argument $w$ takes place inside the
functions $\overline{ W}^{ j, \alpha}$ and $\overline{ Z}^{ k,
\alpha}$ appearing in the first line. Thus, one is led, for an
arbitrary converging 
holomorphic power series $\overline{ A} = \overline{ A} ( w) =
\sum_{ \gamma\in \N^d}\, 
\frac{ \partial^{ \vert \gamma \vert} \overline{ A}}{
\partial w^\gamma} ( 0)\, w^\gamma$, to apply
the well known basic infinite Taylor series formula under the
following slightly
artificial form:
\[
\small
\aligned
\overline{A}\big(-w+\overline{\Xi}\big)
&
=
\overline{A}\big(w+(-2w+\overline{\Xi})\big)
\\
&
=
\sum_{\gamma\in\N^d}\,
\frac{\partial^{\vert\gamma\vert}\overline{A}}{
\partial w^\gamma}(w)\,
\frac{1}{\gamma!}\,
\big(-2w+\overline{\Xi}(\underline{z},z,w)\big)^\gamma.
\endaligned
\]
When one does this, one transforms the first
lines of the previous $d$ equations as follows:
\begin{equation}
\label{before-expansion-z-bar-z}
\footnotesize
\aligned
0
&
\equiv
\sum_{\alpha\in\N^n}\,
\sum_{\gamma\in\N^d}\,
\frac{1}{\gamma!}\,
\underline{z}^\alpha\,
\big(-2w+\overline{\Xi}(\underline{z},z,w)\big)^\gamma\,
\frac{\partial^{\vert\gamma\vert}\overline{W}^{j,\alpha}}{
\partial w^\gamma}(w)
-
\\
&
\ \ \ \ \
-
\sum_{k=1}^n\,\sum_{\alpha\in\N^n}\,\sum_{\gamma\in\N^d}\,
\frac{1}{\gamma!}\,
\underline{z}^\alpha\,
\big(-2w+\overline{\Xi}(\underline{z},z,w)\big)^\gamma\,
\frac{\partial^{\vert\gamma\vert}\overline{Z}^{k,\alpha}}{
\partial w^\gamma}(w)
+
\\
&
\ \ \ \ \
+
\sum_{\beta\in\N^n}\,
z^\beta\,W^{j,\beta}(w)
-
\sum_{k=1}^n\,\sum_{\beta\in\N^n}\,
z^\beta\,Z^{k,\beta}(w)\,
\frac{\partial\overline{\Xi}_j}{\partial z_k}
(\underline{z},z,w)
-
\sum_{l=1}^d\,\sum_{\beta\in\N^n}\,
z^\beta\,W^{l,\beta}(w)\,
\frac{\partial\overline{\Xi}_j}{\partial w_l}
(\underline{z},z,w)
\\
&
\ \ \ \ \ \ \ \ \ \ \ \ \ \ \ \ \ \ \ \ \ \ \ \ \ \ 
\ \ \ \ \ \ \ \ \ \ \ \ \ \ \ \ \ \ \ \ \ \ \ \ \ \
{\scriptstyle{(j\,=\,1\,\cdots\,d)}}.
\endaligned
\end{equation}
But still, we must expand and reorganize everything in terms of the powers
$\underline{ z}^\alpha\, z^\beta$ of $(\underline{ z}, z)$. At
first, we must do this for the multipowers 
$\big(-2w + \overline{ \Xi} ( \underline{ z}, z, w)\big)^\gamma$. 

\subsection{Expansion, reorganization and 
associated linear {\sc pde} system}
To begin with, 
let us denote the $(\underline{ z}, z)$-power series expansion
of $-2w_j + \overline{ \Xi}_j$ by:
\[
-2w_j
+
\overline{\Xi}_j(\underline{z},z,w)
=
\sum_{\alpha\in\N^n}\,\sum_{\beta\in\N^n}\,
\underline{z}^\alpha\,z^\beta\,
\overline{\Xi}_{j,\alpha,\beta}^\sim(w)
\ \ \ \ \ \ \ \ \ \ \ \ \ {\scriptstyle{(j\,=\,1\,\cdots\,d)}},
\]
with the understanding that the coefficients of the expansion of
$\overline{ \Xi}_j$ would be denoted plainly $\overline{ \Xi}_{ j,
\alpha, \beta} ( w)$, without $\sim$ sign. Hence, as $\overline{
\Xi}_j$ was assumed to be an ${\rm O} ( 2)$ at the beginning, we adopt
the convention that in this right-hand side, the $\overline{ \Xi}_{
j,\alpha, \beta}^\sim ( w)$ for $\alpha = \beta = 0$ comes not from
$\overline{ \Xi}_j$ itself, but from the supplementary 
first-order term $-2\, w_j$.

Thus, denoting:
\[
\gamma
=
(\gamma_1,\gamma_2,\dots,\gamma_d)
\in
\N^d,
\]
we may expand explicitly the exponentiated product under
consideration, and the intermediate, detailed computations read as follows:
\[
\footnotesize
\aligned
&
\prod_{j=1}^d\,
\big(-2\,w_j
+
\overline{\Xi}_j(\underline{z},z,w)\big)^{\gamma_j}
=
\\
&
=
\prod_{j=1}^d\,
\bigg(
\sum_{\alpha\in\N^n}\,\sum_{\beta\in\N^n}\,
\underline{z}^\alpha\,z^\beta\,
\overline{\Xi}_{j,\alpha,\beta}^\sim(w)
\bigg)^{\gamma_j}
\\
&
=
\prod_{j=1}^d\,
\bigg[
\sum_{\alpha\in\N^n}\,\sum_{\beta\in\N^n}\,
\underline{z}^\alpha\,z^\beta\,
\bigg(
\sum_{\alpha_1+\cdots+\alpha_{\gamma_j}=\alpha
\atop
\beta_1+\cdots+\beta_{\gamma_j}=\beta}\,
\overline{\Xi}_{j,\alpha_1,\beta_1}^\sim(w)
\cdots\,
\overline{\Xi}_{j,\alpha_{\gamma_j},\beta_{\gamma_j}}^\sim(w)
\bigg)
\bigg]
\endaligned
\]
\[
\footnotesize
\aligned
&
=
\sum_{\alpha\in\N^n}\,\sum_{\beta\in\N^n}\,
\underline{z}^\alpha\,z^\beta
\bigg[
\sum_{\alpha^1+\cdots+\alpha^d=\alpha
\atop
\beta_1+\cdots+\beta^d=\beta}\,
\sum_{\alpha_1^1+\cdots+\alpha_{\gamma_1}^1=\alpha^1
\atop
\beta_1^1+\cdots+\beta_{\gamma_1}^1=\beta^1}\,
\cdots\,
\sum_{\alpha_1^d+\cdots+\alpha_{\gamma_d}^d=\alpha^d
\atop
\beta_1^d+\cdots+\beta_{\gamma_d}^d=\beta_d}
\\
&
\ \ \ \ \ \ \ \ \ \ \ \ \ \ \ \ \ \ \ \ \ \ \ \ \ \ \ \ \ \ \ \ \ \ \
\overline{\Xi}_{1,\alpha_1^1,\beta_1^1}^\sim(w)
\cdots
\overline{\Xi}_{1,\alpha_{\gamma_1}^1,\beta_{\gamma_1}^1}^\sim(w)
\cdots\cdots\,
\overline{\Xi}_{d,\alpha_1^d,\beta_1^d}^\sim(w)
\cdots
\overline{\Xi}_{d,\alpha_{\gamma_d}^d,\beta_{\gamma_d}^d}^\sim(w)
\bigg]
\\
&
=:
\sum_{\alpha\in\N^n}\,\sum_{\beta\in\N^n}\,
\underline{z}^\alpha\,z^\beta\,
\mathcal{A}_{\alpha,\beta,\gamma}
\Big(
\big\{
\overline{\Xi}_{\widehat{j},\widehat{\alpha},\widehat{\beta}}^\sim(w)
\big\}_{\widehat{j}\in\N,\widehat{\alpha}\in\N^n,\widehat{\beta}\in\N^n}
\Big),
\endaligned
\]
where we introduce a collection of certain polynomial functions
$\mathcal{ A}_{ \alpha, \beta, \gamma}$ of all the $\overline{ \Xi}_{
\widehat{ j}, \widehat{ \alpha}, \widehat{ \beta}}^\sim ( w)$ that
appear naturally in the large brackets of the penultimate equality,
namely where we set:
\[
\small
\aligned
\mathcal{A}_{\alpha,\beta,\gamma}
\Big(
\big\{
\overline{\Xi}_{\widehat{j},\widehat{\alpha},\widehat{\beta}}^\sim(w)
&
\big\}_{\widehat{j}\in\N,\widehat{\alpha}\in\N^n,\widehat{\beta}\in\N^n}
\Big)
:=
\sum_{\alpha^1+\cdots+\alpha^d=\alpha
\atop
\beta_1+\cdots+\beta^d=\beta}\,
\sum_{\alpha_1^1+\cdots+\alpha_{\gamma_1}^1=\alpha^1
\atop
\beta_1^1+\cdots+\beta_{\gamma_1}^1=\beta^1}\,
\cdots\,
\sum_{\alpha_1^d+\cdots+\alpha_{\gamma_d}^d=\alpha^d
\atop
\beta_1^d+\cdots+\beta_{\gamma_d}^d=\beta_d}
\\
&
\overline{\Xi}_{1,\alpha_1^1,\beta_1^1}^\sim(w)
\cdots
\overline{\Xi}_{1,\alpha_{\gamma_1}^1,\beta_{\gamma_1}^1}^\sim(w)
\cdots\cdots\,
\overline{\Xi}_{d,\alpha_1^d,\beta_1^d}^\sim(w)
\cdots
\overline{\Xi}_{d,\alpha_{\gamma_d}^d,\beta_{\gamma_d}^d}^\sim(w).
\endaligned
\] 
At present, coming back to the $d$ equations~\thetag{
\ref{before-expansion-z-bar-z}} we left momentarily untouched, we see
that in them, five sums are extant and we now want to expand and to
reorganize properly each one of these terms as a $(\underline{ z},
z)$-power series $\sum_{ \alpha \in \N^n} \, \sum_{ \beta \in \N^n}\,
\underline{ z}^\alpha \, z^\beta \big( \cdots \big)$. For the sum
in~\thetag{
\ref{before-expansion-z-bar-z}}, 
we therefore compute, changing in advance the index $\alpha$ to
$\alpha'$:
\begin{equation}
\label{I}
\footnotesize
\aligned
&
\sum_{\alpha'\in\N^n}\,\sum_{\gamma\in\N^d}\,
\frac{1}{\gamma!}\,\underline{z}^{\alpha'}\,
\frac{\partial^{\vert\gamma\vert}\overline{W}^{j,\alpha'}}{
\partial w^\gamma}(w)
\big(-2w+\overline{\Xi}(\underline{z},z,w)\big)^\gamma\,
=
\\
&
=
\sum_{\alpha'\in\N^n}\,\sum_{\beta\in\N^n}\,
\frac{1}{\gamma!}\,\underline{z}^{\alpha'}\,
\frac{\partial^{\vert\gamma\vert}
\overline{W}^{j,\alpha'}}{\partial w^\gamma}(w)
\sum_{\alpha''\in\N^n}\,\sum_{\beta\in\N^n}\,
\underline{z}^{\alpha''}\,z^\beta\,
\mathcal{A}_{\alpha'',\beta,\gamma}
\Big(\big\{
\overline{\Xi}_{\widehat{j},\widehat{\alpha},\widehat{\beta}}^\sim(w)
\big\}\Big)
\\
&
=
\sum_{\alpha\in\N^n}\,\sum_{\beta\in\N^n}\,
\underline{z}^\alpha\,z^\beta
\bigg[
\sum_{\gamma\in\N^d}\,
\sum_{\alpha=\alpha'+\alpha''}\,
\frac{1}{\gamma!}\,
\mathcal{A}_{\alpha'',\beta,\gamma}
\Big(\big\{
\overline{\Xi}_{\widehat{j},\widehat{\alpha},\widehat{\beta}}^\sim(w)
\big\}\Big)
\cdot
\frac{\partial^{\vert\gamma\vert}\overline{W}^{j,\alpha'}}{
\partial w^\gamma}(w)
\bigg].
\endaligned
\end{equation}
The computations for the second sum in~\thetag{
\ref{before-expansion-z-bar-z}} are essentially exactly the same:
\begin{equation}
\label{II}
\footnotesize
\aligned
&
\ \ \ \ \
-\sum_{k=1}^n\,\sum_{\alpha'\in\N^n}\,\sum_{\gamma\in\N^d}\,
\frac{1}{\gamma!}\,\underline{z}^{\alpha'}\,
\frac{\partial^{\vert\gamma\vert}\overline{Z}^{k,\alpha'}}{
\partial w^\gamma}(w)
\big(-2w+\overline{\Xi}(\underline{z},z,w)\big)^\gamma\,
\\
&
=
-\sum_{k=1}^n\,\sum_{\alpha'\in\N^n}\,\sum_{\gamma\in\N^d}\,
\frac{1}{\gamma!}\,\underline{z}^{\alpha'}\,
\frac{\partial^{\vert\gamma\vert}\overline{Z}^{k,\alpha'}}{
\partial w^\gamma}(w)\,
\sum_{\alpha''\in\N^n}\,\sum_{\beta\in\N^n}\,
\underline{z}^{\alpha''}\,z^\beta\,
\mathcal{A}_{\alpha'',\beta,\gamma}
\Big(\big\{
\overline{\Xi}_{\widehat{j},\widehat{\alpha},\widehat{\beta}}^\sim(w)
\big\}\Big)
\\
&
=
\sum_{\alpha\in\N^n}\,\sum_{\beta\in\N^n}\,
\underline{z}^\alpha\,z^\beta
\bigg[
-\sum_{k=1}^n\,\sum_{\gamma\in\N^d}\,
\sum_{\alpha'+\alpha''=\alpha}\,
\frac{1}{\gamma!}\,
\mathcal{A}_{\alpha'',\beta,\gamma}
\Big(\big\{
\overline{\Xi}_{\widehat{j},\widehat{\alpha},\widehat{\beta}}^\sim(w)
\big\}\Big)
\cdot
\frac{\partial^{\vert\gamma\vert}\overline{Z}^{k,\alpha'}}{
\partial w^\gamma}(w)
\bigg].
\endaligned
\end{equation}
The third sum in~\thetag{
\ref{before-expansion-z-bar-z}} is already almost well written, for we
indeed have, if we denote by ${\bf 0} = (0, \dots, 0) 
\in \N^n$ the zero-multiindex:
\begin{equation}
\label{III}
\small
\aligned
\sum_{\beta\in\N^n}\,
z^\beta\,W^{j,\beta}(w)
=
\sum_{\alpha\in\N^n}\,\sum_{\beta\in\N^n}\,
\underline{z}^\alpha\,z^\beta\,
\big[
\delta_\alpha^{\bf 0}
\cdot
W^{j,\beta}(w)
\big],
\endaligned
\end{equation}
where $\delta_{\sf a}^{\sf b} = 0$ if ${\sf a} \neq {\sf b}$ 
and $1$ if ${\sf a} = {\sf b}$.
To transform the fourth sum in~\thetag{
\ref{before-expansion-z-bar-z}}, we must at first compute, for each $k =
1, \dots, n$ (and for each $j = 1, \dots, d$), the first-order partial
derivatives $\frac{ \partial \overline{ \Xi}_j}{ \partial z_k}$, which
gives, if we denote simply by ${\bf 1}_k$ the multiindex $(0, \dots,
1, \dots, 0)$ of $\N^n$ with $1$ at the $k$-th place and zero
elsewhere:
\[
\small
\aligned
\frac{\partial\overline{\Xi}_j}{\partial z_k}
(\underline{z},z,w)
&
=
\sum_{\alpha\in\N^n}\,\sum_{\beta\in\N^n\atop \beta_k\geqslant 1}\,
\underline{z}^\alpha\,\beta_k\,z^{\beta-{\bf 1}_k}\,
\overline{\Xi}_{j,\alpha,\beta}^\sim(w)
\\
&
=
\sum_{\alpha\in\N^n}\,\sum_{\beta\in\N^n}\,
\underline{z}^\alpha\,z^\beta\,
(\beta_k+1)\,
\overline{\Xi}_{j,\alpha,\beta+{\bf 1}_k}^\sim(w).
\endaligned
\]
Thanks to this, the fourth sum in~\thetag{
\ref{before-expansion-z-bar-z}} may be reorganized as wanted:
\begin{equation}
\label{IV}
\footnotesize
\aligned
&
-\sum_{k=1}^n\,\sum_{\beta'\in\N^n}\,
z^{\beta'}\,Z^{k,\beta'}(w)\,
\frac{\partial\overline{\Xi}_j}{\partial z_k}
(\underline{z},z,w)
=
\\
&
=
-\sum_{k=1}^n\,\sum_{\beta'\in\N^n}\,
z^{\beta'}\,Z^{k,\beta'}(w)\,
\sum_{\alpha\in\N^n}\,\sum_{\beta''\in\N^n}\,
\underline{z}^\alpha\,z^{\beta''}\,
(1+\beta_k'')\,
\overline{\Xi}_{j,\alpha,\beta''+{\bf 1}_k}^\sim(w)
\\
&
=
\sum_{\alpha\in\N^n}\,\sum_{\beta\in\N^n}\,
\underline{z}^\alpha\,z^\beta
\bigg[
-\sum_{k=1}^n\,
\sum_{\beta'+\beta''=\beta}\,
(\beta_k''+1)\,
\overline{\Xi}_{j,\alpha,\beta''+{\bf 1}_k}^\sim(w)
\cdot
Z^{k,\beta'}(w)
\bigg].
\endaligned
\end{equation}
Lastly, in order to transform the fifth sum in~\thetag{
\ref{before-expansion-z-bar-z}}, we must at first
compute, for each $l = 1, \dots, d$ (and for each $j = 1, \dots, d$),
the first-order partial derivatives $\frac{ \partial \overline{
\Xi}_j}{ \partial w_l}$, and to this aim, we start by rewriting:
\[
\overline{\Xi}_j(\underline{z},z,w)
=
2w_j
+
\sum_{\alpha\in\N^n}\,\sum_{\beta\in\N^n}\,
\underline{z}^\alpha\,z^\beta\,
\overline{\Xi}_{j,\alpha,\beta}^\sim(w),
\]
whence it immediately follows:
\[
\small
\aligned
\frac{\partial\overline{\Xi}_j}{\partial w_l}
(\underline{z},z,w)
=
2\,\delta_j^l
+
\sum_{\alpha\in\N^n}\,\sum_{\beta\in\N^n}\,
\underline{z}^\alpha\,z^\beta\,
\big(\partial\overline{\Xi}_{j,\alpha,\beta}^\sim(w)\big/\partial w_l\big).
\endaligned
\]
Thanks to this, the fifth sum in~\thetag{
\ref{before-expansion-z-bar-z}}, too, may be reorganized appropriately:
\begin{equation}
\label{V}
\footnotesize
\aligned
&
-\sum_{l=1}^d\,\sum_{\beta'\in\N^n}\,
z^{\beta'}\,W^{l,\beta'}(w)\,
\frac{\partial\overline{\Xi}_j}{\partial w_l}
(\underline{z},z,w)
=
\\
&
=
-\sum_{l=1}^d\,\sum_{\beta'\in\N^n}\,
z^{\beta'}\,W^{l,\beta'}(w)\,
\bigg[
2\,\delta_j^l
+
\sum_{\alpha\in\N^n}\,\sum_{\beta''\in\N^n}\,
\underline{z}^\alpha\,z^{\beta''}\,
\big(
\partial\overline{\Xi}_{j,\alpha,\beta''}^\sim(w)
\big/
\partial w_l
\big)
\bigg]
\\
&
=
\sum_{\alpha\in\N^n}\,\sum_{\beta\in\N^n}\,
\underline{z}^\alpha\,z^\beta\,
\bigg[
-2\,\delta_\alpha^{\bf 0}
\cdot
W^{j,\beta}(w)
-
\sum_{l=1}^d\,
\sum_{\beta'+\beta''=\beta}\,
\big(
\partial\overline{\Xi}_{j,\alpha,\beta''}^\sim(w)
\big/
\partial w_l
\big)
\cdot
W^{l,\beta'}(w)
\bigg].
\endaligned
\end{equation}
Summing up these five reorganized sums appearing in~\thetag{
\ref{before-expansion-z-bar-z}} as a double sum $\sum_{ \alpha} \,
\sum_\beta\, \underline{ z}^\alpha\, z^\beta \, \big( {\sf coeff}_{ j,
\alpha, \beta} \big)$, and equating to zero all the obtained
coefficients~\thetag{ \ref{I}}, \thetag{ \ref{II}}, \thetag{
\ref{III}}, \thetag{ \ref{IV}} and~\thetag{ \ref{V}}, we deduce the
following fundamental statement.

\begin{Theorem}
Let $M$ be a generic real analytic CR-submanifold of $\C^{ n+d}$
having positive codimension $d \geqslant 1$ and positive CR dimension
$n \geqslant 1$ which is represented, in local holomorphic coordinates
$(z, w) = ( z_1, \dots, z_n, w_1, \dots, w_d)$ by $d$ complex defining
equations of the shape:
\[
\overline{w}_j+w_j
=
\overline{\Xi}_j(\overline{z},z,w)
\ \ \ \ \ \ \ \ \ \ \ \ \ {\scriptstyle{(j\,=\,1\,\cdots\,d)}},
\]
denote by $(\underline{ z}, \underline{ w})$ the extrinsic
complexifications of the antiholomorphic variables $(\overline{ z},
\overline{ w})$ and introduce the power series expansion with respect
to the variables $( \underline{ z}, z)$:
\[
-2w_j
+
\overline{\Xi}_j(\underline{z},z,w)
=:
\sum_{\alpha\in\N^n}\,\sum_{\beta\in\N^n}\,
\underline{z}^\alpha\,z^\beta\,
\overline{\Xi}_{j,\alpha,\beta}^\sim(w)
\ \ \ \ \ \ \ \ \ \ \ \ \ {\scriptstyle{(j\,=\,1\,\cdots\,d)}}.
\]
For every multiindex $\alpha \in \N^n$, 
every multiindex $\beta \in \N^n$ and every multiindex
$\gamma \in \N^d$, introduce also the explicit universal polynomial:
\[
\small
\aligned
\mathcal{A}_{\alpha,\beta,\gamma}
\Big(
\big\{
\overline{\Xi}_{\widehat{j},\widehat{\alpha},\widehat{\beta}}^\sim(w)
&
\big\}_{\widehat{j}\in\N,\widehat{\alpha}\in\N^n,\widehat{\beta}\in\N^n}
\Big)
:=
\sum_{\alpha^1+\cdots+\alpha^d=\alpha
\atop
\beta_1+\cdots+\beta^d=\beta}\,
\sum_{\alpha_1^1+\cdots+\alpha_{\gamma_1}^1=\alpha^1
\atop
\beta_1^1+\cdots+\beta_{\gamma_1}^1=\beta^1}\,
\cdots\,
\sum_{\alpha_1^d+\cdots+\alpha_{\gamma_d}^d=\alpha^d
\atop
\beta_1^d+\cdots+\beta_{\gamma_d}^d=\beta_d}
\\
&
\overline{\Xi}_{1,\alpha_1^1,\beta_1^1}^\sim(w)
\cdots
\overline{\Xi}_{1,\alpha_{\gamma_1}^1,\beta_{\gamma_1}^1}^\sim(w)
\cdots\cdots\,
\overline{\Xi}_{d,\alpha_1^d,\beta_1^d}^\sim(w)
\cdots
\overline{\Xi}_{d,\alpha_{\gamma_d}^d,\beta_{\gamma_d}^d}^\sim(w).
\endaligned
\] 
Then a general holomorphic vector field:
\[
{\sf X}
=
\sum_{k=1}^n\,Z^k(z,w)\,
\frac{\partial}{\partial z_k}
+
\sum_{l=1}^d\,W^l(z,w)\,
\frac{\partial}{\partial w_l}
\]
is an infinitesimal CR-automorphism of $M$ belonging to $\mathfrak{
hol} ( M)$, namely it has the property that $\overline{\sf X} 
+ {\sf
X}$ is tangent to $M$ {\em if and only if}, for every $j = 1, \dots,
d$, for every $\alpha \in \N^n$ and for every $\beta \in \N^n$, the
following linear holomorphic partial differential equation:
\[
\boxed{
\footnotesize
\aligned
0
&
\equiv
\sum_{\gamma\in\N^d}\,
\sum_{\alpha=\alpha'+\alpha''}\,
\frac{1}{\gamma!}\,
\mathcal{A}_{\alpha'',\beta,\gamma}
\Big(
\big\{
\overline{\Xi}_{\widehat{j},\widehat{\alpha},\widehat{\beta}}^\sim(w)
\big\}_{\widehat{j}\in\N,\widehat{\alpha}\in\N^n,\widehat{\beta}\in\N^n}
\Big)
\cdot
\frac{\partial^{\vert\gamma\vert}\overline{W}^{j,\alpha'}}{
\partial w^\gamma}(w)
-
\\
&
-\sum_{k=1}^n\,\sum_{\gamma\in\N^d}\,
\sum_{\alpha'+\alpha''=\alpha}\,
\frac{1}{\gamma!}\,
\mathcal{A}_{\alpha'',\beta,\gamma}
\Big(\big\{
\overline{\Xi}_{\widehat{j},\widehat{\alpha},\widehat{\beta}}^\sim(w)
\big\}\Big)
\cdot
\frac{\partial^{\vert\gamma\vert}\overline{Z}^{k,\alpha'}}{
\partial w^\gamma}(w)
+
\\
&
+
\delta_\alpha^{\bf 0}
\cdot
W^{j,\beta}(w)
-
\\
&
-\sum_{k=1}^n\,
\sum_{\beta'+\beta''=\beta}\,
(\beta_k''+1)\,
\overline{\Xi}_{j,\alpha,\beta''+{\bf 1}_k}^\sim(w)
\cdot
Z^{k,\beta'}(w)
-
\\
&
-2\,\delta_\alpha^{\bf 0}
\cdot
W^{j,\beta}(w)
-
\sum_{l=1}^d\,
\sum_{\beta'+\beta''=\beta}\,
\big(
\partial\overline{\Xi}_{j,\alpha,\beta''}^\sim(w)
\big/
\partial w_l
\big)
\cdot
W^{l,\beta'}(w)
\endaligned}
\]
which is {\em linear} with respect to the partial derivatives:
\[
\footnotesize
\aligned
\frac{\partial^{\vert\gamma\vert}Z^{k,\alpha}}{
\partial w^\gamma}(w),
\ \ \ \ \ \
\frac{\partial^{\vert\gamma'\vert}Z^{k',\alpha'}}{
\partial w^{\gamma'}}(w),
\ \ \ \ \ \
\frac{\partial^{\vert\gamma''\vert}W^{l,\alpha''}}{
\partial w^{\gamma''}}(w),
\ \ \ \ \ \
\frac{\partial^{\vert\gamma'''\vert}W^{l',\alpha'''}}{
\partial w^{\gamma'''}}(w)
\endaligned
\]
is satisfied identically in $\C\{ w\}$
by the four families of functions:
\[
Z^{k,\alpha}(w),
\ \ \ \ \ \
Z^{k',\alpha'}(w),
\ \ \ \ \ \
W^{l,\alpha'''}(w),
\ \ \ \ \ \
W^{l',\alpha'''}(w).
\]
depending only upon the $d$ holomorphic variables
$(w_1, \dots, w_d)$. 
\end{Theorem}

Then the resolution of this linear system of holomorphic partial
differential equations (having nonconstant coefficients in general) is
often delicate when dealing with concrete, specific functions
$\overline{ \Xi}_j ( \underline{ z}, z, w)$. Of course, starting from
an $M$ of equation $v = \varphi( x, y, u)$ with $T_0 M = \{ {\rm Im}\,
w = 0\}$, instead of $T_0 M = \{ {\rm Re}\, w = 0\}$, the same process
of extracting linear partial differential equations providing (after
resolution) access to all ${\sf X} \in \mathfrak{ hol} ( M)$ may be
conducted quite similarly, the only difference being that an
$i$-factor comes regularly into play, for the complex defining
equations of $M$ must then be thought to be of the general form:
\[
w_j
=
\overline{w}_j
+i\,\Xi_j(z,\overline{z},\overline{w})
\ \ \ \ \ \ \ \ \ \ \ \ \ {\scriptstyle{(j\,=\,1\,\cdots\,d)}},
\]
because each $w_j - \overline{ w}_j = v_j$ is purely real, or because
a reality condition like~\thetag{ \ref{reality-Theta}} must hold true.
The presence of the $i$-factor is especially visible in the case where
the $M$ under consideration is of the particular (and quite
convenient) form, sometimes called {\sl rigid} in the literature,
where the right-hand side functions $\Xi_j$ are {\em completely
independent} of the variable $w \in \C^d$, since in this case if one
writes:
\[
w_j
=
\overline{w}_j
+
i\,\Xi_j(z,\overline{z})
\ \ \ \ \ \ \ \ \ \ \ \ \ {\scriptstyle{(j\,=\,1\,\cdots\,d)}},
\]
it is clear that each right-hand side function $\Xi_j ( z, \overline{
z})$ must be purely real, namely must satisfy:
\[
\Xi_j(z,\underline{z})
\equiv
\overline{\Xi}_j(\underline{z},z)
\ \ \ \ \ \ \ \ \ \ \ \ \ {\scriptstyle{(j\,=\,1\,\cdots\,d)}}
\]
identically in $\C\{ z, \underline{ z} \}$. A concrete
example is on.

\section{CR symmetries of the Heisenberg Sphere $\mathbb
H^3\subseteq\mathbb C^2$}
\label{Infinitesimal-CR-section}

\HEAD{CR symmetries of the Heisenberg Sphere $\mathbb
H^3\subseteq\mathbb C^2$}{
Mansour Aghasi, Joël Merker, and Masoud Sabzevari}

\subsection{Infinitesimal CR automorphisms of $\mathbb H^3$}
\label{Infinitesimal-CR-subsection}
We now consider the Heisenberg sphere $\mathbb H^3$ in $\mathbb
C^2$, equipped with coordinates $(z, w)$, of equation:
\[
0 
= 
w-\overline{w}-2i\,z\overline{z}.
\]
A local $(1, 0)$ vector field defined in a neighborhood of the origin:
\[
{\sf X} 
= 
Z(z,w)\,\frac{\partial}{\partial z} 
+ 
W(z,w)\,\frac{\partial}{\partial w}
\]
having {\em holomorphic} coefficients $Z ( z, w)$ and $W ( z, w)$ is
an {\sl infinitesimal CR automorphism} of $M$ if and only $X +
\underline{ X}$ is tangent to the extrinsic complexification $M^{
e_c}$, that is to say, if and only if the following equation:
\[
0 
\equiv 
\big[W-2i\underline{z}\,Z-\overline{W}-2iz\,\overline{Z}
\big]_{w=\underline{w}+2i\,z\underline{z}}
\]
holds identically in $\mathbb C \{ z, \underline{ z}, \underline{ w}
\}$, that is to say if again, and only if:
\begin{equation}
\label{fund-eq} 0 \equiv W(z,\,\underline{w}+2i\,z\underline{z}) -
2i\underline{z}\,Z(z,\,\underline{w}+2i\,z\underline{z}) -
\overline{W}(\underline{z},\,\underline{w}) -
2iz\,\overline{Z}(\underline{z},\underline{w}).
\end{equation}
Since the two coefficients $Z$ and $W$ of $L$ are analytic, we may
expand them with respect to the powers of $z$:
\[
Z(z,w) 
= 
\sum_{k\in\mathbb N}\,z^k\,Z_k(w) 
\ \ \ \ \ \ \ \ \ \ \ 
\text{\rm and} 
\ \ \ \ \ \ \ \ \ \
\ W(z,w) = \sum_{k\in\mathbb N}\,z^k\,W_k(z,w),
\]
and the fundamental 
equation~\thetag{ \ref{fund-eq}} just written becomes:
\[
\aligned 0 & \equiv \sum_{k\in\mathbb N}\, z^k\,W_k(\underline{w}+2i\,z\underline{z}) -
2i\underline{z}\, \sum_{k\in\mathbb N}\, z^k\,Z_k(\underline{w}+2i\,z\underline{z}) \,-
\\
& \ \ \ \ \ -\, \sum_{k\in\mathbb N}\,
\underline{z}^k\,\overline{W}^k(\underline{w}) - 2iz\,
\sum_{k\in\mathbb N}\, \underline{z}^k\,\overline{Z}_k(\underline{w}).
\endaligned
\]
Furthermore, if $A = A ( \underline{ w}) = \sum_{ l \in \mathbb N} \,
A_{ \underline{ w}^l}(0)\, \frac{ 1}{ l!}\, \underline{ w}^l$ is a
function holomorphic with respect to $\underline{ w}$ near the origin,
where $A_{\underline{ w}}$, $A_{ \underline{ w}^2}$, \dots, $A_{
\underline{ w}^l}$ denote (partial) derivatives, we may yet expand:
\begin{equation}
\label{expansion-A} 
A(\underline{w}+2i\,z\underline{z}) 
= 
\sum_{l\in\mathbb N}\,
A_{\underline{w}^l}(\underline{w})\, 
(2i\,z\underline{z})^l\, \frac{1}{l!},
\end{equation}
and here, this gives us:
\[
\aligned 0 & \equiv 
\sum_{k\in\mathbb N}\,\sum_{l\in\mathbb N}\, \bigg( z^k\,
(2i\,z\underline{z})^l\, \frac{1}{l!}\, W_{k,\underline{w}^l}(\underline{w}) -
2i\,\underline{z}\,
z^k\, (2i\,z\underline{z})^l\, 
\frac{1}{l!}\, Z_{k,\underline{w}^l}(\underline{w})
\bigg) -
\\
& \ \ \ \ \ 
- 
\sum_{k\in\mathbb N}\, \underline{z}^k\, \Big(
\overline{W}_k(\underline{w}) +
2i\,z\,\overline{Z}_k(\underline{w}) \Big).
\endaligned
\]
In this equation, the coefficients of the monomials $\underline{ z}^k$
for all $k \geqslant 2$ and the coefficients of the monomials $z \,
\underline{ z}^{ k'}$ for all $k' \geqslant 3$ must vanish, and this
simply yields:
\[
0 \equiv \overline{W}_k(\underline{w}) \ \ \ \text{\rm for all} \ \ k\geqslant 2 \ \
\ \ \ \ \
\text{\rm and} \ \ \ \ \ \ \ 0 \equiv \overline{Z}_{k'}(\underline{w}) \ \ \
\text{\rm for all} \ \
k'\geqslant 3.
\]
Consequently, the two coefficients $Z$ and $W$ of our infinitesimal CR automorphism
greatly
simplify, and they receive the (truncated) form:
\[
\left[ \aligned Z(z,w) & = Z_0(w)+z\,Z_1(w)+z^2\,Z_2(w)
\\
W(z,w) & = W_0(w)+z\,W_1(w).
\endaligned\right.
\]
After this simplification, the fundamental equation~\thetag{ \ref{fund-eq}} becomes:
\[
\aligned 0 & \equiv W_0(\underline{w}+2i\,z\underline{z}) +
z\,W_1(\underline{w}+2i\,z\underline{z}) -
\\
& \ \ \ \ \ - 2i\underline{z}\, Z_0(\underline{w}+2i\,z\underline{z}) -
2iz\underline{z}\,
Z_1(\underline{w}+2i\,z\underline{z}) - 2iz^2\underline{z}\,
Z_2(\underline{w}+2i\,z\underline{z})
-
\\
& \ \ \ \ \ - \overline{W}_0(\underline{w}) -
\underline{z}\,\overline{W}_1(\underline{w}) -
2iz\,\overline{Z}_0(\underline{w}) - 2iz\underline{z}\,\overline{Z}(\underline{w}) -
2iz\underline{z}^2\,\overline{Z}_2(\underline{w}).
\endaligned
\]
We now expand all the series $A ( \underline{ w} + 2i\, z \underline{ z})$ appearing
in the first
two lines, using~\thetag{ \ref{expansion-A}}:
\[
\aligned 0 & \equiv W_0(\underline{w}) +
2iz\underline{z}\,W_{0,\underline{w}}(\underline{w}) -
4z^2\underline{z}^2\, {\textstyle{\frac{1}{2!}}}\,
W_{0,\underline{w}^2}(\underline{w}) -
8iz^3\underline{z}^3\, {\textstyle{\frac{1}{3!}}}\,
W_{0,\underline{w}^3}(\underline{w}) +\cdots
\\
& \ \ \ \ \ + z\,W_1(\underline{w}) +
2iz^2\underline{z}\,W_{1,\underline{w}}(\underline{w}) -
4z^3\underline{z}^2\, {\textstyle{\frac{1}{2!}}}\,
W_{1,\underline{w}^2}(\underline{w}) -\cdots
\\
& \ \ \ \ \ -\,2i\underline{z}\, Z_0(\underline{w}) + 4\,z\underline{z}^2\,
Z_{0,\underline{w}}(\underline{w}) + 8iz^2\underline{z}^3\,
{\textstyle{\frac{1}{2!}}}\,
Z_{0,\underline{w}^2}(\underline{w}) +\cdots
\\
& \ \ \ \ \ -\,2i\underline{z}z\, Z_1(\underline{w}) + 4z^2\underline{z}^2\,
Z_{1,\underline{w}}(\underline{w}) + 8iz^3\underline{z}^3 {\textstyle{\frac{1}{2!}}}\,
Z_{1,\underline{w}^2}(\underline{w}) +\cdots
\\
& \ \ \ \ \ -\,2i\underline{z}z^2 Z_2(\underline{w}) +
4z^3\underline{z}^2\,Z_{2,\underline{w}}
+\cdots
\\
& \ \ \ \ \ -\,\overline{W}_0(\underline{w}) - \underline{z}\,
\overline{W}_1(\underline{w}) -
2iz\, \overline{Z}_0(\underline{w}) - 2iz\underline{z}\,
\overline{Z}_1(\underline{w}) -
2iz\underline{z}^2\, \overline{Z}_2(\underline{w}).
\endaligned
\]
Now, we extract the coefficients of the monomials $z^\mu \underline{
z}^\nu$ for small values of $\mu$ and $\nu$, and these coefficients
must all vanish identically in $\mathbb C \{ \underline{ w} \}$. What
is left out in the cdots will not be useful to us.

First of all, for $(\mu, \nu)$ equal to $(0, 0)$ and to $(1, 0)$, we
get two equations:
\begin{eqnarray}
\label{0-0} && 0 \equiv W_0(\underline{w}) - \overline{W}_0(\underline{w})
\\
\label{1-0} && 0 \equiv W_1(\underline{w}) - 2i\,\overline{Z}_0(\underline{w}),
\end{eqnarray}
holding identically in $\mathbb C \{ \underline{ w} \}$, while for
$(\mu, \nu) = (0, 1)$, we get $0 \equiv - \, \overline{ W}_1
(\underline{w}) - 2i\, Z_0 ( \underline{ w})$ which is fully
equivalent to~\thetag{ \ref{1-0}}, after conjugation and replacement
of the variable $w$ by $\underline{ w}$ (a power series $\varphi (
\underline{ w} )$ is identically zero if and only if $\varphi ( w)$ is
identically zero). Next, for $(\mu, \nu) = (2, 0)$, nothing comes,
while for $(\mu, \nu) = (1, 1)$, we obtain:
\begin{equation}
\label{1-1} 0 \equiv 2i\,W_{0,\underline{w}}(\underline{w}) - 2i\,Z_1(\underline{w}) -
2i\,\overline{Z}_1(\underline{w}).
\end{equation}
Next, for $(2, 1)$ and for $(1, 2)$ we obtain:
\begin{equation}
\label{2-1} 0 \equiv 2i\,W_{1,\underline{w}}(\underline{w}) - 2i\,Z_2(\underline{w})
\end{equation}
and: $0 \equiv 4\, Z_{ 0, w}(\underline{w}) - 2i\, \overline{
Z}_2(\underline{w})$, but this second equation visibly follows from
the ones already obtained, hence will be disregarded. Next, for $(2,
2)$, for $(2, 3)$, for $(3, 2)$ and for $(3, 3)$, we obtain:
\begin{eqnarray}
\label{2-2} && 0 \equiv -4\, {\textstyle{\frac{1}{2!}}}\,
W_{0,\underline{w}^2}(\underline{w}) +
4\,Z_{1,\underline{w}}(\underline{w})
\\
\label{2-3} && 0 \equiv 8i\, {\textstyle{\frac{1}{2!}}}\,
Z_{0,\underline{w}^2}(\underline{w})
\\
\label{3-2} && 0 \equiv -4\, {\textstyle{\frac{1}{2!}}}\,
W_{1,\underline{w}^2}(\underline{w}) +
4\,Z_{2,\underline{w}}(\underline{w})
\\
\label{3-3} && 0 \equiv -8i\, {\textstyle{\frac{1}{3!}}}\,
W_{0,\underline{w}^3}(\underline{w}) -
8i\, {\textstyle{\frac{1}{2!}}}\, Z_{1,\underline{w}^2}(\underline{w}).
\end{eqnarray}
Clearly, \thetag{ \ref{2-3}} yields that $Z_0$ is affine:
\begin{equation}
Z_0(w) = {\sf z}_{0,0} + {\sf z}_{0,1}\,w,
\end{equation}
where ${\sf z}_{ 0, 0} = {\sf x}_{ 0, 0} + i\, {\sf y}_{ 0, 0}$ and
${\sf z}_{ 0, 1} = {\sf x}_{ 0, 1} + i \, {\sf y}_{ 0, 1}$ are {\em
complex} constants in $\mathbb C$. From~\thetag{ \ref{1-0}}, it then
follows immediately that:
\begin{equation}
\label{w-1} W_1(w) = 2i\,\overline{\sf z}_{0,0} + 2i\,\overline{\sf z}_{0,1}\,w.
\end{equation}
Next, differentiating~\thetag{ \ref{2-2}} once with respect to $\underline{ w}$ and
comparing
to~\thetag{ \ref{3-3}}, we get:
\[
0 \equiv W_{0,\underline{w}^3}(\underline{w}) \ \ \ \ \ \ \ \ \text{\rm and} \ \ \ \
\ \ \ \ 0
\equiv Z_{1,\underline{w}}(\underline{w}).
\]
It follows firstly that $W_0$ is quadratic:
\begin{equation}
W_0(w) = {\sf u}_{0,0}+{\sf u}_{0,1}\,w+{\sf u}_{0,2}\,w^2,
\end{equation}
but taking account of~\thetag{ \ref{0-0}}, we see that the three
appearing coefficients ${\sf u}_{ 0, 0}$, ${\sf u}_{ 0, 1}$, ${\sf
u}_{ 0, 2}$ must even all be {\em real}. Secondly, it follows that
$Z_1 (w) = {\sf z}_{ 1, 0} + {\sf z}_{ 1, 1} \, w$ is affine, but
moreover, taking in addition account of~\thetag{ \ref{1-1}} and
of~\thetag{ \ref{2-2}}, we see that:
\begin{equation}
Z_1(w) 
= 
{\textstyle{\frac{1}{2}}}\,{\sf u}_{0,1}
+
i\,{\sf y}_{1,0} + {\sf u}_{0,2}\,w.
\end{equation}
Finally, \thetag{ \ref{2-1}} and~\thetag{ \ref{w-1}} 
give that $Z_2 ( w)$ is constant:
\[
Z_2(w) = 2\,{\sf y}_{0,1} + 2i\,{\sf x}_{0,1}.
\]

\subsection{Solution}
\label{Solution}
The eight real constants found in this way:
\[
{\sf x}_{0,0},\ \ \ {\sf y}_{0,0},\ \ \ {\sf x}_{0,1},\ \ \ {\sf y}_{0,1},\ \ \ {\sf
u}_{0,0},\ \ \
{\sf u}_{0,1},\ \ \ {\sf u}_{0,2},\ \ \ {\sf y}_{1,0}
\]
give us eight $\mathbb R$-linearly independent infinitesimal
automorphisms of the Heisenberg sphere, when one sets one of these
constants equal to $1$, while the 7 remaining constants are set equal
to $0$:
\[
\aligned &
\partial_z+2iz\,\partial_w
\\
& i\,\partial_z+2z\,\partial_w
\\
& (w+2iz^2)\,\partial_z+2izw\,\partial_w
\\
& (iw+2z^2)\,\partial_z+2zw\,\partial_w
\\
&
\partial_w
\\
& {\textstyle{\frac{1}{2}}}z\,\partial_z+w\,\partial_w
\\
& zw\,\partial_z+w^2\,\partial_w
\\
& iz\,\partial_z.
\endaligned
\]
By straightforward computations, one verifies that indeed the real
part of ${\sf X} + \underline{\sf X}$, where ${\sf X}$ is any of these
eight holomorphic vector fields, is tangent to $M^{ e_c}$; for
instance, for the third vector field, we get:
\[
\aligned & 4z^2\underline{z} - 4z\underline{z}^2 - 2i\underline{z}w -
2iz\underline{w} + 2izw +
2i\underline{z}\underline{w} \\
& = (w-\underline{w}-2iz\underline{z})[2iz-2i\underline{z}],
\endaligned
\]
which identifies to the equation of $M^{ e_c}$ multiplied by a factor,
hence vanishes on $M^{ e_c}$.

\subsection{ Homogeneities and graded structure}
\label{Homogeneities}
The anisotropic real dilation $(z, w) \longmapsto (c \, z, \, c^2 \,
w)$ with $c \in \mathbb R$ visibly stabilizes $\mathbb H^3$, hence it
is natural to ascribe homogeneity $1$ to the variable $z$ and
homogeneity $2$ to the variable $w$. Accordingly, $\partial_z$ and
$\partial_w$ have homogeneity $-1$ and $-2$, respectively, and a
holomorphic field like $zw\, \partial_w$, for instance, has
homogeneity $1 + 2 - 2 = 1$. One may thus list the eight generators
found above according to their homogeneities, which take the values
$-2$, $-1$, $0$, $1$ and $2$ and this conducts us to represent:
\[
\mathfrak{hol}(\mathbb H^3) 
= 
\mathfrak{h}_{-2} \oplus \mathfrak{h}_{-1} \oplus
\mathfrak{h}_0
\oplus \mathfrak{h}_1 \oplus \mathfrak{h}_2
\]
as a direct sum of five components of dimensions $1$, $2$, $2$, $2$,
$1$ defined by:
\[
\aligned \mathfrak{h}_{-2} & = \mathbb R\,{\sf T}, 
\ \ \ \ \ 
\mathfrak{h}_{-1} 
=
\mathbb R\,{\sf
H}_1 \oplus \mathbb R\,{\sf H}_2
\\
\mathfrak{h}_0 & = \mathbb R\,{\sf D} \oplus \mathbb R\,{\sf R},
\\
\mathfrak{h}_1 & = \mathbb R\,{\sf I}_1 \oplus \mathbb R\,{\sf I}_2, \ \ \ \ \
\mathfrak{h}_2 
=
\mathbb R\,{\sf J}.
\endaligned
\]
where:
\[
\aligned 
& 
\mathfrak{h}_{-2}\colon \left\{ {\sf T} 
:=
\partial_w
\right. 
\ \ \ \ \ \ \ \ \ \ \ 
\mathfrak{h}_{-1}\colon \left\{ \aligned {\sf H}_1 &
:=
\partial_z
+ 2iz\,\partial_w
\\
{\sf H}_2 
& 
:= 
i\,\partial_z + 2z\,\partial_w
\endaligned
\right.
\\
& \mathfrak{h}_0\colon 
\left\{ 
\aligned 
{\sf D} & 
:= 
z\,\partial_z + 2w\,\partial_w
\\
{\sf R} & := iz\,\partial_z
\\
\endaligned\right.
\\
& \mathfrak{h}_1\colon \left\{ \aligned {\sf I}_1 & := (w+2iz^2)\,\partial_z +
2izw\,\partial_w
\\
{\sf I}_2 & := (iw+2z^2)\,\partial_z + 2zw\,\partial_w
\endaligned\right.
\ \ \ \ \ \ \ \ \ \ \ \mathfrak{h}_2\colon \left\{ {\sf J} := zw\,\partial_z +
w^2\,\partial_w.
\right.
\endaligned
\]
Here with $t \in \mathbb R$, the flow $(z, w) \mapsto (z, \, w + t)$
of ${\sf T}$ is \underline{\sf t}ransversal to the complex tangent
bundle $T^cM$, spanned by ${\rm Re}\, \mathcal{ L}$ and by ${\rm Im}\,
\mathcal{ L}$, where:
\[
\mathcal{L} 
:=
\partial_z
+ 2i\overline{z}\,\partial_w;
\]
the flows of ${\sf H}_1$ and ${\sf H}_2$, namely $(z + t, \, w+ 2izt +
it^2) $ and $(z + it , \, w + 2zt + it^2)$ are somewhat \underline{\sf
h}orizontal; the flow of ${\sf D}$ is just the \underline{\sf
d}ilation $(e^tz, \, e^{ 2t} \, w)$; the flow of ${\sf R}$ is just the
imaginary \underline{\sf r}otation of the $z$-coordinate $(e^{ it} z,
\, w)$. On the other hand, it is known since Poincar\'e~\cite{
Poincare1907} that any holomorphic automorphism of the Heisenberg
sphere fixing the origin is a fractional linear transformation of the
general form:
\[
(z,\,w) 
\longmapsto 
\Big({\textstyle{\frac{c(z+aw)}{1-2i\overline{a}z 
-
(r+ia\overline{a})w}}},
\ \ 
{\textstyle{\frac{\rho w}{1-2i\overline{a}z-(r+ia\overline{a})w}}}\Big),
\]
where $c \in \mathbb C$, $a \in \mathbb C$, $r \in \mathbb R$ and
$\rho \in \mathbb R$. Such a general expression may be recovered by
concatenating the eight flows, after a change of parameters. However,
as understood originally by Lie himself (\cite{ EngelLie, Merker7}),
except in some specific situations, it is essentially useless to
explicit the finite equations of a local transformation group, because
the infinitesimal description shows better the structures.

With the convention that $\mathfrak{ h}_k = \{ 0\}$ for either $k
\leqslant -3$ or $k \geqslant 3$, one may then verify the property
that:
\[
\big[\mathfrak{h}_{k_1},\,\mathfrak{h}_{k_2}\big] \subset \mathfrak{h}_{k_1+k_2}
\]
for any two $k_1 \leqslant k_2 \in \mathbb Z$, and more precisely,
this fact follows by inspecting the full commutator table between
these eight generators of $\mathfrak{ hol}( \mathbb H^3)$:
\medskip
\begin{center}
\begin{tabular} [t] { l | l l l l l l l l }
& ${\sf T}$ & ${\sf H}_1$ & ${\sf H}_2$ & ${\sf D}$ & ${\sf R}$ & ${\sf I}_1$ &
${\sf I}_2$ & ${\sf
J}$
\\
\hline ${\sf T}$ & $0$ & $0$ & $0$ & $2\,{\sf T}$ & $0$ & ${\sf H}_1$ & ${\sf H}_2$
& ${\sf D}$
\\
${\sf H}_1$ & $*$ & $0$ & $4\,{\sf T}$ & ${\sf H}_1$ & ${\sf H}_2$ & $6\,{\sf R}$ &
$2\,{\sf D}$ &
${\sf I}_1$
\\
${\sf H}_2$ & $*$ & $*$ & $0$ & ${\sf H}_2$ & $-{\sf H}_1$ & $-2\,{\sf D}$ &
$6\,{\sf R}$ & ${\sf
I}_2$
\\
${\sf D}$ & $*$ & $*$ & $*$ & $0$ & $0$ & ${\sf I}_1$ & ${\sf I}_2$ & $2\,{\sf J}$
\\
${\sf R}$ & $*$ & $*$ & $*$ & $*$ & $0$ & $-{\sf I}_2$ & ${\sf I}_1$ & $0$
\\
${\sf I}_1$ & $*$ & $*$ & $*$ & $*$ & $*$ & $0$ & $4\,{\sf J}$ & $0$
\\
${\sf I}_2$ & $*$ & $*$ & $*$ & $*$ & $*$ & $*$ & $0$ & $0$
\\
${\sf J}$ & $*$ & $*$ & $*$ & $*$ & $*$ & $*$ & $*$ & $0$ %
\end{tabular}
\end{center}

Clearly also, the isotropy algebra of the origin is just the
nonnegative part of the sum:
\[
\mathfrak{hol}(\mathbb H^3,0) 
= 
\mathfrak{h}_0 \oplus \mathfrak{h}_1 
\oplus
\mathfrak{h}_2.
\]

\section{Tanaka Prolongation}
\label{Tanaka-Prolongation}

\HEAD{Tanaka Prolongation}{
Mansour Aghasi, Joël Merker, and Masoud Sabzevari}

\subsection{The prolongation procedure in the CR context}
\label{Prolongation-definition}
Consider a finite-dimensional graded real Lie algebra indexed by
negative integers:
\[
\mathfrak{g}_- 
= 
\mathfrak{g}_{-\mu} \oplus \cdots \oplus \mathfrak{g}_{-2} \oplus
\mathfrak{g}_{-1},
\]
satisfying $[ \mathfrak{ g}_{ - l_1}, \, \mathfrak{ g}_{ -l_2} 
] \subset \mathfrak{ g}_{ -l_1 - l_2}$ with the convention that
$\mathfrak{ g}_k = 0$ for $k \leqslant - \mu - 1$. Following~\cite{
Tanaka1}, $\mathfrak{ g}_-$ will be said to be {\sl of $\mu$-th
kind}. Assume that there is a complex structure $J \colon \mathfrak{
g}_{ -1} \to \mathfrak{ g}_{ -1}$ such that $J^2 = - {\rm Id}$, whence
$\mathfrak{ g}_{ -1}$ is even-dimensional and bears a natural
structure of a complex vector space. Tanaka's prolongation of
$\mathfrak{ g}_-$ is an algebraic procedure which generates a certain
larger graded Lie algebra:
\[
\mathfrak{g} = \mathfrak{g}_- \oplus \mathfrak{g}_0 \oplus \mathfrak{g}_1 \oplus
\mathfrak{g}_2
\oplus \cdots
\]
in the following way.

By definition, the order-zero component $\mathfrak{ g}_0$ consists of
all linear endomorphisms ${\sf d} \colon \mathfrak{ g}_- \to
\mathfrak{ g}_-$ which preserve gradation: ${\sf d} ( \mathfrak{ g}_k)
\subset \mathfrak{ g}_k$, which respect the complex structure: ${\sf
d} ( J \, {\sf x}) = J {\sf d} ( {\sf x})$ for all ${\sf x} \in
\mathfrak{ g}_{ -1}$ and which are {\em derivations}, namely satisfy
${\sf d} ( [ {\sf y}, \, {\sf z}]) = [ {\sf d} ( {\sf x}), \, {\sf y}]
+ [{\sf x}, \, {\sf d} ( {\sf y}) ]$ for every ${\sf y}, {\sf z} \in
\mathfrak{ g}_-$. Then the bracket between a ${\sf d} \in \mathfrak{
g}_0$ and an ${\sf x} \in \mathfrak{ g}_-$ is simply defined by $[
{\sf d}, \, {\sf x}] := {\sf d} ( {\sf x})$, while the bracket between
{\em two} elements ${\sf d}', {\sf d}'' \in \mathfrak{ g}_0$ is
defined to be the commutator ${\sf d}'\circ {\sf d}'' - {\sf d}''
\circ {\sf d}'$ between endomorphisms. One checks at once that Jacobi
relations hold, hence $\mathfrak{ g}_- \oplus \mathfrak{ g}_0$ becomes
a true Lie algebra.

By contrast, for any $l \geqslant 1$, no constraint with respect to
$J$ is required. Assuming that the components $\mathfrak{ g}_{ l'}$
are already constructed for any $l' \leqslant l - 1$, the $l$-th
component $\mathfrak{ g}_l$ of the prolongation consists of
$l$-shifted graded linear morphisms $\mathfrak{ g}_- \to \mathfrak{
g}_- \oplus \mathfrak{ g}_0 \oplus \mathfrak{ g}_1 \oplus \cdots
\oplus \mathfrak{ g}_{ l-1}$ that are derivations, namely:
\begin{equation}
\label{definition-prolongation} \mathfrak{g}_l = \Big\{ {\sf d} \in
\bigoplus_{k\leqslant-1}\, {\rm
Lin}(\mathfrak{g}_k,\,\mathfrak{g}_{k+l}) \colon {\sf d}([{\sf y},\,{\sf z}]) =
[{\sf d}({\sf
y}),\,{\sf z}] + [{\sf y},\,{\sf d}({\sf z})], \ \ \ \ \ \forall\, {\sf y},\,{\sf
z}\in\mathfrak{g}_- \Big\}.
\end{equation}
Now, for ${\sf d} \in \mathfrak{ g}_k$ and ${\sf e} \in \mathfrak{
g}_l$, by induction on the integer $k + l \geqslant 0$, one defines
the bracket $[ {\sf d}, \, {\sf e} ] \in \mathfrak{ g}_{ k
+ l} \otimes \mathfrak{ g}_-^*$ by:
\begin{equation}
\label{d-e-x} [{\sf d},\,{\sf e}]({\sf x}) = \big[[{\sf d},\,{\sf x}],\,{\sf e}\big]
+ \big[{\sf
d},\,[{\sf e},\,{\sf x}]\big] \ \ \ \ \ \ \ \ \ \ \text{\rm for}\ \ {\sf
x}\in\mathfrak{g}_-.
\end{equation}
One notes that, for $k = l = 0$, this definition coincides with the
above one for $[ \mathfrak{ g}_0, \, \mathfrak{ g}_0]$. It follows by
induction (\cite{ Tanaka1, Yamaguchi}) that $[ {\sf d}, \, {\sf e} ]
\in \mathfrak{ g}_{ k+l}$ and that with this bracket, the sum
$\mathfrak{ g}_- \bigoplus_{ k \geqslant 1}\, \mathfrak{ g}_k$ becomes
a graded Lie algebra, because the general Jacobi identity:
\[
0 
= 
\big[[{\sf d},\,{\sf e}],\,{\sf f}\big] 
+ 
\big[[{\sf f},\,{\sf d}],\,{\sf
e}\big] + \big[[{\sf
e},\,{\sf f}],\,{\sf d}\big]
\]
for ${\sf d} \in \mathfrak{ g}_k$, ${\sf e} \in \mathfrak{ g}_l$ and
${\sf f} \in \mathfrak{ g}_m$ follows by definition when one of $k$,
$l$, $m$ is negative, and can be shown by induction on the integer $k
+ l + m \geqslant 0$ when all of $k$, $l$, $m$ are nonnegative.

\subsection{The Heisenberg algebra}
\label{Heisenberg-algebra}
The symbol Lie algebra $\mathfrak{ g}_- := \mathfrak{ g}_{-2} \oplus
\mathfrak{ g}_{ -1}$ associated to any Levi nondegenerate CR manifold
$M^3 \subset \mathbb C^2$ equipped with the distribution $T^cM$ is
three-dimensional, with $\mathfrak{ g}_{-2} = \mathbb R\, {\sf x}_1$,
with $\mathfrak{ g}_{-1} = \mathbb R \, {\sf x}_2 \oplus \mathbb R
{\sf x}_3$ with ${\sf x}_3 = J({\sf x}_2)$ and with only nonzero Lie
bracket $[{\sf x}_2,\, {\sf x}_3] =4 \, {\sf x}_1$. If one disregards
$J$, such a $\mathfrak{ g}_-$ is the unique irreducible
three-dimensional nilpotent Lie algebra, denoted $\mathfrak{ n}_3^1$
in~\cite{Goze}. Now, what is the Tanaka prolongation of $\mathfrak{
g}_-$?

By definition, an element of $\mathfrak{ g}_0$ is a derivation ${\sf
d} \in (\mathfrak{ g}_{ -2} \otimes \mathfrak{ g}_{ -2}^*) \oplus
(\mathfrak{ g}_{ -1} \otimes \mathfrak{ g}_{ -1}^*)$, hence it writes:
\[
{\sf d}({\sf x}_1)=k{\sf x}_1, 
\ \ \ \ \ \ \ 
{\sf d}({\sf x}_2)=r_1{\sf x}_2+r_2{\sf
x}_3, 
\ \ \ \ \ \ \ 
{\sf d}({\sf x}_3)=r_3{\sf x}_2+r_4{\sf x}_3,
\]
for some five real, unknown constants. But because ${\sf d}$ preserves
the complex structure $J$ on $\mathfrak{ g}_{-1}$, one also has ${\sf
d} (J{\sf x}_2) = J{\sf d} ({\sf x}_2)$, i.e. $r_3{\sf x}_2+r_4{\sf
x}_3=r_1{\sf x}_3-r_2{\sf x}_2$, that is to say: $r_1=r_4$ and
$r_2=-r_3$. Moreover, applying the derivation property of ${\sf d}$,
one must have:
\[
\aligned 4k\,{\sf x}_1 
= 
4{\sf d}({\sf x}_1) = {\sf d}([{\sf x}_2,{\sf x}_3]) 
& 
=
[{\sf d}({\sf
x}_2),{\sf x}_3] + [{\sf x}_2,{\sf d}({\sf x}_3)]
\\
& = [r_1{\sf x}_2+r_2{\sf x}_3,{\sf x}_3] 
+ 
[{\sf x}_2,r_3{\sf x}_2+r_4{\sf x}_3]
\\
& = 4r_1{\sf x}_1+4r_4{\sf x}_1,
\endaligned
\]
which yields that $k=r_1+r_4$. These three linear equations solve as
$r_3 = -r_2$, $r_4 = r_1$ and $k = 2 r_1$ with free $r_1$ and
$r_2$. It follows that $\mathfrak{ g}_0$ is two-dimensional and
generated over $\mathbb R$ by two derivations (corresponding to the
two choices: $r_1 = -1$, $r_2 0$ and $r_1 = 0$, $r_2 = -1$) that we
will denote ${\sf x}_4$ and ${\sf x}_5$ and which are defined
by:
\[
\aligned {\sf x}_4 & \colon \ \ \ \ \ {\sf x}_1 
\mapsto -2\,{\sf x}_1, \ \ \ \ \ \ \
{\sf x}_2
\mapsto -{\sf x}_2, \ \ \ \ \ \ \ {\sf x}_3 \mapsto -{\sf x}_3,
\\
{\sf x}_5 & \colon \ \ \ \ \ {\sf x}_1 \mapsto 0, \ \ \ \ \ \ \ \ \ \ \ \ \ \ \ {\sf
x}_2 \mapsto
-{\sf x}_3, \ \ \ \ \ \ \ {\sf x}_3 \mapsto {\sf x}_2.
\endaligned
\]
Then the commutator ${\sf x}_4 \circ {\sf x}_5 - {\sf x}_5 \circ {\sf
x}_4 = 0$ vanishes, and at this stage, the Lie brackets between the
obtained ${\sf x}_k$ read as follows, if listed by increasing
homogeneity:
\[
\aligned 
\boxed{{\scriptstyle{-3}}} & \colon\,\, \big\{\,\, [{\sf x}_1,\,{\sf
x}_2]=0, \ \ \ \ \ \
\ [{\sf x}_1,\,{\sf x}_3]=0,
\\
\boxed{{\scriptstyle{-2}}} & \colon\,\, 
\big\{\,\, [{\sf x}_2,\,{\sf x}_3]=4\,{\sf
x}_1, \ \ \ 
[{\sf x}_1,\,{\sf x}_4]=2\,{\sf x}_1, 
\ \ \ \ \ \ \ [{\sf x}_1,\,{\sf x}_5]=0,
\\
\boxed{{\scriptstyle{-1}}} & \colon\,\, 
\big\{\,\, [{\sf x}_2,\,{\sf x}_4]={\sf
x}_2, \ \ \ \ \ \
[{\sf x}_3,\,{\sf x}_4]
=
{\sf x}_3, \ \ \ \ \ \ \ \ \
[{\sf x}_2,\,{\sf x}_5]={\sf x}_3,
\ \ \ \ \ \ \
[{\sf x}_3,\,{\sf x}_5]=-{\sf x}_2
\\
\boxed{{\scriptstyle{0}}} & 
\colon\,\, \big\{\,\, [{\sf x}_4,\,{\sf x}_5]=0.
\endaligned
\]
Here, we see that $\mathfrak{ g}_{ -2} \oplus \mathfrak{ g}_{ -1}
\oplus \mathfrak{ g}_0$ is a Lie algebra in itself. Moreover, we
observe that $\mathfrak{ g}_{ -2} \oplus \mathfrak{ g}_{ -1} \oplus
\mathfrak{ g}_0$ is isomorphic to the isotropy subalgebra $\mathfrak{
h}_{ -2} \oplus \mathfrak{ h}_{ -1} \oplus \mathfrak{ h}_0$ of the
Heisenberg sphere through the simple map:
\[
\aligned {\sf T}\mapsto{\sf x}_1, \ \ \ \ \ \ \ {\sf H}_1\mapsto{\sf x}_2, \ \ \ \ \
\ \ {\sf
H}_2\mapsto{\sf x}_3, \ \ \ \ \ \ \ {\sf D}\mapsto{\sf x}_4, \ \ \ \ \ \ \ {\sf
R}\mapsto{\sf x}_5.
\endaligned
\]

Next, we compute $\mathfrak{ g}_1$. An element ${\sf d}$ of
$\mathfrak{ g}_1$ belongs to $(\mathfrak{ g}_{-1} \otimes \mathfrak{
g}_{ -2}^*) \oplus (\mathfrak{ g}_0 \otimes \mathfrak{ g}_{ -1}^*)$,
hence it writes:
\[
{\sf d}({\sf x}_1)
=
k{\sf x}_2+l{\sf x}_3, 
\ \ \ \ \ \ \ 
{\sf d}({\sf x}_2)=m{\sf
x}_4+n{\sf x}_5, \
\ \ \ \ \ \ {\sf d}({\sf x}_3)=p{\sf x}_4+q{\sf x}_5,
\]
for some six real, unknown constants. But the condition that ${\sf
d}$ be a derivation gives us exactly three constraints, firstly:
\[
0 
= 
{\sf d}([{\sf x}_1,\,{\sf x}_2]) 
=
[k{\sf x}_2+l{\sf x}_3,\,{\sf x}_2] + [{\sf
x}_2,\,m{\sf
x}_4+n{\sf x}_5] = -4l{\sf x}_1+2m{\sf x}_1
\]
that is to say: $0 = -2l + m$; secondly:
\[
0 
= 
{\sf d}([{\sf x}_1,\,{\sf x}_3]) 
= 
[k{\sf x}_2+l{\sf x}_3,\,{\sf x}_3] + [{\sf
x}_1,\,p{\sf
x}_4+q{\sf x}_5] = 4k{\sf x}_1+2p{\sf x}_1,
\]
that is to say: $0 = 2k+p$; thirdly and lastly:
\[
\aligned
4k{\sf x}_2+4l{\sf x}_3 
= 
4{\sf d}({\sf x}_1) 
= 
{\sf d}([{\sf x}_1,\,{\sf x}_2]) 
&
=
[m{\sf
x}_4+n{\sf x}_5,\,{\sf x}_3] 
+ 
[{\sf x}_2,\,p{\sf x}_4+q{\sf x}_5] 
=
\\
&
=
-m{\sf x}_3+n{\sf x}_2+p{\sf
x}_2+q{\sf x}_3,
\endaligned
\]
that is to say: $4k = n+p$ and $4l = -m+q$. These four linear
equations solve as $m = 2l$, $p -2k$, $n = 6k$, $q = 6l$ with free $k$
and $l$. It follows that $\mathfrak{ g}_1$ is two-dimensional and
generated over $\mathbb R$ by two derivations (corresponding to the
two choices: $k = -1$, $l = 0$ and $k = 0$, $l = -1$):
\[
\aligned {\sf x}_6 & \colon \ \ \ \ \ 
{\sf x}_1 \mapsto -{\sf x}_2, \ \ \ \ \ \ \
{\sf x}_2 \mapsto
-6\,{\sf x}_5, \ \ \ \ \ \ \ 
{\sf x}_3 \mapsto 2\,{\sf x}_4,
\\
{\sf x}_7 & \colon \ \ \ \ \ 
{\sf x}_1 \mapsto -{\sf x}_3, \ \ \ \ \ \ \ {\sf x}_2
\mapsto -2\,{\sf
x}_4, \ \ \ \ \ \ \ {\sf x}_3 \mapsto -6\,{\sf x}_5.
\endaligned
\]
We still need to know the brackets structures $[ \mathfrak{ g}_1, \,
\mathfrak{ g}_0 ]$ and $[ \mathfrak{ g}_1, \, \mathfrak{ g}_1]$. At
this stage in fact, we can only determine $[ \mathfrak{ g}_1, \,
\mathfrak{ g}_0 ]$. By definition, with ${\sf d} \in \mathfrak{ g}_1$
and ${\sf e} \in \mathfrak{ g}_0$, the bracket $[ {\sf d}, \, {\sf e}]
\in \mathfrak{ g}_1 \subset ( \mathfrak{ g}_{ -1} \otimes \mathfrak{
g}_{ -2}^*) \oplus ( \mathfrak{ g}_0 \otimes \mathfrak{ g}_{ -1}^*)$
is determined by his action on the three vectors ${\sf x}_1$, ${\sf
x}_2$, ${\sf x}_3$ generating $\mathfrak{ g}_-$ through the
formula~\thetag{ \ref{d-e-x}}, hence we compute three times at once:
\[
\aligned{} [{\sf x}_4,\,{\sf x}_6] \Big( \underset{{\sf
x}_3}{\overset{\rule[-3pt]{0pt}{11pt}{\sf
x}_1}{{\scriptstyle{{\sf x}_2}}}} \Big) := &\, \Big[ \Big[{\sf x}_4,\, \underset{{\sf
x}_3}{\overset{\rule[-3pt]{0pt}{11pt}{\sf x}_1}{{\scriptstyle{{\sf x}_2}}}} \Big],\,
{\sf x}_6
\Big] + \Big[ {\sf x}_4,\, \Big[ {\sf x}_6,\, \underset{{\sf
x}_3}{\overset{\rule[-3pt]{0pt}{11pt}{\sf x}_1}{{\scriptstyle{{\sf x}_2}}}} \Big]
\Big] = \Big[
\underset{-{\sf x}_3}{\overset{\rule[-3pt]{0pt}{11pt}{-2\,\sf
x}_1}{{\scriptstyle{-{\sf x}_2}}}},\,
{\sf x}_6 \Big] + \Big[ {\sf x}_4,\, \underset{-2\,{\sf
x}_4}{\overset{\rule[-3pt]{0pt}{11pt}-{\sf
x}_2}{{\scriptstyle{-6\,{\sf x}_5}}}} \Big]
\\
= &\, \Big( \underset{2\,{\sf x}_4}{\overset{\rule[-3pt]{0pt}{11pt}-2\,{\sf
x}_2}{{\scriptstyle{-6\,{\sf x}_5}}}} \Big) + \Big(
\underset{0}{\overset{\rule[-3pt]{0pt}{11pt}{\sf x}_2}{{\scriptstyle{0}}}} \Big) =
\Big(
\underset{2\,{\sf x}_4}{\overset{\rule[-3pt]{0pt}{11pt}-{\sf
x}_2}{{\scriptstyle{-6\,{\sf x}_5}}}}
\Big) = \text{\rm map}\,{\sf x}_6 \Big( \underset{{\sf
x}_3}{\overset{\rule[-3pt]{0pt}{11pt}{\sf
x}_1}{{\scriptstyle{{\sf x}_2}}}} \Big).
\endaligned
\]
and we recognize ${\sf x}_6$, as a linear map, whence $[ {\sf x}_4, \,
{\sf x}_6 ] = {\sf x}_6$. In a similar way, one may compute the three
remaining brackets. In summary, one obtains the following
supplementary brackets, listed by increasing homogeneity:
\[
\aligned
\boxed{{\scriptstyle{-1}}} & \colon\,\, \big\{\,\, [{\sf x}_1,\,{\sf
x}_6]={\sf x}_2, \ \
\ \ \ \ \ [{\sf x}_1,\,{\sf x}_7]={\sf x}_3,
\\
\boxed{{\scriptstyle{0}}} & \colon\,\, 
\big\{\,\, [{\sf x}_2,\,{\sf x}_6]=6\,{\sf
x}_5, 
\ \ \ \ 
[{\sf x}_3,\,{\sf x}_6]=-2\,{\sf x}_4, 
\ \ \ \ 
[{\sf x}_2,\,{\sf x}_7]=2\,{\sf x}_4, 
\ \ \ \ 
[{\sf x}_3,\,{\sf x}_7]=6\,{\sf x}_5,
\\
\boxed{{\scriptstyle{1}}} & \colon\,\, \big\{\,\, [{\sf x}_4,\,{\sf x}_6]={\sf x}_6,
\ \ \ \ \ \ \
[{\sf x}_4,\,{\sf x}_7]={\sf x}_7, \ \ \ \ \ \ \ \ \ \,
[{\sf x}_5,\,{\sf x}_6]=-{\sf x}_7,
\ \ \ \,
[{\sf x}_5,\,{\sf x}_7]={\sf x}_6.
\endaligned
\]

Now, we are in a position to compute $\mathfrak{ g}_2$. An element
${\sf d}$ of $\mathfrak{ g}_2$ belongs to $( \mathfrak{ g}_0 \otimes
\mathfrak{ g}_{ -2}^*) \oplus ( \mathfrak{ g}_{ -1} \oplus \mathfrak{
g}_{ -1}^*)$, hence it writes:
\[
{\sf d}({\sf x}_1)=k{\sf x}_4+l{\sf x}_5, 
\ \ \ \ \ \ \ {\sf d}({\sf x}_2)=m{\sf
x}_6+n{\sf x}_7, \
\ \ \ \ \ \ {\sf d}({\sf x}_3)=p{\sf x}_6+q{\sf x}_7,
\]
for some six real, unknown constants. Again, the condition that ${\sf
d}$ be a derivation gives exactly three constraints, firstly:
\[
0 
= 
{\sf d}([{\sf x}_1,\,{\sf x}_2]) 
= 
[k{\sf x}_4+l{\sf x}_5,\,{\sf x}_2] + [{\sf
x}_1,\,m{\sf
x}_6+n{\sf x}_7] 
= 
-k{\sf x}_2-l{\sf x}_3+m{\sf x}_2+n{\sf x}_3,
\]
that is to say: $0 = -k + m$ and $0 = - l + n$; secondly:
\[
0 
= 
{\sf d}([{\sf x}_1,\,{\sf x}_3]) 
= 
[k{\sf x}_4+l{\sf x}_5,\,{\sf x}_3] 
+
[{\sf x}_1,\,p{\sf
x}_6+q{\sf x}_7] 
= 
-k{\sf x}_3+l{\sf x}_2+p{\sf x}_2+q{\sf x}_3,
\]
that is to say: $0 = l + p$ and $0 = - k + q$; thirdly and lastly:
\[
\aligned
4k{\sf x}_4+4l{\sf x}_5 
= 
4{\sf d}({\sf x}_1) = {\sf d}([{\sf x}_2,\,{\sf x}_3]) 
&
=
[m{\sf
x}_6+n{\sf x}_7,\,{\sf x}_3] 
+ 
[{\sf x}_2,\,p{\sf x}_6+q{\sf x}_7] 
\\
&
= 
2m{\sf
x}_4-6n{\sf x}_5+6p{\sf
x}_5+2q{\sf x}_4,
\endaligned
\]
that is to say: $4k = 2m+2q$ and $4l = -n + p$. These four linear
equations solve, up to a dilation factor, as $0 = l = n = p$ and $m =
q = k = - 1$, whence it follows that $\mathfrak{ g}_2$ is
one-dimensional and generated over $\mathbb R$ by the single
derivation:
\[
{\sf x}_8 \colon \ \ \ \ \ {\sf x}_1\mapsto-{\sf x}_4, \ \ \ \ \ \ \ {\sf
x}_2\mapsto-{\sf x}_6, \
\ \ \ \ \ \ {\sf x}_3\mapsto-{\sf x}_7.
\]
We still need to know the bracket structures $[ \mathfrak{ g}_0, \,
\mathfrak{ g}_2]$ and $[ \mathfrak{ g}_1, \, \mathfrak{ g}_1]$. After
a few computations using the already known brackets, one obtains the
supplementary brackets:
\[
\aligned 
\boxed{
{\scriptstyle{0}}} 
& 
\colon\,\, \big\{\,\, [{\sf x}_1,\,{\sf
x}_8]={\sf x}_4,
\\
\boxed{{\scriptstyle{1}}} 
& 
\colon\,\, \big\{\,\, [{\sf x}_2,\,{\sf x}_8]={\sf x}_6,
\ \ \ \ \ \ \ \ \ \,
[{\sf x}_3,\,{\sf x}_8]={\sf x}_7,
\\
\boxed{{\scriptstyle{2}}} & \colon\,\, 
\big\{\,\, [{\sf x}_4,\,{\sf x}_8]=2\,{\sf
x}_8, \ \ \ \ \ \
\ [{\sf x}_5,\,{\sf x}_8]=0, 
\ \ \ \ \ \ \ 
[{\sf x}_6,\,{\sf x}_7]=4\,{\sf x}_8,
\endaligned
\]
Finally, the prolongation stops, for one may easily verify that
$\mathfrak{ g}_3 = \{ 0\}$, while it is known that $\mathfrak{ g}_k \{
0\}$ for some $k \geqslant 0$ implies $\mathfrak{ g}_l = \{ 0 \}$ for
every $l \geqslant k$ (\cite{Yamaguchi}, p.~433). Also, $0 = [
\mathfrak{ g}_1, \, \mathfrak{ g}_2] = [ \mathfrak{ g}_2, \,
\mathfrak{ g}_2]$, which ends up the process, all brackets being
computed and:
\[
\mathfrak{g} 
:= 
\mathfrak{g}_{-2} 
\oplus 
\mathfrak{g}_{-1} \oplus \mathfrak{g}_0 \oplus
\mathfrak{g}_1 \oplus \mathfrak{g}_2
\]
being a graded $1 + 2 + 2 + 2 + 1 =$ eight-dimensional Lie algebra.

\subsection{Fundamental isomorphism}
Now by inspecting all brackets just obtained, we observe that this
$J$-compatible Tanaka prolongation $\mathfrak{ g}$ (of the above
three-dimensional Heisenberg nilpotent Lie algebra $\mathfrak{
n}_3^1$) {\em isomorphically coincides} with the Lie algebra
$\mathfrak{ hol} (\mathbb{ H}^3)$ of CR-automorphisms of the
Heisenberg sphere through the plain identifications:
\[
{\sf x}_1\leftrightarrow{\sf T},\ \ \
{\sf x}_2\leftrightarrow{\sf H}_1,\ \ \
{\sf x}_3\leftrightarrow{\sf H}_2,\ \ \
{\sf x}_4\leftrightarrow{\sf D},\ \ \
{\sf x}_5\leftrightarrow{\sf R},\ \ \
{\sf x}_6\leftrightarrow{\sf I}_1,\ \ \
{\sf x}_7\leftrightarrow{\sf I}_2,\ \ \
{\sf x}_8\leftrightarrow{\sf J}.
\]
In fact, this coincidence comes from Tanaka's general (\cite{
Yamaguchi}) theorem that the prolongation $\mathfrak{ g}$ can
naturally be identified with the Lie algebra of all $J$-compatible
infinitesimal automorphisms of the unique connected simply connected
three-dimensional Lie group with the Lie algebra $[{\sf x}_2, \, {\sf
x}_3] = 4\, {\sf x}_1$, ${\sf x}_3 = J ({\sf x}_2)$. For more
clarity, we will employ from now on the letters ${\sf t}$, ${\sf
h}_1$, ${\sf h}_2$, ${\sf d}$, ${\sf r}$, ${\sf i}_1$, ${\sf i}_2$,
${\sf j}$ instead of ${\sf x}_1$, ${\sf x}_2$, ${\sf x}_3$, ${\sf
x}_4$, ${\sf x}_5$, ${\sf x}_6$, ${\sf x}_7$, ${\sf x}_8$ as
generators of $\mathfrak{ g}$. For later use, we draw up the full
commutator table between these eight generators of the abstract Lie
algebra $\mathfrak{ g}$:

\medskip
\begin{center}
\begin{tabular} [t] { l | l l l l l l l l }
\label{structure-g}
& ${\sf t}$ & ${\sf h}_1$ & ${\sf h}_2$ & ${\sf d}$ & ${\sf r}$ &
${\sf i}_1$ & ${\sf i}_2$ & ${\sf j}$
\\
\hline ${\sf t}$ & $0$ & $0$ & $0$ & $2\,{\sf t}$ 
& $0$ & ${\sf h}_1$ & ${\sf h}_2$
& ${\sf d}$
\\
${\sf h}_1$ & $*$ & $0$ & $4\,{\sf t}$ 
& ${\sf h}_1$ & ${\sf h}_2$ & $6\,{\sf r}$ &
$2\,{\sf d}$ &
${\sf i}_1$
\\
${\sf h}_2$ & $*$ & $*$ & $0$ & ${\sf h}_2$ 
& $-{\sf h}_1$ & $-2\,{\sf d}$ &
$6\,{\sf r}$ & ${\sf
i}_2$
\\
${\sf d}$ & $*$ & $*$ & $*$ & $0$ 
& $0$ & ${\sf i}_1$ & ${\sf i}_2$ & $2\,{\sf j}$
\\
${\sf r}$ & $*$ & $*$ & $*$ & $*$ 
& $0$ & $-{\sf i}_2$ & ${\sf i}_1$ & $0$
\\
${\sf i}_1$ & $*$ & $*$ & $*$ & $*$ 
& $*$ & $0$ & $4\,{\sf j}$ & $0$
\\
${\sf i}_2$ & $*$ & $*$ & $*$ & $*$ & $*$ & $*$ & $0$ & $0$
\\
${\sf j}$ & $*$ & $*$ & $*$ & $*$ & $*$ & $*$ & $*$ & $0$
\end{tabular}
\end{center}

\section{Second Cohomology of the Heisenberg Lie Algebra}
\label{Cohomology-section}

\HEAD{Second Cohomology of the Heisenberg Lie Algebra}{
Mansour Aghasi, Joël Merker, and Masoud Sabzevari}

Specifically for a Cartan geometry modeled on a pair $(\mathfrak{
g}_-, \mathfrak{ g})$ of Lie algebras $\mathfrak{ g}_- \subset
\mathfrak{ g}$, the second cohomology $H^2 ( \mathfrak{ g}_-,
\mathfrak{ g})$ {\em a priori} provides some useful algebraic
information about the number of functionally independent Cartan
curvatures. Similarly as for Tanaka prolongations, the computations
for $H^2 ( \mathfrak{ g}_-, \mathfrak{ g})$ are of purely algebraic
nature, without any differentialo-geometric invariant coming into the
picture hence more elementary.
In Section~\ref{second-cohomology-graded}, 
general considerations and formulas
about (second) cohomologies of (graded) Lie algebra
are set up. 

Thus, let $\frak g$ be an $r$-dimensional Lie algebra and
let $\frak g_-$ be an $n$-dimensional ($1\leqslant n \leqslant r - 1$)
subalgebra of $\frak g$. For any $k \geqslant 1$, the space $\mathcal{
C}^k ( \mathfrak{ g}_-, \mathfrak{ g})$ of {\em k-cochains} consists
by definition (\cite{ Goze}, Chap.~3) of the space of linear maps from
$\Lambda^k \mathfrak{ g}_-$ to $\mathfrak{ g}$, that is to say:
\[
\mathcal{C}^k(\mathfrak{g}_-,\mathfrak{g})
=
{\rm Lin}\big(\Lambda^k\mathfrak{g}_-,\,\mathfrak{g}\big).
\]
When $\mathfrak{ g}$ is equipped with the structure of a graded Lie
algebra:
\[
\mathfrak{g}
=
\underbrace{
\mathfrak{g}_{-\mu}
\oplus\cdots\oplus
\mathfrak{g}_{-1}}_{\mathfrak{g}_-}
\oplus
\mathfrak{g}_0
\oplus
\mathfrak{g}_1
\oplus\cdots\oplus
\mathfrak{g}_\nu,
\]
with $\big[ \mathfrak{ g}_{ k_1} , \, \mathfrak{ g}_{ k_2} \big]
\subset \mathfrak{ g}_{ k_1+k_2}$ for any $k_1, k_2 \in \Z$ (with the
convention that $\mathfrak{ g}_k = \{ 0\}$ whenever $k \leqslant -\mu
- 1$ or $k \geqslant \nu + 1$), each vector space
$\mathcal{C}^k(\mathfrak{g}_-,\mathfrak{g})$ naturally splits into a
direct sum of so-called homogeneous cochains as follows: a $k$-cochain
$\Phi \in \mathcal{ C}^k (\mathfrak{g}_-,\mathfrak{g})$ is said to be {\em
of homogeneity} $h \in \mathbb Z$ whenever for arbitrary elements:
\[
{\sf z}_{i_1}\in\mathfrak{g}_{i_1},
\ldots\ldots,
{\sf z}_{i_k}\in\mathfrak{g}_{i_k}
\]
belonging to certain arbitrary determined $\mathfrak{ g}$-component,
its value:
\[
\Phi({\sf z}_{i_1},\ldots{\sf z}_{i_k})
\in
\mathfrak{g}_{i_1+\cdots+i_k+h}
\]
belongs to the $(i_1+\ldots +i_k+i)$-th component of $\frak g$. In
fact, one easily convinces oneself that any $k$-cochain $\Phi$ splits
as a direct sum of $k$-cochains of fixed homogeneity:
\[
\Phi
=
\cdots+
\Phi^{[h-1]}
+
\Phi^{[h]}
+
\Phi^{[h+1]}
+\cdots
\]
where we denote the $h$-th 
component of $\Phi$ just by $\Phi^{ [h]}$.

For each $k = 0, 1, \dots, n$,
the {\sl differential operator} (\cite{Goze, Cap, EMS}):
\[
\partial^k\colon\ \ \
\mathcal{C}^k\big(\mathfrak{g}_-,\,\mathfrak{g}\big) \longrightarrow
\mathcal{C}^{k+1}\big(\mathfrak{g}_-,\,\mathfrak{g}\big)
\]
assigns to a $k$-cochain $\Phi \in \mathcal{ C}^k \big( \mathfrak{
g}_-, \mathfrak{ g}\big)$ the $k+1$-cochain $\partial^k \Phi$ whose
value on any collection of $k+1$ vectors ${\sf z}_0, {\sf z}_1, \dots,
{\sf z}_k$ is defined through the formula:
\begin{equation}
\label{cohomology}
\small
\aligned
&
\
(\partial^k\Phi) 
\big({\sf z}_0,{\sf z}_1,\dots,{\sf z}_k\big) 
:= 
\sum_{i=0}^k\,(-1)^i
\big[{\sf{z}}_i,\,\Phi({\sf z}_0,\dots,
\widehat{\sf{z}}_i,\dots,{\sf z}_k)\big]_{\mathfrak{g}}
+
\\
&
\ \ \ \ \
+
\sum_{0\leqslant
i<j\leqslant k}\, (-1)^{i+j}\, 
\Phi
\big([{\sf z}_i,{\sf z}_j]_{\mathfrak{g}},{\sf
z}_0,\dots,\widehat{\sf{z}}_i,
\dots, \widehat{\sf{z}}_j,\dots,{\sf z}_k\big),
\endaligned
\end{equation}
where $\widehat{\sf{z}}_l$ means removal of the term
${\sf{z}}_l$. This action is linear with respect
to each argument ${\sf{z}}_i$, $i= 0, 1, \dots, k$. One can check that
the composition $\partial^{ k+1}
\circ \partial^k$ vanishes for any
$k\in\mathbb N$ and we have the following {\em cochain complex}:
\[
0 \overset{\partial^0}{\longrightarrow} \mathcal{C}^1
\overset{\partial^1}{\longrightarrow}
\mathcal{C}^2 \overset{\partial^2}{\longrightarrow} \cdots
\overset{\partial^{n-2}}{\longrightarrow} \mathcal{C}^{n-1}
\overset{\partial^{n-1}}{\longrightarrow} \mathcal{C}^n
\overset{\partial^n}{\longrightarrow} 0.
\]
The $k$-th cohomological space $H^k (\frak g_-,\frak g)$ is then
defined as being the following quotient:
\[
H^k\big(\mathfrak{g}_-,\mathfrak{g}\big)
=
\frac{{\rm ker}(\partial^k)}{{\rm im}(\partial^{k-1})}.
\]

Before entering specific computations, let us briefly motivate and
anticipate. Only the second cohomology space $H^2(\frak g_-,\frak g)$
will be of interest to us. In
Section~\ref{Cartan-connections-coordinates} ({\em cf.} also~\cite{
Cap, EMS}), we will indeed see that the the curvature function
(Definition~\ref{definition-curvature-function}) associated to a
Cartan connection takes its image in the set of 2-cochains ${\mathcal
C}^2(\frak g_-,\frak g)$. Also, the curvature function naturally
splits in homogeneous components. As explained in~\cite{ EMS}, the
so-called {\sl Bianchi-Tanaka identities} stated by
Proposition~\ref{Bianchi-Tanaka} below entail in particular that the
lowest order nonvanishing curvature must be $\partial$-closed, and
more generally, any homogeneous curvature component is determined by
the lower components up to a $\partial$-closed component ({\em see}
also Proposition~\ref{cohomology-result} below). Hence, some of the
linear-like properties of Cartan curvatures rely on probing the
corresponding second cohomological spaces.

\subsection{Cochains $\mathcal{ C}^2(\frak g_-,\frak g)$ for
the prolonged Heisenberg Lie algebra}
\label{homogeneity-Psi-0-5}
Now, let $\mathfrak{ g} = \mathfrak{
g}_{ -2} \oplus \mathfrak{ g}_{ -1} \oplus \mathfrak{ g}_0 \oplus
\mathfrak{ g}_1 \oplus \mathfrak{ g}_2$ be the
eight-dimensional (abstract) graded Lie algebra under
study with:
\[
\mathfrak{g}_{-2}=\mathbb{R}\,{\sf t},\ \ \
\mathfrak{g}_{-1}=\mathbb{R}\,{\sf h}_1\oplus\R\,{\sf h}_2,\ \ \
\mathfrak{g}_0=\mathbb{R}\,{\sf d}\oplus\R\,{\sf r},\ \ \
\mathfrak{g}_1=\mathbb{R}\,{\sf i}_1\oplus\R\,{\sf i}_2,\ \ \
\mathfrak{g}_2=\mathbb{R}\,{\sf j},
\]
its commutator table being shown above. Let $\{ {\sf t}^*, \, {\sf
h}_1^*, \, {\sf h}_2^* \}$ be the dual basis of $\frak g_- = \R \,
{\sf t} \oplus \R\, {\sf h}_1 \oplus \R\, {\sf h}_2$. Then with such
bases, a general $2$-cochain $\Phi \in \Lambda^2 \mathfrak{ g}_-^*
\otimes \mathfrak{ g}$ writes under the explicit expanded form:
\[
\scriptsize
\aligned
\Phi
&
=
{\sf t}^*\wedge{\sf h}_1^*
\otimes
\Big(
\phi_t^{th_1}\,{\sf t}
+
\phi_{h_1}^{th_1}\,{\sf h}_1
+
\phi_{h_2}^{th_1}\,{\sf h}_2
+
\phi_d^{th_1}\,{\sf d}
+
\phi_r^{th_1}\,{\sf r}
+
\phi_{i_1}^{th_1}\,{\sf i}_1
+
\phi_{i_2}^{th_1}\,{\sf i}_2
+
\phi_j^{th_1}\,{\sf j}
\Big)
+
\\
&
+
{\sf t}^*\wedge{\sf h}_2^*
\otimes
\Big(
\phi_t^{th_2}\,{\sf t}
+
\phi_{h_1}^{th_2}\,{\sf h}_1
+
\phi_{h_2}^{th_2}\,{\sf h}_2
+
\phi_d^{th_2}\,{\sf d}
+
\phi_r^{th_2}\,{\sf r}
+
\phi_{i_1}^{th_2}\,{\sf i}_1
+
\phi_{i_2}^{th_2}\,{\sf i}_2
+
\phi_j^{th_2}\,{\sf j}
\Big)
+
\\
&
+
{\sf h}_1^*\wedge{\sf h}_2^*
\otimes
\Big(
\phi_t^{h_1h_2}\,{\sf t}
+
\phi_{h_1}^{h_1h_2}\,{\sf h}_1
+
\phi_{h_2}^{h_1h_2}\,{\sf h}_2
+
\phi_d^{h_1h_2}\,{\sf d}
+
\phi_r^{h_1h_2}\,{\sf r}
+
\phi_{i_1}^{h_1h_2}\,{\sf i}_1
+
\phi_{i_2}^{h_1h_2}\,{\sf i}_2
+
\phi_j^{h_1h_2}\,{\sf j}
\Big),
\endaligned
\]
where the 24 real coefficients $\phi_t^{ th_1}$, \dots,
$\phi_j^{ h_1h_2}$ are arbitrary. This $2$-cochain could also
be written ({\em cf.}~\cite{ BES})
under a condensed symbolic form as:
\[
\Phi
=
\sum_{x<y}\,
\sum_v\,
\phi^{xy}_{v}\,
{\sf x}^\ast\wedge{\sf y}^\ast\otimes{\sf v},
\]
for ${\sf x}, {\sf y} \in \{{\sf h}^\ast_1, {\sf h}^\ast_2, {\sf
t}^\ast\}$ and for ${\sf v}\in \{ {\sf t}, {\sf h}_1, {\sf h}_2, {\sf
d}, {\sf r}, {\sf i}_1, {\sf i}_2, {\sf j} \}$, but the first,
complete writing is much more suited to effective calculations. Also,
it is useful to reorganize the 24 components of $\Phi$ by collecting,
in one and a single line, all those which have the same homogeneity:
\[
\scriptsize
\aligned
\Phi
&
=
\phi_t^{h_1h_2}\,{\sf h}_1^*\!\!\wedge{\sf h}_2^*\otimes{\sf t}
+
\ \ \ \ \
\fbox{\tiny 0}
\\
\fbox{\tiny 1}
&
\ \ \ \ \
+
\phi_t^{th_1}\,{\sf t}^*\!\!\wedge{\sf h}_1^*\otimes{\sf t}
+
\phi_t^{th_2}\,{\sf t}^*\!\!\wedge{\sf h}_2^*\otimes{\sf t}
+
\phi_{h_1}^{h_1h_2}\,{\sf h}_1^*\!\!\wedge{\sf h}_2^*\otimes{\sf t}
+
\phi_{h_2}^{h_1h_2}\,{\sf h}_1^*\!\!\wedge{\sf h}_2^*\otimes{\sf h}_2
+
\\
\fbox{\tiny 2}
&
\ \ \ \ \
+
\phi_{h_1}^{th_1}\,{\sf t}^*\!\!\wedge{\sf h}_1^*\otimes{\sf h}_1
+
\phi_{h_2}^{th_1}\,{\sf t}^*\!\!\wedge{\sf h}_1^*\otimes{\sf h}_2
+
\phi_{h_1}^{th_2}\,{\sf t}^*\!\!\wedge{\sf h}_2^*\otimes{\sf h}_1
+
\phi_{h_2}^{th_2}\,{\sf t}^*\!\!\wedge{\sf h}_2^*\otimes{\sf h}_2
+
\\
&
\ \ \ \ \ \ \ \ \ \ \ \ \ \ \ \ \ \ \ \ \ \ \ \ \ \ \ \ \ \ \ \ \ \
\ \ \ \ \ \ \ \ \ \ \ \ \ \ \ \ \ \ \ \ \ \ \ \ \ \ \ \ \ \ \ \ \ \ 
\ \ \ \
+
\phi_d^{h_1h_2}\,{\sf h}_1^*\!\!\wedge{\sf h}_2^*\otimes{\sf d}
+
\phi_r^{h_1h_2}\,{\sf h}_1^*\!\!\wedge{\sf h}_2^*\otimes{\sf r}
+
\\
\fbox{\tiny 3}
&
\ \ \ \ \
+
\phi_d^{th_1}\,{\sf t}^*\!\!\wedge{\sf h}_1^*\otimes{\sf d}
+
\phi_r^{th_1}\,{\sf t}^*\!\!\wedge{\sf h}_1^*\otimes{\sf r}
+
\phi_d^{th_2}\,{\sf t}^*\!\!\wedge{\sf h}_2^*\otimes{\sf d}
+
\phi_r^{th_2}\,{\sf t}^*\!\!\wedge{\sf h}_2^*\otimes{\sf r}
\\
&
\ \ \ \ \ \ \ \ \ \ \ \ \ \ \ \ \ \ \ \ \ \ \ \ \ \ \ \ \ \ \ \ \ \
\ \ \ \ \ \ \ \ \ \ \ \ \ \ \ \ \ \ \ \ \ \ \ \ \ \ \ \ \ \ \ \ 
+
\phi_{i_1}^{h_1h_2}\,{\sf h}_1^*\!\!\wedge{\sf h_2}^*\otimes{\sf i}_1
+
\phi_{i_2}^{h_1h_2}\,{\sf h}_1^*\!\!\wedge{\sf h}_2^*\otimes{\sf i}_2
+
\\
\fbox{\tiny 4}
&
\ \ \ \ \
+
\phi_{i_1}^{th_1}\,{\sf t}^*\!\!\wedge{\sf h}_1^*\otimes{\sf i}_1
+
\phi_{i_2}^{th_1}\,{\sf t}^*\!\!\wedge{\sf h}_1^*\otimes{\sf i}_2
+
\phi_{i_1}^{th_2}\,{\sf t}^*\!\!\wedge{\sf h}_2^*\otimes{\sf i}_1
+
\phi_{i_2}^{th_2}\,{\sf t}^*\!\!\wedge{\sf h}_2^*\otimes{\sf i}_2
\\
&
\ \ \ \ \ \ \ \ \ \ \ \ \ \ \ \ \ \ \ \ \ \ \ \ \ \ \ \ \ \ \ \ \ \
\ \ \ \ \ \ \ \ \ \ \ \ \ \ \ \ \ \ \ \ \ \ \ \ \ \ \ \ \ \ \ \ \ \ 
\ \ \ \ \ \ \ \ \ \ \ \ \ \ \ \ \ \ \ \ \ \ \ \ \ \ \ \ \ \ \ \ \, 
+
\phi_j^{h_1h_2}\,{\sf h}_1^*\!\!\wedge{\sf h}_2^*\otimes{\sf j}
+
\\
\fbox{\tiny 5}
&
\ \ \ \ \
+
\phi_j^{th_1}\,{\sf t}^*\!\!\wedge{\sf h}_1^*\otimes{\sf j}
+
\phi_j^{th_2}\,{\sf t}^*\!\!\wedge{\sf h}_1^*\otimes{\sf j},
\endaligned
\]
starting from homogeneity $0$ (first line) up to homogeneity $5$ (last
line). Thus, the graded dimensions of $\mathcal{ C}^2_0$ $\mathcal{
C}_1^2$, $\mathcal{ C}_2^2$, $\mathcal{ C}_3^2$, $\mathcal{ C}_4^2$,
$\mathcal{ C}_5^2$ are equal, respectively, to: $1$, $4$, $6$, $6$,
$5$, $2$, {\em cf.} also the summarizing table at the end of this
section.

\subsection{Computations of $\mathcal Z^2(\frak g_-,\frak g)$}
Now, such a general $2$-cochain $\Phi$ belongs to ${\mathcal Z}^2$ if
and only if the value of $\partial \Phi$ on each antisymmetric
$3$-vector of $\Lambda^3 \mathfrak{ g}_-$ vanishes. But here,
$\Lambda^3 \mathfrak{ g}_-$ is one-dimensional, generated by just
${\sf t} \wedge {\sf h}_1 \wedge {\sf h}_2$. Consequently, applying
the definition~\thetag{ \ref{cohomology}}, the cocycle condition
amounts to the single equation:
\[
\small
\aligned
0
&
=
\partial\Phi({\sf t},{\sf h}_1,{\sf h}_2)
\\
&
=
\big[{\sf t},\,\Phi({\sf h}_1,{\sf h}_2)\big]_{\mathfrak{g}}
-
\big[{\sf h}_1,\,\Phi({\sf t},{\sf h}_2)\big]_{\mathfrak{g}}
+
\big[{\sf h}_2,\,\Phi({\sf t},{\sf h}_1)\big]_{\mathfrak{g}}
-
\\
&
\ \ \ \ \
-
\Phi\big([{\sf t},{\sf h}_1]_{\mathfrak{g}},\,{\sf h}_2\big)
+
\Phi\big([{\sf t},{\sf h}_2]_{\mathfrak{g}},\,{\sf h}_1\big)
-
\Phi\big([{\sf h}_1,{\sf h}_2]_{\mathfrak{g}},\,{\sf t}\big),
\endaligned
\]
and then, after substituting all corresponding values of
$\Phi( \cdot, \cdot)$, we get:
\[
\footnotesize
\aligned
0
&
=
\big[{\sf t},\,
\phi_t^{h_1h_2}\,{\sf t}
+
\phi_{h_1}^{h_1h_2}\,{\sf h}_1
+
\phi_{h_2}^{h_1h_2}\,{\sf h}_2
+
\phi_d^{h_1h_2}\,{\sf d}
+
\phi_r^{h_1h_2}\,{\sf r}
+
\phi_{i_1}^{h_1h_2}\,{\sf i}_1
+
\phi_{i_2}^{h_1h_2}\,{\sf i}_2
+
\phi_j^{h_1h_2}\,{\sf j}
\big]
-
\\
&
\ \ \ \ \
-
\big[{\sf h}_1,\,
\phi_t^{th_2}\,{\sf t}
+
\phi_{h_1}^{th_2}\,{\sf h}_1
+
\phi_{h_2}^{th_2}\,{\sf h}_2
+
\phi_d^{th_2}\,{\sf d}
+
\phi_r^{th_2}\,{\sf r}
+
\phi_{i_1}^{th_2}\,{\sf i}_1
+
\phi_{i_2}^{th_2}\,{\sf i}_2
+
\phi_j^{th_2}\,{\sf j}
\big]
+
\\
&
\ \ \ \ \
+
\big[{\sf h}_2,\,
\phi_t^{th_1}\,{\sf t}
+
\phi_{h_1}^{th_1}\,{\sf h}_1
+
\phi_{h_2}^{th_1}\,{\sf h}_2
+
\phi_d^{th_1}\,{\sf d}
+
\phi_r^{th_1}\,{\sf r}
+
\phi_{i_1}^{th_1}\,{\sf i}_1
+
\phi_{i_2}^{th_1}\,{\sf i}_2
+
\phi_j^{th_1}\,{\sf j}
\big]
-
\\
&
\ \ \ \ \
-
\zero{\Phi(0,{\sf h}_2)}
+
\zero{\Phi(0,{\sf h}_1)}
-
\zero{\Phi(4{\sf t},{\sf t})}.
\endaligned
\]
Using the commutator table, we may then replace each appearing
Lie bracket:
\[
\footnotesize
\aligned
0
&
=
2\phi_d^{h_1h_2}\,{\sf t}
+
\phi_{i_1}^{h_1h_2}\,{\sf h}_1
+
\phi_{i_2}^{h_1h_2}\,{\sf h}_2
+
\phi_j^{h_1h_2}\,{\sf d}
-
\\
&
\ \ \ \ \
-4\phi_{h_2}^{th_2}\,{\sf t}
-
\phi_d^{th_2}\,{\sf h}_1
-
\phi_r^{th_2}\,{\sf h}_2
-
6\phi_{i_1}^{th_2}\,{\sf r}
-
2\phi_{i_2}^{th_2}\,{\sf d}
-
\phi_j^{th_2}\,{\sf i}_1
-
\\
&
\ \ \ \ \
-4\phi_{h_1}^{th_1}\,{\sf t}
+
\phi_d^{th_1}\,{\sf h}_2
-
\phi_r^{th_1}\,{\sf h}_1
-
2\phi_{i_1}^{th_1}\,{\sf d}
+
6\phi_{i_2}^{th_1}\,{\sf r}
+
\phi_j^{th_1}\,{\sf i}_2,
\endaligned
\]
and lastly, gather the coefficients of the appearing
vectors ${\sf t}$, \dots, ${\sf i}_2$:
\[
\footnotesize
\aligned
0
&
=
(2\phi_d^{h_1h_2}-4\phi_{h_2}^{th_2}-4\phi_{h_1}^{th_1})\,{\sf t}
+
(\phi_{i_1}^{h_1h_2}-\phi_d^{th_2}-\phi_r^{th_1})\,{\sf h}_1
+
(\phi_{i_2}^{h_1h_2}-\phi_r^{th_2}+\phi_d^{th_1})\,{\sf h}_2
+
\\
&
\ \ \ \ \
+
(\phi_j^{h_1h_2}-2\phi_{i_2}^{th_2}-2\phi_{i_1}^{th_1})\,{\sf d}
+
(-6\phi_{i_1}^{th_2}+6\phi_{i_2}^{th_1})\,{\sf r}
+
(-\phi_j^{th_2})\,{\sf i}_1
+
(\phi_j^{th_1})\,{\sf i}_2.
\endaligned
\]
Thus, a $2$-cochain $\Phi$ is a $2$-cocycle if and only if its 24
coefficients satisfy the following seven linear equations, ordered by
increasing homogeneity:
\[
\footnotesize
\aligned
&
\fbox{\tiny 2}
\ \ \ \ \ \ \ \ \ \
0
=
2\phi_d^{h_1h_2}-4\phi_{h_2}^{th_2}-4\phi_{h_1}^{th_1},
\\
&
\fbox{\tiny 3}
\ \ \ \ \ \ \ \ \ \
0
=
\phi_{i_1}^{h_1h_2}-\phi_d^{th_2}-\phi_r^{th_1},
\ \ \ \ \ \ \ \ \ \ \ \ \ \
0
=
\phi_{i_2}^{h_1h_2}-\phi_r^{th_2}+\phi_d^{th_1},
\\
&
\fbox{\tiny 4}
\ \ \ \ \ \ \ \ \ \
0
=
\phi_j^{h_1h_2}-2\phi_{i_2}^{th_2}-2\phi_{i_1}^{th_1},
\ \ \ \ \ \ \ \ \ \
0
=
-6\phi_{i_1}^{th_2}+6\phi_{i_2}^{th_1},
\\
&
\fbox{\tiny 5}
\ \ \ \ \ \ \ \ \ \
0
=
-\phi_j^{th_2},
\ \ \ \ \ \ \ \ \ \
0
=
\phi_j^{th_1}.
\endaligned
\]
All these equations are visibly linearly independent, and we deduce that
the homogeneous components $\mathcal{ Z}_{[h]}^2$ of:
\[
\mathcal{Z}^2
=
\mathcal{Z}_{[0]}^2
\oplus
\mathcal{Z}_{[1]}^2
\oplus
\mathcal{Z}_{[2]}^2
\oplus
\mathcal{Z}_{[3]}^2
\oplus
\mathcal{Z}_{[4]}^2
\oplus
\mathcal{Z}_{[5]}^2
\]
have codimensions within $\mathcal{ C}_{[h]}^2$ equal to $0$, $0$, $1$,
$2$, $2$, $2$, hence are of dimensions equal to $1$, $4$, $5$, $4$,
$3$, $0$ respectively,

\subsection{Determination of $\mathcal B^2(\frak g_-,\frak g)$}
\label{B2} Now, let $\Psi \in \Lambda^1 \mathfrak{ g}_-^*
\otimes \mathfrak{ g}$ be a general $1$-cochain. In terms of the bases
$\{ {\sf t}^*, \, {\sf h}_1^*, \, {\sf h}_2^* \}$ of $\mathfrak{
g}_-^*$ and $\{ {\sf t}, \, {\sf h}_1, \, {\sf h}_2, \, {\sf d}, \,
{\sf r}, \, {\sf i}_1, \, {\sf i}_2,\, {\sf j} \}$ of $\mathfrak{ g}$,
it writes under the explicit expanded form:
\[
\footnotesize
\aligned
\Psi
&
=
{\sf t}^*\otimes
\Big(
\psi_t^t\,{\sf t}
+
\psi_{h_1}^t\,{\sf h}_1
+
\psi_{h_2}^t\,{\sf h}_2
+
\psi_d^t\,{\sf d}
+
\psi_r^t\,{\sf r}
+
\psi_{i_1}^t\,{\sf i}_1
+
\psi_{i_2}^t\,{\sf i}_2
+
\psi_j^t\,{\sf j}
\Big)
+
\\
&
\ \ \ \ \
+
{\sf h}_1^*
\otimes
\Big(
\psi_t^{h_1}\,{\sf t}
+
\psi_{h_1}^{h_1}\,{\sf h}_1
+
\psi_{h_2}^{h_1}\,{\sf h}_2
+
\psi_d^{h_1}\,{\sf d}
+
\psi_r^{h_1}\,{\sf r}
+
\psi_{i_1}^{h_1}\,{\sf i}_1
+
\psi_{i_2}^{h_1}\,{\sf i}_2
+
\psi_j^{h_1}\,{\sf j}
\Big)
+
\\
&
\ \ \ \ \
+
{\sf h}_1^*
\otimes
\Big(
\psi_t^{h_2}\,{\sf t}
+
\psi_{h_1}^{h_2}\,{\sf h}_1
+
\psi_{h_2}^{h_2}\,{\sf h}_2
+
\psi_d^{h_2}\,{\sf d}
+
\psi_r^{h_2}\,{\sf r}
+
\psi_{i_1}^{h_2}\,{\sf i}_1
+
\psi_{i_2}^{h_2}\,{\sf i}_2
+
\psi_j^{h_2}\,{\sf j}
\Big),
\endaligned
\]
where the 24 real coefficients $\psi_t^t$, \dots, $\psi_j^{j_2}$ are
arbitrary. Equivalently, by collecting in one and
a single line all components having equal homogeneity, such
a general $1$-cochains writes:
\[
\footnotesize
\aligned
\Psi
&
=
\psi_t^{h_1}\,{\sf h}_1^*\otimes{\sf t}
+
\psi_t^{h_2}\,{\sf h}_2^*\otimes{\sf t}
+
\ \ \ \ \ \ \ \ \ \ \
\fbox{\tiny -1}
\\
\fbox{\tiny 0}
&
\ \ \ \ \
+
\psi_t^t\,{\sf t}^*\otimes{\sf t}
+
\psi_{h_1}^{h_1}\,{\sf h}_1^*\otimes{\sf h}_1
+
\psi_{h_2}^{h_1}\,{\sf h}_1^*\otimes{\sf h}_2
+
\psi_{h_1}^{h_2}\,{\sf h}_2^*\otimes{\sf h}_1
+
\psi_{h_2}^{h_2}\,{\sf h}_2^*\otimes{\sf h}_2
+
\\
\fbox{\tiny 1}
&
\ \ \ \ \
+
\psi_{h_1}^t\,{\sf t}^*\otimes{\sf h}_1
+
\psi_{h_2}^t\,{\sf t}^*\otimes{\sf h}_2
+
\psi_d^{h_1}\,{\sf h}_1^*\otimes{\sf d}
+
\psi_r^{h_1}\,{\sf h}_1^*\otimes{\sf r}
+
\psi_d^{h_2}\,{\sf h}_2^*\otimes{\sf d}
+
\psi_r^{h_2}\,{\sf h}_2^*\,\otimes{\sf r}
+
\\
\fbox{\tiny 2}
&
\ \ \ \ \
+
\psi_d^t\,{\sf t}^*\otimes{\sf d}
+
\psi_r^t\,{\sf t}^*\otimes{\sf r}
+
\psi_{i_1}^{h_1}\,{\sf h}_1^*\otimes{\sf i}_1
+
\psi_{i_2}^{h_1}\,{\sf h}_1^*\otimes{\sf i}_2
+
\psi_{i_1}^{h_2}\,{\sf h}_2^*\otimes{\sf i}_1
+
\psi_{i_2}^{h_2}\,{\sf h}_2^*\otimes{\sf i}_2
+
\\
\fbox{\tiny 3}
&
\ \ \ \ \
+
\psi_{i_1}^t\,{\sf t}^*\otimes{\sf i}_1
+
\psi_{i_2}^t\,{\sf t}^*\otimes{\sf i}_2
+
\psi_j^{h_1}\,{\sf h}_1^*\otimes{\sf j}
+
\psi_j^{h_2}\,{\sf h}_2^*\otimes{\sf j}
+
\\
\fbox{\tiny 4}
&
\ \ \ \ \
+
\psi_j^t\,{\sf t}^*\otimes{\sf j}.
\endaligned
\]
In order to characterize when a $2$-cochain $\Phi$ is of the form
$\partial \Psi$ namely is a coboundary, applying the
definition~\thetag{ \ref{cohomology}}, we at first compute the values
of $\partial \Psi$ on each of the three antisymmetric $2$-vectors
${\sf t} \wedge {\sf h}_1$, ${\sf t} \wedge {\sf h}_2$, ${\sf h}_1
\wedge {\sf h}_2$ which make up a natural basis for $\Lambda^2
\mathfrak{ g}_-$, and with intermediate details, we obtain:
\[
\footnotesize
\aligned
(\partial\Psi)({\sf t},{\sf h}_1)
&
=
\big[{\sf t},\Psi({\sf h}_1)\big]
-
\big[{\sf h}_1,\Psi({\sf t})\big]
-
\Psi\big([{\sf t},{\sf h}_1]\big)
\\
&
=
\big[{\sf t},\,
\psi_t^{h_1}\,{\sf t}
+
\psi_{h_1}^{h_1}\,{\sf h}_1
+
\psi_{h_2}^{h_1}\,{\sf h}_2
+
\psi_d^{h_1}\,{\sf d}
+
\psi_r^{h_1}\,{\sf r}
+
\psi_{i_1}^{h_1}\,{\sf i}_1
+
\psi_{i_2}^{h_1}\,{\sf i}_2
+
\psi_j^{h_1}\,{\sf j}
\big]
-
\\
&
\ \ \ \ \
-
\big[{\sf h}_1,\,
\psi_t^t\,{\sf t}
+
\psi_{h_1}^t\,{\sf h}_1
+
\psi_{h_2}^t\,{\sf h}_2
+
\psi_d^t\,{\sf d}
+
\psi_r^t\,{\sf r}
+
\psi_{i_1}^t\,{\sf i}_1
+
\psi_{i_2}^t\,{\sf i}_2
+
\psi_j^t\,{\sf j}
\big]
-
\\
&
\ \ \ \ \
-\zero{\Psi(0)}
\\
&
=
2\psi_d^{h_1}\,{\sf t}
+
\psi_{i_1}^{h_1}\,{\sf h}_1
+
\psi_{i_2}^{h_1}\,{\sf h}_2
+
\psi_j^{h_1}\,{\sf d}
-
\\
&
\ \ \ \ \
-4\psi_{h_2}^t\,{\sf t}
-
\psi_d^t\,{\sf h}_1
-
\psi_r^t\,{\sf h}_2
-
6\psi_{i_1}^t\,{\sf r}
-
2\psi_{i_2}^t\,{\sf d}
-
\psi_j^t\,{\sf i}_1
\\
&
=
(2\psi_d^{h_1}-4\psi_{h_2}^t)\,{\sf t}
+
(\psi_{i_1}^{h_1}-\psi_d^t)\,{\sf h}_1
+
(\psi_{i_2}^{h_1}-\psi_r^t)\,{\sf h}_2
+
(\psi_j^{h_1}-2\psi_{i_2}^t)\,{\sf d}
+
\\
&
\ \ \ \ \
+
(-6\psi_{i_1}^t)\,{\sf r}
+
(-\psi_j^t)\,{\sf i}_1
+
0\,{\sf i}_2
+
0\,{\sf j},
\endaligned
\]
\[
\footnotesize
\aligned
(\partial\Psi)({\sf t},{\sf h}_2)
&
=
\big[{\sf t},\Psi({\sf h}_2)\big]
-
\big[{\sf h}_2,\Psi({\sf t})\big]
-
\Psi\big([{\sf t},{\sf h}_2]\big)
\\
&
=
\big[{\sf t},\,
\psi_t^{h_2}\,{\sf t}
+
\psi_{h_1}^{h_2}\,{\sf h}_1
+
\psi_{h_2}^{h_2}\,{\sf h}_2
+
\psi_d^{h_2}\,{\sf d}
+
\psi_r^{h_2}\,{\sf r}
+
\psi_{i_1}^{h_2}\,{\sf i}_1
+
\psi_{i_2}^{h_2}\,{\sf i}_2
+
\psi_j^{h_2}\,{\sf j}
\big]
-
\\
&
\ \ \ \ \
-
\big[{\sf h}_2,\,
\psi_t^t\,{\sf t}
+
\psi_{h_1}^t\,{\sf h}_1
+
\psi_{h_2}^t\,{\sf h}_2
+
\psi_d^t\,{\sf d}
+
\psi_r^t\,{\sf r}
+
\psi_{i_1}^t\,{\sf i}_1
+
\psi_{i_2}^t\,{\sf i}_2
+
\psi_j^t\,{\sf j}
\big]
-
\\
&
\ \ \ \ \
-
\zero{\Psi(0)}
\\
&
=
2\psi_d^{h_1}\,{\sf t}
+
\psi_{i_1}^{h_2}\,{\sf h}_1
+
\psi_{i_2}^{h_2}\,{\sf h}_2
+
\psi_j^{h_2}\,{\sf d}
+
\\
&
\ \ \ \ \
+
4\psi_{h_1}^t\,{\sf t}
-
\psi_d^t\,{\sf h}_2
+
\psi_r^t\,{\sf h}_1
+
2\psi_{i_1}^t\,{\sf d}
-
6\psi_{i_2}^r\,{\sf r}
-
\psi_j^t\,{\sf i}_2
\\
&
=
(2\psi_d^{h_2}+4\psi_{h_1}^t)\,{\sf t}
+
(\psi_{i_1}^{h_2}+\psi_r^t)\,{\sf h}_1
+
(\psi_{i_2}^{h_2}-\psi_d^t)\,{\sf h}_2
+
(\psi_j^{h_2}+2\psi_{i_1}^t)\,{\sf d}
+
\\
&
\ \ \ \ \
+
(-6\psi_{i_2}^t)\,{\sf r}
+
0\,{\sf i}_1
+
(-\psi_j^t)\,{\sf i}_2
+
0\,{\sf j},
\endaligned
\]
\[
\footnotesize
\aligned
(\partial\Psi)({\sf h}_1,{\sf h}_2)
&
=
\big[{\sf h}_1,\Psi({\sf h}_2)\big]
-
\big[{\sf h}_2,\Psi({\sf h}_1)\big]
-
\Psi\big([{\sf h}_1,{\sf h}_2]\big)
\\
&
=
\big[{\sf h}_1,\,
\psi_t^{h_2}\,{\sf t}
+
\psi_{h_1}^{h_2}\,{\sf h}_1
+
\psi_{h_2}^{h_2}\,{\sf h}_2
+
\psi_d^{h_2}\,{\sf d}
+
\psi_r^{h_2}\,{\sf r}
+
\psi_{i_1}^{h_2}\,{\sf i}_1
+
\psi_{i_2}^{h_2}\,{\sf i}_2
+
\psi_j^{h_2}\,{\sf j}
\big]
-
\\
&
\ \ \ \ \
-
\big[{\sf h}_2,\,
\psi_t^{h_1}\,{\sf t}
+
\psi_{h_1}^{h_1}\,{\sf h}_1
+
\psi_{h_2}^{h_1}\,{\sf h}_2
+
\psi_d^{h_1}\,{\sf d}
+
\psi_r^{h_1}\,{\sf r}
+
\psi_{i_1}^{h_1}\,{\sf i}_1
+
\psi_{i_2}^{h_1}\,{\sf i}_2
+
\psi_j^{h_1}\,{\sf j}
\big]
-
\\
&
\ \ \ \ \
-
4\psi_t^t\,{\sf t}
-
4\psi_{h_1}^t\,{\sf h}_1
-
4\psi_{h_2}^t\,{\sf h}_2
-
4\psi_d^t\,{\sf d}
-
4\psi_r^t\,{\sf r}
-
4\psi_{i_1}^t\,{\sf i}_1
-
4\psi_{i_2}^t\,{\sf i}_2
-
4\psi_j^t\,{\sf j}
\\
&
=
4\psi_{h_2}^{h_2}\,{\sf t}
+
\psi_d^{h_2}\,{\sf h}_1
+
\psi_r^{h_2}\,{\sf h}_2
+
6\psi_{i_1}^{h_2}\,{\sf r}
+
2\psi_{i_2}^{h_2}\,{\sf d}
+
\psi_j^{h_2}\,{\sf i}_1
+
\\
&
\ \ \ \ \
+
4\psi_{h_1}^{h_1}\,{\sf t}
-
\psi_d^{h_1}\,{\sf h}_2
+
\psi_r^{h_1}\,{\sf h}_1
+
2\psi_{i_1}^{h_1}\,{\sf d}
-
6\psi_{i_2}^{h_1}\,{\sf r}
-
\psi_j^{h_1}\,{\sf i}_2
-
\\
&
\ \ \ \ \ \
-
4\psi_t^t\,{\sf t}
-
4\psi_{h_1}^t\,{\sf h}_1
-
4\psi_{h_2}^t\,{\sf h}_2
-
4\psi_d^t\,{\sf d}
-
4\psi_r^t\,{\sf r}
-
4\psi_{i_1}^t\,{\sf i}_1
-
4\psi_{i_2}^t\,{\sf i}_2
-
4\psi_j^t\,{\sf j}
\\
&
=
(4\psi_{h_2}^{h_2}+4\psi_{h_1}^{h_1}-4\psi_t^t)\,{\sf t}
+
(\psi_d^{h_2}+\psi_r^{h_1}-4\psi_{h_1}^t)\,{\sf h}_1
+
(\psi_r^{h_2}-\psi_d^{h_1}-4\psi_{h_2}^t)\,{\sf h}_2
+
\\
&
\ \ \ \ \
+
(2\psi_{i_2}^{h_2}+2\psi_{i_1}^{h_1}-4\psi_d^t)\,{\sf d}
+
(6\psi_{i_1}^{h_2}-6\psi_{i_2}^{h_1}-4\psi_r^t)\,{\sf r}
+
(\psi_j^{h_2}-4\psi_{i_1}^t)\,{\sf i}_1
+
\\
&
\ \ \ \ \
+
(-\psi_j^{h_1}-4\psi_{i_2}^t)\,{\sf i}_2
+
(-4\psi_j^t)\,{\sf j}.
\endaligned
\]
As a result, a $2$-cochain $\Phi$ of the general form written above
equals the coboundary $\partial \Psi$ of a $1$-cochain if and only if
there exist 24 quantities $\psi_\cdot^\cdot$ such that $\Psi$'s
three families of eight coefficients $\phi_\cdot^{ th_1}$,
$\phi_\cdot^{ th_2}$, $\phi_\cdot^{ h_1h_2}$ are equal, respectively,
to the three collections of eight coefficients just found:
\[
\footnotesize
\aligned
\fbox{\tiny 1}\ \ \ \ \
\phi_t^{th_1}
&
=
2\psi_d^{h_1}-4\psi_{h_2}^t
\\
\fbox{\tiny 2}\ \ \ \ \
\phi_{h_1}^{th_1}
&
=
\psi_{i_1}^{h_1}-\psi_d^t
\\
\fbox{\tiny 2}\ \ \ \ \
\phi_{h_2}^{th_1}
&
=
\psi_{i_2}^{h_1}-\psi_r^t
\\
\fbox{\tiny 3}\ \ \ \ \
\phi_d^{th_1}
&
=
\psi_j^{h_1}-2\psi_{i_2}^t
\\
\fbox{\tiny 3}\ \ \ \ \
\phi_r^{th_1}
&
=
-6\psi_{i_1}^t
\\
\fbox{\tiny 4}\ \ \ \ \
\phi_{i_1}^{th_1}
&
=
-\psi_j^t
\\
\fbox{\tiny 4}\ \ \ \ \
\phi_{i_2}^{th_1}
&
=
0
\\
\fbox{\tiny 5}\ \ \ \ \
\phi_j^{th_1}
&
=
0
\endaligned
\ \ \ \ \
\aligned
\fbox{\tiny 1}\ \ \ \ \
\phi_t^{th_2}
&
=
2\psi_d^{h_2}+4\psi_{h_1}^t
\\
\fbox{\tiny 2}\ \ \ \ \
\phi_{h_1}^{th_2}
&
=
\psi_{i_1}^{h_2}+\psi_r^t
\\
\fbox{\tiny 2}\ \ \ \ \
\phi_{h_2}^{th_2}
&
=
\psi_{i_2}^{h_2}-\psi_d^t
\\
\fbox{\tiny 3}\ \ \ \ \
\phi_d^{th_2}
&
=
\psi_j^{h_2}+2\psi_{i_1}^t
\\
\fbox{\tiny 3}\ \ \ \ \
\phi_r^{th_2}
&
=
-6\psi_{i_2}^t
\\
\fbox{\tiny 4}\ \ \ \ \
\phi_{i_1}^{th_2}
&
=
0
\\
\fbox{\tiny 4}\ \ \ \ \
\phi_{i_2}^{th_2}
&
=
-\psi_j^t
\\
\fbox{\tiny 5}\ \ \ \ \
\phi_j^{th_2}
&
=
0
\endaligned
\ \ \ \ \
\aligned
\fbox{\tiny 0}\ \ \ \ \
\phi_t^{h_1h_2}
&
=
4\psi_{h_2}^{h_2}+4\psi_{h_1}^{h_1}-4\psi_t^t
\\
\fbox{\tiny 1}\ \ \ \ \
\phi_{h_1}^{h_1h_2}
&
=
\psi_d^{h_2}+\psi_r^{h_1}-4\psi_{h_1}^t
\\
\fbox{\tiny 1}\ \ \ \ \
\phi_{h_2}^{h_1h_2}
&
=
\psi_r^{h_2}-\psi_d^{h_1}+4\psi_{h_2}^t
\\
\fbox{\tiny 2}\ \ \ \ \
\phi_d^{h_1h_2}
&
=
2\psi_{i_2}^{h_2}+2\psi_{i_1}^{h_1}-4\psi_d^t
\\
\fbox{\tiny 2}\ \ \ \ \
\phi_r^{h_1h_2}
&
=
6\psi_{i_1}^{h_2}-6\psi_{i_2}^{h_1}-4\psi_r^t
\\
\fbox{\tiny 3}\ \ \ \ \
\phi_{i_1}^{h_1h_2}
&
=
\psi_j^{h_2}-4\psi_{i_1}^t
\\
\fbox{\tiny 3}\ \ \ \ \
\phi_{i_2}^{h_1h_2}
&
=
-\psi_j^{h_1}-4\psi_{i_2}^t
\\
\fbox{\tiny 4}\ \ \ \ \
\phi_j^{h_1h_2}
&
=
-4\psi_j^t.
\endaligned
\]

\subsection{Graded computation of $H^2(\frak g_-,\frak g)$}
\label{H2} Now, the map $\Psi \mapsto \partial \Psi =: \Phi$
so obtained explicitly is visibly linear $(\psi_\cdot^\cdot) \longmapsto
(\phi_\cdot^{ \cdot, \cdot})$, and furthermore, because $\partial \Psi
\in \mathcal{ Z}^2$ and because cochains naturally split in
homogeneous components, this map happens to be a direct sum of
six linear maps:
\[
\mathcal{C}_{[0]}^1\to\mathcal{Z}_{[0]}^2,\ \ \
\mathcal{C}_{[1]}^1\to\mathcal{Z}_{[1]}^2,\ \ \
\mathcal{C}_{[2]}^1\to\mathcal{Z}_{[2]}^2,\ \ \
\mathcal{C}_{[3]}^1\to\mathcal{Z}_{[3]}^2,\ \ \
\mathcal{C}_{[4]}^1\to\mathcal{Z}_{[4]}^2,\ \ \
\mathcal{C}_{[5]}^1\to\mathcal{Z}_{[5]}^2,
\]
the last one being just $\{ 0\} \to \{ 0\}$, that is to say a direct
sum of the following five explicit nonzero linear maps:
\[
\footnotesize
\aligned
\partial_{[0]}\colon\ \
\big(\psi_t^t,\psi_{h_1}^{h_1},\psi_{h_2}^{h_1},
\psi_{h_1}^{h_2},\psi_{h_2}^{h_2}\big)
\longmapsto
&\,
\big(4\psi_{h_2}^{h_2}+4\psi_{h_1}^{h_1}-4\psi_t^t\big)
\\
&\,
=
\big(\phi_t^{h_1h_2}\big)
\endaligned
\]
\[
\footnotesize
\aligned
\partial_{[1]}\colon\ \
\big(\psi_{h_1}^t,\psi_{h_2}^t,\psi_d^{h_1},\psi_r^{h_1},\psi_d^{h_2},
\psi_r^{h_2}\big)
\longmapsto
&\,
\big(2\psi_d^{h_1}-4\psi_{h_2}^t,2\psi_d^{h_2}+4\psi_{h_1}^t,
\psi_d^{h_2}+\psi_r^{h_1}-4\psi_{h_1}^t,
\\
&
\ \ \ \ \ \ \ \ \ \ \ \ \ \ \ \ \ \ \ \ \ \ \ \ \ \ \ \ \ \ \ \ \ \ \
\ \ \ \ \ \ \ \ \ \ \ \ \
\psi_r^{h_2}-\psi_d^{h_1}+4\psi_{h_2}^t
\big)
\\
&\,
=
\big(\phi_t^{th_1},\phi_t^{th_2},\phi_{h_1}^{h_1h_2},
\phi_{h_2}^{h_1h_2}\big)
\endaligned
\]
\[
\footnotesize
\aligned
\partial_{[2]}\colon\ \
\big(\psi_d^t,\psi_r^t,\psi_{i_1}^{h_1},\psi_{i_2}^{h_1},
\psi_{i_1}^{h_2},\psi_{i_2}^{h_2}\big)
\longmapsto
&\,
\big(
\psi_{i_1}^{h_1}-\psi_d^t,\psi_{i_2}^{h_1}-\psi_r^t,
\psi_{i_1}^{h_2}+\psi_r^t,\psi_{i_2}^{h_2}-\psi_d^t,
\\
&
\ \ \ \ \ \ \ \ \ \ \ \ \ \ \ \ \ \
2\psi_{i_2}^{h_2}+2\psi_{i_1}^{h_1}-4\psi_d^t,
6\psi_{i_1}^{h_2}-6\psi_{i_2}^{h_1}-4\psi_r^t
\big)
\\
&\,
=
\big(\phi_{h_1}^{th_1},\phi_{h_2}^{th_1},\phi_{h_1}^{th_2},
\phi_{h_2}^{th_2},\phi_d^{h_1h_2},\phi_r^{h_1h_2}\big)
\endaligned
\]
\[
\footnotesize
\aligned
\partial_{[3]}\colon\ \
\big(\psi_{i_1}^t,\psi_{i_2}^t,\psi_j^{h_1},\psi_j^{h_2}\big)
\longmapsto
&\,
\big(
\psi_j^{h_1}-2\psi_{i_2}^t,-6\psi_{i_1}^t,\psi_j^{h_2}+2\psi_{i_1}^t,
-6\psi_{i_2}^t,\psi_j^{h_2}-4\psi_{i_1}^t,
\\
&
\ \ \ \ \ \ \ \ \ \ \ \ \ \ \ \ \ \ \ \ \ \ \ \ \ \ \ \ \ \ \ \ \ \ \
\ \ \ \ \ \ \ \ \ \ \ \ \ \ \ \ \ \ \ \ \ \ \ \ \ \ \ \ \ \
-\psi_j^{h_1}-4\psi_{i_2}^t
\big)
\\
&\,
=
\big(\phi_d^{th_1},\phi_r^{th_1},\phi_d^{th_2},\phi_r^{th_2},
\phi_{i_1}^{h_1h_2},\phi_{i_2}^{h_1h_2}\big)
\endaligned
\]
\[
\footnotesize
\aligned
\partial_{[4]}\colon\ \
\big(\psi_j^t)
\longmapsto
&\,
\big(
-\psi_j^t,0,0,-\psi_j^t,-4\psi_j^t
\big)
\\
&\,
=
\big(\phi_{i_1}^{th_1},\phi_{i_2}^{th_1},\phi_{i_1}^{th_2},
\phi_{i_2}^{th_2},\phi_j^{h_1h_2}\big).
\endaligned
\]
One checks easily that the images of $\partial_{[2]}$, $\partial_{[3]}$,
$\partial_{[4]}$ satisfy the equations of $\mathcal{ Z}_{[2]}^2$, 
$\mathcal{
Z}_{[3]}^2$, $\mathcal{ Z}_{[4]}^2$ 
found above. Now, it is easy to view the
dimensions of the homogeneous components of the second cohomology
space:
\[
H^2(\mathfrak{g}_-,\mathfrak{g})
=
\bigoplus_{h\in\Z}\,
H_{[h]}^2(\mathfrak{g}_-,\mathfrak{g})
\ \ \ \ \
\text{\rm with}
\ \ \ \ \
H_{[h]}^2(\mathfrak{g}_-,\mathfrak{g})
:=
\frac{
\mathcal{Z}_{[h]}^2(\mathfrak{g}_-,\mathfrak{g})}{
\mathcal{B}_{[h]}^2(\mathfrak{g}_-,\mathfrak{g})}.
\]
Remind that $\mathcal{ Z}_{[0]}^2 \simeq \R^1$, $\mathcal{ Z}_{[1]}^2
\simeq \R^4$, $\mathcal{ Z}_{[2]}^2 \simeq \R^5$, $\mathcal{
Z}_{[3]}^2 \simeq \R^4$, $\mathcal{ Z}_{[4]}^2 \simeq \R^3$. Clearly,
$\partial_{[0]} \colon \R^5 \to \mathcal{ Z}_{[0]}^2 \simeq \R$ is
onto, whence $H_{[0]}^2 = \{ 0\}$. Similarly, one easily convinces
oneself with almost no computations that $\partial_{[1]}$ is of rank
$4$, that $\partial_{[2]}$ is of rank $5$, that $\partial_{[3]}$ is of
rank $4$ and that $\partial_{[5]}$ is of rank $1$. It follows that
$H_{[1]}^2 = \{ 0\}$, that $H_2^2 = \{ 0 \}$, that $H_{[3]}^2 = \{
0\}$, the only nonzero cohomology space being $H_{[4]}^2 \simeq \R^2$
which is $2$-dimensional. Finally, one also sees that $H^2(\frak
g_-,\frak g) = H_{[4]}^2$ is generated by the following two
independent $2$-cochains:
\[
\boxed{
\aligned
&
{\sf t}^*\wedge{\sf h}_2^*
\otimes
{\sf i}_2
-
2{\sf h}_1^*\wedge
{\sf h}_2^*
\otimes
{\sf j}
\\
\text{\rm and:}
\ \ \ \ \
&
{\sf t}^*\wedge{\sf h}_2^*
\otimes
{\sf i}_1
-
{\sf t}^*\wedge{\sf h}_1^*
\otimes
{\sf i}_2.
\endaligned}
\]
In conclusion, let us summarize the results obtained
by means of a dimensional table obtained in~\cite{ EMS}
that we recover here:
\begin{center}
\label{dimensional-cohomologies}
\begin{tabular}{|c|c|c|c|c|}
\hline\vspace{-10pt} &&&\\
\text{\rm Homogeneity} & $\dim\mathcal{C}^2$ & $\dim\mathcal{Z}^2$
& $\dim\mathcal{B}^2$ & $\dim H^2$
\\
\hline
0 & 1 & 1 & 1 & 0 \\
1 & 4 & 4 & 4 & 0 \\
2 & 6 & 5 & 5 & 0 \\
3 & 6 & 4 & 4 & 0 \\
4 & 5 & 3 & 1 & 2 \\
5 & 2 & 0 & 0 & 0 \\
\hline
\end{tabular}
\end{center}

\subsection{Codifferential}

When the Lie algebra $\frak g$ is semi-simple, there
exists another, 
degree-decreasing linear operator on the space of cochains:
\[
\partial^{*k}
\colon\ \ \
\mathcal{C}^{k+1}\big(\mathfrak{g}_-,\,\mathfrak{g}\big) \longrightarrow
\mathcal{C}^k\big(\mathfrak{g}_-,\,\mathfrak{g}\big),
\] 
called the {\em codifferential operator}, which is defined as
follows. For a Lie algebra $\frak g$ defined over a commutative field
$\mathbb K$, recall that ${\rm ad} \colon \frak g \rightarrow {\rm
End}_{\mathbb K}\frak g$ denotes the adjoint action of $\mathfrak{ g}$
on its space of endomorphisms:
\[
({\rm ad}({\sf x}))({\sf y})
:=
[{\sf x},{\sf y}]_{\mathfrak{g}}
\ \ \ \ \ \ \ \ \ \ \ \ \ 
{\scriptstyle{({\sf x},\,{\sf y}\,\in\,\mathfrak{g})}}.
\]
Recall also (\cite{ Knapp})
that the Killing form $B \colon \frak g \times\frak g\rightarrow
\mathbb K$ is defined as being the symmetric bilinear form:
\[
B({\sf x},{\sf y})
:=
{\rm Tr}\big({\rm ad}({\sf x})\circ{\rm ad}{(\sf y)}\big)),
\]
and that its nondegeneracy is equivalent to the semi-simplicity of
$\mathfrak{ g}$. Furthermore (\cite{ Cap}), if $\frak g = \frak
g_{-\mu} \oplus \cdots \oplus \frak g_\mu$ is a graded semi-simple Lie
algebra, then $B$ induces an isomorphism $\frak g_i^\ast\cong \frak
g_{-i}$ of $\frak g_0$-modules for $i = 1, \ldots,\mu$. If we denote
$\frak g_1 \oplus \cdots \oplus \frak g_\mu$ by $\frak g_+$, then each
space $\mathcal C^k (\frak g_-,\frak g)\cong \wedge^k \frak g_-^\ast
\otimes \frak g$ can be identified with the dual space of the space
$\wedge^k \frak g_+^\ast \otimes \frak g\cong \mathcal C^k( \frak g_+,
\frak g)$. In particular, the negative of the dual map of $ \partial^k
\colon\ \ \ \mathcal{ C}^k \big( \mathfrak{ g}_+,\, \mathfrak{ g}
\big) \longrightarrow \mathcal{ C}^{ k+1} \big( \mathfrak{
g}_+,\,\mathfrak{ g}\big)$ can be viewed as a linear map which is
exactly the codifferential operator $\partial^{*k} \colon\ \ \
\mathcal{ C}^{ k+1}\big(\mathfrak{g}_-,\,\mathfrak{g}\big)
\longrightarrow \mathcal{ C}^k \big(\mathfrak{ g}_-,\,\mathfrak{
g}\big)$. From this definition, it immediately follows that
$\partial^{*(k-1)} \circ \partial^{*k} = 0$, whence one has a second
{\em cochain complex}:
\[
0 \overset{\partial^{*n}}{\longrightarrow} 
\mathcal{C}^n
\overset{\partial^{*(n-1)}}{\longrightarrow}
\mathcal{C}^{n-1} 
\overset{\partial^{*(n-2)}}{\longrightarrow} 
\cdots
\overset{\partial^{*2}}{\longrightarrow} 
\mathcal{C}^2
\overset{\partial^{*1}}{\longrightarrow} 
\mathcal{C}^1
\overset{\partial^{*0}}{\longrightarrow} 0.
\]
Lastly (\cite{Sato}), for any $k+1$ elements ${\sf
z}_1, \ldots, {\sf z}_k$ of $\frak g_-$ and for any $(k+1)$-cochain
$\Psi \in \mathcal C^{k+1} (\frak g_-, \frak g)$, the expression of
$\partial^{ *k} \Psi$ realizes as follows:
\[
\small
\aligned
&
(\partial^{*k}\Psi)
\big({\sf z}_1,\dots,{\sf z}_k\big)
:=
\sum_{i=1}^n\big[{\sf v}^\ast_i,\Psi\big({\sf v}_i,{\sf z}_1,
\dots,{\sf z}_k\big)\big]_{\mathfrak{g}}
+
\\
&
\ \ \ \ \ \
+
\frac{1}{2}\,
\sum_{i=1}^n\,
\sum_{j=1}^k\,
(-1)^{j+1}\,
\Psi
\Big({\sf proj}_{\mathfrak{g}_-}\big(
\big[{\sf v}^\ast_i,{\sf z}_j\big]_{\mathfrak{g}}\big),\,\,
{\sf v}_i,\,{\sf z}_1,
\dots,
\widehat{{\sf z}}_j,\ldots,{\sf z}_k
\Big),
\endaligned
\]
where ${\sf v}_1,\ldots,{\sf v}_n$ are independent basis elements of
$\frak g_-$, where ${\sf v}^\ast_i$ ($i = 1, \dots, n$) is the dual of
${\sf v}_i$ with respect to the Killing form, where ${\sf
proj}_{\mathfrak{g}_-} \big( \big[{\sf v}^\ast_i,{\sf
z}_j\big]_{\mathfrak{g}}\big)$ denotes the $\frak g_-$-component of
$\big[{\sf v}^\ast_i,{\sf z}_j\big]_{ \mathfrak{ g}}$ with respect to
the decomposition $\mathfrak{ g} = \mathfrak{ g}_- \oplus \mathfrak{
p}$, and where $\mathfrak{ p} := \mathfrak{ g}_0 \oplus \mathfrak{ g}_+$.

\section{Initial Frame on a Strongly Pseudoconvex $M^3 \subset \C^2$}
\label{initial-frame}

\HEAD{Initial Frame on a Strongly Pseudoconvex $M^3 \subset \C^2$}{
Mansour Aghasi, Joël Merker, and Masoud Sabzevari}

\subsection{Explicit CR structure}
\label{Explicit-Cr-structure}
Let $M$ be a real $\mathcal{ C}^1$-smooth 
hypersurface of $\mathbb C^2$, represented by:
\[
v 
= 
\varphi(x,y,u)
\]
in coordinates $(z, w) = (x+iy, u+iv)$. After a linear straightening,
we may assume $0 \in M$ and $T_0 M = \{ {\rm Im}\, w = 0\}$, so that
$\varphi ( 0) = \varphi_x ( 0) = \varphi_y ( 0) = \varphi_u ( 0) = 0$.
A $(0, 1)$ vector field of the form:
\[
\overline{\mathcal{L}}
= 
\frac{\partial}{\partial\overline{z}} 
+
{\tt A}\,\frac{\partial}{\partial\overline{w}}
\]
is tangent to $M$ if and only if its coefficient ${\tt A}$ satisfies:
\[
0 = \frac{{\tt A}}{2i} + \varphi_{\overline{z}} 
+
\frac{{\tt A}}{2}\,\varphi_u,
\]
or equivalently:
\[
\boxed{{\tt A} 
= 
\frac{2\,\varphi_{\overline{z}}}{i-\varphi_u}}\,.
\]
Consequently the vector field:
\[
\overline{\mathcal{L}} = \frac{1}{2}\,\frac{\partial}{\partial x} +
\frac{i}{2}\,\frac{\partial}{\partial y} + \bigg(
\frac{2\,\varphi_{\overline{z}}}{i-\varphi_u}
\bigg) \bigg( \frac{1}{2}\,\frac{\partial}{\partial u} +
\frac{i}{2}\,\frac{\partial}{\partial v}
\bigg)
\]
generates $T^{ 0, 1}M$ in a neighborhood of the origin, since $T^{ 0,
1}M$ is obviously 
of rank $\dim \mathbb C^2 - {\rm CRdim}\, M = 2 - 1 = 1$.

We notice that this $\overline{ \mathcal{ L}}$ is written here {\em
extrinsically}, namely it involves the extra coordinate $v$ and it
lives in a neighborhood of $M$, in $\mathbb C^2$, while $M$ itself,
which is three-dimensional, is naturally equipped with the three real
coordinates $(x, y, u)$. Since we want two
intrinsic sections of:
\[
T^cM 
= 
{\rm Re}\,\big(T^{0,1}M\big),
\]
we need at first to pullback this $\overline{ \mathcal{ L}}$ to $M$,
which simply means dropping the basic vector field $\frac{\partial}{
\partial v}$ and replacing $v$ by $\varphi ( x, y, u)$ in the
coefficient functions (in fact here, no $v$ appears), and we get the
following section:
\[
\overline{\mathcal{L}}\big\vert_M 
=
\frac{1}{2}\,\frac{\partial}{\partial x} +
\frac{i}{2}\,\frac{\partial}{\partial y} + \bigg(
\frac{2\,\varphi_{\overline{z}}}{i-\varphi_u}
\bigg) \bigg( \frac{1}{2}\,\frac{\partial}{\partial u} \bigg)
\]
which generates $T^{ 0, 1}M$, {\em intrinsically} ({\em see} also the
basic first chapters of~\cite{BER, Boggess, Jacobowitz}).

So it remains only to decompose $\overline{ \mathcal{ L}} \big\vert_M$
in real and imaginary parts, and at first, we do this for the
coefficient:
\[
{\tt A} 
= 
\frac{\varphi_x+i\,\varphi_y}{i-\varphi_u}
=
\frac{\varphi_y-\varphi_x\,\varphi_u}{1+\varphi_u^2} 
+ 
i\,\frac{-\varphi_x-\varphi_y\,\varphi_u}{1+\varphi_u^2}.
\]
Hence we can provide an explicit representation of two independent
real vector fields that are generators for $T^cM$ near the origin,
namely $2 \, {\rm Re}\, \big( \overline{ \mathcal{ L}}
\big\vert_M\big)$ and $2 \, {\rm Im}\, \big( \overline{ \mathcal{ L}}
\big\vert_M\big)$, multiplying by a factor $2$ to simplify a bit.

\begin{Lemma}
For any local $\mathcal{ C}^1$-smooth
real hypersurface $M^3$
of $\mathbb C^2$ which is represented as a graph:
\[
v 
= 
\varphi(x,y,u)
\]
in coordinates $(z, w) = ( x + iy, \, u + iv)$, the complex tangent
bundle $T^cM = {\rm Re}\, T^{ 0, 1}M$ is generated by the two explicit
vector fields:
\[
\left\{ 
\aligned 
H_1 
& 
:= 
\frac{\partial}{\partial x} 
+ 
\bigg(\frac{\varphi_y-\varphi_x\,\varphi_u}{1+\varphi_u^2}\bigg)
\frac{\partial}{\partial u}
\\
H_2
&
:= 
\frac{\partial}{\partial y} + \bigg(
\frac{-\varphi_x-\varphi_y\,\varphi_u}{1+\varphi_u^2}
\bigg) \frac{\partial}{\partial u}.
\endaligned\right.
\]
\end{Lemma}

In fact, as one easily verifies, one does not need that $\varphi ( 0)
= \varphi_x ( 0) = \varphi_y ( 0) = \varphi_u ( 0) = 0$ for the lemma
to hold true (but we will always assume that such an affine
normalization is done in advance, since it is free).

Some further notation will be useful. If we set:
\[
\Delta
:=
1+\varphi_u^2,
\ \ \ \ \ \ \ \ \ \
\Lambda_1
:=
\varphi_y-\varphi_x\,\varphi_u,
\ \ \ \ \ \ \ \ \ \
\Lambda_2
:=
-\,\varphi_x-\varphi_y\,\varphi_u,
\]
our two intrinsic $T^cM$-tangent vector fields become:
\[
H_1
=
\frac{\partial}{\partial x}
+
\frac{\Lambda_1}{\Delta}\,
\frac{\partial}{\partial u}
\ \ \ \ \ \ \ \ \ \ 
\text{\rm and}
\ \ \ \ \ \ \ \ \ \ 
H_2
=
\frac{\partial}{\partial y}
+
\frac{\Lambda_2}{\Delta}\,
\frac{\partial}{\partial u}.
\]

\subsection{Levi nondegeneracy assumption}
\label{Levi-nondegeneracy}

Now, we assume that $M$ is Levi nondegenerate at the origin, so that
second order terms can be assumed to be normalized as:
\[
v 
= 
\varphi(x,y,u) 
= 
x^2+y^2+{\rm O}(3).
\]
We may therefore compute the bracket $[ H_1, \, H_2]$ using these
notations, and realize that two terms underlined cancel:
\begin{equation*}
\small 
\aligned 
{} 
\big[H_1,\,H_2\big]
& 
= 
\Big[{\textstyle{\frac{\partial}{\partial x}}} 
+
\big({\textstyle{\frac{\Lambda_1}{\Delta}}}\big)
{\textstyle{\frac{\partial}{\partial u}}}, 
\ \
{\textstyle{\frac{\partial}{\partial y}}} 
+
\big({\textstyle{\frac{\Lambda_2}{\Delta}}}\big)
{\textstyle{\frac{\partial}{\partial u}}}\Big]
\\
& 
= 
\Big( {\textstyle{\frac{\Lambda_{2,x}}{\Delta}}} 
- 
\Lambda_2\,{\textstyle{\frac{\Delta_x}{\Delta^2}}} 
+ 
{\textstyle{\frac{\Lambda_1}{\Delta}}}
{\textstyle{\frac{\Lambda_{2,u}}{\Delta}}} 
-
\zero{\textstyle{\frac{\Lambda_1}{\Delta}}\Lambda_2
{\textstyle{\frac{\Delta_u}{\Delta^2}}}}
-
\\
&
\ \ \ \ \ \ \ 
-\,
{\textstyle{\frac{\Lambda_{1,y}}{\Delta}}} +
\Lambda_1 {\textstyle{\frac{\Delta_y}{\Delta^2}}} -
{\textstyle{\frac{\Lambda_2}{\Delta}}}
{\textstyle{\frac{\Lambda_{1,u}}{\Delta}}} +
\zero{\textstyle{\frac{\Lambda_2}{\Delta}}\Lambda_1
{\textstyle{\frac{\Delta_u}{\Delta^2}}}}\Big)\,
{\textstyle{\frac{\partial}{\partial u}}},
\endaligned
\end{equation*}
so that the common denominator is {\em not} equal to $\Delta^3$ as one
would have expected, but is equal to $\Delta^2$. Expanding the partial
derivatives and simplifying either by hand or with a computer
(\cite{AMSMaple}), we therefore
get:
\[
\boxed{ 
\small 
\aligned 
{} 
\big[H_1,\,H_2\big] 
& 
= 
\Big[
{\textstyle{\frac{\partial}{\partial x}}} 
+
{\textstyle{\frac{\varphi_y-\varphi_x\,\varphi_u}{1+\varphi_u^2}}}\,
{\textstyle{\frac{\partial}{\partial u}}}, 
\ \ 
{\textstyle{\frac{\partial}{\partial y}}} 
+
{\textstyle{\frac{-\varphi_x+\varphi_y\,\varphi_u}{1+\varphi_u^2}}}\,
{\textstyle{\frac{\partial}{\partial u}}}\Big]
\\
& 
= 
\Big( {\textstyle{\frac{1}{(1+\varphi_u^2)^2}}} 
\big\{
-\varphi_{xx}-\varphi_{yy} -
2\,\varphi_y\,\varphi_{xu}-\varphi_x^2\,\varphi_{uu}
+
2\,\varphi_x\,\varphi_{yu}-\varphi_y^2\,\varphi_{uu} 
+
\\
& \ \ \ \ \ \ \ \ \ \ \ \ \ \ \ \ \ \ \ \ \ \ \ \ +
2\,\varphi_y\,\varphi_u\,\varphi_{yu} +
2\,\varphi_x\,\varphi_u\,\varphi_{xu} - \varphi_u^2\,\varphi_{xx} -
\varphi_u^2\,\varphi_{yy}
\big\} \Big)\, \frac{\partial}{\partial u}.
\endaligned}
\]
Equivalently, as we want for later use to specify the numerator, 
if we set:
\[
\aligned
\Upsilon
&
:=
-\varphi_{xx}-\varphi_{yy} -
2\,\varphi_y\,\varphi_{xu}-\varphi_x^2\,\varphi_{uu}
+
2\,\varphi_x\,\varphi_{yu}-\varphi_y^2\,\varphi_{uu} 
+
\\
& \ \ \ \ \ 
+
2\,\varphi_y\,\varphi_u\,\varphi_{yu} +
2\,\varphi_x\,\varphi_u\,\varphi_{xu} - \varphi_u^2\,\varphi_{xx} -
\varphi_u^2\,\varphi_{yy},
\endaligned
\]
we can write shortly:
\[
\big[H_1,\,H_2\big]
= 
\frac{\Upsilon}{\Delta^2}\,
\frac{\partial}{\partial u}.
\]
Now, thanks to the Levi-nondegeneracy assumption and because of the
normalizations $0 = \varphi ( 0) = \varphi_x ( 0) = \varphi_y ( 0) =
\varphi_u ( 0)$, we have $\Upsilon ( 0) = - 4$ (notice the minus sign),
that is to say:
\[
[H_1,H_2] \big\vert_0 
= 
-\,4\, 
{\textstyle{\frac{\partial}{\partial u}}} \big\vert_0.
\]
So, if we introduce the vector field (we choose a 
plus sign in the factor $\frac{ 1}{ 4}$):
\[
T 
:= 
\frac{1}{4}\,\frac{\Upsilon}{\Delta^2}\, 
\frac{\partial}{\partial u},
\]
we may rewrite:
\[
\boxed{ [H_1,\,H_2] = 4\,T}\,.
\]

\subsection{Length-three brackets}
\label{length-three-brackets}
At the next step, we must compute the two brackets $4\, 
[ H_1, T]$ and $4\, [
H_2, T]$, for instance the first one, in which we see how the 
denominator grows and of which we extract the numerator:
\[
\aligned
{} 
4\,\big[H_1,\,T\big] 
& 
=
\Big[{\textstyle{\frac{\partial}{\partial x}}} 
+
\big({\textstyle{\frac{\Lambda_1}{\Delta}}}\big)
{\textstyle{\frac{\partial}{\partial u}}},\ \ 
\big({\textstyle{\frac{\Upsilon}{\Delta^2}}}\big)
{\textstyle{\frac{\partial}{\partial u}}}
\Big]
\\
& 
= 
\Big( {\textstyle{\frac{\Upsilon_x}{\Delta^2}}} 
- 
2\,\Upsilon\,
{\textstyle{\frac{\Delta_x}{\Delta^3}}} 
+ 
{\textstyle{\frac{\Lambda_1}{\Delta}}}\,
{\textstyle{\frac{\Upsilon_u}{\Delta^2}}} 
- 
2\,{\textstyle{\frac{\Lambda_1}{\Delta}}}\,
\Upsilon\, {\textstyle{\frac{\Delta_u}{\Delta^3}}} 
-
{\textstyle{\frac{\Upsilon}{\Delta^2}}}\,
{\textstyle{\frac{\Lambda_{1,u}}{\Delta}}} 
+
{\textstyle{\frac{\Upsilon}{\Delta^2}}}\,\Lambda_1\,
{\textstyle{\frac{\Delta_u}{\Delta^2}}}\Big)
\frac{\partial}{\partial u}
\\
& 
= 
\bigg(
\frac{
\Delta^2\big[\Upsilon_x\big]
+ 
\Delta\big[-2\,\Upsilon\,\Delta_x+\Lambda_1\,\Upsilon_u 
-
\Upsilon\,\Lambda_{1,u}\big] -
\Lambda_1\,\Upsilon\,\Delta_u
}{\Delta^4}
\bigg)
\frac{\partial}{\partial u}.
\endaligned
\]
Exchanging $H_1$ with $H_2$, which means replacing $\Lambda_1$ by
$\Lambda_2$ and $\frac{ \partial}{ \partial x}$ by $\frac{ \partial }{
\partial y}$, we get similarly and without any computation:
\[
4\,\big[H_1,\,T\big]
=
\bigg(
\frac{
\Delta^2\big[\Upsilon_y\big]
+ 
\Delta\big[-2\,\Upsilon\,\Delta_y+\Lambda_2\,\Upsilon_u 
-
\Upsilon\,\Lambda_{2,u}\big] -
\Lambda_2\,\Upsilon\,\Delta_u
}{\Delta^4}
\bigg)
\frac{\partial}{\partial u}.
\] 
Let us therefore introduce two new summarizing names:
\[
\aligned
A_1
&
:=
\Delta^2\big[\Upsilon_x\big]
+ 
\Delta\big[-2\,\Upsilon\,\Delta_x+\Lambda_1\,\Upsilon_u 
-
\Upsilon\,\Lambda_{1,u}\big] -
\Lambda_1\,\Upsilon\,\Delta_u,
\\
A_2
&
:=
\Delta^2\big[\Upsilon_y\big]
+ 
\Delta\big[-2\,\Upsilon\,\Delta_y+\Lambda_2\,\Upsilon_u 
-
\Upsilon\,\Lambda_{2,u}\big] -
\Lambda_2\,\Upsilon\,\Delta_u,
\endaligned
\]
for the two appearing numerators. Now, for later use, 
we want to re-express these two 
brackets $[ H_1, T]$ and $[ H_2, T]$ 
in terms of the third field $T$ transverse
to $T^cM$, and for this, it suffices to 
simply replace the basic field:
\[
\frac{\partial}{\partial u}
=
\frac{4\,\Delta^2}{\Upsilon}\,\,T,
\]
so that doing this just yields expressions of the two
supplementary brackets
\[
\aligned
\big[H_1,\,T\big]
&
=
\frac{1}{4}\,
\frac{A_1}{\Delta^4}\,
\frac{4\,\Delta^2}{\Upsilon}\,\,T
=
\frac{A_1}{\Delta^2\,\Upsilon}\,\,T,
\\
\big[H_2,\,T\big]
&
=
\frac{1}{4}\,
\frac{A_2}{\Delta^4}\,
\frac{4\,\Delta^2}{\Upsilon}\,\,T
=
\frac{A_2}{\Delta^2\,\Upsilon}\,\,T.
\endaligned
\] 
However, these two numerators $A_1$ and $A_2$ are not yet expanded
as explicit polynomials in the third-order jet $J_{ x, y, u}^3
\varphi$ of the graphing function $\varphi ( x, y, u)$ for $M$. This
can be done either by hand or using a computer (\cite{AMSMaple}),
hence we directly summarize the fundamental result fully describing a
useful initial frame for $TM$ which is naturally produced by $T^cM$.

\begin{Proposition}
\label{explicit-J-3-varphi}
If $M^3$ is an arbitrary local 
$\mathcal{ C}^3$-smooth Levi nondegenerate real
hypersurface of $\mathbb C^2$ represented in coordinates $(z, w) = (x
+ iy, \, u + iv)$ as a graph:
\[
v 
= 
\varphi(x,y,u) 
= 
x^2+y^2+{\rm O}(3),
\]
and whose complex tangent bundle $T^cM$ is generated by the two
explicit vector fields:
\[
H_1 := {\textstyle{\frac{\partial}{\partial x}}} + \big(
{\textstyle{\frac{\varphi_y-\varphi_x\,\varphi_u}{1+\varphi_u^2}}} \big)
{\textstyle{\frac{\partial}{\partial u}}} 
\ \ \ \ \ \ \ \ \ \ \text{\em and} \ \ \ \ \ \ \ \ \ \
H_2 := {\textstyle{\frac{\partial}{\partial y}}} 
+ 
\big({\textstyle{\frac{-\varphi_x-\varphi_y\,\varphi_u}{1+\varphi_u^2}}} 
\big)
{\textstyle{\frac{\partial}{\partial u}}},
\]
satisfying $H_1 \vert_0 = \frac{ \partial}{ \partial x} \big\vert_0$
and $H_2 \vert_0 = \frac{ \partial}{ \partial y} \big\vert_0$, then
the third, bracketed vector field:
\[
\aligned 
T 
:= 
&\, 
{\textstyle{\frac{1}{4}}}\, [H_1,H_2]
\\
= 
&\,
\Big({\textstyle{\frac{1}{4}}}\,
{\textstyle{\frac{1}{(1+\varphi_u^2)^2}}}
\big\{
-
\varphi_{xx}-\varphi_{yy}-2\,\varphi_y\,\varphi_{xu}
- 
\varphi_x^2\,\varphi_{uu} 
+
2\,\varphi_x\,\varphi_{yu} 
- 
\varphi_y^2\,\varphi_{uu} 
+
\\
& \ \ \ \ \ \ \ \ \ \ \ \ \ \ \ \ \ \ \ \ \ \ \ \ +
2\,\varphi_y\,\varphi_u\,\varphi_{yu} +
2\,\varphi_x\,\varphi_u\,\varphi_{xu} - \varphi_u^2\,\varphi_{xx} -
\varphi_u^2\,\varphi_{yy}
\big\} \Big)\, \frac{\partial}{\partial u}
\\
=: 
&\, 
\Big( {\textstyle{\frac{1}{4}}}\, 
{\textstyle{\frac{1}{(\Delta)^2}}}\,
\big\{\Upsilon\big\}\, \Big) 
\frac{\partial}{\partial u}
\endaligned
\]
satisfying $T \big\vert_0 = - \frac{ \partial}{ \partial u}
\big\vert_0$ produces, jointly with $H_1$ and $H_2$ of which it is
locally linearly independent, a frame for $TM$ in a neighborhood of
the origin. Furthermore, the remaining Lie bracket structure of this
frame:
\[
\boxed{
[H_1,T]
=
\Phi_1\,T 
\ \ \ \ \ \ 
\text{\rm and} 
\ \ \ \ \ \ 
[H_2,T]
=
\Phi_2\,T}\,,
\]
involves two further rational functions:
\[
\Phi_1 
= 
\frac{A_1}{\Delta^2\,\Upsilon}
\ \ \ \ \ \
\text{\rm and}
\ \ \ \ \ \ 
\Phi_2 
= 
\frac{A_2}{\Delta^2\,\Upsilon}
\]
having common denominator equal to $\Delta^2\, \Upsilon$ and whose two
numerators $A_1$ and $A_2$, depending both upon the third-order jet
$J_{ x, y, u}^3 \varphi$, read explicitly as follows:
\[
\footnotesize
\aligned A_1 
& 
= 
-\varphi_{xxx} -\varphi_{xyy} +2\,\varphi_x\,\varphi_{xyu}
-3\,\varphi_y\,\varphi_{xxu} -\varphi_y\,\varphi_{yyu} -3\,\varphi_y^2\,\varphi_{xuu}
+
\\
& \ \ \ \ 
+
2\,\varphi_x\,\varphi_y\,\varphi_{yuu} 
-\varphi_x^2\,\varphi_{xuu} -\varphi_x^2\,\varphi_y\,\varphi_{uuu}
-\varphi_y^3\,\varphi_{uuu} 
-
2\,\varphi_x\,\varphi_y\,\varphi_{xu}\,\varphi_{uu}
-
\\
& \ \ \ \ \ 
-
3\,\varphi_x\varphi_{xx}\,\varphi_{uu} 
+\varphi_y^2\,\varphi_{uu}\,\varphi_{yu}
-2\,\varphi_y\,\varphi_{uu}\,\varphi_{xy}
+3\,\varphi_x^2\,\varphi_{yu}\,\varphi_{uu}
-
\\
&
\ \ \ \ \
-\varphi_x\,\varphi_{yy}\,\varphi_{uu}
-\varphi_x\,\varphi_y^2\,\varphi_{uu}^2 
+
4\,\varphi_y\,\varphi_{xu}\,\varphi_{yu} -\varphi_x^3\,\varphi_{uu}^2
+\,\varphi_{yu}\,\varphi_{yy} 
+
\\
& \ \ \ \ 
+
3\,\varphi_{xx}\,\varphi_{yu}
-2\,\varphi_{xu}\,\varphi_{xy} +
2\,\varphi_x\,\varphi_{xu}^2 -2\,\varphi_x\,\varphi_{yu}^2 +
\endaligned
\]
\[
\footnotesize
\aligned & + \varphi_u \Big( 3\,\varphi_x\,\varphi_{xxu} +2\,\varphi_y\,\varphi_{xyu}
+\varphi_x\,\varphi_{yyu} +4\,\varphi_x\,\varphi_y\,\varphi_{xuu}
+2\,\varphi_y^2\,\varphi_{yuu} -
\\
& \ \ \ \ \ \ \ \ \ \ -2\,\varphi_x^2\,\varphi_{yuu}
+\varphi_x\,\varphi_y^2\,\varphi_{uuu}
+\varphi_x^2\,\varphi_{uuu} +2\,\varphi_x^2\,\varphi_{uu}^2\,\varphi_y
+5\,\varphi_{uu}\,\varphi_{xu}\,\varphi_x^2 -
\\
& \ \ \ \ \ \ \ \ \ \ -8\,\varphi_x\,\varphi_{xu}\,\varphi_{yu}
+7\,\varphi_y^2\,\varphi_{xu}\,\varphi_{uu} +\varphi_{yy}\,\varphi_{xu}
+2\,\varphi_y^3\,\varphi_{uu}^2 +3\,\varphi_{xx}\,\varphi_{xu} +
\\
& \ \ \ \ \ \ \ \ \ \ +8\,\varphi_y\,\varphi_{xu}^2 +2\,\varphi_{xy}\,\varphi_{yu}
-2\,\varphi_x\,\varphi_y\,\varphi_{yu}\,\varphi_{uu} \Big) +
\endaligned
\]
\[
\footnotesize
\aligned & + \varphi_u^2 \Big( -3\,\varphi_{xxx} -3\,\varphi_{xyy}
-6\,\varphi_y\,\varphi_{xxu}
-2\,\varphi_y\,\varphi_{yyu} +4\,\varphi_x\,\varphi_{xyu}
-4\,\varphi_y^2\,\varphi_{xuu} -
\\
& \ \ \ \ \ \ \ \ \ \ -4\,\varphi_x^2\,\varphi_{xuu} -\varphi_y^3\,\varphi_{uuu}
-\varphi_y\,\varphi_x^2\,\varphi_{uuu} -2\,\varphi_x\,\varphi_{uu}\,\varphi_{yy}
+7\,\varphi_x^2\,\varphi_{yu}\,\varphi_{uu} -
\\
& \ \ \ \ \ \ \ \ \ \ -6\,\varphi_x\,\varphi_{uu}\,\varphi_{xx}
-4\,\varphi_y\,\varphi_{uu}\,\varphi_{xy} -3\,\varphi_y^2\,\varphi_{uu}^2\,\varphi_x
-3\,\varphi_y^2\,\varphi_{uu}\,\varphi_{yu} -4\,\varphi_{xu}\,\varphi_{xy} -
\\
& \ \ \ \ \ \ \ \ \ \ -3\,\varphi_x^3\,\varphi_{uu}^2 +6\,\varphi_{xx}\,\varphi_{yu}
-4\,\varphi_x\,\varphi_{yu}^2 -4\,\varphi_x\,\varphi_{xu}^2
+2\,\varphi_{yu}\,\varphi_{yy}
-10\,\varphi_x\,\varphi_y\,\varphi_{xu}\,\varphi_{uu} \Big) +
\endaligned
\]
\[
\footnotesize
\aligned & + \varphi_u^3 \Big( 6\,\varphi_x\,\varphi_{xxu} +4\,\varphi_y\,\varphi_{xyu}
+2\,\varphi_x\,\varphi_{yyu} +4\,\varphi_x\,\varphi_y\,\varphi_{xuu}
-2\,\varphi_x^2\,\varphi_{yuu}
+2\,\varphi_y^2\,\varphi_{yuu} +
\\
& \ \ \ \ \ \ \ \ \ \ +\varphi_x^3\,\varphi_{uuu}
+\varphi_x\,\varphi_y^2\,\varphi_{uuu}
+3\,\varphi_y^2\,\varphi_{xu}\,\varphi_{uu} -8\,\varphi_{xu}\,\varphi_{yu}\,\varphi_x
+9\,\varphi_{uu}\,\varphi_{xu}\,\varphi_x^2 +
\\
& \ \ \ \ \ \ \ \ \ \ +4\,\varphi_{xy}\,\varphi_{yu} +8\,\varphi_y\,\varphi_{xy}^2
+2\,\varphi_{yy}\,\varphi_{xu} +6\,\varphi_{xx}\,\varphi_{xu}
+6\,\varphi_x\,\varphi_y\,\varphi_{yu}\,\varphi_{uu} \Big) +
\endaligned
\]
\[
\footnotesize
\aligned & + \varphi_u^4 \Big( -3\,\varphi_{xxx} -3\,\varphi_{xyy}
+2\,\varphi_x\,\varphi_{xyu}
-\varphi_y\,\varphi_{yyu} -3\,\varphi_y\,\varphi_{xxu}
-3\,\varphi_x^2\,\varphi_{xuuu} -
\\
& \ \ \ \ \ \ \ \ \ \ -2\,\varphi_x\,\varphi_y\,\varphi_{yuu}
-\varphi_y^2\,\varphi_{xuu}
-3\,\varphi_x\,\varphi_{uu}\,\varphi_{xx} -\varphi_x\,\varphi_{uu}\,\varphi_{yy}
-6\,\varphi_x\,\varphi_{xu}^2 -
\\
& \ \ \ \ \ \ \ \ \ \ -2\,\varphi_x\,\varphi_{yu}^2
-2\,\varphi_y\,\varphi_{uu}\,\varphi_{xy}
-4\,\varphi_y\,\varphi_{xu}\,\varphi_{yu} -2\,\varphi_{xu}\,\varphi_{xy}
+3\,\varphi_{xx}\,\varphi_{yu} +\varphi_{yu}\,\varphi_{yy} \Big) +
\endaligned
\]
\[
\footnotesize
\aligned & + \varphi_u^5 \Big( \varphi_x\,\varphi_{yyu} +2\,\varphi_y\,\varphi_{xyu}
+3\,\varphi_x\,\varphi_{xxu} +\varphi_{yy}\,\varphi_{xu} +3\,\varphi_{xx}\,\varphi_{xu}
+2\,\varphi_{xy}\,\varphi_{yu} \Big) +
\ \ \ \ \ \ \ \ \ \ \
\\
\ \ \ \ \ \ \ \ \ & + \varphi_u^6 \Big( -\varphi_{xxx} -\varphi_{xyy} \Big).
\endaligned
\]
while $A_2$ is obtained from $A_1$ by just exchanging $x$ and $y$.
\end{Proposition}

\subsection{Abstract shape of the initial frame on $M^3 \subset \mathbb C^2$}
\label{Abstract-shape}

From now on, we shall restart from the beginning by assuming that we
are given an initial frame $(H_1, H_2, T)$ for $TM$ made of certain
two linearly independent vector fields which generate $T^cM$ locally:
\[
H_1\in\Gamma(T^cM) 
\ \ \ \ \ \ \ 
\text{\rm and} 
\ \ \ \ \ \ \ 
H_2\in\Gamma(T^cM),
\]
together with their bracket:
\[
T 
:= 
{\textstyle{\frac{1}{4}}}\, [H_1,H_2] \in\Gamma(TM)
\]
enjoying the following commutator relations:
\[
[H_1,T]
=
\Phi_1\,T 
\ \ \ \ \ \ 
\text{\rm and} 
\ \ \ \ \ \ 
[H_2,T]=\Phi_2\,T,
\]
were we now consider the two $\mathcal{ C}^\infty$ functions $\Phi_1
\colon M \to \mathbb R$ and $\Phi_2 \colon M \to \mathbb R$ as basic
data, without it to be necessary to know that they both depend
explicitly on some local graphing function $\varphi$ for $M$, as was
stated by the preceding proposition. In the subsequent section, we
shall construct a Cartan connection just in terms of $\Phi_1$ and
$\Phi_2$, not trying to express the newly constructed functions and
curvatures explicitly in terms of the graphing function $\varphi$, for
the sizes of such expressions might explode dramatically. Only at the
very end, after all the computations in terms of just $\Phi_1$ and
$\Phi_2$ are finalized, will we give it to a computer to expand the
gained curvatures in terms of the sixth-order jet $J_{ x, y, u}^6
\varphi$.

At least at the moment, it is useful two explore in advance what
relations come out when one takes iterated brackets:
\[
\small
\aligned
\big[H_i,\,T\big]
&
=
{\textstyle{\frac{1}{4}}}\,
\big[H_i,
[H_1,H_2]\big]
\\
\big[H_i,\big[H_j,\,T\big]\big]
&
=
{\textstyle{\frac{1}{4}}}\,
\big[H_i,\big[H_j,\,
[H_1,H_2]\big]\big]
\\
\big[H_i,\big[H_j,\big[H_k,\,T\big]\big]\big]
&
=
{\textstyle{\frac{1}{4}}}\,
\big[H_i,\big[H_j,\big[H_k,\,
[H_1,H_2]\big]\big]\big]
\\
\big[H_i,\big[H_j,\big[H_k,\big[H_l,\,T\big]\big]\big]\big]
&
=
{\textstyle{\frac{1}{4}}}\,
\big[H_i,\big[H_j,\big[H_k,\big[H_l,\,
[H_1,H_2]\big]\big]\big]\big]
\endaligned
\] 
up to length $\leqslant 6$, where $i, j, k, l = 1, 2$. A first
observation is as follows, but a more systematic exploration of higher
order relations will be achieved in the next section.

\begin{Lemma}
\label{H1-Phi2-H2-Phi1}
The two functions $H_1(\Phi_2)$ and $H_2(\Phi_1)$ are equal.
\end{Lemma}

\proof 
By what has been seen at the moment, we have by definition:
\[
[H_1,T]
=
\Phi_1\,T, 
\ \ \ \ \ \ \ \ \ \ \ \
[H_2,T]
=
\Phi_2\,T,
\]
whence, by bracketing the second 
(resp. first) equation with $[ H_1, \cdot ]$
(resp. $[ H_2, \cdot]$):
\[
\aligned
\big[H_1,[H_2,T]\big]
&
=
\big[H_1,\Phi_2\,T\big]
=
H_1(\Phi_2)\,T
+
\Phi_2\Phi_1\,T,
\\
\big[H_2,[H_1,T]]
&
=
\big[H_2,\Phi_1\,T\big]
=
H_2(\Phi_1)\,T
+
\Phi_1\Phi_2\,T.
\endaligned
\]
On the other hand, the Jacobi identity enables us to realize that these
two iterated brackets of length 3 are in fact equal:
\[
\big[H_1,[H_2,T]]-[H_2,[H_1,T]\big]
=
-
\big[T,[H_1,H_2]\big]
=
\big[T,\,4\,T\big]
=
0,
\]
so that we deduce at once:
\[
H_1(\Phi_2)=H_2(\Phi_1),
\]
as was claimed.
\endproof

\section{Free Lie Algebras of Rank Two 
\\
and Relations Between Brackets of Length $\leqslant 6$}
\label{free-Lie-algebras}

\HEAD{Free Lie Algebras of Rank Two
and Relations Between Brackets of Length $\leqslant 6$}{
Mansour Aghasi, Joël Merker, and Masoud Sabzevari}

\subsection{Free Lie algebras of rank two}
To reach higher order relations, one must at first count the maximal
number of iterated Lie brackets between $H_1$ and $H_2$ which are
linearly independent modulo skew-symmetry and Jacobi identity, just
abstractly, without using $[ H_i, \, T] = \Phi_i\, T$. For this, one
calls to the concept of free Lie algebra of rank $2$, {\em
cf.} the
reference~\cite{ MerkerPorten}, 
pp.~9--11 of which we borrow the notations.

Let $h_1, h_2$ be two linearly independent elements of a certain
vector space over $\R$. The {\sl free Lie algebra} $\mathcal{ F}$ of
rank $2$ is the smallest (non-commutative, non-associative)
$\R$-algebra having $h_1, h_2$ as elements, with bilinear
multiplication: 
\[
(h,h') 
\longmapsto 
[h, h']
\in
\mathcal{F} 
\ \ \ \ \ \ \ \ \ \ \ \ \ 
{\scriptstyle{(h,\,h'\,\in\,\mathcal{F})}}
\]
satisfying skew-symmetry: 
\[
0 
=
[h,\,h']
+ 
[h',\,h]
\ \ \ \ \ \ \ \ \ \ \ \ \ 
{\scriptstyle{(h,\,h'\,\in\,\mathcal{F})}}
\] 
and a general Jacobi-like identity: 
\[
0
= 
\big[h,\,[h',\,h'']\big] 
+ 
\big[h'',\,[h,\,h']\big]
+ 
\big[h',\,[h'',\,h]\big]
\ \ \ \ \ \ \ \ \ \ \ \ \ 
{\scriptstyle{(h,\,h'\,h''\,\in\,\mathcal{F})}}. 
\]
Such an algebra
$\mathcal{ F}$ is unique up to isomorphism. Thus, the
multiplication in $\mathcal{ F}$ plays the role of the concrete Lie
bracket between vector fields. But importantly, {\it no linear
relation exists between iterated multiplications}, {\it i.e.} between
iterated Lie brackets, {\it except those generated
just by antisymmetry and by
Jacobi identity}: this is {\sl freeness} of the algebra.

Since the bracket multiplication is not associative, one must
carefully write down the occurring
brackets, for instance:
\[
\big[[h_1,\,h_2],\,h_2\big], 
\ \ \ \ \ \ \ \ \
\big[ h_1,\,\big[h_2,\,[h_1,\,h_2]\big],
\ \ \ \ \ \ \ \ \
\big[ [h_1,h_2],\,\big[ h_1,\,[h_1,\,h_2]\big]\big]. 
\]
Writing all such words only with the alphabet $\{ h_1,
h_2\}$, we define the {\sl length} of a word $\mathbf{ h}$ to be the
number of elements $h_{ i_\alpha }$ in it, $i_\alpha = 1, 2$. 
For $\ell \in \N$ with
$\ell \geqslant 1$, let $\mathcal{ W}^\ell$ denote the set of words of
length equal to $\ell$ and let $\mathcal{ W} = \bigcup_{ \ell
\geqslant 1} \mathcal{ W}^\ell$ be the set of all words.

Define $\mathcal{F}_1$ to be the $\R$-vector space generated by $h_1,
h_2$ and for $\ell \geqslant 2$, define $\mathcal{F}_\ell$ to be the
$\R$-vector space generated by all words of length $\leqslant \ell$.
In this way, $\mathcal{F} = \bigcup_{ \ell \geqslant 1}\, \mathcal{
F}_\ell$ naturally becomes a graded Lie algebra, because by applying
inductively the Jacobi identity, one may rather easily establish by
induction that (but see also explicit examples below):
\[
\big[\mathcal{F}_{\ell_1},\,\mathcal{F}_{ \ell_2}\big]
\subset
\mathcal{F}_{\ell_1+\ell_2}. 
\] 
Again by an induction based on the Jacobi identity, 
it also follows that $\mathcal{ 
F}_\ell$ is generated, as an $\R$-vector space, 
by only those words that of the form:
\[
\big[
h_{i_1},\,\big[h_{i_2},\,\big[\dots\big[h_{i_{\ell'-1}},\,h_{i_{\ell'}}
\big]\dots\big]\big],
\]
and which are called {\sl simple}, where $\ell' \leqslant \ell$ and where
$1 \leqslant i_1, i_2, \dots, i_{ \ell' -1}, i_{ \ell'} \leqslant
2$. For instance, the non-simple word 
$\big[ [h_1,\, h_2],\, \big[ h_1,\, [ h_1,\, h_2]\big]\big]$
may be written as a certain 
linear combination of simple words of this kind having length
$5$, as we will see quite explicitly in a while. Let us denote by:
\[
\mathcal{ SW} 
= 
\bigcup_{\ell\geqslant 1}\
\mathcal{ SW}^\ell 
\]
the set of all the simple words, where $\mathcal{ SW}^\ell$ denotes
the set of simple words of length $\ell$. Thus, a rough induction
argument based on Jacobi shows that $\mathcal{ SW}$ generates
$\mathcal{F}$ as a vector space over $\R$, but there are further
linear dependence relations between simple words, as is known and as
will be visible in examples.

\subsection{All relations up to length 5}
Thus, we are interested in words, namely in abstract-free Lie
brackets, just up to length $6$, and this,
unexpectedly, will happen to already be a bit
not straightforward. According to a known theorem ({\em
see} {\em e.g.}~\cite{ MerkerPorten}, p.~11), the dimensions $n_\ell -
n_{ \ell - 1}$ of $\mathcal{ F}_\ell / \mathcal{ F}_{ \ell - 1} $
satisfy the induction relations:
\[
n_\ell
-
n_{\ell-1}
=
\frac{1}{\ell}\,
\sum_{d\ {\rm divides}\ \ell}
\mu(d)\,2^{\frac{\ell}{d}},
\] 
where $\mu$ is the M\"obius function:
\[
\mu(d)
=
\left\{
\aligned
&
1,\ {\rm if}\ d=1;
\\
&
0,\ {\rm if}\ d\
\text{\rm contains square integer factors};
\\
&
(-1)^\nu,\ {\rm if}\ d=p_1\cdots p_\nu\
\text{\rm is the product of}\ \nu \
\text{\rm distinct prime numbers}.
\endaligned
\right.
\]
Thus, a direct application of this general formula yields:
\[
\small
\aligned
n_2-n_1
&
=
{\textstyle{\frac{1}{2}}}\,
\big(
\mu(1)\,2^{\frac{2}{1}}
+
\mu(2)\,2^{\frac{2}{2}}
\big)
=
{\textstyle{\frac{1}{2}}}\,
\big(
2^2-2
\big)
=
1,
\\
n_3
-
n_2
&
=
{\textstyle{\frac{1}{3}}}\,
\big(
\mu(1)\,2^{\frac{3}{1}}
+
\mu(3)\,2^{\frac{3}{3}}
\big)
=
{\textstyle{\frac{1}{3}}}\,
\big(8-2\big)
=
2,
\\
n_4
-
n_3
&
=
{\textstyle{\frac{1}{4}}}\,
\big(
\mu(1)\,2^{\frac{4}{1}}
+
\mu(2)\,2^{\frac{4}{2}}
+
\mu(4)\,2^{\frac{4}{4}}
\big)
=
{\textstyle{\frac{1}{4}}}\,
\big(16-4+0\big)
=
3,
\\
n_5
-
n_4
&
=
{\textstyle{\frac{1}{5}}}\,
\big(
\mu(1)\,2^{\frac{5}{1}}
+
\mu(5)\,2^{\frac{5}{5}}
\big)
=
{\textstyle{\frac{1}{5}}}\,
\big(32-2\big)
=
6,
\\
n_6
-
n_5
&
=
{\textstyle{\frac{1}{6}}}\,
\big(
\mu(1)\,2^{\frac{6}{1}}
+
\mu(2)\,2^{\frac{6}{2}}
+
\mu(3)\,2^{\frac{6}{3}}
+
\mu(6)\,2^{\frac{6}{6}}
\big)
=
{\textstyle{\frac{1}{6}}}\,
\big(64-8-4+2\big)
=
9.
\endaligned
\]
Now, in length $\ell = 2$ it is clear that there is, up to
skew-symmetry, only {\em one} simple word: 
\[
[h_1,\,h_2],
\]
confirming $n_2 - n_1 = 1$ while $n_1 = 2$ of course, because $h_1$
and $h_2$ are two independent simple words of length $1$.

Next, in length $\ell = 3$, it is again clear that up to skew-symmetry,
there are only two simple words:
\[
\big[h_1,\,[h_1,\,h_2]\big]
\ \ \ \ \ \
\text{\rm and}
\ \ \ \ \ \
\big[h_2,\,[h_1,\,h_2]\big],
\] 
while no word is not simple.

It is only in length $\ell = 4$ that nontrivial relations come
out. Indeed, again up to the skew-symmetry inside the `core' $[h_1, \,
h_2]$, there are {\em a priori} $4$ distinct simple words generating
$\mathcal{ SW}^4$, namely:
\[
\footnotesize
\aligned
\big[h_1,\,\big[h_1,\,[h_1,\,h_2]\big]\big],
\ \ \ \ \
\big[h_1,\,\big[h_2,\,[h_1,\,h_2]\big]\big],
\ \ \ \ \
\big[h_2,\,\big[h_1,\,[h_1,\,h_2]\big]\big],
\ \ \ \ \
\big[h_2,\,\big[h_2,\,[h_1,\,h_2]\big]\big],
\endaligned
\]
but an obvious Jacobi identity provides one linear relation
between simple words\footnote{\,
For arbitrary words $h$, $h'$, $h''$ of length $\geqslant 1$, our
convention for writing out any Jacobi identity under either
one or the other form:
\[
\aligned
0
= 
\big[h,\,[h',\,h'']\big] 
+ 
\big[h'',\,[h,\,h']\big]
+ 
\big[h',\,[h'',\,h]\big]
\\
0
= 
\big[[h,\,h'],\,h''\big] 
+ 
\big[[h'',\,h],\,h'\big] 
+ 
\big[[h',\,h''],\,h\big] 
\endaligned
\] 
consists in subjecting the terms to a circular permutation, the last
term being brought back to the first position while other terms are
simultaneously shifted (pushed) from left to right.
}: 
\[
0
=
\big[h_1,\,\big[h_2,\,
[h_1,\,h_2]\big]\big]
+
\zero{\big[[h_1,\,h_2],\,\big[h_1,\,h_2\big]\big]}
+
\big[[h_1,\,h_2],\,\big[h_1,\,h_2\big]\big],
\] 
where the central term trivially vanishes, and this is coherent with
$n_4 - n_3 = 3$. Furthermore, one easily convinces oneself that, up to
skew-symmetry, the Jacobi identity cannot produce any other nontrivial
relation, for instance:
\[
0
=
\big[h_1,\,\big[h_1,\,
[h_1,\,h_2]\big]\big]
+
\big[[h_1,\,h_2],\,\zero{\big[h_1,\,h_1\big]}\big]
+
\big[h_1,\,\big[[h_1,\,h_2],\,h_1,\big]\big],
\]
is a trivial relation, it gives nothing. In fact, one realizes that
{\em all} brackets between two words of length $\ell = 2$ vanish. As
a basis for $\mathcal{ SW}^4$, let us therefore choose the three 
simple words:
\[
\aligned
\big[h_1,\,\big[h_1,\,[h_1,\,h_2]\big]\big],\ \ \ \ \ \ \ \ \
\big[h_1,\,\big[h_2,\,[h_1,\,h_2]\big]\big],\ \ \ \ \ \ \ \ \
\big[h_2,\,\big[h_2,\,[h_1,\,h_2]\big]\big],
\endaligned
\]
remembering that the fourth simple word is simply given by:
\begin{equation}
\label{relation-4}
\big[h_2,\,\big[h_1,\,[h_1,\,h_2]\big]\big]
=
\big[h_1,\,\big[h_2,\,[h_1,\,h_2]\big]\big].
\end{equation}

Next, in length $\ell = 5$, applying $[ h_i, \cdot]$, $i = 1, 2$, 
to these three simple words, we deduce that
$\mathcal{ SW}^5$ is generated by the following six
simple words:
\[
\small
\aligned
&
\big[h_1,\big[h_1,\big[h_1,[h_1,h_2]\big]\big]\big],
\ \ \ \ \ \ \ \ \
\big[h_1,\big[h_1,\big[h_2,[h_1,h_2]\big]\big]\big],
\ \ \ \ \ \ \ \ \
\big[h_1,\big[h_2,\big[h_2,[h_1,h_2]\big]\big]\big],
\\
&
\big[h_2,\big[h_1,\big[h_1,[h_1,h_2]\big]\big]\big],
\ \ \ \ \ \ \ \ \
\big[h_2,\big[h_1,\big[h_2,[h_1,h_2]\big]\big]\big],
\ \ \ \ \ \ \ \ \
\big[h_2,\big[h_2,\big[h_2,[h_1,h_2]\big]\big]\big].
\endaligned
\]
Are there other linear dependence relations between these six simple
words? Certainly not, because of $n_5 - n_4 = 6$; alternatively,
one could also
realize this by trying to apply Jacobi to all possible triples of
words, the sum-length of which equals $5$. 
In addition, it is also important,
for later use, to explicitly represent all length-$5$ multiple
iterated brackets as specific linear combinations between simple
brackets. For instance, there are exactly two brackets between two
basic words of lengths $2$ and $3$, and the Jacobi identity
gives\footnote{\,
For clarity, we underline the three terms that are subjected to
a circular permutation.
}: 
\[
\footnotesize
\aligned
0
&
=
\big[[\underline{h_1},\,\underline{h_2}],\,
\underline{\big[h_1,[h_1,h_2]\big]}\big]
+
\big[\big[\big[h_1,\,[h_1,\,h_2]\big],\,h_1\big],\,h_2\big]
+
\big[\big[h_2,\,\big[h_1,[h_1,h_2]\big]\big],\,h_1\big],
\\
0
&
=
\big[[\underline{h_1},\,\underline{h_2}],\,
\underline{\big[h_2,[h_1,h_2]\big]}\big]
+
\big[\big[\big[h_2,\,[h_1,\,h_2]\big],\,h_1\big],\,h_2\big]
+
\big[\big[h_2,\,\big[h_2,[h_1,h_2]\big]\big],\,h_1\big].
\endaligned
\]
Here, in each one of the two lines, the last two words happen, thanks
to skew-symmetry, to all be simple, whence (using~\thetag{
\ref{relation-4}} for the last term of the first line):
\begin{equation}
\label{relations-5}
\footnotesize
\aligned
\big[[h_1,\,h_2],\,
\big[h_1,[h_1,h_2]\big]\big]
&
=
-\,
\big[h_2,\,\big[h_1,\,\big[h_1,\,[h_1,\,h_2]\big]\big]\big]
+
\big[h_1,\,\big[h_1,\,\big[h_2,\,[h_1,\,h_2]\big]\big]\big],
\\
\big[[h_1,\,h_2],\,\big[h_2,[h_1,h_2]\big]\big]
&
=
-
\big[h_2,\,\big[h_1,\,\big[h_2,\,[h_1,\,h_2]\big]\big]\big]
+
\big[h_1,\,\big[h_2,\,\big[h_2,\,[h_1,\,h_2]\big]\big]\big].
\endaligned
\end{equation}

\subsection{All relations in length 6}
In the next length $\ell = 6$, more complexity occurs. By applying $[
h_i, \cdot]$, $i = 1, 2$, to the above six linearly independent simple
words of length $5$, we at first get the following twelve simple
words:
\[
\scriptsize
\aligned
&
\big[h_1,\big[h_1,\big[h_1,\big[h_1,[h_1,h_2]
\big]\big]\big]\big],
\ \ \ \ \ 
\big[h_1,\big[h_1,\big[h_1,\big[h_2,[h_1,h_2]
\big]\big]\big]\big],
\ \ \ \ \ 
\big[h_1,\big[h_1,\big[h_2,\big[h_2,[h_1,h_2]
\big]\big]\big]\big],
\\
&
\big[h_1,\big[h_2,\big[h_1,\big[h_1,[h_1,h_2]
\big]\big]\big]\big],
\ \ \ \ \ 
\big[h_1,\big[h_2,\big[h_1,\big[h_2,[h_1,h_2]
\big]\big]\big]\big],
\ \ \ \ \ 
\big[h_1,\big[h_2,\big[h_2,\big[h_2,[h_1,h_2]
\big]\big]\big]\big],
\\
&
\big[h_2,\big[h_1,\big[h_1,\big[h_1,[h_1,h_2]
\big]\big]\big]\big],
\ \ \ \ \ 
\big[h_2,\big[h_1,\big[h_1,\big[h_2,[h_1,h_2]
\big]\big]\big]\big],
\ \ \ \ \ 
\big[h_2,\big[h_1,\big[h_2,\big[h_2,[h_1,h_2]
\big]\big]\big]\big],
\\
&
\big[h_2,\big[h_2,\big[h_1,\big[h_1,[h_1,h_2]
\big]\big]\big]\big],
\ \ \ \ \ 
\big[h_2,\big[h_2,\big[h_1,\big[h_2,[h_1,h_2]
\big]\big]\big]\big],
\ \ \ \ \ 
\big[h_2,\big[h_2,\big[h_2,\big[h_2,[h_1,h_2]
\big]\big]\big]\big].
\endaligned
\]
However, according to the dimensional count made above, $n_6 - n_5 =
9$, so there must exist three independent linear relations between
these twelve simple words.

We begin by exploring Lie brackets between two words of length
$3$. Such words are automatically simple. Since there are only two
words of length $3$, only one bracket exists, and the Jacobi identity
can give only two different relations. The first relation is:
\[
\scriptsize
\aligned
0
&
=
\big[\underline{\big[h_1,[h_1,h_2]\big]},
\big[\underline{h_2},\underline{[h_1,h_2]}\big]\big]
+
\big[[h_1,h_2],
\big[\big[h_1,[h_1,h_2]\big],h_2\big]\big]
+
\big[h_2,\big[[h_1,h_2],
\big[h_1,[h_1,h_2]\big]\big]\big]
\\
&
=
\big[\big[h_1,[h_1,h_2]\big],
\big[h_2,[h_1,h_2]\big]\big]
-
\big[[h_1,h_2],
\big[h_2,\big[h_1,[h_1,h_2]\big]\big]\big]
+
\big[h_2,\big[[h_1,h_2],
\big[h_1,[h_1,h_2]\big]\big]\big]
\\
&
=
\big[\big[h_1,[h_1,h_2]\big],
\big[h_2,[h_1,h_2]\big]\big]
-
\big[[h_1,h_2],
\big[h_1,\big[h_2,[h_1,h_2]\big]\big]\big]
-
\big[h_2,\big[h_2,\big[h_1,\big[h_1,
[h_1,h_2]\big]\big]\big]\big]
+
\\
&
\ \ \ \ \ \ \ \ \ \ \ \ \ \ \ \ \ \ \ \ \ \ \ \ \ \ \ \ \ \ \ \ \ \ \
\ \ \ \ \ \ \ \ \ \ \ \ \ \ \ \ \ \ \ \ \ \ \ \ \ \ \ \ \ \ \ \ \ \ \
\ \ \ \ \ \ \ \ \ \ \ \ \ \ \ \ \ \ \ \ \ \ \ \ \ \ \ \ \ \ \ \ \ \ \
\ \ \ \ \ 
+
\big[h_2,\big[h_1,\big[h_1,\big[h_2,
[h_1,h_2]\big]\big]\big]\big].
\endaligned
\]
where we leave untouched the first term, where we apply~\thetag{
\ref{relation-4}} to normalize the second term, and where the third
term expresses as a linear combination of two simple words thanks to
the first relation~\thetag{ \ref{relations-5}}. The second relation,
just with a different underlining for applying Jacobi, is:
\[
\scriptsize
\aligned
0
&
=
\big[\big[\underline{h_1},\underline{[h_1,h_2]}\big],
\underline{\big[h_2,[h_1,h_2]\big]}\big]
+
\big[\big[\big[h_2,[h_1,h_2]\big],h_1\big],[h_1,h_2]\big]
+
\big[\big[[h_1,h_2],
\big[h_2,[h_1,h_2]\big]\big],h_1\big]
\\
&
=
\big[\big[h_1,[h_1,h_2]\big],
\big[h_2,[h_1,h_2]\big]\big]
+
\big[[h_1,h_2],
\big[h_1,\big[h_2,[h_1,h_2]\big]\big]\big]
-
\big[h_1,\big[[h_1,h_2],
\big[h_2,[h_1,h_2]\big]\big]\big]
\\
&
=
\big[\big[h_1,[h_1,h_2]\big],
\big[h_2,[h_1,h_2]\big]\big]
+
\big[[h_1,h_2],
\big[h_1,\big[h_2,[h_1,h_2]\big]\big]\big]
+
\big[h_1,\big[h_2,\big[h_1,\big[h_2,
[h_1,h_2]\big]\big]\big]\big]
-
\\
&
\ \ \ \ \ \ \ \ \ \ \ \ \ \ \ \ \ \ \ \ \ \ \ \ \ \ \ \ \ \ \ \ \ \ \
\ \ \ \ \ \ \ \ \ \ \ \ \ \ \ \ \ \ \ \ \ \ \ \ \ \ \ \ \ \ \ \ \ \ \
\ \ \ \ \ \ \ \ \ \ \ \ \ \ \ \ \ \ \ \ \ \ \ \ \ \ \ \ \ \ \ \ \ \ \
\ \ \ \ 
-
\big[h_1,\big[h_1,\big[h_2,\big[h_2,
[h_1,h_2]\big]\big]\big]\big].
\endaligned
\]

Next, there are exactly six Jacobi identities for Lie brackets between
two words having lengths $2$ and $4$. The first pair is:
\[
\scriptsize
\aligned
0
&
=
\big[\underline{[h_1,h_2]},\big[
\underline{h_1},\underline{\big[h_1,[h_1,h_2]\big]}\big]\big]
+
\zero{\big[\big[h_1,[h_1,h_2]\big],
\big[[h_1,h_2],h_1\big]\big]}
+
\big[h_1,\big[\big[h_1,[h_1,h_2]\big],
[h_1,h_2]\big]\big]
\\
&
=
\big[[h_1,h_2],\big[h_1,\big[h_1,[h_1,h_2]\big]\big]\big]
-
\big[h_1,\big[[h_1,h_2],
\big[h_1,[h_1,h_2]\big]\big]\big]
\\
&
=
\big[[h_1,h_2],\big[h_1,\big[h_1,[h_1,h_2]\big]\big]\big]
+
\big[h_1,\big[h_2,\big[h_1,\big[h_1,
[h_1,h_2]\big]\big]\big]\big]
-
\big[h_1,\big[h_1,\big[h_1,\big[h_2,
[h_1,h_2]\big]\big]\big]\big];
\endaligned
\]
\[
\scriptsize
\aligned
0
&
=
\big[[\underline{h_1},\underline{h_2}],
\underline{\big[h_1,\big[h_1,[h_1,h_2]\big]\big]}\big]
+
\big[\big[\big[h_1,\big[h_1,
[h_1,h_2]\big]\big],h_1\big],h_2\big]
+
\big[\big[h_2,\big[h_1,\big[h_1,[h_1,h_2]\big]\big]\big],
h_1\big]
\\
&
=
\big[[h_1,h_2],\big[h_1,\big[h_1,[h_1,h_2]\big]\big]\big]
+
\big[h_2,\big[h_1,\big[h_1,\big[h_1,
[h_1,h_2]\big]\big]\big]\big]
-
\big[h_1,\big[h_2,\big[h_1,\big[h_1,
[h_1,h_2]\big]\big]\big]\big].
\endaligned
\]
The second pair is:
\[
\scriptsize
\aligned
0
&
=
\big[\underline{[h_1,h_2]},\big[
\underline{h_1},\underline{\big[h_2,[h_1,h_2]\big]}\big]\big]
+
\big[\big[h_2,[h_1,h_2]\big],
\big[[h_1,h_2],h_1\big]\big]
+
\big[h_1,\big[\big[h_2,[h_1,h_2]\big],
[h_1,h_2]\big]\big]
\\
&
=
\big[[h_1,h_2],\big[h_1,\big[h_2,[h_1,h_2]\big]\big]\big]
+
\big[\big[h_1,[h_1,h_2]\big],
\big[h_2,[h_1,h_2]\big]\big]
-
\big[h_1,\big[[h_1,h_2],
\big[h_2,[h_1,h_2]\big]\big]\big]
\\
&
=
\big[[h_1,h_2],\big[h_1,\big[h_2,[h_1,h_2]\big]\big]\big]
+
\big[\big[h_1,[h_1,h_2]\big],
\big[h_2,[h_1,h_2]\big]\big]
+
\big[h_1,\big[h_2,\big[h_1,\big[h_2,
[h_1,h_2]\big]\big]\big]\big]
-
\\
&
\ \ \ \ \ \ \ \ \ \ \ \ \ \ \ \ \ \ \ \ \ \ \ \ \ \ \ \ \ \ \ \ \ \ \
\ \ \ \ \ \ \ \ \ \ \ \ \ \ \ \ \ \ \ \ \ \ \ \ \ \ \ \ \ \ \ \ \ \ \
\ \ \ \ \ \ \ \ \ \ \ \ \ \ \ \ \ \ \ \ \ \ \ \ \ \ \ \ \ \ \ \ \ \ \
\ \ \ \
-
\big[h_1,\big[h_1,\big[h_2,\big[h_2,
[h_1,h_2]\big]\big]\big]\big];
\endaligned
\]
\[
\scriptsize
\aligned
0
&
=
\big[[\underline{h_1},\underline{h_2}],
\underline{\big[h_1,\big[h_2,[h_1,h_2]\big]\big]}\big]
+
\big[\big[\big[h_1,\big[h_2,
[h_1,h_2]\big]\big],h_1\big],h_2\big]
+
\big[\big[h_2,\big[h_1,\big[h_2,[h_1,h_2]\big]\big]\big],
h_1\big]
\\
&
=
\big[[h_1,h_2],\big[h_1,\big[h_2,[h_1,h_2]\big]\big]\big]
+
\big[h_2,\big[h_1,\big[h_1,\big[h_2,
[h_1,h_2]\big]\big]\big]\big]
-
\big[h_1,\big[h_2,\big[h_1,\big[h_2,
[h_1,h_2]\big]\big]\big]\big].
\endaligned
\]
The third pair is:
\[
\scriptsize
\aligned
0
&
=
\big[\underline{[h_1,h_2]},\big[
\underline{h_2},\underline{\big[h_2,[h_1,h_2]\big]}\big]\big]
+
\zero{\big[\big[h_2,[h_1,h_2]\big],
\big[[h_1,h_2],h_2\big]\big]}
+
\big[h_2,\big[\big[h_2,[h_1,h_2]\big],
[h_1,h_2]\big]\big]
\\
&
=
\big[[h_1,h_2],\big[h_2,\big[h_2,[h_1,h_2]\big]\big]\big]
-
\big[h_2,\big[[h_1,h_2],
\big[h_2,[h_1,h_2]\big]\big]\big]
\\
&
=
\big[[h_1,h_2],\big[h_2,\big[h_2,[h_1,h_2]\big]\big]\big]
+
\big[h_2,\big[h_2,\big[h_1,\big[h_2,
[h_1,h_2]\big]\big]\big]\big]
-
\big[h_2,\big[h_1,\big[h_2,\big[h_2,
[h_1,h_2]\big]\big]\big]\big];
\endaligned
\]
\[
\scriptsize
\aligned
0
&
=
\big[[\underline{h_1},\underline{h_2}],
\underline{\big[h_2,\big[h_2,[h_1,h_2]\big]\big]}\big]
+
\big[\big[\big[h_2,\big[h_2,
[h_1,h_2]\big]\big],h_1\big],h_2\big]
+
\big[\big[h_2,\big[h_2,\big[h_2,[h_1,h_2]\big]\big]\big],
h_1\big]
\\
&
=
\big[[h_1,h_2],\big[h_2,\big[h_2,[h_1,h_2]\big]\big]\big]
+
\big[h_2,\big[h_1,\big[h_2,\big[h_2,
[h_1,h_2]\big]\big]\big]\big]
-
\big[h_1,\big[h_2,\big[h_2,\big[h_2,
[h_1,h_2]\big]\big]\big]\big].
\endaligned
\]
We thus have got eight relations involving simple words: the single
bracket between two words of length $3$ and the six brackets between
a word of length $2$ and a word of length $4$.
We number these eight equations and,
for easier readability, we underbrace
the four non-simple words ${\tt t}$, ${\tt u}$, ${\tt v}$, ${\tt w}$:
\[
\scriptsize
\aligned
0
&
\overset{1}{=}
\underbrace{\big[\big[h_1,[h_1,h_2]\big],
\big[h_2,[h_1,h_2]\big]\big]}_{\tt t}
-
\underbrace{\big[[h_1,h_2],
\big[h_1,\big[h_2,[h_1,h_2]\big]\big]\big]}_{\tt v}
-
\big[h_2,\big[h_2,\big[h_1,\big[h_1,
[h_1,h_2]\big]\big]\big]\big]
+
\\
&
\ \ \ \ \ \ \ \ \ \ \ \ \ \ \ \ \ \ \ \ \ \ \ \ \ \ \ \ \ \ \ \ \ \ \
\ \ \ \ \ \ \ \ \ \ \ \ \ \ \ \ \ \ \ \ \ \ \ \ \ \ \ \ \ \ \ \ \ \ \
\ \ \ \ \ \ \ \ \ \ \ \ \ \ \ \ \ \ \ \ \ \ \ \ \ \ \ \ \ \ \ \ \ \ \
\ \ \ \ 
+
\big[h_2,\big[h_1,\big[h_1,\big[h_2,
[h_1,h_2]\big]\big]\big]\big]
\\
0
&
\overset{2}{=}
\underbrace{\big[\big[h_1,[h_1,h_2]\big],
\big[h_2,[h_1,h_2]\big]\big]}_{\tt t}
+
\underbrace{\big[[h_1,h_2],
\big[h_1,\big[h_2,[h_1,h_2]\big]\big]\big]}_{\tt v}
+
\big[h_1,\big[h_2,\big[h_1,\big[h_2,
[h_1,h_2]\big]\big]\big]\big]
-
\\
&
\ \ \ \ \ \ \ \ \ \ \ \ \ \ \ \ \ \ \ \ \ \ \ \ \ \ \ \ \ \ \ \ \ \ \
\ \ \ \ \ \ \ \ \ \ \ \ \ \ \ \ \ \ \ \ \ \ \ \ \ \ \ \ \ \ \ \ \ \ \
\ \ \ \ \ \ \ \ \ \ \ \ \ \ \ \ \ \ \ \ \ \ \ \ \ \ \ \ \ \ \ \ \ \ \
\ \ \ \
-
\big[h_1,\big[h_1,\big[h_2,\big[h_2,
[h_1,h_2]\big]\big]\big]\big],
\endaligned
\]
\[
\scriptsize
\aligned
0
&
\overset{3}{=}
\underbrace{
\big[[h_1,h_2],\big[h_1,\big[h_1,[h_1,h_2]\big]\big]\big]}_{
\tt u}
+
\big[h_1,\big[h_2,\big[h_1,\big[h_1,
[h_1,h_2]\big]\big]\big]\big]
-
\big[h_1,\big[h_1,\big[h_1,\big[h_2,
[h_1,h_2]\big]\big]\big]\big],
\\
0
&
\overset{4}{=}
\underbrace{
\big[[h_1,h_2],\big[h_1,\big[h_1,[h_1,h_2]\big]\big]\big]}_{
\tt u}
+
\big[h_2,\big[h_1,\big[h_1,\big[h_1,
[h_1,h_2]\big]\big]\big]\big]
-
\big[h_1,\big[h_2,\big[h_1,\big[h_1,
[h_1,h_2]\big]\big]\big]\big],
\endaligned
\]
\[
\scriptsize
\aligned
0
&
\overset{5}{=}
\underbrace{
\big[[h_1,h_2],\big[h_1,\big[h_2,[h_1,h_2]\big]\big]\big]}_{
\tt v}
+
\underbrace{
\big[\big[h_1,[h_1,h_2]\big],
\big[h_2,[h_1,h_2]\big]\big]}_{\tt t}
+
\big[h_1,\big[h_2,\big[h_1,\big[h_2,
[h_1,h_2]\big]\big]\big]\big]
-
\\
&
\ \ \ \ \ \ \ \ \ \ \ \ \ \ \ \ \ \ \ \ \ \ \ \ \ \ \ \ \ \ \ \ \ \ \
\ \ \ \ \ \ \ \ \ \ \ \ \ \ \ \ \ \ \ \ \ \ \ \ \ \ \ \ \ \ \ \ \ \ \
\ \ \ \ \ \ \ \ \ \ \ \ \ \ \ \ \ \ \ \ \ \ \ \ \ \ \ \ \ \ \ \ \ \ \
\ \ \ \
-
\big[h_1,\big[h_1,\big[h_2,\big[h_2,
[h_1,h_2]\big]\big]\big]\big],
\\
0
&
\overset{6}{=}
\underbrace{
\big[[h_1,h_2],\big[h_1,\big[h_2,[h_1,h_2]\big]\big]\big]}_{
\tt v}
+
\big[h_2,\big[h_1,\big[h_1,\big[h_2,
[h_1,h_2]\big]\big]\big]\big]
-
\big[h_1,\big[h_2,\big[h_1,\big[h_2,
[h_1,h_2]\big]\big]\big]\big].
\endaligned
\]
\[
\scriptsize
\aligned
0
&
\overset{7}{=}
\underbrace{
\big[[h_1,h_2],\big[h_2,\big[h_2,[h_1,h_2]\big]\big]\big]}_{
\tt w}
+
\big[h_2,\big[h_2,\big[h_1,\big[h_2,
[h_1,h_2]\big]\big]\big]\big]
-
\big[h_2,\big[h_1,\big[h_2,\big[h_2,
[h_1,h_2]\big]\big]\big]\big],
\\
0
&
\overset{8}{=}
\underbrace{
\big[[h_1,h_2],\big[h_2,\big[h_2,[h_1,h_2]\big]\big]\big]}_{
\tt w}
+
\big[h_2,\big[h_1,\big[h_2,\big[h_2,
[h_1,h_2]\big]\big]\big]\big]
-
\big[h_1,\big[h_2,\big[h_2,\big[h_2,
[h_1,h_2]\big]\big]\big]\big].
\endaligned
\]
Visibly, the fifth equation coincides with the second one. There
remain seven equations. Since four non-simple words are involved, one
may expect to see here the three linearly independent relations
between simple words that we are looking for. Firstly, subtracting the
third equation to the fourth, we get a first relation of this kind:
\[
\scriptsize
\aligned
0
&
\overset{9}{=}
\big[h_1,\big[h_1,\big[h_1,\big[h_2,
[h_1,h_2]\big]\big]\big]\big]
-
2\times
\big[h_1,\big[h_2,\big[h_1,\big[h_1,
[h_1,h_2]\big]\big]\big]\big]
+
\\
&
\ \ \ \ \
+
\big[h_2,\big[h_1,\big[h_1,\big[h_1,
[h_1,h_2]\big]\big]\big]\big].
\endaligned
\]
Secondly, subtracting the eighth equation to the
seventh, we get a second, visibly independent relation:
\[
\scriptsize
\aligned
0
&
\overset{10}{=}
\big[h_2,\big[h_2,\big[h_1,\big[h_2,
[h_1,h_2]\big]\big]\big]\big]
-
2\times
\big[h_2,\big[h_1,\big[h_2,\big[h_2,
[h_1,h_2]\big]\big]\big]\big]
+
\\
&
\ \ \ \ \ 
+
\big[h_1,\big[h_2,\big[h_2,\big[h_2,
[h_1,h_2]\big]\big]\big]\big].
\endaligned
\]
Thirdly and lastly, adding the sixth equation multiplied by $2$
to the first one and subtracting the 
second one, we get a third independent
relation between simple words:
\[
\scriptsize
\aligned
0
&
\overset{11}{=}
\big[h_1,\big[h_1,\big[h_2,\big[h_2,
[h_1,h_2]\big]\big]\big]\big]
-
3\times
\big[h_1,\big[h_2,\big[h_1,\big[h_2,
[h_1,h_2]\big]\big]\big]\big]
+
\\
&
\ \ \ \ \ \
+
3\times
\big[h_2,\big[h_1,\big[h_1,\big[h_2,
[h_1,h_2]\big]\big]\big]\big]
-
\big[h_2,\big[h_2,\big[h_1,\big[h_1,
[h_1,h_2]\big]\big]\big]\big].
\endaligned
\]

\subsection{Iterated brackets of $H_1$ and $H_2$ on $M$}
Now, we come back to our two vector fields $H_1$ and $H_2$
on $M$ satisfying the two specific relations:
\[
\big[H_1,[H_1,H_2]\big]
=
\Phi_1[H_1,H_2]
\ \ \ \ \ \ \ \
\text{\rm and}
\ \ \ \ \ \ \ \
\big[H_2,[H_1,H_2]\big]
=
\Phi_2[H_1,H_2],
\]
for certain two functions $\Phi_1$ and $\Phi_2$ on $M$ whose explicit
expressions (not needed here), in terms of the third-order jet $J_{
x, y, u}^3 \varphi$ of the graphing function for $M$, have already been
shown in Proposition~\ref{explicit-J-3-varphi}. These two relations
show well that $H_1$ and $H_2$ do {\em not} behave as the two abstract
{\em totally free} elements $h_1$ and $h_2$ considered above. In fact, a
straightforward induction argument shows that every iterated
simple-word bracket:
\[
\big[H_{i_1},\big[H_{i_2},
\big[\cdots,
\big[H_{i_{\ell-1}},H_{i_\ell}\big],\cdots\big]\big]\big]
=
\Phi_{i_1,i_2,\dots,i_{\ell-1},i_\ell}
[H_1,H_2]
\]
of arbitrary length $\ell \geqslant 2$, where $i_1, \dots, i_\ell = 1,
2$, must always be a multiple of $[ H_1, \, H_2]$ by means of a
certain function $\Phi_{ i_1, \dots, i_\ell}$ which depends on
$\Phi_1$ and $\Phi_2$, but whose explicit expression in terms of
$\Phi_1$ and $\Phi_2$ is not immediate. For $\ell = 4, 5, 6$, we must
now compute all these $\Phi_{ i_1, \dots, i_\ell}$ so that the abstract
relations between iterated brackets of the free $h_1$ and $h_2$
computed above provide us with interesting relations that will be useful
later.

At first, applying the basic formula 
$[ fX,\, gY] = fX(g)\,Y - g\,Y(f)\, X + fg\,[X, Y]$, 
we have in length $\ell = 4$:
\[
\small
\aligned
\big[H_1,\,\big[H_1,\,[H_1,\,H_2]\big]\big]
&
=
H_1(\Phi_1)\,[H_1,\,H_2]
+
\Phi_1\,\big[H_1,\,[H_1,\,H_2]\big]
\\
&
\overset{1}{=}
\big(H_1(\Phi_1)+(\Phi_1)^2\big)\,[H_1,\,H_2],
\endaligned
\]
\[
\small
\aligned
\big[H_1,\,\big[H_2,\,[H_1,\,H_2]\big]\big]
&
=
H_1(\Phi_2)\,[H_1,\,H_2]
+
\Phi_2\,\big[H_2,\,[H_1,\,H_2]\big]
\\
&
\overset{2}{=}
\big(H_1(\Phi_2)+\Phi_2\Phi_1\big)\,[H_1,\,H_2],
\endaligned
\]
\[
\small
\aligned
\big[H_2,\,\big[H_1,\,[H_1,\,H_2]\big]\big]
&
=
H_2(\Phi_1)\,[H_1,\,H_2]
+
\Phi_1\,\big[H_2,\,[H_1,\,H_2]\big]
\\
&
\overset{3}{=}
\big(H_2(\Phi_1)+\Phi_1\Phi_2\big)\,[H_1,\,H_2],
\endaligned
\]
\[
\small
\aligned
\big[H_2,\,\big[H_2,\,[H_1,\,H_2]\big]\big]
&
=
H_2(\Phi_2)\,[H_1,\,H_2]
+
\Phi_2\,\big[H_2,\,[H_1,\,H_2]\big]
\\
&
\overset{4}{=}
\big(H_2(\Phi_2)+(\Phi_2)^2\big)\,[H_1,\,H_2].
\endaligned
\]
But the known relation $\big[ h_2, \, \big[ h_1, \, [ h_1, \,
h_2]\big] \big] = \big[ h_1, \, \big[ h_2, \, [ h_1, \, h_2]\big]
\big]$ between free elements imposes here:
\[
H_2(\Phi_1)+\Phi_1\Phi_2
=
H_1(\Phi_2)+\Phi_2\Phi_1,
\]
a relation already seen in Lemma~\ref{H1-Phi2-H2-Phi1}.

Next, setting aside the consideration of $\big[ H_2, \, \big[ H_1, \,
[H_1,\, H_2] \big] \big]$, we compute as follows the six simple
iterated brackets of length $\ell = 5$ (we replace $H_2( \Phi_1)$ by
$H_1 ( \Phi_2)$ wherever it occurs):
\[
\scriptsize
\aligned
\big[H_1,\big[H_1,\big[H_1,[H_1,H_2]\big]\big]\big]
&
=
\big(H_1(H_1(\Phi_1))+2\Phi_1H_1(\Phi_1)
+
\Phi_1H_1(\Phi_1)+(\Phi_1)^3\big)
[H_1,H_2]
\\
&
\overset{1}{=}
\big(
H_1(H_1(\Phi_1))+3\Phi_1H_1(\Phi_1)+(\Phi_1)^3
\big)
[H_1,H_2],
\endaligned
\]
\[
\label{length-5-2}
\scriptsize
\aligned
\big[H_1,\big[H_1,\big[H_2,[H_1,H_2]\big]\big]\big]
&
=
\big(
H_1(H_1(\Phi_2))+\Phi_1H_1(\Phi_2)
+
\Phi_2H_1(\Phi_1)+\Phi_1H_1(\Phi_2)
+
(\Phi_1)^2\Phi_2
\big)
[H_1,H_2]
\\
&
\overset{2}{=}
\big(
H_1(H_1(\Phi_2))+2\Phi_1H_1(\Phi_2)
+
\Phi_2H_1(\Phi_1)+(\Phi_1)^2\Phi_2
\big)
[H_1,H_2],
\endaligned
\]
\[
\scriptsize
\aligned
\big[H_1,\big[H_2,\big[H_2,[H_1,H_2]\big]\big]\big]
&
\overset{3}{=}
\big(
H_1(H_2(\Phi_2))+2\Phi_2H_1(\Phi_2)
+
\Phi_1H_2(\Phi_2)+\Phi_1(\Phi_2)^2
\big)
[H_1,H_2],
\endaligned
\]
\[
\scriptsize
\aligned
\big[H_2,\big[H_1,\big[H_1,[H_1,H_2]\big]\big]\big]
&
=
\big(
H_2(H_1(\Phi_1))+2\Phi_1H_2(\Phi_1)
+
\Phi_2H_1(\Phi_1)+(\Phi_1)^2\Phi_2
\big)
[H_1,H_2]
\\
&
\overset{4}{=}
\big(
H_2(H_1(\Phi_1))+2\Phi_1H_1(\Phi_2)
+
\Phi_2H_1(\Phi_1)+(\Phi_1)^2\Phi_2
\big)
[H_1,H_2],
\endaligned
\]
\[
\scriptsize
\aligned
\big[H_2,\big[H_1,\big[H_2,[H_1,H_2]\big]\big]\big]
&
=
\big(
H_2(H_1(\Phi_2))+\Phi_1H_2(\Phi_2)
+
\Phi_2H_2(\Phi_1)+\Phi_2H_1(\Phi_2)
+
\Phi_1(\Phi_2)^2
\big)
[H_1,H_2]
\\
&
\overset{5}{=}
\big(
H_2(H_1(\Phi_2))+\Phi_1H_2(\Phi_2)
+
2\Phi_2H_1(\Phi_2)
+
\Phi_1(\Phi_2)^2
\big)
[H_1,H_2],
\endaligned
\]
\[
\scriptsize
\aligned
\big[H_2,\big[H_2,\big[H_2,[H_1,H_2]\big]\big]\big]
&
\overset{6}{=}
\big(
H_2(H_2(\Phi_2))+3\Phi_2H_2(\Phi_2)+(\Phi_2)^3
\big)
[H_1,H_2].
\endaligned
\]
Also, one may compute the two brackets between the unique simple word
of length $2$ and the two simple words of length $3$, expanding $[
H_1, H_2] ( \Psi)$ just as $H_1 ( H_2 ( \Psi)) - H_2 ( H_1 (
\Psi))$:
\[
\scriptsize
\aligned
\big[[H_1,H_2],\big[H_1,[H_1,H_2]\big]\big]
&
=
\big(H_1(H_2(\Phi_1))-H_2(H_1(\Phi_1))\big)[H_1,H_2],
\\
\big[[H_1,H_2],\big[H_2,[H_1,H_2]\big]\big]
&
=
\big(H_1(H_2(\Phi_2))-H_2(H_1(\Phi_2))\big)[H_1,H_2].
\endaligned
\]
Unexpectedly, if one looks at the two relations~\thetag{
\ref{relations-5}}, one only gets twice the trivial relation $0 = 0$. 
Only in length $\ell = 6$ will one find new nontrivial
relations.

Now, here are the twelve (not
linearly independent) simple iterated brackets of length $\ell = 6$:
\[
\scriptsize
\aligned
\big[H_1,\big[H_1,\big[H_1,\big[H_1,
[H_1,H_2]\big]\big]\big]\big]
&
\overset{1}{=}
\big(
H_1(H_1(H_1(\Phi_1)))
+
4\Phi_1H_1(H_1(\Phi_1))
+
\\
&
\ \ \ \ \ 
+
3H_1(\Phi_1)H_1(\Phi_1)
+
6(\Phi_1)^2H_1(\Phi_1)
+
(\Phi_1)^4
\big)
[H_1,H_2],
\endaligned
\]
\[
\scriptsize
\aligned
\big[H_1,\big[H_1,\big[H_1,\big[H_2,
[H_1,H_2]\big]\big]\big]\big]
&
\overset{2}{=}
\big(
H_1(H_1(H_1(\Phi_2)))
+
3\Phi_1H_1(H_1(\Phi_2))
+
\\
&
\ \ \ \ \
+
\Phi_2H_1(H_1(\Phi_1))
+
3H_1(\Phi_1)H_1(\Phi_2)
+
\\
&
\ \ \ \ \ 
+
3\Phi_1\Phi_2H_1(\Phi_1)
+
3(\Phi_1)^2H_1(\Phi_2)
+
(\Phi_1)^3\Phi_2
\big)
[H_1,H_2],
\endaligned
\]
\[
\scriptsize
\aligned
\big[H_1,\big[H_1,\big[H_2,\big[H_2,
[H_1,H_2]\big]\big]\big]\big]
&
\overset{3}{=}
\big(
H_1(H_1(H_2(\Phi_2)))
+
2\Phi_1H_1(H_2(\Phi_2))
+
\\
&
\ \ \ \ \
+
2\Phi_2H_1(H_1(\Phi_2))
+
H_1(\Phi_1)H_2(\Phi_2)
+
2H_1(\Phi_2)H_1(\Phi_2)
+
\\
&
\ \ \ \ \ 
+
(\Phi_1)^2H_2(\Phi_2)
+
(\Phi_2)^2H_1(\Phi_1)
+
4\Phi_1\Phi_2H_1(\Phi_2)
+
(\Phi_1)^2(\Phi_2)^2
\big)
[H_1,H_2],
\endaligned
\]
\[
\scriptsize
\aligned
\big[H_1,\big[H_2,\big[H_1,\big[H_1,
[H_1,H_2]\big]\big]\big]\big]
&
\overset{4}{=}
\big(
H_1(H_2(H_1(\Phi_1)))
+
2\Phi_1H_1(H_1(\Phi_2))
+
\\
&
\ \ \ \ \
+
\Phi_2H_1(H_1(\Phi_1))
+
\Phi_1H_2(H_1(\Phi_1))
+
3H_1(\Phi_1)H_2(\Phi_2)
+
\\
&
\ \ \ \ \ 
+
3\Phi_1\Phi_2H_1(\Phi_1)
+
3(\Phi_1)^2H_1(\Phi_2)
+
(\Phi_1)^3\Phi_2
\big)
[H_1,H_2],
\endaligned
\]
\[
\scriptsize
\aligned
\big[H_1,\big[H_2,\big[H_1,\big[H_2,
[H_1,H_2]\big]\big]\big]\big]
&
\overset{5}{=}
\big(
H_1(H_2(H_1(\Phi_2)))
+
\Phi_1H_1(H_2(\Phi_2))
+
2\Phi_2H_1(H_1(\Phi_2))
+
\\
&
\ \ \ \ \
+
\Phi_1H_2(H_1(\Phi_2))
+
H_1(\Phi_1)H_2(\Phi_2)
+
2H_1(\Phi_2)H_1(\Phi_2)
+
\\
&
\ \ \ \ \ 
+
4\Phi_1\Phi_2H_1(\Phi_2)
+
(\Phi_1)^2H_2(\Phi_2)
+
(\Phi_2)^2H_1(\Phi_1)
+
(\Phi_1)^2(\Phi_2)^2
\big)
[H_1,H_2],
\endaligned
\]
\[
\scriptsize
\aligned
\big[H_1,\big[H_2,\big[H_2,\big[H_2,
[H_1,H_2]\big]\big]\big]\big]
&
\overset{6}{=}
\big(
H_1(H_2(H_2(\Phi_2)))
+
3\Phi_2H_1(H_2(\Phi_2))
+
\\
&
\ \ \ \ \
+
\Phi_1H_2(H_2(\Phi_2))
+
3H_1(\Phi_2)H_2(\Phi_2)
+
\\
&
\ \ \ \ \ 
+
3(\Phi_2)^2H_1(\Phi_2)
+
3\Phi_1\Phi_2H_2(\Phi_2)
+
\Phi_1(\Phi_2)^3
\big)
[H_1,H_2],
\endaligned
\]
\[
\scriptsize
\aligned
\big[H_2,\big[H_1,\big[H_1,\big[H_1,
[H_1,H_2]\big]\big]\big]\big]
&
\overset{7}{=}
\big(
H_2(H_1(H_1(\Phi_1)))
+
3\Phi_1H_2(H_1(\Phi_1))
+
\\
&
\ \ \ \ \
+
\Phi_2H_1(H_1(\Phi_1))
+
3H_1(\Phi_1)H_1(\Phi_2)
+
\\
&
\ \ \ \ \ 
+
3(\Phi_1)^2H_1(\Phi_2)
+
3\Phi_1\Phi_2H_1(\Phi_1)
+
(\Phi_1)^3\Phi_2
\big)
[H_1,H_2],
\endaligned
\]
\[
\scriptsize
\aligned
\big[H_2,\big[H_1,\big[H_1,\big[H_2,
[H_1,H_2]\big]\big]\big]\big]
&
\overset{8}{=}
\big(
H_2(H_1(H_1(\Phi_2)))
+
2\Phi_1H_2(H_1(\Phi_2))
+
\Phi_2H_2(H_1(\Phi_1))
+
\\
&
\ \ \ \ \
+
\Phi_2H_1(H_1(\Phi_2))
+
2H_1(\Phi_2)H_1(\Phi_2)
+
H_2(\Phi_2)H_1(\Phi_1)
+
\\
&
\ \ \ \ \ 
+
4\Phi_1\Phi_2H_1(\Phi_2)
+
(\Phi_1)^2H_2(\Phi_2)
+
(\Phi_2)^2H_1(\Phi_1)
+
(\Phi_1)^2(\Phi_2)^2
\big)
[H_1,H_2],
\endaligned
\]
\[
\scriptsize
\aligned
\big[H_2,\big[H_1,\big[H_2,\big[H_2,
[H_1,H_2]\big]\big]\big]\big]
&
\overset{9}{=}
\big(
H_2(H_1(H_2(\Phi_2)))
+
2\Phi_2H_2(H_1(\Phi_2))
+
\\
&
\ \ \ \ \
+
\Phi_1H_2(H_2(\Phi_2))
+
\Phi_2H_1(H_2(\Phi_2))
+
3H_1(\Phi_2)H_2(\Phi_2)
+
\\
&
\ \ \ \ \ 
+
3(\Phi_2)^2H_1(\Phi_2)
+
3\Phi_1\Phi_2H_2(\Phi_2)
+
\Phi_1(\Phi_2)^3
\big)
[H_1,H_2],
\endaligned
\]
\[
\scriptsize
\aligned
\big[H_2,\big[H_2,\big[H_1,\big[H_1,
[H_1,H_2]\big]\big]\big]\big]
&
\overset{10}{=}
\big(
H_2(H_2(H_1(\Phi_1)))
+
2\Phi_1H_2(H_1(\Phi_2))
+
\\
&
\ \ \ \ \
+
2\Phi_2H_2(H_1(\Phi_1))
+
2H_1(\Phi_2)H_1(\Phi_2)
+
H_1(\Phi_1)H_2(\Phi_2)
+
\\
&
\ \ \ \ \ 
+
4\Phi_1\Phi_2H_1(\Phi_2)
+
(\Phi_1)^2H_2(\Phi_2)
+
(\Phi_2)^2H_1(\Phi_1)
+
(\Phi_1)^2(\Phi_2)^2
\big)
[H_1,H_2],
\endaligned
\]
\[
\scriptsize
\aligned
\big[H_2,\big[H_2,\big[H_1,\big[H_2,
[H_1,H_2]\big]\big]\big]\big]
&
\overset{11}{=}
\big(
H_2(H_2(H_1(\Phi_2)))
+
\Phi_1H_2(H_2(\Phi_2))
+
\\
&
\ \ \ \ \
+
3\Phi_2H_2(H_1(\Phi_2))
+
3H_1(\Phi_2)H_2(\Phi_2)
+
\\
&
\ \ \ \ \ 
+
3(\Phi_2)^2H_1(\Phi_2)
+
3\Phi_1\Phi_2H_2(\Phi_2)
+
\Phi_1(\Phi_2)^3
\big)
[H_1,H_2],
\endaligned
\]
\[
\scriptsize
\aligned
\big[H_2,\big[H_2,\big[H_2,\big[H_2,
[H_1,H_2]\big]\big]\big]\big]
&
\overset{12}{=}
\big(
H_2(H_2(H_2(\Phi_2)))
+
4\Phi_2H_2(H_2(\Phi_2))
+
\\
&
\ \ \ \ \ 
+
3H_2(\Phi_2)H_2(\Phi_2)
+
6(\Phi_2)^2H_2(\Phi_2)
+
(\Phi_2)^4
\big)
[H_1,H_2].
\endaligned
\]
Also, one may compute the single Lie bracket between
two simple words of length $3$ and the three 
Lie brackets between the single simple word of length
$2$ and the three simple words of length $3$:
\[
\scriptsize
\aligned
\big[
\big[H_1,[H_1,H_2]\big],
\big[H_2,[H_1,H_2]\big]
\big]
&
\overset{13}{=}
\big(
\Phi_1H_1(H_2(\Phi_2))
-
\Phi_1H_2(H_1(\Phi_2))
-
\\
&
\ \ \ \ \
-
\Phi_2H_1(H_1(\Phi_2))
+
\Phi_2H_2(H_1(\Phi_1))
\big)[H_1,H_2],
\endaligned
\]
\[
\scriptsize
\aligned
\big[[H_1,H_2],\big[H_1,[H_1,H_2]\big]\big]
&
\overset{14}{=}
\big(
H_1(H_2(H_1(\Phi_2)))
-
H_2(H_1(H_1(\Phi_1)))
+
\\
&
\ \ \ \ \ 
+
2\Phi_1H_1(H_1(\Phi_2))
-
2\Phi_1H_2(H_1(\Phi_1))
\big)[H_1,H_2],
\endaligned
\]
\[
\scriptsize
\aligned
\big[[H_1,H_2],\big[H_1,[H_2,H_2]\big]\big]
&
\overset{15}{=}
\big(H_1(H_2(H_1(\Phi_2)))
-
H_2(H_1(H_1(\Phi_2)))
+
\\
&
\ \ \ \ \
+
\Phi_2H_1(H_2(\Phi_1))
-
\Phi_2H_2(H_1(\Phi_1))
+
\\
&
\ \ \ \ \
+
\Phi_1H_1(H_2(\Phi_2))
-
\Phi_1H_2(H_1(\Phi_2))
\big)[H_1,H_2],
\endaligned
\]
\[
\scriptsize
\aligned
\big[[H_1,H_2],\big[H_2,[H_2,H_2]\big]\big]
&
\overset{16}{=}
\big(
H_1(H_2(H_2(\Phi_2)))
-
H_2(H_1(H_2(\Phi_2)))
+
\\
&
\ \ \ \ \ 
+
2\Phi_2H_1(H_2(\Phi_2))
-
2\Phi_2H_2(H_1(\Phi_2))
\big)[H_1,H_2].
\endaligned
\]

\begin{Proposition}
\label{relations-Phi-H-I-V}
The two functions $\Phi_1$ and $\Phi_2$ identically satisfy:
\[ 
\small
\aligned
H_2(\Phi_1) 
\equiv 
H_1(\Phi_2)
\endaligned
\]
together with the following five third-order relations:
\[
\small
\aligned
0
&
\overset{\rm I}{\equiv}
-\,H_1(H_2(H_1(\Phi_2)))
+
2\,H_2(H_1(H_1(\Phi_2)))
-
H_2(H_2(H_1(\Phi_1)))
-
\\
&
\ \ \ \ \ 
-\,\Phi_2\,H_1(H_2(\Phi_1))
+
\Phi_2\,H_2(H_1(\Phi_1)),
\endaligned
\]
\[
\small
\aligned
0
&
\overset{\rm II}{\equiv}
-\,H_2(H_1(H_1(\Phi_2)))
+
2\,H_1(H_2(H_1(\Phi_2)))
-
H_1(H_1(H_2(\Phi_2)))
-
\\
&
\ \ \ \ \ 
-\,\Phi_1\,H_2(H_1(\Phi_2))
+
\Phi_1\,H_1(H_2(\Phi_2)),
\endaligned
\]
\[
\small
\aligned
0
&
\overset{\rm III}{\equiv}
-\,H_1(H_1(H_1(\Phi_2)))
+
2\,H_1(H_2(H_1(\Phi_1)))
-
H_2(H_1(H_1(\Phi_1)))
+
\\
&
\ \ \ \ \ 
+\Phi_1\,H_1(H_1(\Phi_2))
-
\Phi_1\,H_2(H_1(\Phi_1)),
\endaligned
\]
\[
\small
\aligned
0
&
\overset{\rm IV}{\equiv}
H_2(H_2(H_1(\Phi_2)))
-
2\,H_2(H_1(H_2(\Phi_2)))
+
H_1(H_2(H_2(\Phi_2)))
-
\\
&
\ \ \ \ \
-\,\Phi_2\,H_2(H_1(\Phi_2))
+
\Phi_2\,H_1(H_2(\Phi_2)),
\endaligned
\]
\[
\footnotesize
\aligned
0
&
\overset{\rm V}{\equiv}
H_1(H_1(H_2(\Phi_2)))
-
3\,H_1(H_2(H_1(\Phi_2)))
+
3\,H_2(H_1(H_1(\Phi_2)))
-
H_2(H_2(H_1(\Phi_1)))
-
\\
&
\ \ \ \ \
-\,
\Phi_2\,H_1(H_2(\Phi_1))
+
\Phi_2\,H_2(H_1(\Phi_1))
-
\Phi_1\,H_1(H_2(\Phi_2))
+
\Phi_1\,H_2(H_1(\Phi_2)),
\endaligned
\]
the first four being linearly independent, while the fifth
coincides with ${\rm I} - {\rm II}$. 
\end{Proposition}

\proof
Using the representations $\overset{ 1}{ =}$, \dots, $\overset{ 16}{
=}$ of the iterated brackets between $H_1$ and $H_2$ of lengths $\ell
= 6$, one may substitute them in the eleven free-Lie relations
$\overset{ 1}{ =}$, \dots, $\overset{ 11}{ =}$
Some of the obtained equations are redundant, and
some reduce to the trivial identity $0 = 0$.
\endproof

\begin{Corollary}
\label{0-Delta2-Delta3-2-Delta4}
The following two quantities are identically zero:
\[
\footnotesize
\aligned
0
&
\overset{a}{\equiv}
-\,H_2(H_2(H_1(\Phi_1)))
+
H_2(H_1(H_1(\Phi_2)))
+
H_1(H_2(H_1(\Phi_2)))
-
H_1(H_1(H_2(\Phi_2)))
+
\\
&
\ \ \ \ \
+
\Phi_1\,H_1(H_2(\Phi_2))
-
\Phi_1\,H_2(H_1(\Phi_2))
-
\Phi_2\,H_1(H_1(\Phi_2))
+
\Phi_2\,H_2(H_1(\Phi_1)),
\endaligned
\]
\[
\footnotesize
\aligned
0
&
\overset{b}{\equiv}
H_1(H_2(H_2(\Phi_2)))
-
2\,H_2(H_1(H_2(\Phi_2)))
+
H_2(H_2(H_1(\Phi_2)))
-
\\
&
\ \ \ \ \
-
H_2(H_1(H_1(\Phi_1)))
+
2\,H_1(H_2(H_1(\Phi_1)))
-
H_1(H_1(H_1(\Phi_2)))
+
\\
&
\ \ \ \ \
\Phi_1\,H_1(H_1(\Phi_2))
+
\Phi_2\,H_1(H_2(\Phi_2))
-
\Phi_2\,H_2(H_1(\Phi_2))
-
\Phi_1\,H_2(H_1(\Phi_1)).
\endaligned
\]
\end{Corollary}

\proof
Indeed, using the proposition, the first identity is just ${\rm I} +
{\rm II}$, while the second is just ${\rm III} + {\rm IV}$.
\endproof

\section{Cartan Connections in Terms of Coordinates and Bases}
\label{Cartan-connections-coordinates}

\HEAD{Cartan Connections in Terms of Coordinates and Bases}{
Mansour Aghasi, Joël Merker, and Masoud Sabzevari}

\subsection{Definition of Cartan connections
``à la Ehresmann'':}
Let $\K$ be either $\C$ or $\R$. Let $G$ be a local Lie group and let
$H$ be a local Lie subgroup of $G$. Denote by $\mathfrak{ g}$ and
$\mathfrak{ h}$ their respective Lie algebras which are $\K$-vector
spaces, with $\big[ \mathfrak{ h}, \mathfrak{ h} \big]_{\mathfrak{ g}}
\subset \mathfrak{ h}$. In order to set up a clear notational
distinction, we shall write $[ X, Y]$ for (Lie) brackets between vector
fields, and $[ {\sf x}, {\sf y}]_{ \mathfrak{ g}}$ for (Lie) brackets
between vectors of an abstract Lie algebra $\mathfrak{ g}$.
Following~\cite{ Kobayashi, Sharpe, Le}, we start with
a definition formulated independently of any coordinate
system or of any basis for $\mathfrak{ g}$.

\begin{Definition}
A {\sl Cartan geometry} on a $\mathcal{ C}^\infty$ manifold $M$ (over
$\K$) modeled on $(\mathfrak{ g}, \mathfrak{ h})$ consists of the
following data:

\medskip\noindent$\bullet$ 
a principal $H$ bundle $\mathcal{P} \to M$ over $M$, the
right action of $H$
on $\mathcal{ P}$
being on the right:
\[
R_h
\colon\
p\mapsto ph 
\ \ \ \ \ \ \ \ \ \ \ \ \ 
{\scriptstyle{(p\,\in\,\mathcal{P})}};
\]

\medskip\noindent$\bullet$ 
a $\mathfrak{ g}$-valued $1$-form $\omega$ on 
$\mathcal{ P}$ which enjoys the following three properties:

\begin{itemize}

\smallskip\item[{\bf (i)}]
for every point $p \in \mathcal{ P}$, the linear map:
\[
\omega_p\colon 
T_p\mathcal{P}
\longrightarrow
\mathfrak{g}
\]
is an isomorphism; 

\smallskip\item[{\bf (ii)}]
if, for every element ${\sf y} \in \mathfrak{ h}$, one defines the
fundamental vector field $Y^\dag$ on $\mathcal{ P}$ by differentiating the
action of $\exp (t{\sf y}) \in H$ on $\mathcal{ P}$:
\[
Y^\dag\vert_p
:=
\frac{d}{dt}\Big\vert_0
\Big(
p\exp{(t{\sf y})}
\Big)
\ \ \ \ \ \ \ \ \ \ \ \ \ 
{\scriptstyle{(p\,\in\,\mathcal{P})}},
\]
then $\omega$ satisfies:
\[
\omega(Y^\dag)
=
{\sf y};
\] 

\smallskip\item[{\bf (iii)}]
again with the right translation $R_h \colon p \mapsto p \, h$ on
$\mathcal{ P}$ by means of an element $h \in H$, the $\mathfrak{
g}$-valued $1$-form $\omega$ satisfies:
\[
\omega_{ph}\big(R_{h*}(v_p)\big)
=
{\rm Ad}(h^{-1})
\big[\omega_p(v_p)\big],
\]
for every tangent vector $v_p \in T_p \mathcal{ P}$ at every point $p
\in \mathcal{ P}$.

\end{itemize}

\end{Definition}

Assuming that a Cartan connection $1$-form
$\omega \colon T\mathcal{ P} \to \mathfrak{ g}$ is given, our main aim
in the next paragraphs will be to express more concretely its
properties in terms of a certain local coordinate systems on $M$, $H$,
$\mathcal{ P}$, and in terms of a fixed basis for $\mathfrak{ g}$.

\subsection{First consequences}
The way $\omega$ behaves through right translations {\bf (iii)} may
also be abbreviated without arguments as:
\[
R_h^*(\omega)
=
{\rm Ad}(h^{-1})\circ \omega,
\]
where the composition is:
\[
T\mathcal{P}
\overset{\omega}{\longrightarrow}
\mathfrak{g}
\overset{{\rm Ad}(h^{-1})}{\longrightarrow}
\mathfrak{g}.
\]

Of course, the principal bundle $\mathcal{P}$ is foliated by copies of
$H$. It will be convenient to denote by:
\[
\mathcal{H}_p
:=
\{ph\colon h\in H\}
\] 
the $H$-orbit ($\simeq H$) of an arbitrary point $p$ in the
$H$-principal bundle $\mathcal{P}$, which is a $\mathcal{ C}^\infty$
submanifold of $\mathcal{ P}$. Also, we shall denote by $\mathcal{
H}$ the whole foliation of $\mathcal{ P}$ by these copies of $H$. Then
property {\bf (ii)} means that each:
\[
\omega_p
\colon
T_p\mathcal{H}_p
\overset{\simeq}{\longrightarrow}
\mathfrak{h}
\ \ \ \ \ \ \ \ \ \ \ \ \
{\scriptstyle{(p\,\in\,\mathcal{P})}}
\]
is the {\em identity} isomorphism, if one interprets $T_p \mathcal{
H}_p$ as the tangent space to the Lie group copy $\mathcal{ H}_p
\simeq H$.

\subsection{Curvature $2$-form and curvature
function} For reasons of clear
notational distinction of objects that live in
different spaces, we shall
always denote by $\widehat{ X}, \widehat{ Y}, \widehat{ Z}$ or by
$\widetilde{ X}$, $\widetilde{ Y}$, $\widetilde{ Z}$ vector fields on
the bundle $\mathcal{ P}$, while vector fields $X, Y, Z$ on $M$ will
be systematically denoted without any hat or tilde.

Notably, for any vector ${\sf x} \in \mathfrak{ g}$, the inverse
image of ${\sf x}$ through $\omega$ at any point $p \in \mathcal{ P}$,
namely:
\[
\widehat{X}_p
:=
\omega_p^{-1}({\sf x})
\ \ \ \ \ \ \ \ \ \ \ \ \
{\scriptstyle{(p\,\in\,\mathcal{P})}}
\]
provides, as the point $p$ varies all over $\mathcal{ P}$, a
well-defined $\mathcal{ C}^\infty$ vector field which is sometimes
called the {\sl constant vector field}: 
\[
\widehat{X}
:=
\omega^{-1}({\sf x})
\]
associated to ${\sf x}$, the slight abuse of notation ``$\omega^{ -1}
( {\sf x})$'' being admissible here. Also, because all the $\omega_p
\colon T_p \mathcal{ P} \to \mathfrak{ g}$ are isomorphisms, for any
choice of a vector space basis $({\sf x}_k)_{ 1 \leqslant k \leqslant
\dim_\K \mathfrak{ g}}$ for $\mathfrak{ g}$, the collection of the $r$
vector fields:
\[
\widehat{X}_k
:=
\omega_p^{-1}({\sf x}_k)
\ \ \ \ \ \ \ \ \ \ \ \ \
{\scriptstyle{(k\,=\,1\,\cdots\,r)}}
\] 
visibly makes a global frame on $\mathcal{ P}$, that is to say, at
every point $p \in \mathcal{ P}$, the vectors $\widehat{ X}_1\vert_p$,
\dots, $\widehat{ X}_r \vert_p$ make a basis of $T_p \mathcal{
P}$. Furthermore, property {\bf (ii)} that the $\mathfrak{ g}$-valued
Cartan-connection $1$-form $\omega$ should satisfy implies that for
every element ${\sf y} \in \mathfrak{ h}$, one has at every point $p
\in \mathcal{ P}$:
\[
Y_p^\dag
=
\omega_p^{-1}({\sf y})
=
\widehat{Y}_p,
\]
so that (only) for these elements ${\sf y} \in \mathfrak{ h}$, 
one has the coincidence of vector fields on $\mathcal{ P}$:
\[
Y^\dag
=
\widehat{Y}
\ \ \ \ \ \ \ \ \ \ \ \ \
{\scriptstyle{({\sf y}\,\in\,\mathfrak{h})}}.
\]

\begin{Definition}
\label{definition-curvature-form} 
The {\sl curvature form} of the Cartan connection is the $2$-form on
$\mathcal{ P}$ which acts on pairs of vectors $(\widetilde{ X}_p,
\widetilde{ Y}_p)$ based at an arbitrary point $p \in \mathcal{ P}$
through the formula:
\[
\Omega_p
\big(\widetilde{X}_p,\,\widetilde{Y}_p\big)
:=
d\omega_p\big(\widetilde{X},\,\widetilde{Y}\big)
+
\big[\omega_p\big(\widetilde{X}\big),\,
\omega_p\big(\widetilde{Y}\big)\big]_{\mathfrak{g}},
\]
where $\widetilde{ X}$ and $\widetilde{ Y}$ denote any two local
$\mathcal{ C}^\infty$ extensions near $p$ satisfying $\widetilde{
X}\big\vert_p = \widetilde{ X}_p$ and $\widetilde{ Y}\big\vert_p =
\widetilde{ Y}_p$; the obtained value $\Omega_p
\big(\widetilde{X}_p,\, \widetilde{Y}_p\big)$ is easily seen to be
independent of these extensions, and also to be skew-symmetric.
\end{Definition}

\begin{Lemma}
The curvature vanishes as soon as at least one of its two arguments is
tangent to the $H$-principal bundle foliation $\mathcal{ H}$, namely:
\[
0
=
\Omega_p
\big(
\widetilde{X}_p,\,\widetilde{Y}_p
\big)
\]
whenever either $\widetilde{ X}_p \in T_p \mathcal{ H}_p$ or
$\widetilde{ Y}_p \in T_p \mathcal{ H}_p$.
\end{Lemma}

\proof
Since the $\omega_p^{ -1} ( {\sf y})$ span $T_p \mathcal{ H}_p$ when
${\sf y}$ varies in $\mathfrak{ h}$, and since $\Omega$ is
skew-symmetric, it suffices in fact to establish this vanishing
curvature property:
\[
0
=
\Omega\big(\widehat{X},\widehat{Y}\big),
\]
for any constant vector field $\widehat{ X} = \omega^{ -1} ({\sf x})$
with ${\sf x} \in \mathfrak{ g}$ arbitrary and any vertical constant
vector field $\widehat{ Y} = \omega^{ -1} ( {\sf y})$ with ${\sf y}
\in \mathfrak{ h}$ arbitrary. But then, this is a known consequence
of {\bf (iii)}.
\endproof

Now, at each point $p \in \mathcal{ P}$, for every ${\sf x} \in
\mathfrak{ g}$, the inverse image through $\omega_p$ of the element
${\sf x} + \mathfrak{ h}$ of the quotient $\mathfrak{ g} / \mathfrak{
h}$:
\[
\omega_p^{-1}
\big({\sf x}+\mathfrak{h}\big)
=
\omega_p^{-1}({\sf x})
+
\omega_p^{-1}(\mathfrak{h})
=
\omega_p^{-1}({\sf x})+T_p\mathcal{H}_p
\]
is defined exactly up to the $H$-principal bundle tangent space. 
By bilinearity, it
then follows immediately from the above proposition
that
\[
\Omega_p
\big(
\omega_p^{-1}({\sf x}'+\mathfrak{h}),\,
\omega_p^{-1}({\sf x}''+\mathfrak{h})
\big)
=
\Omega_p
\big(\omega_p^{-1}({\sf x}'),\,
\omega_p^{-1}({\sf x}'')\big).
\] 
In other words, the curvature $2$-form acts in fact on the quotient $T
\mathcal{ P} \big/ T \mathcal{ H}$ by the vertical $H$-foliation.
This observation yields a path to a more concrete access to Cartan
curvature which is quite useful for effective computations.

\begin{Definition}
\label{definition-curvature-function} 
The {\sl curvature function} of a Cartan connection $\omega \colon
T\mathcal{ P} \to \mathfrak{ g}$ is the map:
\[
\kappa
\in
\mathcal{C}^\infty
\big(
\mathcal{P},\,
{\rm End}
\big(\Lambda^2(\mathfrak{g}/\mathfrak{h}),\,\mathfrak{g}\big)
\big)
\]
which sends every point $p \in \mathcal{ P}$ to the following
$\mathfrak{ g}$-valued alternating bilinear
map:
\[
\aligned
(\mathfrak{g}/\mathfrak{h})
\wedge
(\mathfrak{g}/\mathfrak{h})
&
\longrightarrow
\mathfrak{g}
\\
\kappa(p)
\colon
\ \ \ \ \ \ \ \ \ \ \ \ \ \ \ \ \ 
\big({\sf x}'\,{\rm mod}\,\mathfrak{h}\big)
\wedge
\big({\sf x}''\,{\rm mod}\,\mathfrak{h}\big)
&
\longmapsto
\Omega_p
\big(
\omega_p^{-1}({\sf x}'),\,
\omega_p^{-1}({\sf x}'')
\big)
\\
&
\ \ \ \,
=
\Omega_p\big(\widehat{X}_p',\widehat{X}_p''\big).
\endaligned
\]
\end{Definition}

Conversely, 
one easily convinces oneself that the curvature function 
determines the curvature form in a unique way.

\begin{Lemma}
\label{first-curvature-function}
For any two representatives ${\sf x}', {\sf x}'' \in \mathfrak{ g}$ of
two elements of $\mathfrak{ g} / \mathfrak{ h}$, one has:
\[
\kappa(p)({\sf x}',{\sf x}'')
=
[{\sf x}',{\sf x}'']_{\mathfrak{g}}
-
\omega_p
\big(
\big[\omega^{-1}({\sf x}'),\,\omega^{-1}({\sf x}'')\big]
\big).
\]
\end{Lemma}

\proof
Starting from the definition just stated and applying the so-called
{\sl Cartan formula}\footnote{\,
{\em See}~\cite{ Sharpe}, p.~58: if $\omega$ is an arbitrary
$1$-form on a smooth manifold $N$ valued in some $\K$-vector space
$V$, then $d\omega (X, Y) = X ( \omega (Y)) - Y ( \omega (X)) - \omega
([X, Y])$ for any two vector fields $X$ and $Y$ on $M$.
} 
to expand the $d\omega_p$-term, we easily get:
\[
\footnotesize
\aligned
\kappa(p)({\sf x}',{\sf x}'')
&
=
\Omega_p(\widehat{X}',\,\widehat{X}'')
\\
&
=
d\omega_p\big(\widehat{X}',\widehat{X}'')
+
\big[
\omega_p(\widehat{X}'),\,
\omega_p(\widehat{X}'')
\big]_{\mathfrak{g}}
\\
&
=
\widehat{X}'
\big(\omega_p(\widehat{X}'')\big)
-
\widehat{X}''
\big(\omega_p(\widehat{X}')\big)
-
\omega_p\big(
[\widehat{X}',\,\widehat{X}'']\big)
+
[{\sf x}',{\sf x}'']_{\mathfrak{g}}
\\
&
=
\zero{\widehat{X}'
\big({\sf x}''\big)}
-
\zero{\widehat{X}''
\big({\sf x}'\big)}
-
\omega_p\big(
\big[\omega^{-1}({\sf x}'),\,
\omega^{-1}({\sf x}'')\big]\big)
+
[{\sf x}',{\sf x}'']_{\mathfrak{g}},
\endaligned
\]
where underlined terms vanish as do all differentiated
constants.
\endproof

\begin{Proposition}
\label{equivariancy-curvature-function}
For any $h \in H$, the curvature function enjoys the 
${\rm Ad}$-equivariancy property:
\[
\kappa\big(ph\big)
\big({\sf x}'\,{\rm mod}\,\mathfrak{h},\,
{\sf x}''\,{\rm mod}\,\mathfrak{h}\big)
=
{\rm Ad}(h^{-1})
\big[
\kappa(p)
\big({\rm Ad}(h)\,{\sf x}',\,{\rm Ad}(h)\,{\sf x}''\big)
\big],
\]
where ${\sf x}'\, {\rm mod}\, \mathfrak{ h}$ and ${\sf x}'\, {\rm
mod}\, \mathfrak{ h}$ are any two elements of $\mathfrak{ g} /
\mathfrak{ h}$. Furthermore, for any fundamental field $Y^\dag =
\frac{ d}{ dt}\big\vert_0 R_{ \exp(t{\sf y})}$ on $\mathcal{ P}$
associated to an arbitrary ${\sf y} \in \mathfrak{ h}$, one has:
\[
(Y^\dag\kappa)(p)({\sf x}',{\sf x}'')
=
-\big[{\sf y},\,\kappa(p)({\sf x}',{\sf x}'')\big]_{\mathfrak{g}}
+
\kappa(p)\big([{\sf y},{\sf x}']_{\mathfrak{g}},\,
{\sf x}''\big)
+
\kappa(p)\big({\sf x}',\,
[{\sf y},{\sf x}'']_{\mathfrak{g}}\big),
\]
where the two arguments ${\sf x}'$ and
${\sf x}''$ of the curvature function
$\kappa (p)$ are (implicitly) taken ${\rm mod}\, 
\mathfrak{ h}$.
\end{Proposition}

\proof
First of all, the right-equivariancy of the connection
form and of the curvature may be read as two equations:
\begin{equation}
\label{reformulation-double-equivariancy}
\aligned
\omega_{ph}\big((R_h)_{\ast p}(\widetilde{X}_p)\big)
&
=
{\rm Ad}(h^{-1})
\big[
\omega_p\big(\widetilde{X}_p\big)
\big],
\\
\Omega_{ph}\big(
(R_h)_{*p}(\widetilde{X}_p'),\,\,
(R_h)_{*p}(\widetilde{X}_p'')
\big)
&
=
{\rm Ad}(h^{-1})
\big[
\Omega_p\big(\widetilde{X}_p',\,\widetilde{X}_p''\big)
\big],
\endaligned
\end{equation}
valid for any $h \in H$ and for any three (vector) fields $\widetilde{
X}_p, \widetilde{ X}_p', \widetilde{ X}_p'' \in T_p \mathcal{ P}$
based at an arbitrary point $p \in \mathcal{ P}$. So, let us apply
the first equivariancy in which we replace $\widetilde{ X}_p$ by the
constant vector (field) $\widehat{ X}_p' := \omega_p^{-1} ({\sf x}')$
associated to an arbitrary ${\sf x}' \in \mathfrak{ g}$:
\[
\omega_{ph}\big((R_h)_{*p}\big(\omega_p^{-1}({\sf x}')\big)\big)
=
\omega_{ph}\big((R_h)_{*p}\big(\widehat{X}_p'\big)\big)
=
{\rm Ad}(h^{-1})
\big[
\omega_p\big(\widehat{X}_p'\big)
\big]
=
{\rm Ad}(h^{-1})
[{\sf x}'],
\]
whence equivalently if we apply $\omega_{ ph}^{-1}$ to both extreme sides:
\[
(R_h)_{*p}
\big(\omega_p^{-1}({\sf x}')\big)
=
\omega_{ph}^{-1}
\big({\rm Ad}(h^{-1})\,{\sf x}'\big).
\]
If we now replace ${\sf x}'$ by just ${\sf Ad} (h)\, {\sf x}'$
here, we get the useful formula:
\begin{equation}
\label{useful-R-omega}
(R_h)_{*p}
\big(
\omega_p^{-1}\big({\rm Ad}(h)\,{\sf x}'\big)
\big)
=
\omega_{ph}^{-1}({\sf x}').
\end{equation}
Thanks to this preliminary, we may now compute:
\[
\aligned
\kappa(ph)({\sf x}',{\sf x}'')
&
\overset{\rm def}{=}
\Omega_{ph}\big(\omega_{ph}^{-1}({\sf x}'),\,
\omega_{ph}^{-1}({\sf x}'')\big)
\\
\explain{apply~\thetag{\ref{useful-R-omega}}}
\ \ \ \ \ \ \ \ 
&
=
\Omega_{ph}\Big(
(R_h)_{*p}\big(\omega_p^{-1}\big({\rm Ad}(h)\,{\sf x}'
\big)\big),\,\,
(R_h)_{*p}\big(\omega_p^{-1}\big({\rm Ad}(h)\,{\sf x}''\big)\big)
\Big)
\\
\explain{use~\thetag{\ref{reformulation-double-equivariancy}}}
\ \ \ \ \ \ \ \ \ \ \,
&
=
{\rm Ad}(h^{-1})
\big[
\Omega_p\big(
\omega_p^{-1}\big({\rm Ad}(h)\,{\sf x}'\big),\,\,
\omega_p^{-1}\big({\rm Ad}(h)\,{\sf x}''\big)\big)
\big]
\\
&
\overset{\rm def}{=}
{\rm Ad}(h^{-1})
\big[
\kappa(p)\big(
{\rm Ad}(h)\,{\sf x}',\,
{\rm Ad}(h)\,{\sf x}''\big)
\big],
\endaligned
\]
which completes the verification of the first formula claimed in the
proposition.

Next, let us assume that the translational element $h \in H$ is of the
form $h_t := \exp ( t{\sf y})$, for some ${\sf y} \in \mathfrak{ h}$,
where the parameter $t \in \K$ is small, so that we have, after
multiplying both sides by ${\rm Ad}(h_t)$:
\[
{\rm Ad}(h_t)\big[\kappa(ph_t)({\sf x}',{\sf x}'')\big]
=
\kappa(p)\big({\rm Ad}(h_t)\,{\sf x}',\,{\rm Ad}(h_t)\,{\sf x}''\big).
\]
From the very definition of
the vector field $Y^\dag = \frac{ d}{ dt} \big\vert_0 {\sf R}_{ \exp (
t {\sf y})}$ on $\mathcal{ P}$, it classically follows that for any
function $f \colon \mathcal{ P} \to \C$, one has:
\[
\frac{d}{dt}\Big\vert_0
\Big(
f\big(
R_{\exp(t{\sf y})}(p)\big)
\Big)
=
\frac{d}{dt}\Big\vert_0
\Big(
f\big(
p\exp(t{\sf y})\big)
\Big)
=
(Y^\dag f)(p).
\]
Consequently, if we apply $\frac{ d}{ dt}\big\vert_0$ to the above
identity, using ${\rm Ad}(h_0) = {\rm Id}$,
using Leibniz' formula and using bilinearity of the
curvature function, we get:
\[
\aligned
&
\frac{d}{dt}\Big\vert_0\,{\rm Ad}(\exp(t{\sf y}))
\big[\kappa(p)({\sf x}',{\sf x}'')\big]
+
(Y^\dag\kappa)(p)({\sf x}',{\sf x}'')
=
\\
&
\ \ \ \ \ \ \ \ \ \ \ \ \ \ \ \ \
=
\kappa(p)
\big(
{\textstyle{\frac{d}{dt}}}\big\vert_0\,{\rm Ad}(\exp(t{\sf y}))\,{\sf x}',\,\,
{\sf x}''
\big)
+
\kappa(p)
\big(
{\sf x}',\,
{\textstyle{\frac{d}{dt}}}\big\vert_0\,{\rm Ad}(\exp(t{\sf y}))\,{\sf x}''
\big).
\endaligned
\]
Applying then twice within the right hand side the classical identity:
\[
\frac{d}{dt}\Big\vert_0
{\rm Ad}\big(\exp(t{\sf y})\big)[{\sf z}]
=
[{\sf y},{\sf z}]_{\mathfrak{g}}
=
{\rm ad}({\sf y})({\sf z}),
\] 
which can also be thought of as defining the adjoint representation
(restricted to the subalgebra ${\sf h} \subset \mathfrak{ g}$):
\[
{\rm ad}\colon\mathfrak{h} 
\longrightarrow {\rm End} ( \mathfrak{ g}), 
\ \ \ \ \ \ \ \
{\sf y}
\longmapsto
\big({\sf z}\mapsto[{\sf y},{\sf z}]_{\mathfrak{g}}\big),
\]
and putting $(Y^\dag \kappa)$ single on the left, we finally obtain
the second identity claimed in the proposition:
\[
(Y^\dag\kappa)(p)({\sf x}',{\sf x}'')
=
-\big[{\sf y},\,\kappa(p)({\sf x}',{\sf x}'')\big]_{\mathfrak{g}}
+
\kappa(p)\big([{\sf y},{\sf x}']_{\mathfrak{g}},\,
{\sf x}''\big)
+
\kappa(p)\big({\sf x}',\,
[{\sf y},{\sf x}']_{\mathfrak{g}}\big).
\]
In the case where $H$ and $G$ are connected and simply connected, a
natural converse holds by integrating these differential relations,
but as we will work with local Lie groups, we will not need such a
converse.
\endproof

Computing this curvature function explicitly in 
coordinates and with a basis of $\mathfrak{ g}$ will 
exhibit much more the algebraic complexity it can possess.

\subsection{Lie algebra bases}
In order to achieve the construction of Cartan connections for
specific new geometric structures, one must work and compute 
effectively with some concrete
bases for $\mathfrak{ g}$ and $\mathfrak{ h}$. Denote then by:
\[
r:=\dim_\K\,\mathfrak{g},
\ \ \ \ \
n:=\dim_\K\,\big(\mathfrak{g}/\mathfrak{h}\big)
\ \ \ \ \
\text{\rm and}
\ \ \ \ \
n-r=\dim_\K\,\mathfrak{h},
\]
suppose $r \geqslant 2$, $n \geqslant 1$, $n-r \geqslant 1$ so that
$\mathfrak{ g}$, $\mathfrak{ g} / \mathfrak{ h}$ and $\mathfrak{ h}$
are all nonzero. For an adapted basis $({\sf x}_k)_{1\leqslant
k\leqslant r}$ which may be thought of, in concrete examples, to enjoy
some useful normalizations/simplification (think of root systems for
semi-simple Lie algebras), we then have:
\[
\aligned
\mathfrak{g}
&
=
{\rm Span}_\K
\big({\sf x}_1,\dots,{\sf x}_n,{\sf x}_{n+1},\dots,{\sf x}_r\big),
\\
\mathfrak{h}
&
=
\ \ \ \ \ \ \ \ \ \ \ \ \ \ \ \ \ \ 
{\rm Span}_\K
\big({\sf x}_{n+1},\dots,{\sf x}_r\big).
\endaligned
\]
Also, we may suppose that $(n-r)$ representatives in $\mathfrak{ g}
/ \mathfrak{ h}$ are just:
\[
{\sf y}_1
:=
{\sf x}_{n+1},\
\dots\dots,\
{\sf y}_{n-r}
:=
{\sf x}_r,
\]
and we shall have to remember the notational coincidences:
\[
{\sf y}_j
\equiv
{\sf x}_{n+j}
\ \ \ \ \ \ \ \ \ \ \ \ \
{\scriptstyle{(j\,=\,1\,\cdots\,r\,-\,n)}}.
\]

Next, let $\mathfrak{ g}^* = {\rm Lin} ( \mathfrak{ g}, \C)$ denote
the dual of the Lie algebra $\mathfrak{ g}$, viewed as a plain vector
space (it has no natural Lie bracket structure). If we introduce the
basis of $\mathfrak{ g}^*$:
\[
{\sf x}_1^*,\dots,{\sf x}_n^*,{\sf x}_{ n+1}^*,\dots,{\sf x}_r^*
\] 
which is dual to the previously fixed basis ${\sf x}_1, \dots, {\sf
x}_n, {\sf x}_{ n+1}, \dots, {\sf x}_r$ of $\mathfrak{ g}$, then by
definition, for any $i, j = 1, \dots, n, n+1, \dots, r$, we have:
\[
{\sf x}_i^*({\sf x}_j)
=
\delta_j^i
:=
\left\{
\aligned
&
1\ \ 
\text{\rm if}\ \
i=j,
\\
&
0\ \ 
\text{\rm otherwise}.
\endaligned\right. 
\]
Of course, $\mathfrak{ h}^*$ is then spanned by ${\sf x}_1^*, \dots,
{\sf x}_{r-n}^*$. 

For later use, we remind also that if $E$ is a finite-dimensional
$\K$-vector space and if $\omega^*$, $\pi^*$ are one-forms
belonging to its dual $E^* = 
{\rm Lin} ( E, \C)$, then the two-form $\omega^* \wedge \pi^*$
acts on pairs $(e, f) \in E^2$ by definition as:
\[
\omega^*\wedge\pi^*(e,f)
\overset{\rm def}{=}
\omega^*(e)\,\pi^*(f)-\omega^*(f)\,\pi^*(e).
\] 
In particular, for any $i_1$, $i_2$ with 
$i_1 < i_2$ and for any $j_1$, $j_2$ without restriction, we have: 
\begin{equation}
\label{i-1-i-2-j-1-j-2}
\aligned
{\sf x}_{i_1}^*\wedge{\sf x}_{i_2}^*
({\sf x}_{j_1},{\sf x}_{j_2})
&
=
{\sf x}_{i_1}^*({\sf x}_{j_1})\,
{\sf x}_{i_2}^*({\sf x}_{j_2})
-
{\sf x}_{i_1}^*({\sf x}_{j_2})\,
{\sf x}_{i_2}^*({\sf x}_{j_1})
\\
&
=
\delta_{j_1}^{i_1}\,\delta_{j_2}^{i_2}
-
\delta_{j_2}^{i_1}\,\delta_{j_1}^{i_2}.
\endaligned
\end{equation}

With any such choice of a basis, brackets of the Lie algebra
$\mathfrak{ g}$ are encoded by the so-called {\sl structure
constants}:
\[
\aligned
&
\big[{\sf x}_{k_1},\,{\sf x}_{k_2}\big]_{
\mathfrak{g}}
=
\sum_{s=1}^r\,
c_{k_1,k_2}^s\,{\sf x}_s
\\
&
\ \ \ \ \ \
{\scriptstyle{(k_1,\,k_2\,=\,1\,\cdots\,n,\,n\,+\,1,\dots,r)}}.
\endaligned
\]
Of course, the skew-symmetry and the Jacobi identity:
\begin{equation}
\label{skew-Jacobi}
\aligned
0
&
=
[{\sf x}_{k_1},{\sf x}_{k_2}]_{\mathfrak{g}}
+
[{\sf x}_{k_1},{\sf x}_{k_2}]_{\mathfrak{g}}
\\
0
&
=
\big[
[{\sf x}_{k_1},{\sf x}_{k_2}]_{\mathfrak{g}},
{\sf x}_{k_3}\big]_{\mathfrak{g}}
+
\big[
[{\sf x}_{k_3},{\sf x}_{k_1}]_{\mathfrak{g}},
{\sf x}_{k_2}\big]_{\mathfrak{g}}
+
\big[
[{\sf x}_{k_2},{\sf x}_{k_3}]_{\mathfrak{g}},
{\sf x}_{k_1}\big]_{\mathfrak{g}}
\endaligned
\end{equation}
read at the level of structure constants as:
\[
\aligned
0
&
=
c_{k_1,k_2}^s
+
c_{k_2,k_1}^s
\\
0
&
=
\sum_{s=1}^r\,
\big(
c_{k_1,k_2}^s\,c_{s,k_3}^l
+
c_{k_3,k_1}^s\,c_{s,k_2}^l
+
c_{k_2,k_3}^s\,c_{s,k_1}^l
\big)
\\
&
\ \ \ \ \ \ \ \ \ \ \ \ \ \ \ \ \ \ \ \ \ \ \
{\scriptstyle{(k_1,\,\,k_2,\,\,k_3,\,\,l\,=\,1\,\cdots\,r)}}.
\endaligned
\]

Since $\mathfrak{ h}$ is a Lie subalgebra of $\mathfrak{ g}$, one has
$c_{ k_1, k_2}^s = 0$ whenever $n+1
\leqslant k_1, k_2 \geqslant r$, for any $s
\leqslant n$. In terms of the generators ${\sf y}_j = {\sf x}_{ n+j}$
of $\mathfrak{ h}$, we may also write the Lie algebra structure of
$\mathfrak{ h}$ under the form:
\[
\big[{\sf y}_{j_1},\,{\sf y}_{j_2}\big]_{\mathfrak{g}}
=
\sum_{t=1}^{r-n}\,
d_{j_1,j_2}^t\,{\sf y}_t
\ \ \ \ \ \ \ \ \ \ \ \ \
{\scriptstyle{(j_1,\,j_2,\,t\,=\,1\,\cdots\,r\,-\,n)}},
\]
where we admit the notational coincidences:
\[
d_{j_1,j_2}^t
=
c_{n+j_1,n+j_2}^{n+t}
\ \ \ \ \ \ \ \ \ \ \ \ \
{\scriptstyle{(j_1,\,j_2,\,t\,=\,1\,\cdots\,r\,-\,n)}}.
\]

\subsection{Expansion of the $\mathfrak{ g}$-valued Cartan connection
form $\omega$} In terms of the basis $({\sf x}_k)_{ 1 \leqslant k
\leqslant r}$ of $\mathfrak{ g}$, the $\mathfrak{ g}$-valued $1$-form
$\omega$ expands under the shape:
\[
\omega
=
\sum_{k=1}^r\,\omega^k\,{\sf x}_k,
\] 
where the $\omega^k$ are {\em scalar}, {\em i.e.} $\K$-valued,
$1$-forms on $\mathcal{P}$. If $\varpi$ is another $\mathfrak{
g}$-valued $1$-form that one wants to ``wedge'' with $\omega$, it is
natural to simultaneously ``wedge'' the $1$-forms and let the Lie
bracket act on the basis elements ${\sf x}_k$ as follows:
\[
\aligned
\omega\wedge_{\mathfrak{g}}\varpi
:=
&\,
\Big(
\sum_{k_1=1}^r\,\omega^{k_1}\,{\sf x}_{k_1}
\Big)
\wedge_{\mathfrak{g}}
\Big(
\sum_{k_2=1}^r\,\varpi^{k_2}\,{\sf x}_{k_2}
\Big)
\\
=
&\,
\sum_{k_1=1}^r\,\sum_{k_2=1}^r\,
\omega^{k_1}\wedge\varpi^{k_2}\,
[{\sf x}_{k_1},\,{\sf x}_{k_2}]_{\mathfrak{g}}.
\endaligned
\]
We employ here the intuitively clear notation ``$\wedge_{ \mathfrak{
g}}$'' in order to expressly notify that both a wedge product and a
Lie bracket act simultaneously. With this natural definition, it then
follows from the classical formula:
\[
(\omega^{k_1}\wedge\varpi^{k_2})
\big(\widetilde{X},\widetilde{Y}\big)
=
\omega^{k_1}(\widetilde{X})\,
\varpi^{k_2}(\widetilde{Y})
-
\omega^{k_1}(\widetilde{Y})\,
\varpi^{k_2}(\widetilde{X})
\]
that $\omega \wedge_{ \mathfrak{ g}} \varpi$ acts on pairs 
$(\widetilde{ X}, \widetilde{ Y})$ of vector
fields on $\mathcal{ P}$ just as:
\[
\footnotesize
\aligned
(\omega\wedge_{\mathfrak{g}}\varpi)
\big(\widetilde{X},\widetilde{Y}\big)
&
=
\sum_{k_1,k_2=1}^r\,
\big(
\omega^{k_1}(\widetilde{X})\,\varpi^{k_2}(\widetilde{Y})
-
\omega^{k_2}(\widetilde{X})\,\varpi^{k_1}(\widetilde{Y})
\big)\,
[{\sf x}_{k_1},{\sf x}_{k_2}]_{\mathfrak{g}}
\\
&
=
\big[
\omega(\widetilde{X}),\,
\varpi(\widetilde{Y})
\big]_{\mathfrak{g}}
-
\big[
\omega(\widetilde{Y}),\,
\varpi(\widetilde{X})
\big]_{\mathfrak{g}}.
\endaligned
\]
In particular the connection form $\omega$
wedged with itself acts as:
\[
\{\omega\wedge_{\mathfrak{g}}\omega\}
\big(\widetilde{X},\widetilde{Y}\big)
=
2\,
\big[
\omega(\widetilde{X}),\,
\omega(\widetilde{Y})
\big]_{\mathfrak{g}},
\]
from which we deduce the following alternative, equivalent 
expression of the curvature.

\begin{Lemma}
The curvature $2$-form $\Omega \in \Lambda^2 \big(
T^* \mathcal{ P}, \mathfrak{ g}\big)$ writes shortly as:
\[
\Omega
=
d\omega
+
{\textstyle{\frac{1}{2}}}\,
\omega\wedge_{\mathfrak{g}}\omega,
\]
and also in terms of a basis $({\sf x}_k)_{ 1\leqslant k \leqslant r}$
for $\mathfrak{ g}$, just as:
\[
\Omega
=
\sum_{k=1}^r\,
\Big\{
d\omega^k
+
\sum_{1\leqslant l_1<l_2\leqslant r}\,
c_{l_1,l_2}^k\,\omega^{l_1}\wedge\omega^{l_2}
\Big\}\,
{\sf x}_k.
\]
\end{Lemma}

\proof
We compute $\Omega$ and reorganize its expression using the natural
linearly independent frame $\omega^{ l_1} \wedge \omega^{ l_2}$, $l_1
< l_2$, on the space of $2$-forms on $\mathcal{ P}$:
\[
\aligned
\Omega
&
=
d\omega
+
{\textstyle{\frac{1}{2}}}\,
\omega\wedge_{\mathfrak{g}}\omega
\\
&
=
\sum_{k=1}^r\,d\omega^k\,{\sf x}_k
+
{\textstyle{\frac{1}{2}}}\,
\sum_{l_1=1}^r\,\sum_{l_2=1}^r\,
\omega^{l_1}\wedge\omega^{l_2}\,
[{\sf x}_{l_1},{\sf x}_{l_2}]_{\mathfrak{g}}
\\
&
=
\sum_{k=1}^r\,
\Big\{
d\omega^k
+
{\textstyle{\frac{1}{2}}}\,
\sum_{l_1=1}^r\,\sum_{l_2=1}^r\,
c_{l_1,l_2}^k\,\omega^{l_1}\wedge\omega^{l_2}
\Big\}\,
{\sf x}_k
\\
\explain{use $c_{l_2,l_1}^k = - c_{l_1, l_2}^k$}
\ \ \ \ \ \ \
&
=
\sum_{k=1}^r\,
\Big\{
d\omega^k
+
\sum_{1\leqslant l_1<l_2\leqslant r}\,
c_{l_1,l_2}^k\,\omega^{l_1}\wedge\omega^{l_2}
\Big\}\,
{\sf x}_k,
\endaligned
\]
which is what was claimed.
\endproof

{\em However}, this representation of the curvature
is not at all finalized for at least two reasons:

\smallskip\noindent$\bullet$
each $d\omega^k$ should in principle be re-expressed in terms the
basis $\big( \omega^{ l_1} \wedge \omega^{ l_2}
\big)_{ 1\leqslant l_1 < l_2 \leqslant r}$ when
one really wants to perform explicit computations; 

\smallskip\noindent$\bullet$
more deeply, any $2$-form as the curvature $\Omega$ should be
expressed as a $2$-form in the coordinates $(x, a)$ on $\mathcal{P}$
in order to view explicitly how the curvature depends upon some basic
initial geometric data, namely upon the coefficients $\eta_{i, i'}
(x)$ of a certain initial frame $(H_1, \dots, H_n)$ on $M$ ({\em see}
below) and also upon the coefficients $\upsilon_{ j, j'} (a)$ of some
explicit frame $(Y_1, \dots, Y_{ n-r})$ of left-invariant vector
fields on the (local) Lie group $H$.

\smallskip
Before entering further the construction of effective Cartan
connections, we proceed to a preliminary investigation of what basic
relations come at the level of only linear algebra.

\subsection{Back to the curvature function}
By definition, at any fixed point $p \in \mathcal{ P}$, the curvature
function $\kappa (p)$ happens to be a linear map from $\Lambda^2 (
\mathfrak{ g} / \mathfrak{ h})$ into $\mathfrak{ g}$, whence, thanks
to the canonical, classical identification ${\rm Lin} ( E, F) \simeq
E^* \otimes F$ which reads here as:
\[
\aligned
{\rm Lin}
\big(\Lambda^2(\mathfrak{g}/\mathfrak{h}),\,
\mathfrak{g}\big)
&
\simeq
\big(\Lambda^2\mathfrak{g}/\mathfrak{h}\big)^*
\otimes
\mathfrak{g}
\\
&
\simeq
\Lambda^2(\mathfrak{g}^*/\mathfrak{h}^*)
\otimes
\mathfrak{g},
\endaligned
\]
we may express as follows the curvature function in terms of the basis
elements ${\sf x}_1, \dots, {\sf x}_n, {\sf x}_{ n+1}, \dots, {\sf
x}_r$ for $\mathfrak{ g}$ and in terms of the basis (representatives)
elements ${\sf x}_1^*, \dots, {\sf x}_n^*$ of $\mathfrak{ g}^* /
\mathfrak{ h}^*$:
\[
\boxed{
\kappa(p)
=
\sum_{1\leqslant i_1<i_2\leqslant n}\,
\sum_{k=1}^r\,
\kappa_{i_1,i_2}^k(p)\,\,
{\sf x}_{i_1}^*\wedge{\sf x}_{i_2}^*
\otimes
{\sf x}_k
}\,,
\]
where the $p \mapsto \kappa_{i_1, i_2}^k ( p)$ are a certain
collection of $\K$-valued function on $\mathcal{ P}$.

Now, we remind that Proposition~\ref{equivariancy-curvature-function} on
p.~\pageref{equivariancy-curvature-function} showed us that:
\[
(Y^\dag\kappa)(p)({\sf x}',{\sf x}'')
=
-\big[{\sf y},\,\kappa(p)({\sf x}',{\sf x}'')\big]_{\mathfrak{g}}
+
\kappa(p)\big([{\sf y},{\sf x}']_{\mathfrak{g}},\,
{\sf x}''\big)
+
\kappa(p)\big({\sf x}',\,
[{\sf y},{\sf x}'']_{\mathfrak{g}}\big),
\]
for any fundamental field $Y^\dag = \frac{ d}{ dt} \big\vert_0 R_{
\exp (t {\sf y})}$ on $\mathcal{ P}$ associated to an arbitrary ${\sf
y} \in \mathfrak{ h}$. In terms of Lie algebra bases, it suffices to
inspect what this formula gives when ${\sf x}' = {\sf x}_{ i_1}$, when
${\sf x}'' = {\sf x}_{ i_2}$ for some $1 \leqslant i_1 < i_2 \leqslant
n$, and when ${\sf y} = {\sf y}_j = {\sf x}_{ n+j}$ for some $j$ with
$1 \leqslant j \leqslant r - n$. Plugging in then these values, 
we receive the equation:
\[
\aligned
\sum_{k=1}^r\,\big(Y^\dag\kappa_{i_1,i_2}^k\big)(p)\,{\sf x}_k
&
=
-
\sum_{k=1}^r\,\kappa_{i_1,i_2}^k(p)\,
[{\sf y}_j,\,{\sf x}_k]_{\mathfrak{g}}
+
\\
&
\ \ \ \ \
+
\kappa(p)\big({\sf x}_{i_1},[{\sf y}_j,{\sf x}_{i_2}]_{\mathfrak{g}}\big)
-
\kappa(p)\big({\sf x}_{i_2},[{\sf y}_j,{\sf x}_{i_1}]_{\mathfrak{g}}\big).
\endaligned
\]
On the first hand, using structure constants,
we may easily compute the first (among three)
terms appearing in the right-hand side: 
\[
-\sum_{k=1}^r\,\kappa_{i_1,i_2}^k(p)\,
[{\sf x}_{n+j},{\sf x}_k]_{\mathfrak{g}}
=
-
\sum_{k=1}^r\,
\sum_{s=1}^r\,
\kappa_{i_1,i_2}^k(p)\,
c_{n+j,k}^s\,{\sf x}_s.
\]
On the other hand, the two brackets appearing in the second line
express themselves by means of structure constants just as follows
modulo $\mathfrak{ h}$:
\[
\aligned
{}
[{\sf y}_j,{\sf x}_{i_1}]_{\mathfrak{g}}
=
[{\sf x}_{n+j},{\sf x}_{i_1}]_{\mathfrak{g}}
=
\sum_{k=1}^r\,c_{n+j,i_1}^k\,{\sf x}_k
&
\equiv
\sum_{i_3=1}^n\,c_{n+j,i_1}^{i_3}\,{\sf x}_{i_3}\,\,\,
{\rm mod}\,\mathfrak{h},
\\
{}
[{\sf y}_j,{\sf x}_{i_2}]_{\mathfrak{g}}
=
[{\sf x}_{n+j},{\sf x}_{i_2}]_{\mathfrak{g}}
=
\sum_{k=1}^r\,c_{n+j,i_2}^k\,{\sf x}_k
&
\equiv
\sum_{i_3=1}^n\,c_{n+j,i_2}^{i_3}\,{\sf x}_{i_3}\,\,\,
{\rm mod}\,\mathfrak{h}.
\endaligned
\]
We may hence derive that the entire second line above
equals:
\[
\sum_{i_3=1}^n\,c_{n+j,i_2}^{i_3}\,\kappa(p)({\sf x}_{i_1},{\sf x}_{i_3})
-
\sum_{i_3=1}^n\,c_{n+j,i_1}^{i_3}\,\kappa(p)({\sf x}_{i_2},{\sf x}_{i_3}).
\]
At this point however, the computation is not yet finished. Indeed,
reminding the basic formula~\thetag{ \ref{i-1-i-2-j-1-j-2}} on
p.~\pageref{i-1-i-2-j-1-j-2}, we may pursue for instance the expansion
of $\kappa(p) ( {\sf x}_{ i_1}, {\sf x}_{ i_3})$ in the first sum:
\[
\aligned
\kappa(p)({\sf x}_{i_1},{\sf x}_{i_3})
&
=
\sum_{1\leqslant i_1'<i_2'\leqslant n}\,
\sum_{k=1}^r\,
\kappa_{i_1',i_2'}^k(p)\,
{\sf x}_{i_1'}^*\wedge{\sf x}_{i_2'}^*({\sf x}_{i_1},{\sf x}_{i_3})
\otimes
{\sf x}_k
\\
&
=
\sum_{1\leqslant i_1'<i_2'\leqslant n}\,
\sum_{k=1}^r\,
\kappa_{i_1',i_2'}^k(p)\,
\big[
\delta_{i_1}^{i_1'}\delta_{i_3}^{i_2'}
-
\delta_{i_3}^{i_1'}\delta_{i_1}^{i_2'}
\big]\,
{\sf x}_k.
\endaligned
\]
Thanks to this simple formula, the first $\sum_{i_3}$ continues as
(complete explanations follow just after):
\[
\footnotesize
\aligned
\sum_{i_3=1}^n\,c_{n+j,i_2}^{i_3}\,\kappa(p)({\sf x}_{i_1},{\sf x}_{i_3})
&
=
\sum_{k=1}^r\,
\bigg(
\sum_{i_3=1}^n\,\sum_{1\leqslant i_1'<i_2'\leqslant n}\,
c_{n+j,i_2}^{i_3}\,\kappa_{i_1',i_2'}^k(p)
\big[
\delta_{i_1}^{i_1'}\delta_{i_3}^{i_2'}
-
\delta_{i_3}^{i_1'}\delta_{i_1}^{i_2'}
\big]
\bigg)\,{\sf x}_k
\\
&
=
\sum_{k=1}^r\,
\bigg(
\sum_{1\leqslant i_1'<i_2'\leqslant n}\,
c_{n+j,i_2}^{i_2'}\,\kappa_{i_1',i_2'}^k(p)\,\delta_{i_1}^{i_1'}
-
\sum_{1\leqslant i_1'<i_2'\leqslant n}\,
c_{n+j,i_2}^{i_1'}\,\kappa_{i_1',i_2'}^k(p)\,\delta_{i_1}^{i_2'}
\bigg)\,{\sf x}_k
\\
&
=
\sum_{k=1}^r\,
\bigg(
\sum_{i_2'=i_1+1}^n\,c_{n+j,i_2}^{i_2'}\,\kappa_{i_1,i_2'}^k(p)
-
\sum_{i_1'=1}^{i_1}\,c_{n+j,i_2}^{i_1'}\,\kappa_{i_1',i_1}^k(p)
\bigg)\,{\sf x}_k
\\
&
=
\sum_{k=1}^r\,
\bigg(
-
\sum_{i'=1}^{i_1}\,c_{n+j,i_2}^{i'}\,\kappa_{i',i_1}^k(p)
+
\sum_{i'=i_1+1}^n\,c_{n+j,i_2}^{i'}\,\kappa_{i_1,i'}^k(p)
\bigg)\,{\sf x}_k.
\endaligned
\]
From the first to the second line, the $\sum_{ i_3}$ is automatically
killed because of the presence of the two $\delta_{ i_3}^{ i_2'}$ and
$- \delta_{ i_3}^{ i_1'}$. Then the two $\sum_{ i_1' < i_2'}$ both
reduce to a single sum, because of the presence of the two remaining
$\delta_{ i_1}^{ i_1'}$ and $-\delta_{ i_1}^{ i_2'}$. Finally, the
order of appearance of the two sums $\sum_{ i_2'}$ and $- \sum_{
i_1'}$ is interchanged and both summation indices are 
given the same name $i'$.
Summing up all terms, we have shown the following

\begin{Proposition}
For every ${\sf y} \in \mathfrak{ h}$ and every
$1 \leqslant i_1 < i_2 \leqslant n$, one has:
\[
\footnotesize
\aligned
\sum_{k=1}^r\,\big(Y^\dag\kappa_{i_1,i_2}^k\big)(p)\,{\sf x}_k
&
=
\sum_{k=1}^r\,
\bigg(
\sum_{s=1}^r\,c_{n+j,s}^k\,\kappa_{i_1,i_2}^s(p)
-
\\
&
\ \ \ \ \ \ \ \ \ \ \ \ \ \ \
-
\sum_{i'=1}^{i_1}\,c_{n+j,i_2}^{i'}\,\kappa_{i',i_1}^k(p)
+
\sum_{i'=i_1+1}^n\,c_{n+j,i_2}^{i'}\,\kappa_{i_1,i'}^k(p)
\\
&
\ \ \ \ \ \ \ \ \ \ \ \ \ \ \
+
\sum_{i'=1}^{i_2}\,c_{n+j,i_1}^{i'}\,\kappa_{i',i_2}^k(p)
-
\sum_{i'=i_2+1}^n\,c_{n+j,i_1}^{i'}\,\kappa_{i_2,i'}^k(p)
\bigg)\,{\sf x}_k,
\endaligned
\]
where $Y^\dag = \widehat{ Y} = \omega^{ -1} ( {\sf y})$ is
the constant field associated to ${\sf y}$.
\qed
\end{Proposition}

\subsection{Bianchi identity}
Looking at $\Omega = d\omega + \frac{ 1}{ 2}\, \omega \wedge_{
\mathfrak{ g}} \omega$, if
one applies once more the exterior differential operator $d$,
the term $d d\omega$ drops.
This observation leads to Bianchi-type identities which were first
understood by Christoffel and Lipschitz 
in the context of Riemannian geometry.

\begin{Theorem}
The curvature $2$-form:
\[
\Omega
=
\sum_{k=1}^r\,
\Omega^k\,{\sf x}_k
=
\sum_{k=1}^r\,
\Big\{
d\omega^k
+
\sum_{1\leqslant l_1<l_2\leqslant r}\,
c_{l_1,l_2}^k\,
\omega^{l_1}\wedge\omega^{l_2}
\Big\}\,{\sf x}_k
\]
satisfies the identity:
\[
d\Omega
=
\Omega\wedge_{\mathfrak{g}}\omega,
\]
that is to say in terms of the basis $({\sf x}_k)_{
1\leqslant k\leqslant r}$ of $\mathfrak{ g}$:
\[
d\Omega^k
=
\sum_{i_1,i_2=1}^r\,
c_{i_1,i_2}^k\,
\Omega^{i_1}\wedge\omega^{i_2}
\ \ \ \ \ \ \ \ \ \ \ \ \
{\scriptstyle{(k\,=\,1\,\cdots\,r)}}.
\]
\end{Theorem}

\proof
Applying, as we said, $d$ to the definition $\Omega = d\omega + \frac{
1}{ 2}\, \omega \wedge_{ \mathfrak{ g}} \omega$ of the curvature form,
we get:
\[
\aligned
d\Omega
&
=
\zero{dd\omega}
+
{\textstyle{\frac{1}{2}}}\,
d\omega\wedge_{\mathfrak{g}}\omega
-
{\textstyle{\frac{1}{2}}}\,
\omega\wedge_{\mathfrak{g}}d\omega
\\
&
=
{\textstyle{\frac{1}{2}}}\,
\big(
\Omega
-
{\textstyle{\frac{1}{2}}}\,
\omega\wedge_{\mathfrak{g}}\omega
\big)
\wedge_{\mathfrak{g}}
\omega
-
{\textstyle{\frac{1}{2}}}\,
\omega
\wedge_{\mathfrak{g}}\big(
\Omega
-
{\textstyle{\frac{1}{2}}}\,
\omega\wedge_{\mathfrak{g}}\omega
\big).
\endaligned
\]
Here, one has to mind the fact that, although the wedge product
$\wedge$ between scalar one-forms is associative, the wedge-bracket
product $\wedge_{ \mathfrak{ g}}$ is {\em not} (in general)
associative, just because one does not have in general $\big[ [ {\sf
x}_{ k_1}, {\sf x}_{ k_2} ]_{ \mathfrak{ g}}, \, {\sf x}_{ k_3}
\big]_{ \mathfrak{ g}} = \big[ {\sf x}_{ k_1}, \, [ {\sf x}_{ k_2},
{\sf x}_{ k_3}]_{ \mathfrak{ g}} \big]_{ \mathfrak{ g}}$. Fortunately,
the Jacobi identity guarantees that 
the following is true.

\begin{Lemma}
For any $\mathfrak{ g}$-valued $1$-form $\varpi = 
\sum_{ k=1}^r\, \varpi^k\, {\sf x}_k$, one has:
\[
0
=
\big(\varpi\wedge_{\mathfrak{g}}\varpi\big)\wedge_{\mathfrak{g}}\varpi
\ \ \ \ \ \ \ \ \ \ \
\text{\rm and}
\ \ \ \ \ \ \ \ \ \ \
0
=
\varpi\wedge_{\mathfrak{g}}\big(\varpi\wedge_{\mathfrak{g}}\varpi\big).
\]
\end{Lemma}

\proof
We check the first identity, the second one being similar. 
Expanding, we may write:
\[
\aligned
\big(\varpi\wedge_{\mathfrak{g}}\varpi\big)\wedge_{\mathfrak{g}}\varpi
&
=
\bigg(
\sum_{k_1=1}^r\,\varpi^{k_1}\,{\sf x}_{k_1}
\wedge_{\mathfrak{g}}
\sum_{k_2=1}^r\,\varpi^{k_2}\,{\sf x}_{k_2}
\bigg)
\wedge_{\mathfrak{g}}
\sum_{k_3=1}^r\,\varpi^{k_3}\,{\sf x}_{k_3}
\\
&
=
\bigg(
\sum_{k_1,k_2=1}^r\,
\varpi^{k_1}\wedge\varpi^{k_2}\,
[{\sf x}_{k_1},{\sf x}_{k_2}]_{\mathfrak{g}}
\bigg)
\wedge_{\mathfrak{g}}
\sum_{k_3=1}^r\,\varpi^{k_3}\,{\sf x}_{k_3}
\\
&
=
\sum_{k_1,k_2,k_3=1}^r\,
\varpi^{k_1}\wedge\varpi^{k_2}\wedge\varpi^{k_3}\,
\big[[{\sf x}_{k_1},{\sf x}_{k_2}]_{\mathfrak{g}},\,
{\sf x}_{k_3}\big]_{\mathfrak{g}}.
\endaligned
\]
Clearly, since the three-terms wedge product $\varpi^{ k_1} \wedge
\varpi^{ k_2} \wedge \varpi^{ k_3}$ obviously vanishes as soon as at
least two of the $k_i$'s coincide, the triple sum $\sum_{k_1, k_2,
k_3}$ decomposes plainly in six sub-sums:
\[
\footnotesize
\aligned
\sum_{k_1,k_2,k_3
\atop
\text{\rm pairwise}\neq}
&
=
\sum_{k_1<k_2<k_3}
+
\sum_{k_3<k_1<k_2}
+
\sum_{k_2<k_3<k_1}
+
\\
&
\ \ \ \ \
+
\sum_{k_2<k_1<k_3}
+
\sum_{k_3<k_2<k_1}
+
\sum_{k_1<k_3<k_2}.
\endaligned
\]
Then each triple sum vanishes identically thanks to the Jacobi
identity~\thetag{ \ref{skew-Jacobi}}.
\endproof

Coming back to the $d\Omega$ we left above, the lemma shows that the
two terms involving three times $\omega$ vanish and it remains:
\[
\aligned
d\Omega
&
=
{\textstyle{\frac{1}{2}}}\,
\Omega\wedge_{\mathfrak{g}}\omega
-
{\textstyle{\frac{1}{2}}}\,
\omega\wedge_{\mathfrak{g}}\Omega
\\
&
=
\Omega\wedge_{\mathfrak{g}}\omega,
\endaligned
\]
as claimed by the theorem. Expanding then this identity 
in bases and reorganizing, we get:
\[
\footnotesize
\aligned
\sum_{k=1}^r\,d\Omega^k\,{\sf x}_k
&
=
\Big(\sum_{i_1=1}^r\,\Omega^{i_1}\,{\sf x}_{i_1}\Big)
\wedge_{\mathfrak{g}}
\Big(\sum_{i_2=1}^r\,\Omega^{i_2}\,{\sf x}_{i_2}\Big)
\\
&
=
\sum_{i_1,i_2=1}^r\,
\Omega^{i_1}\wedge\omega^{i_2}\,
[{\sf x}_{i_1},{\sf x}_{i_2}]_{\mathfrak{g}}
\\
&
=
\sum_{k=1}^r\,
\Big(
\sum_{i_1,i_2=1}^r\,c_{i_1,i_2}^k\,
\Omega^{i_1}\wedge\omega^{i_2}
\Big)\,{\sf x}_k
\endaligned
\]
that is to say, for every $k = 1, \dots, r$:
\[
d\Omega^k
=
\sum_{i_1,i_2=1}^r\,
c_{i_1,i_2}^k\,\Omega^{i_1}\wedge\omega^{i_2},
\]
which concludes the proof.
\endproof

The curvature function $p \mapsto \kappa (p)$ is a map $\mathcal{ P}
\to \mathcal{ C}^2 (\mathfrak{ g} / \mathfrak{ h}, \mathfrak{ g})$. A
second version of Bianchi identities is as follows. Remind the
operator $\partial \colon \mathcal{ C}^2 (\mathfrak{ g} / \mathfrak{
h}, \mathfrak{ g}) \longrightarrow \mathcal{ C}^3 ( \mathfrak{ g} /
\mathfrak{ h}, \mathfrak{ g})$.

\begin{Proposition}
\label{Bianchi-Tanaka}
For any three ${\sf x}', {\sf x}'', {\sf x}''' \in \mathfrak{ g}$, 
one has at every point $p \in \mathcal{ P}$:
\[
0
=
(\partial\kappa)(p)({\sf x}',{\sf x}'',{\sf x}''')
+
\sum_{\rm cycl}\,
\kappa(p)\big(
\kappa(p)({\sf x}',{\sf x}''),\,{\sf x}'''\big)
+
\sum_{\rm cycl}\,
\big(\widehat{X}'\kappa\big)(p)({\sf x}'',{\sf x}'''),
\]
that is to say in greater length:
\[
\aligned
0
&
=
\big[{\sf x}',\kappa(p)({\sf x}'',{\sf x}''')\big]_{\mathfrak{g}}
-
\big[{\sf x}'',\kappa(p)({\sf x}',{\sf x}''')\big]_{\mathfrak{g}}
+
\big[{\sf x}''',\kappa(p)({\sf x}',{\sf x}'')\big]_{\mathfrak{g}}
-
\\
&
\ \ \ \ \
-
\kappa(p)\big([{\sf x}',{\sf x}'']_{\mathfrak{g}},\,{\sf x}'''\big)
+
\kappa(p)\big([{\sf x}',{\sf x}''']_{\mathfrak{g}},\,{\sf x}''\big)
-
\kappa(p)\big([{\sf x}'',{\sf x}''']_{\mathfrak{g}},\,{\sf x}'\big)
+
\\
&
\ \ \ \ \
+
\kappa(p)\big(\kappa(p)({\sf x}',{\sf x}''),\,{\sf x}'''\big)
-
\kappa(p)\big(\kappa(p)({\sf x}',{\sf x}'''),\,{\sf x}''\big)
+
\kappa(p)\big(\kappa(p)({\sf x}'',{\sf x}'''),\,{\sf x}''\big)
+
\\
&
\ \ \ \ \
+
\big(\widehat{X}'\kappa\big)(p)({\sf x}'',{\sf x}''')
-
\big(\widehat{X}''\kappa\big)(p)({\sf x}',{\sf x}''')
+
\big(\widehat{X}'''\kappa\big)(p)({\sf x}',{\sf x}'').
\endaligned
\]
\end{Proposition}

\proof
By the Definition~\ref{definition-curvature-function} of the curvature
function and the Definition~\ref{definition-curvature-form} of the
curvature form, for any two ${\sf z}, {\sf x}''' \in \mathfrak{ g}$,
one has:
\[
\aligned
\kappa(p)({\sf z},{\sf x}''')
&
=
\Omega_p\big(\widehat{Z}_p,\widehat{X}_p'''\big)
\\
&
=
d\omega_p\big(\widehat{Z}_p,\widehat{X}_p'''\big)
+
\big[
\omega_p(\widehat{Z}_p),\omega_p(\widehat{X}_p''')
\big]_{\mathfrak{g}}.
\endaligned
\]
Replacing then $\widehat{ Z}_p$ by $\big[ \widehat{ X}_p', \widehat{
X}_p'' \big]$, hence also ${\sf z}$ by $\omega_p \big( \big[ \widehat{
X}_p', \widehat{ X}_p'' \big] \big)$ gives:
\[
\kappa(p)
\big(
\omega_p\big(\big[\widehat{X}_p',\widehat{X}_p''\big]\big),\,\,{\sf x}'''
\big)
=
d\omega_p\big(
\big[\widehat{X}_p',\widehat{X}_p''\big],
\widehat{X}_p'''\big)
+
\big[
\omega_p\big(
\big[\widehat{X}_p',\widehat{X}_p''\big]\big),\,
{\sf x}'''
\big]_{\mathfrak{g}}.
\]
Now, if we apply the Cartan formula to 
expand $d\omega_p$ in the right-hand side, we receive:
\[
\footnotesize
\aligned
&
\kappa(p)
\big(
\omega_p\big(\big[\widehat{X}_p',\widehat{X}_p''\big]\big),\,\,{\sf x}'''
\big)
=
\\
&
\ \ \ \ \ 
=
\zero{\big[\widehat{X}',\widehat{X}''\big]
\big(\omega_p(\widehat{X}''')\big)}
-
\widehat{X}'''\big(\omega_p([\widehat{X}_p',\widehat{X}_p''])\big)
-
\omega_p
\big(\big[[\widehat{X}_p',\widehat{X}_p''],\,\widehat{X}_p'''
\big]\big)
+
\big[
\omega_p\big(
\big[\widehat{X}_p',\widehat{X}_p''\big]\big),\,
{\sf x}'''
\big]_{\mathfrak{g}},
\endaligned
\]
and the first term in the right-hand side vanishes, because $\omega_p
(\widehat{ X}''') = {\sf x}'''$ is constant. Next, both in the
left-hand side and in the right-hand side, we replace thrice from
Lemma~\ref{first-curvature-function}
\[
\omega_p\big(\big[\widehat{X}_p',\widehat{X}_p''\big]\big)
=
[{\sf x}',{\sf x}'']_{\mathfrak{g}}
-
\kappa(p)({\sf x}',{\sf x}''),
\]
and we receive:
\[
\footnotesize
\aligned
&
\kappa(p)\big([{\sf x}',{\sf x}'']_{\mathfrak{g}},\,{\sf x}'''\big)
-
\kappa(p)\big(\kappa(p)({\sf x}',{\sf x}''),\,{\sf x}'''\big)
=
\\
&
\ \ \ \ \ 
=
-
\zero{\widehat{X}'''\big([{\sf x}',{\sf x}'']_{\mathfrak{g}}\big)}
+
\big(\widehat{X}'''\kappa\big)(p)({\sf x}',{\sf x}'')
-
\omega_p
\big(\big[[\widehat{X}_p',\widehat{X}_p''],\,\widehat{X}_p'''
\big]\big)
\\
&
\ \ \ \ \ \ \ \ \ \ 
+
\big[
[{\sf x}',{\sf x}'']_{\mathfrak{g}},\,
{\sf x}'''
\big]_{\mathfrak{g}}
-
\big[
\kappa(p)({\sf x}',{\sf x}''),\,
{\sf x}'''
\big]_{\mathfrak{g}},
\endaligned
\]
Putting all the six remaining terms in the right-hand side and
reorganizing their order of appearance, we get an identity in which:
\[
\aligned
0
&
=
\big[
{\sf x}''',\,
\kappa(p)({\sf x}',{\sf x}'')
\big]_{\mathfrak{g}}
-\kappa(p)\big([{\sf x}',{\sf x}'']_{\mathfrak{g}},\,{\sf x}'''\big)
+
\kappa(p)\big(\kappa(p)({\sf x}',{\sf x}''),\,{\sf x}'''\big)
+
\\
&
\ \ \ \ \ 
+
\big(\widehat{X}'''\kappa\big)(p)({\sf x}',{\sf x}'')
-
\omega_p
\big(\big[[\widehat{X}_p',\widehat{X}_p''],\,\widehat{X}_p'''
\big]\big)
+
\big[
[{\sf x}',{\sf x}'']_{\mathfrak{g}},\,
{\sf x}'''
\big]_{\mathfrak{g}},
\endaligned
\]
when one sums over all cyclic permutations of $\{ {\sf x}', {\sf x}'',
{\sf x}'''\}$, the last two terms disappear thanks to the Jacobi
identity for vector fields (fifth term) and thanks to the Jacobi
identity within $\mathfrak{ g}$ (sixth term) so that we obtain the
formula stated.
\endproof

\begin{Corollary}
For any three ${\sf y}', {\sf y}'', {\sf y}''' \in \mathfrak{ h}$, 
one has:
\[
(\partial\kappa)(p)
\big({\sf x}'+{\sf y}',{\sf x}''+{\sf y}'',{\sf x}'''+{\sf y}'''\big)
=
(\partial\kappa)(p)
\big({\sf x}',{\sf x}'',{\sf x}'''\big)
\]
so that $\partial \kappa$ is well defined in the
space of $3$-cochains $\mathcal{ C}^3
(\mathfrak{ g} / \mathfrak{ h}, \mathfrak{ g})$.
\end{Corollary}

\proof
By symmetry, it suffices to check this when ${\sf y}'' = {\sf y}''' =
0$ with ${\sf y}'$ simply denoted ${\sf y}$. But then, applying twice
the {\em expanded} second formula of the proposition and subtracting,
noticing that $\widehat{ X' + Y} = \widehat{ X}' + Y^\dag$, we obtain:
\[
0
\overset{?}{=}
\big[
{\sf y},\,\kappa(p)({\sf x}'',{\sf x}''')
\big]_{\mathfrak{g}}
-
\kappa(p)\big([{\sf y},{\sf x}'']_{\mathfrak{g}},\,
{\sf x}'''\big)
+
\kappa(p)\big([{\sf y},{\sf x}''']_{\mathfrak{g}},{\sf x}''\big)
+
\big(Y^\dag\kappa\big)(p)({\sf x}'',{\sf x}'''),
\]
which is coherent with the second
formula of Proposition~\ref{equivariancy-curvature-function}.
\endproof

\begin{Proposition}
{\rm (\cite{Cap, EMS})}
\label{Bianchi-Tanaka-graded}
When the Lie algebra 
$\mathfrak{ g} = \mathfrak{ g}_{ -a} \oplus \cdots \oplus
\mathfrak{ g}_b$ is graded as on 
p.~\pageref{form-of-a-graded-Lie-algebra}
with:
\[
\aligned
\mathfrak{g}_-
&
:=
\mathfrak{g}_{-a}
\oplus\cdots\oplus
\mathfrak{g}_{-1},
\\
\mathfrak{h}
&
:=
\mathfrak{g}_0
\oplus\cdots\oplus
\mathfrak{g}_b,
\endaligned
\]
if one decomposes
the curvature function $\kappa = \sum_{ h \in \Z}\, 
\kappa_{ [h]}$ in homogeneous components $\kappa_{ [h]}$, then the
differential of every such a $\kappa_{ [h]}$
is uniquely determined by the lower homogeneous components
$\kappa_{ [h']}$, with $h' \leqslant h -1$, through the specific
formula{\em :}
\[
\aligned
\partial_{[h]}\big(\kappa_{[h]}\big)
({\sf x}',{\sf x}'',{\sf x}''')
&
=
-\,
\sum_{{\rm cycl}}\,
\sum_{h'=1}^{h-1}\,
\Big(
\kappa_{[h-h']}
\big(
{\sf proj}_{\mathfrak{g}_-}
\big(
\kappa_{[h']}({\sf x}',{\sf x}'')
\big),\,
{\sf x}'''
\big)
\Big)
-
\\
&
\ \ \ \ \
-
\sum_{{\rm cycl}}\,
\big(\widehat{X}'\kappa_{[h+\vert{\sf x}'\vert]}\big)
({\sf x}'',{\sf x}'''),
\endaligned
\]
in which ${\sf proj}_{ \mathfrak{ g}_-} ({\sf z})$ denotes the
$\mathfrak{ g}_-$-component of an 
element ${\sf z} \in \mathfrak{ g}$, 
while $\vert {\sf x}' \vert$ denotes the
grade of ${\sf x}'$. 
\end{Proposition}

\proof
Back to the preceding proposition, 
it suffices to pick, in the expression:
\[
-\,
\sum_{\rm cycl}\,
\kappa(p)\big(
\kappa(p)({\sf x}',{\sf x}'')\,\,{\rm mod}\,\mathfrak{h},\,\,
{\sf x}'''\big)
-
\sum_{\rm cycl}\,
\big(\widehat{X}'\kappa\big)(p)({\sf x}'',{\sf x}''')
\]
only the components that are of homogeneous degree $h$. 
\endproof

\subsection{Normality}
We can now present a consequence of the above graded
Bianchi-Tanaka identities that
will be important for the construction of a Cartan connection on
strongly pseudoconvex hypersurfaces $M^3 \subset \C^2$, to be achieved
in an effective way in Section~\ref{Cartan-construction} below. When
$\mathfrak{ g}$ is semi-simple ({\em cf.}~\cite{Cap, Slovak}) with a
grading $\mathfrak{ g} = \mathfrak{ g}_{ - \mu} \oplus \cdots \oplus
\mathfrak{ g}_\mu$ of depth $\mu \geqslant 1$
($a = b$ in this case), it is known ({\em
cf.}~\cite{ Cap}, p.~492) that the curvature function:
\[
\kappa
=
\kappa_{[1]}
+\cdots+
\kappa_{[3\mu]}
\]
has components whose homogeneities $h \geqslant 1$ are all positive.
Then a Cartan connection $\omega \colon T\mathcal{ P} \to \mathfrak{
g}$ is said to be {\sl normal} (in the sense of Tanaka, {\em
cf.}~\cite{ Sato, Cap}) if the codifferential operator annihilates all
the (positive) homogeneous components of the curvature function,
namely if:
\[
0
=
\partial^*
\big(\kappa_{[1]}\big)
=
\cdots
=
\partial^*
\big(\kappa_{[3\mu]}\big).
\]

\begin{Proposition}
\label{cohomology-result}
{\rm (\cite{Cap, EMS})}
If $\mathfrak{ g}$ is semi-simple and if the Cartan connection $\omega
\colon T\mathcal{ P} \to \mathfrak{ g}$ is normal, then the
homogeneity of the first non-zero homogeneous component of the
curvature function $\kappa$ is greater than the homogeneity of the
first non-zero homogeneous component of the second cohomology $H^2
(\frak g_-, \frak g)$.
\end{Proposition}

\proof 
If, as is assumed,
$H_{ [h']}^2 (\frak g_-,\frak g) = 0$ for all $h' \leqslant h$ with
some $h \geqslant 1$, and if, by induction:
\[
0
=
\kappa_{[1]}
=
\cdots
=
\kappa_{[h-1]},
\]
then the graded Bianchi-Tanaka identities of 
Proposition~\ref{Bianchi-Tanaka-graded}
show immediately that:
\[
\partial_{[h]}\big(\kappa_{[h]}\big)
=
0,
\]
which just means that: 
\[
\kappa_{[h]}
\in
{\rm ker}\big(\partial_{[h]}\big).
\]
Furthermore, by normality of the connection, we also have:
\[
\kappa_{[h]}
\in
{\rm ker}\big(\partial_{[h]}^*\big).
\]
In addition, according to Proposition 2.6 of \cite{Cap}, the
semi-simplicity of $\mathfrak{ g}$ insures that for every homogeneous
degree $h \in \Z$, the $h$-homogeneous $2$-cochain space (to which
$\kappa_{ [h]}$ belongs) splits up as:
\[
\aligned
\mathcal{C}_{[h]}^2
&
=
{\rm im}
\big(\partial_{[h]}\big)
\oplus
{\rm ker}
\big(
\partial_{[h]}^*
\big)
\\
&
=
\mathcal{B}_{[h]}^2
\oplus
{\rm ker}\big(\partial_{[h]}^*\big)
\endaligned
\] 
But since $H_{ [h]}^2 = 0$ by assumption, which just means that we can
replace $\mathcal{ B}_{[ h]}^2 = \mathcal{ Z}_{ [h]}^2 = {\rm ker}\,
\partial_{[ h]}$ here, we deduce that:
\[
\kappa_{[h]}
\in
{\rm ker}\big(\partial_{[h]}\big)
\oplus
{\rm ker}\big(\partial_{[h]}^*\big)
\ \ \ \ \
\text{\rm and}
\ \ \ \ \
\kappa_{[h]}
\in
{\rm ker}\big(\partial_{[h]}\big)
\cap
{\rm ker}\big(\partial_{[h]}^*\big),
\]
which trivially implies that $\kappa_{ [ h]} = 0$,
as desired.
\endproof


\subsection{Vector fields and one-forms}
A remark about notations is in order. As did Sophus Lie, we shall most
often prefer to conduct computations in terms of vector fields
(first-order derivations) instead of using \'Elie Cartan's exterior
differential calculus. But as both points of view are complementary,
we shall nonetheless {\em need at the same time} to consider
differential forms, and the best way of remembering the duality
between a frame and a coframe is to put a ``$*$'' as an upper index,
as we already did for the ${\sf x}_k$ and their duals ${\sf x}_k^*$.
This is why we shall not use \'Elie Cartan's classical letters
$\omega^i$ or $\varpi^j$ to denote differential forms, except for $\omega
\colon T\mathcal{ P} \to \mathfrak{ g}$. In fact, in the formalism we
set up presently, Greek letters will be mostly reserved to denote {\em
functions}, namely coefficients of vector fields or of $1$-forms.

\subsection{Left-invariant Lie frame and Maurer coframe}
Now, as we agreed that our constructions shall be local, we must
consider that the local lie group $H$ not only comes with its abstract
Lie algebra $\mathfrak{ h}$, but also with a collection of
left-invariant vector fields near its identity element. So, let us
assume that $H$, of dimension $n$, is coordinatized by means of $(a_1,
\dots, a_{ r - n}) \in \K^{ r - n}$, locally near the origin (its
identity element) and that the Lie algebra of left-invariant vector
fields on $H$ is given by $r - n$ (left-invariant) vector fields:
\[
Y_j
=
\sum_{j'=1}^{r-n}\,
\upsilon_{j,j'}(a_1,\dots,a_{r-n})\,
\frac{\partial}{\partial a_{j'}}
\ \ \ \ \ \ \ \ \ \ \ \ \
{\scriptstyle{(j\,=\,1\,\cdots\,r\,-\,n)}}
\]
having coefficients $\upsilon_{ j, j'} ( a)$ that
are smooth or $\K$-analytic in a neighborhood of the origin
in $\K^{ r - n}$. 

Dually, the Maurer-Cartan coframe on $H$ consists of precisely
the same number $r - n$ of (left-invariant) $1$-forms:
\[
Y_j^*
=
\sum_{j'=1}^{r-n}\,
\upsilon_*^{j,j'}(a)\,da_{j'},
\]
satisfying by definition:
\[
Y_{j_2}^*(Y_{j_1})
=
\delta_{j_1}^{j_2}.
\]
Hence the $(r - n) \times (r - n)$ matrix of functions $\upsilon_*^{j,
j'} (a)$ is the transpose-inverse of the matrix of functions
$\upsilon_{ j, j'} (a)$. 

\subsection{Initial frame on the base manifold}
On the other hand, on the manifold $M$ equipped with local coordinates
$x = (x_1, \dots, x_n)$, we suppose that a certain frame is given
which consists of $n = \dim_\K M$ vector fields:
\[
H_i
=
\sum_{i'=1}^n\,\eta_{i,i'}(x)\,
\frac{\partial}{\partial x_{i'}}
\ \ \ \ \ \ \ \ \ \ \ \ \
{\scriptstyle{(i\,=\,1\,\cdots\,n)}}
\]
having coefficients $\eta_{ i, i'} ( x)$ that are 
smooth or $\K$-analytic in a
neighborhood of the origin in $\K^n$.

Then the dual coframe consists of precisely the same number $n$ of
$1$-forms:
\[
H_i^*
:=
\sum_{i'=1}^n\,\eta_*^{i,i'}(x)\,
dx_{i'}
\ \ \ \ \ \ \ \ \ \ \ \ \
{\scriptstyle{(i\,=\,1\,\cdots\,n)}}
\]
satisfying by definition:
\[
H_{i_2}^*(H_{i_1})
=
\delta_{i_1}^{i_2}.
\]
Hence the $n \times n$ matrix of functions $\eta_*^{i, i'} (x)$ is the
transpose-inverse of the matrix of functions $\eta_{ i, i'} (x)$.

\subsection{Constant frame and coframe on $\mathcal{ P}$}
After all, the coordinates on the total space $\mathcal{ P}$ are just:
\[
(x,a)
=
(x_1,\dots,x_n,a_1,\dots,a_{r-n}).
\]
From now on, in accordance to the convention made a while
ago, we shall denote by:
\[
\widehat{X}_1^*,\dots,\widehat{X}_n^*,
\widehat{X}_{n+1}^*,\dots,\widehat{X}_r^*
\]
the coframe on $\mathcal{ P}$ which is dual to the frame:
\[
\widehat{X}_1,\dots,\widehat{X}_n,
\widehat{X}_{n+1},\dots,\widehat{X}_r
\]
made of all the constant fields $\widehat{ X}_k := 
\omega^{ -1} ( {\sf x}_k)$ on 
$\mathcal{ P}$, $k = 1, \dots, r$, so that
by definition:
\[
\widehat{X}_{k_1}^*
(\widehat{X}_{k_2})
=
\delta_{k_2}^{k_1}
\ \ \ \ \ \ \ \ \ \ \ \ \
{\scriptstyle{(k_1,\,\,k_2\,=\,1\,\cdots\,r)}}.
\]
In fact, because for any $k_1 = 1, \dots, r$, one has:
\[
{\sf x}_{k_2}
=
\omega(\widehat{X}_{k_2})
=
\sum_{k_1=1}^r\,\omega^{k_1}(\widehat{X}_{k_2})\,{\sf x}_{k_1}
\ \ \ \ \ \ \ \
\text{\rm whence}
\ \ \ \ \
\omega^{k_1}(\widehat{X}_{k_2})
=
\delta_{k_2}^{k_1}
\]
also holds, it follows that the $\widehat{ X}_k$ were already known
and that we must admit the coincidence of notation:
\[
\omega^k
\equiv
\widehat{X}_k^*
\ \ \ \ \ \ \ \ \ \ \ \ \
{\scriptstyle{(k\,=\,1\,\cdots\,r)}}.
\]

\begin{Proposition}
\label{representation-curvature-coefficients}
The coefficients $\kappa_{ i_1, i_2}^k (p)$ of the curvature
function:
\[
\kappa(p)
=
\sum_{1\leqslant i_1<i_2\leqslant n}\,
\sum_{k=1}^r\,
\kappa_{i_1,i_2}^k(p)\,\,
{\sf x}_{i_1}^*\wedge{\sf x}_{i_2}^*
\otimes
{\sf x}_k
\]
express explicitly as follows in terms the structure
constants of $\mathfrak{ g}$ and in terms of the constant
frame $\widehat{ X}_i$ and coframe $\widehat{ X}_k^*$:
\[
\boxed{
\kappa_{i_1,i_2}^k(p)
=
c_{i_1,i_2}^k
-
\widehat{X}_k^*
\big(
\big[\widehat{X}_{i_1},\widehat{X}_{i_2}\big]
\big)}\,.
\]
\end{Proposition}

\proof
It suffices to apply $\omega_p^{ -1}$ to the two 
extreme sides of:
\[
\footnotesize
\aligned
\sum_{k=1}^r\,\kappa_{i_1,i_2}^k(p)\,{\sf x}_k
&
=
\kappa(p)({\sf x}_{i_1},{\sf x}_{i_2})
\\
\explain{apply Lemma~\ref{first-curvature-function}}
\ \ \ \ \ \ \ \ \ \ \ 
&
=
[{\sf x}_{i_1},{\sf x}_{i_2}]_{\mathfrak{g}}
-
\omega_p\big([\omega_p^{-1}({\sf x}_{i_1}),\,
\omega_p^{-1}({\sf x}_{i_2})]_{\mathfrak{g}}\big)
\\
&
=
\sum_{k=1}^r\,c_{i_1,i_2}^k\,{\sf x}_k
-
\omega_p\big(\big[\widehat{X}_{i_1},\widehat{X}_{i_2}\big]\big)
\\
&
=
\omega_p
\Big(
\sum_{k=1}^r\,c_{i_1,i_2}^k\,\widehat{X}_k
-
\big[\widehat{X}_{i_1},\widehat{X}_{i_2}\big]
\Big),
\endaligned
\]
which simply yields:
\[
\small
\aligned
\sum_{k=1}^r\,\kappa_{i_1,i_2}^k(p)\,\widehat{X}_k
=
\sum_{k=1}^r\,c_{i_1,i_2}^k\,\widehat{X}_k
-
\big[\widehat{X}_{i_1},\widehat{X}_{i_2}\big].
\endaligned
\]
Finally, the $\widehat{ X}_k$-component of the right-hand side 
is immediately obtained by applying $\widehat{ X}_k^*$ to 
both sides, and this delivers the expression 
boxed in the lemma.
\endproof

\section{Effective Construction of a Normal, Regular 
\\
Cartan-Tanaka Connection}
\label{Cartan-construction}

\HEAD{Effective Construction of a Normal, Regular Cartan-Tanaka 
Connection}{
Mansour Aghasi, Joël Merker, and Masoud Sabzevari}

\subsection{General form of the unknown Cartan connection frame}
\label{General-form}
Again with an $M^3 \subset \mathbb{C}^2$ being a Levi nondegenerate 
$\mathcal{ C}^6$-smooth real
hypersurface which is an arbitrary deformation of the
Heisenberg sphere equipped with three real coordinates $(x, y, u)$ as
above, suppose (in advance) that we may construct a Cartan connection:
\[
\omega_p\colon \ \
T_p\mathcal{P}
\longrightarrow
\mathfrak{g}
\ \ \ \ \ \ \ \ \ \ \ \ \ 
{\scriptstyle{(p\,\in\,\mathcal{G})}}
\]
on a certain (local) $P$-principal bundle $\mathcal{ P}$ locally equal
to the product $M^3 \times P$, where $P$ is the unique (local,
connected) five-dimensional Lie group associated to the {\em abstract}
Lie algebra:
\[
\mathfrak{p}
:=
{\rm Span}_{\mathbb{R}}
\big(
{\sf d},\,{\sf r},\,{\sf i}_1,\,{\sf i}_2,\,{\sf j}
\big)
\]
corresponding to the five generators ${\sf D}$, ${\sf R}$, ${\sf
I}_1$, ${\sf I}_2$, ${\sf J}$ of the Lie isotropy algebra of the
origin $0 \in \mathbb{ H}^3$, namely where the brackets between the
abstract ${\sf d}$, ${\sf r}$, ${\sf i}_1$, ${\sf i}_2$, ${\sf j}$ are
assumed to be exactly the same as those already shown between ${\sf
D}$, ${\sf R}$, ${\sf I}_1$, ${\sf I}_2$, ${\sf J}$, namely $[ {\sf
d}, {\sf i}_1] = {\sf i}_1$, {\em etc.}

At first, one needs to make explicit some corresponding five
left-invariant vector fields on the (local) Lie group $P$. Let $P$ be
equipped with five real coordinates denoted by $(a, b, c, d, e)$. Then
one can simply take exactly the same five left-invariant vector fields
as those shown in~\cite{ EMS}, namely:
\[
\aligned
D
&
:=-
a\,{\textstyle{\frac{\partial}{\partial a}}}
-
b\,{\textstyle{\frac{\partial}{\partial b}}}
-
c\,{\textstyle{\frac{\partial}{\partial c}}}
-
d\,{\textstyle{\frac{\partial}{\partial d}}}
-
2e\,{\textstyle{\frac{\partial}{\partial e}}}
\\
R
&
:=-b\,{\textstyle{\frac{\partial}{\partial a}}}
+
a\,{\textstyle{\frac{\partial}{\partial
b}}}
+
d\,{\textstyle{\frac{\partial}{\partial c}}}
-
c\,{\textstyle{\frac{\partial}{\partial d}}}
\\
I_1
&
:={\textstyle{\frac{\partial}{\partial a}}}
-
b\,{\textstyle{\frac{\partial}{\partial e}}}
\\
I_2
&
:={\textstyle{\frac{\partial}{\partial b}}}
+
a\,{\textstyle{\frac{\partial}{\partial e}}}
\\
J
&
:={\textstyle{\frac{1}{2}}}\,
{\textstyle{\frac{\partial}{\partial e}}},
\endaligned
\]
for one verifies that the Lie brackets between these vector fields are
precisely the same as before, namely:

\medskip
\begin{center}
\label{table-5-restricted}
\begin{tabular}[t]{ l | l l l l l }
& $D$ & $R$ &
$I_1$ & $I_2$ & $J$
\\
\hline
$D$ & $0$ & $0$ & $I_1$ & $I_2$ & $2\,J$
\\
$R$ & $*$ & $0$ & $-I_2$ & $I_1$ & $0$
\\
$I_1$ & $*$ & $*$ & $0$ & $4\,J$ & $0$
\\
$I_2$ & $*$ & $*$ & $*$ & $0$ & $0$
\\
$J$ & $*$ & $*$ & $*$ & $*$ & $0$
\end{tabular}
\end{center}

\noindent
so that they indeed form a basis for the Lie algebra
$\mathfrak{ p}$ of the abstract
Lie group $P$. Then according to property {\bf (ii)} that a Cartan
connection should enjoy on the product $\mathcal{ P} = M^3 \times P$
which is naturally equipped with the eight real coordinates:
\[
(x,y,u,a,b,c,d,e)
=:
\text{\rm an arbitrary point}\,
p\in\mathcal{P},
\]
it follows that some five corresponding generating {\em vertical
constant} fields must necessarily be of the plain form:
\[
\aligned
\widehat{D}\big\vert_p
&
:=
\omega_p^{-1}({\sf d})
\equiv
D
\\
\widehat{R}\big\vert_p
&
:=
\omega_p^{-1}({\sf r})
\equiv
R
\\
\widehat{I}_1\big\vert_p
&
:=
\omega_p^{-1}({\sf i}_1)
\equiv
I_1
\\
\widehat{I}_2\big\vert_p
&
:=
\omega_p^{-1}({\sf i}_2)
\equiv
I_2
\\
\widehat{J}\big\vert_p
&
:=
\omega_p^{-1}({\sf j})
\equiv
J,
\endaligned
\]
where we identify vector fields on the $(x,y,u,a,b,c,d,e)$-space to
vector fields on the $(a,b,c,d,e)$-space. We notice passim and for
coherence of thought that the vanishing of the curvature in vertical
bi-directions just means here that brackets between these five fields
should correspond, through $\omega$, to abstract brackets within the
Lie algebra $\mathfrak{ p}$, say for instance:
\[
\big[\widehat{D},\,\widehat{I}_1\big]
=
\big[\omega^{-1}({\sf d}),\,\omega^{-1}({\sf i}_1)\big]
=
\omega^{-1}
\big(
[{\sf d},\,{\sf i}_1]_{\mathfrak{p}}\big)
=
\omega^{-1}\big(
{\sf i}_1\big)
=
\widehat{I}_1,
\]
and this is clearly the case by the unique choice $\widehat{ D} \equiv
D$ and $\widehat{ I}_1 \equiv I_1$ made at the moment. Only the three
vector fields:
\[
\aligned
\widehat{T}
:=
\omega^{-1}({\sf t}),
\ \ \ \ \
\widehat{H}_1
:=
\omega^{-1}({\sf h}_1),
\ \ \ \ \
\widehat{H}_2
:=
\omega^{-1}({\sf h}_2)
\endaligned
\]
are really unknown, and their determination should be subjected to
specific constraints that shall be presented right below.

In summary, in order to construct the sought Cartan connection, we
must find certain functions $\alpha_{ \cdot \cdot}$ of the eight
variables $(x, y, u, a, b, c, d, e)$ as coefficients for the lifted
horizontal vector fields:
\[
\footnotesize
\left\{ 
\aligned 
\widehat{T}
&
:=
\alpha_{tt}\,T
+
\alpha_{th_1}\,H_1
+
\alpha_{th_2}\,H_2
+
\alpha_{td}\,D
+
\alpha_{tr}\,R
+
\alpha_{ti_1}\,I_1
+
\alpha_{ti_2}\,I_2
+
\alpha_{tj}\,J
\\
\widehat{H}_1
&
:=
\ \ \ \ \ \ \ \ \ \
\alpha_{h_1h_1}\,H_1
+
\alpha_{h_1h_2}\,H_2
+
\alpha_{h_1d}\,D+\alpha_{h_1r}\,R
+
\alpha_{h_1i_1}\,I_1
+
\alpha_{h_1i_2}\,I_2
+
\alpha_{h_1j}\,J
\\
\widehat{H}_2
&
:=
\ \ \ \ \ \ \ \ \ \
\alpha_{h_2h_1}\,H_1
+
\alpha_{h_2h_2}\,H_2
+
\alpha_{h_2d}\,D
+
\alpha_{h_2r}\,R
+
\alpha_{h_2i_1}\,I_1
+
\alpha_{h_2i_2}\,I_2
+
\alpha_{h_2j}\,J,
\endaligned\right.
\]
while, as we said, the lifts of the vertical vector fields identify to
the vertical vector fields themselves:
\[
\left\{ \aligned \widehat{D}
&
\equiv
D
=
-a\,{\textstyle{\frac{\partial}{\partial a}}}
-
b\,{\textstyle{\frac{\partial}{\partial b}}}
-
c\,{\textstyle{\frac{\partial}{\partial c}}}
-
d\,{\textstyle{\frac{\partial}{\partial d}}}
-
2e\,{\textstyle{\frac{\partial}{\partial e}}}
\\
\widehat{R}
&
\equiv
R
=
-b\,{\textstyle{\frac{\partial}{\partial a}}}
+
a\,{\textstyle{\frac{\partial}{\partial b}}}
+
d\,{\textstyle{\frac{\partial}{\partial c}}}
-
c\,{\textstyle{\frac{\partial}{\partial d}}}
\\
\widehat{I}_1
&
\equiv
I_1
=
{\textstyle{\frac{\partial}{\partial a}}}
-
b\,{\textstyle{\frac{\partial}{\partial e}}}
\\
\widehat{I}_2
&
\equiv
I_2
=
{\textstyle{\frac{\partial}{\partial b}}}
+
a\,{\textstyle{\frac{\partial}{\partial e}}}
\\
\widehat{J}
&
\equiv
J
=
{\textstyle{\frac{1}{2}}}\,{\textstyle{\frac{\partial}{\partial e}}}.
\endaligned\right.
\]

\subsection{Conditions for the determination of the Cartan connection}
We have to determine appropriate functions $\alpha_{
{}_\bullet{}_\bullet}$ of the eight coordinates $(x, u, v, a, b, c,
d, e)$ in a neighborhood of the origin in such a way that the
following four conditions are satisfied.

\begin{itemize}
\label{Conditions}

\smallskip\item[{\bf (c1)}]
For any $X = D,R,I_1,I_2,J$ and any $Y=H_1,H_2,T$ with corresponding
${\sf x} = {\sf d}, {\sf r}, {\sf i}_1, {\sf i}_2, {\sf j}$ and ${\sf
y} = {\sf h}_1, {\sf h}_2, {\sf t}$, one should have:
\begin{equation}
\label{zero-curvature-h-v}
\big[
\widehat{X},\widehat{Y}
\big]
=
\widehat{[{\sf x},{\sf y}]_{\mathfrak{g}}},
\end{equation}
or equivalently in terms of the $\mathfrak{ g}$-valued
one-form:
\[
\big[
\omega^{-1}({\sf x}),\,\omega^{-1}({\sf y})
\big]
=
\omega^{-1}\big(
[{\sf x},\,{\sf y}]_{\mathfrak{g}}\big).
\]
As is known ({\em see} {\em e.g.}~\cite{Crampin} page 3), if Lie
groups are assumed to be connected (ours are, because we suppose they
are local), this condition is equivalent to the equivariancy $R_h^* (
\omega) = {\rm Ad} ( h^{ -1}) \circ \omega$ that $\omega$ should enjoy
under right translations by elements $h \in H$
(Section~\ref{Cartan-connections-coordinates}).

\smallskip\item[{\bf (c2)}]
For each $p \in \mathcal{ P}$, the map $\omega_p \colon T_p \mathcal{
P} \to \mathfrak{ g}$ should be an isomorphism. We postpone the
checking of this property to the end of all computations, but at least
here, we may observe that this property is equivalent to the fact that
the eight (local) vector fields $\widehat{ T}$, $\widehat{ H}_1$,
$\widehat{ H}_2$, $\widehat{ D}$, $\widehat{ R}$, $\widehat{ I}_1$,
$\widehat{ I}_2$, $\widehat{ J}$ constitute a {\em frame} near the
origin (linear independency). Furthermore, since $T$, $H_1$, $H_2$
live in the $(x, u, v)$-space and since $D$, $R$, $I_1$, $I_2$, $J$
already make a frame in the $(a, b, c, d, e)$-space, this property is
equivalent to the fact that the $(T, H_1, H_2)$-components of
$\widehat{ T}$, $\widehat{ H}_1$, $\widehat{ H}_2$ are independent,
namely:
\[
\alpha_{tt}\big(
\alpha_{h_1h_1}\alpha_{h_2h_2}
-
\alpha_{h_1h_2}\alpha_{h_2h_1}\big)
\]
should not vanish in a neighborhood of $0$.

\smallskip\item[{\bf (c3)}]
\label{c3-c4-anticipate}
The obtained Cartan connection $\omega$ should be {\em normal}, namely
the codifferential operator $\partial^*$ should annihilate all
homogeneous curvature components, {\em i.e.}:
\[
\partial_{[h]}^\ast\big(\kappa_{[h]}\big)
\equiv 0
\ \ \ \ \ \ \ \ \ \ \ \ \
\text{\rm for}\ \
h=1,\dots,5.
\] 
Similarly as for condition {\bf (c2)}, we
shall also examine this condition at the end of our main calculations.

\smallskip\item[{\bf (c4)}]
The connection should be {\sl regular}, that is to say, all curvatures
of negative homogeneities should vanish. In fact, possible
homogeneities of $2$-cochains are just $0, 1, 2, 3, 4, 5$
(Subsection~\ref{homogeneity-Psi-0-5}), and we will see that making
$\kappa^{[ 0]} = 0$ is the easiest thing. Furthermore, thanks to the
fact that second cohomologies vanish up to homogeneity $h = 3$
according to the table on p.~\pageref{dimensional-cohomologies}), the
Proposition~\ref{cohomology-result} insures that the first nonzero
homogeneous curvature components can only be $\kappa_{ [4]}$, whence
something a bit better than
regularity will in a certain sense hold freely.

\end{itemize}\smallskip

In fact, the process of construction ({\em cf.}~\cite{ EMS}) will
mainly consist in annihilating as many curvatures components as
possible, and without calling to $\partial^*$, we will be able to
annihilate $\kappa_{ [0]}$ (easiest thing), $\kappa_{ [1]}$, $\kappa_{
[2]}$ and $\kappa_{ [3]}$ by an appropriate progressive building of
$\omega$ which requires somewhat hard elimination computations.

\subsection{Explicit (sought) dual coframe}
\label{Dual-section}
Before beginning by carefully inspecting how to fulfill the first,
main condition {\bf (c1)}, we still need further preliminaries.

On the $(x,y,u)$-space, it is clear that the three vector fields
$H_1$, $H_2$ and $T$ make a frame, for we remember that:
\[
H_1\big\vert_0
=
{\textstyle{\frac{\partial}{\partial x}}}\big\vert_0,
\ \ \ \ \
H_2\big\vert_0
=
{\textstyle{\frac{\partial}{\partial y}}}\big\vert_0,
\ \ \ \ \
T\big\vert_0
=
{\textstyle{\frac{\partial}{\partial u}}}\big\vert_0,
\]
and consequently, there exists a well defined {\em dual coframe} which
is composed of three one-forms $H_1^*$, $H_2^*$ and $T^*$
satisfying by definition:
\[
H_1^*(H_1)
=
1,
\ \ \ \ \
H_1^*(H_2)
=
0,
\ \ \ \ \
H_1^*(T)
=
0,
\ \ \ \ \
etc.
\]
At the moment, we shall not consider it to be necessary to
express explicitly $H_1^*$, $H_2^*$ and $T^*$ in terms of ${\rm d} x$,
${\rm d}y$, ${\rm d}u$, leaving such a task to a computer at the very
end of all our constructions.

Now, our eight unknown vector fields $\widehat{ T} = \omega^{ -1} (
{\sf t})$, \dots, $\widehat{ J} = \omega^{ -1} ( {\sf j})$ on the $(x,
y, u, a, b, c, d, e)$-space read as linear combinations:
\[
\footnotesize
\aligned
\widehat{T}
&
=
\alpha_{tt}\,T
+\ \
\alpha_{th_1}\,H_1
+\ \
\alpha_{th_2}\,H_2
+\ \
\alpha_{td}\,D
+\
\alpha_{tr}\,R
+\ \
\alpha_{ti_1}\,I_1
+\ \ \
\alpha_{ti_2}\,I_2
+\ \
\alpha_{tj}\,J
\\
\widehat{H}_1
&
=
\ \ \ \ \ \ \ \ \ \ \ \ \
\alpha_{h_1h_1}\,H_1
+
\alpha_{h_1h_2}\,H_2
+
\alpha_{h_1d}\,D
+
\alpha_{h_1r}\,R
+
\alpha_{h_1i_1}\,I_1
+
\alpha_{h_1i_2}\,I_2
+
\alpha_{h_1j}\,J
\\
\widehat{H}_2
&
=
\ \ \ \ \ \ \ \ \ \ \ \ \
\alpha_{h_2h_1}\,H_1
+
\alpha_{h_2h_2}\,H_2
+
\alpha_{h_2d}\,D
+
\alpha_{h_2r}\,R
+
\alpha_{h_2i_1}\,I_1
+
\alpha_{h_2i_2}\,I_2
+
\alpha_{h_2j}\,J
\\
\widehat{D}
&
=
\ \ \ \ \ \ \ \ \ \ \ \ \ \ \ \ \ \ \ \ \ \ \ \ \ \ \ \ \ \ \ \
\ \ \ \ \ \ \ \ \ \ \ \ \ \ \ \ \ \ \ \ \ \ \ \ \ \ \ \ \ \ 
D
\\
\widehat{R}
&
=
\ \ \ \ \ \ \ \ \ \ \ \ \ \ \ \ \ \ \ \ \ \ \ \ \ \ \ \ \ \ \ \
\ \ \ \ \ \ \ \ \ \ \ \ \ \ \ \ \ \ \ \ \ \ \ \ \ \ \ \ \ \ \ \
\ \ \ \ \ \ \ \ \ \ \ \ \ \ \
R
\\
\widehat{I}_1
&
=
\ \ \ \ \ \ \ \ \ \ \ \ \ \ \ \ \ \ \ \ \ \ \ \ \ \ \ \ \ \ \ \
\ \ \ \ \ \ \ \ \ \ \ \ \ \ \ \ \ \ \ \ \ \ \ \ \ \ \ \ \ \ \ \
\ \ \ \ \ \ \ \ \ \ \ \ \ \ \ \ \ \ \ \ \ \ \ \ \ \ \ \ \ \ \ \ \
I_1
\\
\widehat{I}_2
&
=
\ \ \ \ \ \ \ \ \ \ \ \ \ \ \ \ \ \ \ \ \ \ \ \ \ \ \ \ \ \ \ \
\ \ \ \ \ \ \ \ \ \ \ \ \ \ \ \ \ \ \ \ \ \ \ \ \ \ \ \ \ \ \ \
\ \ \ \ \ \ \ \ \ \ \ \ \ \ \ \ \ \ \ \ \ \ \ \ \ \ \ \ \ \ \ \
\ \ \ \ \ \ \ \ \ \ \ \ \ \ \ \ \ \ \
I_2
\\
\widehat{J}
&
=
\ \ \ \ \ \ \ \ \ \ \ \ \ \ \ \ \ \ \ \ \ \ \ \ \ \ \ \ \ \ \ \
\ \ \ \ \ \ \ \ \ \ \ \ \ \ \ \ \ \ \ \ \ \ \ \ \ \ \ \ \ \ \ \
\ \ \ \ \ \ \ \ \ \ \ \ \ \ \ \ \ \ \ \ \ \ \ \ \ \ \ \ \ \ \ \
\ \ \ \ \ \ \ \ \ \ \ \ \ \ \ \ \ \ \ \ \ \ \ \ \ \ \ \ \ \ \ \ \ \
J
\endaligned
\]
of the eight fields $T$, \dots, $J$ with certain 22 unknown
coefficients $\alpha_{ tt}$, \dots, $\alpha_{ h2j}$, where we use {\em
letters instead of numbers} as lower indices, the logic of indexing
being clearly visible. We need to know explicitly the dual (unknown)
coframe:
\[
\widehat{T}^*,\ \ \
\widehat{H}_1^*,\ \ \
\widehat{H}_2^*,\ \ \
\widehat{D}^*,\ \ \
\widehat{R}^*,\ \ \
\widehat{I}_1^*,\ \ \
\widehat{I}_2^*,\ \ \
\widehat{J}^*,
\]
and the task is simple. Indeed, we recall the elementary fact that, in
a standard vector space ${\rm Span}_{\mathbb{R}} (e_1, e_2, \dots,
e_n)$, the dual of an arbitrary frame:
\[
v_k
:=
{\textstyle{\sum_{i=1}^n}}\,
\alpha_{ki}\,
e_i
\ \ \ \ \ \ \ \ \ \ \ \ \ {\scriptstyle{(k\,=\,1\,\cdots\,n)}}
\]
is a coframe of the form:
\[
v_l^*
:=
{\textstyle{\sum_{j=1}^n}}\,
\beta_{lj}\,e_j^*
\ \ \ \ \ \ \ \ \ \ \ \ \ {\scriptstyle{(l\,=\,1\,\cdots\,n)}},
\]
where the matrix $(\beta_{lj})_{1\leqslant l\leqslant n}^{ 1 \leqslant
j \leqslant n}$ is just the {\em transpose-inverse} of the initial
matrix $(\alpha_{ki})_{1\leqslant k\leqslant n}^{ 1 \leqslant i
\leqslant n}$. Leaving it to a computer to invert and to transpose the
above $8 \times 8$ matrix ({\em see} \cite{ AMSMaple}) whose
determinant is clearly equal to:
\[
\alpha_{tt}
=
\big(
\alpha_{h_1h_1}\alpha_{h_2h_2}
-
\alpha_{h_1h_2}\alpha_{h_2h_1}
\big),
\]
we find without much pain expressions of the form:
\[
\aligned
\widehat{T}^*
&
=
\beta_{tt}\,T^*
\\
\widehat{H}_1^*
&
=
\beta_{h_1t}\,T^*
+
\beta_{h_1h_1}\,H_1^*
+
\beta_{h_1h_2}\,H_2^*
\\
\widehat{H}_2^*
&
=
\beta_{h_2t}\,T^*
+
\beta_{h_2h_1}\,H_1^*
+
\beta_{h_2h_2}\,H_2^*
\\
\widehat{D}^*
&
=
\beta_{dt}\,T^*
+
\beta_{dh_1}\,H_1^*
+
\beta_{dh_2}\,H_2^*
+
D^*
\\
\widehat{R}^*
&
=
\beta_{rt}\,T^*
+
\beta_{rh_1}\,H_1^*
+
\beta_{rh_2}\,H_2^*
+
R^*
\\
\widehat{I}_1^*
&
=
\beta_{i_1t}\,T^*
+
\beta_{i_1h_1}\,H_1^*
+
\beta_{i_1h_2}\,H_2^*
+
I_1^*
\\
\widehat{I}_2^*
&
=
\beta_{i_2t}\,T^*
+
\beta_{i_2h_1}\,H_1^*
+
\beta_{i_2h_2}\,H_2^*
+
I_2^*
\\
\widehat{J}^*
&
=
\beta_{jt}\,T^*
+
\beta_{jh_1}\,H_1^*
+
\beta_{jh_2}\,H_2^*
+
J^*,
\endaligned
\]
where the 22 coefficients $\beta_{tt}$, \dots, $\beta_{ jh_2}$
express themselves rationally and explicitly
in terms of the $\alpha$'s right as follows:
\[
\scriptsize
\aligned
\beta_{tt}
&
:=
{\textstyle{\frac{1}{\alpha_{tt}}}},
\\
\beta_{h_1t}
&
:=
{\textstyle{\frac{
-\alpha_{th_1}\alpha_{h_2h_2}+\alpha_{th_2}\alpha_{h_2h_1}
}{
\alpha_{tt}(\alpha_{h_1h_1}\alpha_{h_2h_2}-
\alpha_{h_1h_2}\alpha_{h_2h_1})}}},
\ \ \ \ \
\beta_{h1h1}
:=
{\textstyle{\frac{
\alpha_{h_2h_2}\alpha_{tt}
}{
\alpha_{tt}(\alpha_{h_1h_1}\alpha_{h_2h_2}-
\alpha_{h_1h_2}\alpha_{h_2h_1})}}},
\\
\beta_{h_1h_2}
&
:=
{\textstyle{\frac{
-\alpha_{h_2h_1}\alpha_{tt}
}{
\alpha_{tt}(\alpha_{h_1h_1}\alpha_{h_2h_2}-
\alpha_{h_1h_2}\alpha_{h_2h_1})}}},
\ \ \ \ \
\beta_{h_2t}
:=
{\textstyle{\frac{
\alpha_{th_1}\alpha_{h_1h_2}-\alpha_{th_2}\alpha_{h_1h1}
}{
\alpha_{tt}(\alpha_{h_1h_1}\alpha_{h_2h_2}-
\alpha_{h_1h_2}\alpha_{h_2h_1})}}},
\\
\beta_{h_2h_1}
&
:=
{\textstyle{\frac{
-\alpha_{h_1h_2}\alpha_{tt}
}{
\alpha_{tt}(\alpha_{h_1h_1}\alpha_{h_2h_2}-
\alpha_{h_1h_2}\alpha_{h_2h_1})}}},
\ \ \ \ \
\beta_{h_2h_2}
:=
{\textstyle{\frac{
\alpha_{h_1h_1}\alpha_{tt}
}{
\alpha_{tt}(\alpha_{h_1h_1}\alpha_{h_2h_2}-
\alpha_{h_1h_2}\alpha_{h_2h_1})}}},
\endaligned
\]
\[
\scriptsize
\aligned
\beta_{dt}
&
:=
{\textstyle{\frac{
-\alpha_{th_1}\alpha_{h_1h_2}\alpha_{h_2d}
+
\alpha_{th_1}\alpha_{h_1d}\alpha_{h_2h_2}
+
\alpha_{th_2}\alpha_{h_2d}\alpha_{h_1h_1}
-
\alpha_{th_2}\alpha_{h_2h_1}\alpha_{h_1d}
-
\alpha_{td}\alpha_{h_1h_1}\alpha_{h_2h_2}
+
\alpha_{td}\alpha_{h_1h_2}\alpha_{h_2h_1}
}{
\alpha_{tt}(\alpha_{h_1h_1}\alpha_{h_2h_2}-
\alpha_{h_1h_2}\alpha_{h_2h_1})}}},
\\
\beta_{dh_1}
&
:=
{\textstyle{\frac{
\alpha_{h_1h_2}\alpha_{h_2d}\alpha_{tt}
-
\alpha_{h_1d}\alpha_{h_2h_2}\alpha_{tt}
}{
\alpha_{tt}(\alpha_{h_1h_1}\alpha_{h_2h_2}-
\alpha_{h_1h_2}\alpha_{h_2h_1})}}},
\ \ \ \ \
\beta_{dh_2}
:=
{\textstyle{\frac{
-\alpha_{h_2d}\alpha_{h_1h_1}\alpha_{tt}
+
\alpha_{h_2h_1}\alpha_{h_1d}\alpha_{tt}
}{
\alpha_{tt}(\alpha_{h_1h_1}\alpha_{h_2h_2}-
\alpha_{h_1h_2}\alpha_{h_2h_1})}}},
\endaligned
\]
\[
\scriptsize
\aligned
\beta_{rt}
&
:=
{\textstyle{\frac{
-\alpha_{th_1}\alpha_{h_1h_2}\alpha_{h_2r}
+
\alpha_{th_1}\alpha_{h_1r}\alpha_{h_2h_2}
+
\alpha_{th_2}\alpha_{h_2r}\alpha_{h_1h_1}
-
\alpha_{th_2}\alpha_{h_2h_1}\alpha_{h_1r}
-
\alpha_{tr}\alpha_{h_1h_1}\alpha_{h_2h_2}
+
\alpha_{tr}\alpha_{h_1h_2}\alpha_{h_2h_1}
}{
\alpha_{tt}(\alpha_{h_1h_1}\alpha_{h_2h_2}-
\alpha_{h_1h_2}\alpha_{h_2h_1})}}},
\\
\beta_{rh_1}
&
:=
{\textstyle{\frac{
\alpha_{h_1h_2}\alpha_{h_2r}\alpha_{tt}
-
\alpha_{h_1r}\alpha_{h_2h_2}\alpha_{tt}
}{
\alpha_{tt}(\alpha_{h_1h_1}\alpha_{h_2h_2}-
\alpha_{h_1h_2}\alpha_{h_2h_1})}}},
\ \ \ \ \
\beta_{rh_2}
:=
{\textstyle{\frac{
-\alpha_{h_2r}\alpha_{h_1h_1}\alpha_{tt}
+
\alpha_{h_2h_1}\alpha_{h_1r}\alpha_{tt}
}{
\alpha_{tt}(\alpha_{h_1h_1}\alpha_{h_2h_2}-
\alpha_{h_1h_2}\alpha_{h_2h_1})}}},
\endaligned
\]
\[
\scriptsize
\aligned
\beta_{i_1t}
&
:=
{\textstyle{\frac{
-\alpha_{th_1}\alpha_{h_1h_2}\alpha_{h_2i_1}
+
\alpha_{th_1}\alpha_{h_1i_1}\alpha_{h_2h_2}
+
\alpha_{th_2}\alpha_{h_2i_1}\alpha_{h_1h_1}
-
\alpha_{th_2}\alpha_{h_2h_1}\alpha_{h_1i_1}
-
\alpha_{ti_1}\alpha_{h_1h_1}\alpha_{h_2h_2}
+
\alpha_{ti_1}\alpha_{h_1h_2}\alpha_{h_2h_1}
}{
\alpha_{tt}(\alpha_{h_1h_1}\alpha_{h_2h_2}-
\alpha_{h_1h_2}\alpha_{h_2h_1})}}},
\\
\beta_{i_1h_1}
&
:=
{\textstyle{\frac{
\alpha_{h_1h_2}\alpha_{h_2i_1}\alpha_{tt}
-
\alpha_{h_1i_1}\alpha_{h_2h_2}\alpha_{tt}
}{
\alpha_{tt}(\alpha_{h_1h_1}\alpha_{h_2h_2}-
\alpha_{h_1h_2}\alpha_{h_2h_1})}}},
\ \ \ \ \
\beta_{i_1h_2}
:=
{\textstyle{\frac{
-\alpha_{h_2i_1}\alpha_{h_1h_1}\alpha_{tt}
+
\alpha_{h_2h_1}\alpha_{h_1i_1}\alpha_{tt}
}{
\alpha_{tt}(\alpha_{h_1h_1}\alpha_{h_2h_2}-
\alpha_{h_1h_2}\alpha_{h_2h_1})}}},
\endaligned
\]
\[
\scriptsize
\aligned
\beta_{i_2t}
&
:=
{\textstyle{\frac{
-\alpha_{th_1}\alpha_{h_1h_2}\alpha_{h_2i_2}
+
\alpha_{th_1}\alpha_{h_1i_2}\alpha_{h_2h_2}
+
\alpha_{th_2}\alpha_{h_2i_2}\alpha_{h_1h_1}
-
\alpha_{th_2}\alpha_{h_2h_1}\alpha_{h_1i_2}
-
\alpha_{ti_2}\alpha_{h_1h_1}\alpha_{h_2h_2}
+
\alpha_{ti_2}\alpha_{h_1h_2}\alpha_{h_2h_1}
}{
\alpha_{tt}(\alpha_{h_1h_1}\alpha_{h_2h_2}-
\alpha_{h_1h_2}\alpha_{h_2h_1})}}},
\\
\beta_{i_2h_1}
&
:=
{\textstyle{\frac{
\alpha_{h_1h_2}\alpha_{h_2i_2}\alpha_{tt}
-
\alpha_{h_1i_2}\alpha_{h_2h_2}\alpha_{tt}
}{
\alpha_{tt}(\alpha_{h_1h_1}\alpha_{h_2h_2}-
\alpha_{h_1h_2}\alpha_{h_2h_1})}}},
\ \ \ \ \
\beta_{i_2h_2}
:=
{\textstyle{\frac{
-\alpha_{h_2i_2}\alpha_{h_1h_1}\alpha_{tt}
+
\alpha_{h_2h_1}\alpha_{h_1i_2}\alpha_{tt}
}{
\alpha_{tt}(\alpha_{h_1h_1}\alpha_{h_2h_2}-
\alpha_{h_1h_2}\alpha_{h_2h_1})}}},
\endaligned
\]
\[
\scriptsize
\aligned
\beta_{jt}
&
:=
{\textstyle{\frac{
-\alpha_{th_1}\alpha_{h_1h_2}\alpha_{h_2j}
+
\alpha_{th_1}\alpha_{h_1j}\alpha_{h_2h_2}
+
\alpha_{th_2}\alpha_{h_2j}\alpha_{h_1h_1}
-
\alpha_{th_2}\alpha_{h_2h_1}\alpha_{h_1j}
-
\alpha_{tj}\alpha_{h_1h_1}\alpha_{h_2h_2}
+
\alpha_{tj}\alpha_{h_1h_2}\alpha_{h_2h_1}
}{
\alpha_{tt}(\alpha_{h_1h_1}\alpha_{h_2h_2}-
\alpha_{h_1h_2}\alpha_{h_2h_1})}}},
\\
\beta_{jh_1}
&
:=
{\textstyle{\frac{
\alpha_{h_1h_2}\alpha_{h_2j}\alpha_{tt}
-
\alpha_{h_1j}\alpha_{h_2h_2}\alpha_{tt}
}{
\alpha_{tt}(\alpha_{h_1h_1}\alpha_{h_2h_2}-
\alpha_{h_1h_2}\alpha_{h_2h_1})}}},
\ \ \ \ \
\beta_{jh_2}
:=
{\textstyle{\frac{
-\alpha_{h_2j}\alpha_{h_1h_1}\alpha_{tt}
+
\alpha_{h_2h_1}\alpha_{h_1j}\alpha_{tt}
}{
\alpha_{tt}(\alpha_{h_1h_1}\alpha_{h_2h_2}-
\alpha_{h_1h_2}\alpha_{h_2h_1})}}},
\endaligned
\]

\subsection{Expressions of the $\alpha_{{}_\bullet {}_\bullet}$
in terms of fiber variables}
Of the mentioned four conditions, we will start by examining
thoroughly {\bf (c1)}, namely~\thetag{
\ref{zero-curvature-h-v}}. 
By expressing, in terms of the basic frame $T$,
$H_1$, $H_2$, $D$, $R$, $I_1$, $I_2$, $J$, each one of the $5 \times 3
= 15$ equalities of the form $[\widehat{X}, \widehat{Y}] = \widehat{[{\sf
x}, {\sf y}]}_{\frak g}$ for ${\sf x}={\sf d,r},{\sf i}_1,{\sf
i}_2,{\sf j}$ and for ${\sf y}={\sf t},{\sf h}_1,{\sf h}_2$, we get at
each time eight scalar equations (sometimes only seven, because one
equation may reduce to just $0 = 0$). To set up these scalar
equations, it suffices to expand any bracket of the form $\big[ \widehat{
 V}, \, \alpha W \big]$ just as $\widehat{ V} ( \alpha) \, W + \alpha\,
\big[ \widehat{ V}, \, \widehat{ W}
\big]$, to remind that any $T$, $H_1$, $H_2$
trivially commutes with any $D$, $R$, $I_1$, $I_2$, $J$ and to expand
any occurring Lie bracket between $D$, $R$, $I_1$, $I_2$, $J$ by means
of the table on p.~\pageref{table-5-restricted}. Then the $15 \times
8$ equations in question read as follows, where we list the
coefficients of $T$, of $H_1$, of $H_2$, of $D$, of $R$, of $I_1$, of
$I_2$, of $J$ rigorously in this order (incidentally, we drop exactly
$2 \times 5 = 10$ trivial equations $0 = 0$ which appear due
to the fact that $\widehat{ H}_1$ and $\widehat{ H}_2$ have zero
$T$-component by our above assumption):

\smallskip\noindent{\bf (1)}\,
$\boxed{[\widehat{D},\,\widehat{T}]+2\widehat{T}=0}:$
\[
\footnotesize
\aligned
&
\widehat{D}(\alpha_{tt})+2\alpha_{tt}=0, \ \ \
\widehat{D}(\alpha_{th_1})+2\alpha_{th_1}=0, \ \ \
\widehat{D}(\alpha_{th_2})+2\alpha_{th_2}=0, \ \ \
\widehat{D}(\alpha_{td})+2\alpha_{td}=0,
\ \ \ \ \ \ \ \ \ \ \ \ \ \ \ \ \ \ \ \ \ \ \ \ \ \ \ \ \ \ \ \ \ \ \ \ \
\\
&
\widehat{D}(\alpha_{tr})+2\alpha_{tr}=0, \ \ \
\widehat{D}(\alpha_{ti_1})+3\alpha_{ti_1}=0, \ \ \
\widehat{D}(\alpha_{ti_2})+3\alpha_{ti_2}=0, \ \ \
\widehat{D}(\alpha_{tj})+4\alpha_{tj}=0.
\endaligned
\]

\noindent{\bf (2)}\,
$\boxed{[\widehat{D},\,\widehat{H}_1]+\widehat{H}_1=0}:$
\[
\scriptsize
\aligned
&
\widehat{D}(\alpha_{h_1h_1})+\alpha_{h_1h_1}=0, \ \ \
\widehat{D}(\alpha_{h_1h_2})+\alpha_{h_1h_2}=0, \ \ \
\widehat{D}(\alpha_{h_1d})+\alpha_{h_1d}=0, \ \ \
\widehat{D}(\alpha_{h_1r})+\alpha_{h_1r}=0,
\ \ \ \ \ \ \ \ \ \ \ \ \ \ \ \ \ \ \ \ \ \ \ \ \ \ \ \ \ \ \ \ \ \ \ \ \
\\
&
\widehat{D}(\alpha_{h_1i_1})+2\alpha_{h_1i_1}=0, \ \ \
\widehat{D}(\alpha_{h_1i_2})+2\alpha_{h_1i_2}=0, \ \ \
\widehat{D}(\alpha_{h_1j})+3\alpha_{h_1j}=0.
\endaligned
\]

\noindent{\bf (3)}\,
$\boxed{[\widehat{D},\,\widehat{H}_2]+\widehat{H}_2=0}:$
\[
\scriptsize
\aligned
&
\widehat{D}(\alpha_{h_2h_1})+\alpha_{h_2h_1}=0, \ \ \
\widehat{D}(\alpha_{h_2h_2})+\alpha_{h_2h_2}=0, \ \ \
\widehat{D}(\alpha_{h_2d})+\alpha_{h_2d}=0, \ \ \
\widehat{D}(\alpha_{h_2r})+\alpha_{h_2r}=0,
\\
&
\widehat{D}(\alpha_{h_2i_1})+2\alpha_{h_2i_1}=0, \ \ \
\widehat{D}(\alpha_{h_2i_2})+2\alpha_{h_2i_2}=0, \ \ \
\widehat{D}(\alpha_{h_2j})+3\alpha_{h_2j}=0.
\endaligned
\]

\noindent{\bf (4)}\,
$\boxed{[\widehat{R},\,\widehat{T}]=0}:$
\[
\footnotesize
\aligned
&
\widehat{R}(\alpha_{tt})=0, \ \ \
\widehat{R}(\alpha_{th_1})=0, \ \ \
\widehat{R}(\alpha_{th_2})=0,\ \ \
\widehat{R}(\alpha_{td})=0,\ \ \
\widehat{R}(\alpha_{tr})=0,
\ \ \ \ \ \ \ \ \ \ \ \ \ \ \ \ \ \ \ \ \ \ \ \ \ \ \ \ \ \ \ \ \ \ \ \ \
\\
&
\widehat{R}(\alpha_{ti_1})+\alpha_{ti_2}=0,\ \ \
\widehat{R}(\alpha_{ti_2})-\alpha_{ti_1}=0,\ \ \
\widehat{R}(\alpha_{tj})=0.
\endaligned
\]

\noindent{\bf (5)}\,
$\boxed{[\widehat{R},\,\widehat{H}_1]+\widehat{H}_2=0}:$
\[
\footnotesize
\aligned
&
\widehat{R}(\alpha_{h_1h_1})+\alpha_{h_2h_1}=0, \ \ \
\widehat{R}(\alpha_{h_1h_2})+\alpha_{h_2h_2}=0,\ \ \
\widehat{R}(\alpha_{h_1d})+\alpha_{h_2d}=0,
\ \ \ \ \ \ \ \ \ \ \ \ \ \ \ \ \ \ \ \ \ \ \ \ \ \ \ \ \ \ \ \ \ \ \ \ \
\\
&
\widehat{R}(\alpha_{h_1r})+\alpha_{h_2r}=0, \ \ \
\widehat{R}(\alpha_{h_1i_1})+\alpha_{h_1i_2}+\alpha_{h_2i_1}=0, \ \ \
\widehat{R}(\alpha_{h_1i_2})-\alpha_{h_1i_1}+\alpha_{h_2i_2}=0,
\\
&
\widehat{R}(\alpha_{h_1j})+\alpha_{h_2j}=0.
\endaligned
\]

\noindent{\bf (6)}\,
$\boxed{[\widehat{R},\,\widehat{H}_2]-\widehat{H}_1=0}:$
\[
\footnotesize
\aligned
&
\widehat{R}(\alpha_{h_2h_1})-\alpha_{h_1h_1}=0,\ \ \
\widehat{R}(\alpha_{h_2h_2})-\alpha_{h_1h_2}=0, \ \ \
\widehat{R}(\alpha_{h_2d})-\alpha_{h_1d}=0,
\ \ \ \ \ \ \ \ \ \ \ \ \ \ \ \ \ \ \ \ \ \ \ \ \ \ \ \ \ \ \ \ \ \ \ \ \
\\
&
\widehat{R}(\alpha_{h_2r})-\alpha_{h_1r}=0, \ \ \
\widehat{R}(\alpha_{h_2i_1})+\alpha_{h_2i_2}-\alpha_{h_1i_1}=0,\ \ \
\widehat{R}(\alpha_{h_2i_2})-\alpha_{h_2i_1}-\alpha_{h_1i_2}=0, \ \ \
\\
&
\widehat{R}(\alpha_{h_2j})-\alpha_{h_1j}=0.
\endaligned
\]

\noindent{\bf (7)}\,
$\boxed{[\widehat{I}_1,\,\widehat{T}]+\widehat{H}_1=0}:$
\[
\footnotesize
\aligned
&
\widehat{I}_1(\alpha_{tt})=0, \ \ \
\widehat{I}_1(\alpha_{th_1})+\alpha_{h_1h_1}=0, \ \ \
\widehat{I}_1(\alpha_{th_2})+\alpha_{h_1h_2}=0, \ \ \
\widehat{I}_1(\alpha_{td})+\alpha_{h_1d}=0,
\ \ \ \ \ \ \ \ \ \ \ \ \ \ \ \ \ \ \ \ \ \ \ \ \ \ \ \ \ \ \ \ \ \ \ \ \
\\
&
\widehat{I}_1(\alpha_{tr})+\alpha_{h_1r}=0, \ \ \
\widehat{I}_1(\alpha_{ti_1})-\alpha_{td}+\alpha_{h_1i_1}=0, \ \ \
\widehat{I}_1(\alpha_{ti_2})+\alpha_{tr}+\alpha_{h_1i_2}=0, \ \ \
\\
&
\widehat{I}_1(\alpha_{tj})+4\alpha_{ti_2}+\alpha_{h_1j}=0.
\endaligned
\]

\noindent{\bf (8)}\,
$\boxed{[\widehat{I}_1,\,\widehat{H}_1]+6\widehat{R}=0}:$
\[
\footnotesize
\aligned
&
\widehat{I}_1(\alpha_{h_1h_1})=0, \ \ \
\widehat{I}_1(\alpha_{h_1h_2})=0, \ \ \
\widehat{I}_1(\alpha_{h_1d})=0, \ \ \
\widehat{I}_1(\alpha_{h_1r})+6=0, \ \ \
\ \ \ \ \ \ \ \ \ \ \ \ \ \ \ \ \ \ \ \ \ \ \ \ \ \ \ \ \ \ \ \ \ \ \ \ \
\\
&
\widehat{I}_1(\alpha_{h_1i_1})-\alpha_{h_1d}=0, \ \ \
\widehat{I}_1(\alpha_{h_1i_2})+\alpha_{h_1r}=0, \ \ \
\widehat{I}_1(\alpha_{h_1j})+4\alpha_{h_1i_2}=0.
\endaligned
\]

\noindent{\bf (9)}\,
$\boxed{[\widehat{I}_1,\,\widehat{H}_2]-2\widehat{D}=0}:$
\[
\footnotesize
\aligned
&
\widehat{I}_1(\alpha_{h_2h_1})=0, \ \ \
\widehat{I}_1(\alpha_{h_2h_2})=0, \ \ \
\widehat{I}_1(\alpha_{h_2d})-2=0, \ \ \
\widehat{I}_1(\alpha_{h_2r})=0, \ \ \
\ \ \ \ \ \ \ \ \ \ \ \ \ \ \ \ \ \ \ \ \ \ \ \ \ \ \ \ \ \ \ \ \ \ \ \ \
\\
&
\widehat{I}_1(\alpha_{h_2i_1})-\alpha_{h_2d}=0, \ \ \
\widehat{I}_1(\alpha_{h_2i_2})+\alpha_{h_2r}=0, \ \ \
\widehat{I}_1(\alpha_{h_2j})+4\alpha_{h_2i_2}=0.
\endaligned
\]

\noindent{\bf (10)}\,
$\boxed{[\widehat{I}_2,\,\widehat{T}]+\widehat{H}_2=0}:$
\[
\footnotesize
\aligned
&
\widehat{I}_2(\alpha_{tt})=0, \ \ \
\widehat{I}_2(\alpha_{th_1})+\alpha_{h_2h_1}=0, \ \ \
\widehat{I}_2(\alpha_{th_2})+\alpha_{h_2h_2}=0, \ \ \
\widehat{I}_2(\alpha_{td})+\alpha_{h_2d}=0,
\ \ \ \ \ \ \ \ \ \ \ \ \ \ \ \ \ \ \ \ \ \ \ \ \ \ \ \ \ \ \ \ \ \ \ \ \
\\
&
\widehat{I}_2(\alpha_{tr})+\alpha_{h_2r}=0,\ \ \
\widehat{I}_2(\alpha_{ti_1})-\alpha_{tr}+\alpha_{h_2i_1}=0, \ \ \
\widehat{I}_2(\alpha_{ti_2})-\alpha_{td}+\alpha_{h_2i_2}=0, \ \ \
\\
&
\widehat{I}_2(\alpha_{tj})-4\alpha_{ti_1}+\alpha_{h_2j}=0.
\endaligned
\]

\noindent{\bf (11)}\,
$\boxed{[\widehat{I}_2,\,\widehat{H}_1]+2\widehat{D}=0}:$
\[
\footnotesize
\aligned
&
\widehat{I}_2(\alpha_{h_1h_1})=0, \ \ \
\widehat{I}_2(\alpha_{h_1h_2})=0, \ \ \
\widehat{I}_2(\alpha_{h_1d})+2=0, \ \ \
\widehat{I}_2(\alpha_{h_1r})=0,
\ \ \ \ \ \ \ \ \ \ \ \ \ \ \ \ \ \ \ \ \ \ \ \ \ \ \ \ \ \ \ \ \ \ \ \ \
\\
&
\widehat{I}_2(\alpha_{h_1i_1})-\alpha_{h_1r}=0, \ \ \
\widehat{I}_2(\alpha_{h_1i_2})-\alpha_{h_1d}=0, \ \ \
\widehat{I}_2(\alpha_{h_1j})-4\alpha_{h_1i_1}=0.
\endaligned
\]

\noindent{\bf (12)}\,
$\boxed{[\widehat{I}_2,\,\widehat{H}_2]+6\widehat{R}=0}:$
\[
\footnotesize
\aligned
&
\widehat{I}_2(\alpha_{h_2h_1})=0, \ \ \
\widehat{I}_2(\alpha_{h_2h_2})=0, \ \ \
\widehat{I}_2(\alpha_{h_2d})=0, \ \ \
\widehat{I}_2(\alpha_{h_2r})+6=0, \ \ \
\ \ \ \ \ \ \ \ \ \ \ \ \ \ \ \ \ \ \ \ \ \ \ \ \ \ \ \ \ \ \ \ \ \ \ \ \
\\
&
\widehat{I}_2(\alpha_{h_2i_1})-\alpha_{h_2r}=0, \ \ \
\widehat{I}_2(\alpha_{h_2i_2})-\alpha_{h_2d}=0, \ \ \
\widehat{I}_2(\alpha_{h_2j})-4\alpha_{h_2i_1}=0.
\endaligned
\]

\noindent{\bf (13)}\,
$\boxed{[\widehat{J},\,\widehat{T}]+\widehat{D}=0}:$
\[
\footnotesize
\aligned
&
\widehat{J}(\alpha_{tt})=0, \ \ \
\widehat{J}(\alpha_{th_1})=0, \ \ \
\widehat{J}(\alpha_{th_2})=0,\ \ \
\widehat{J}(\alpha_{td})+1=0,
\ \ \ \ \ \ \ \ \ \ \ \ \ \ \ \ \ \ \ \ \ \ \ \ \ \ \ \ \ \ \ \ \ \ \ \ \
\\
&
\widehat{J}(\alpha_{tr})=0, \ \ \
\widehat{J}(\alpha_{ti_1})=0, \ \ \
\widehat{J}(\alpha_{ti_2})=0, \ \ \
\widehat{J}(\alpha_{tj})-2\alpha_{td}=0,
\endaligned
\]

\noindent{\bf (14)}\,
$\boxed{[\widehat{J},\,\widehat{H}_1]+\widehat{I}_1=0}:$
\[
\footnotesize
\aligned
&
\widehat{J}(\alpha_{h_1h_1})=0, \ \ \
\widehat{J}(\alpha_{h_1h_2})=0, \ \ \
\widehat{J}(\alpha_{h_1d})=0, \ \ \
\widehat{J}(\alpha_{h_1r})=0, \ \ \
\ \ \ \ \ \ \ \ \ \ \ \ \ \ \ \ \ \ \ \ \ \ \ \ \ \ \ \ \ \ \ \ \ \ \ \ \
\\
&
\widehat{J}(\alpha_{h_1i_1})+1=0, \ \ \
\widehat{J}(\alpha_{h_1i_2})=0, \ \ \
\widehat{J}(\alpha_{h_1j})-2\,\alpha_{h_1d}=0.
\endaligned
\]

\noindent{\bf (15)}\,
$\boxed{[\widehat{J},\,\widehat{H}_2]+\widehat{I}_2=0}:$
\[
\footnotesize
\aligned
&
\widehat{J}(\alpha_{h_2h_1})=0, \ \ \
\widehat{J}(\alpha_{h_2h_2})=0, \ \ \
\widehat{J}(\alpha_{h_2d})=0, \ \ \
\widehat{J}(\alpha_{h_2r})=0, \ \ \
\ \ \ \ \ \ \ \ \ \ \ \ \ \ \ \ \ \ \ \ \ \ \ \ \ \ \ \ \ \ \ \ \ \ \ \ \
\\
&
\widehat{J}(\alpha_{h_2i_1})=0, \ \ \
\widehat{J}(\alpha_{h_2i_2})+1=0, \ \ \
\widehat{J}(\alpha_{h_2j})-2\,\alpha_{h_2d}=0.
\endaligned
\]
We therefore get a system of precisely 110 first-order partial
differential equations having the twenty-two unknowns $\alpha_{ tt}$,
\dots, $\alpha_{ h_2j}$ in the space $(x, u, v, a, b, c, d, e)$, 
and the differentiations only 
involve the five partial derivatives $\frac{ \partial}{ \partial
a}$, $\frac{ \partial}{\partial b}$, $\frac{ \partial}{ \partial c}$,
$\frac{ \partial}{\partial d}$, $\frac{ \partial}{\partial e}$.
Importantly, this system is of first order and is linear, hence this
is the reason why its (large) solution set could be found rather
quickly by means of Maple
(\cite{ AMSMaple}). The 22 appearing functions $\delta_1, \dots,
\delta_{ 22}$ will be determined later: they now constitute the only
remaining unknown part of $\alpha_{tt}, \dots, \alpha_{ h_2j}$.

\begin{Lemma}
\label{fiber-type}
The general solution of the above system {\bf (1)}, {\bf (2)}, \dots,
{\bf (15)} of 110 partial differential equations is polynomial of
degree $\leqslant 4$ with respect to the five vertical variables $a$,
$b$, $c$, $d$, $e$ of the principal bundle $P$, and it involves 22
coefficients $\delta_1 ( x, y, z)$, \dots, $\delta_{ 22} ( x, y, z)$
that are arbitrary (smooth) functions of the horizontal variables
$(x, y, z)$:
\[
\footnotesize
\aligned
\alpha_{tt}
&
=
(c^2+d^2)\,\delta_{22},
\\
\alpha_{th_1}
&
=
-(ad+bc)\,\delta_{13}+(ac-bd)\,\delta_{14}+(c^2+d^2)\,\delta_{21},
\\
\alpha_{th2}
&
=
-(ad+bc)\,\delta_{11}+(ac-bd)\,\delta_{12}+(c^2+d^2)\,\delta_{20}
\\
\alpha_{td}
&
=
(-{\textstyle{\frac{1}{4}}}bc-{\textstyle{\frac{1}{4}}}ad)\,\delta_1
+
({\textstyle{\frac{1}{4}}}ac-{\textstyle{\frac{1}{4}}}bd)\,\delta_2
+
({\textstyle{\frac{1}{4}}}c^2+{\textstyle{\frac{1}{4}}}d^2)\,\delta_{15}
-
2e,
\\
\alpha_{tr}
&
=
({\textstyle{\frac{1}{4}}}c^2+{\textstyle{\frac{1}{4}}}d^2)\,\delta_4
+
({\textstyle{\frac{1}{2}}}ac-{\textstyle{\frac{1}{2}}}bd)\,\delta_7
+
({\textstyle{\frac{1}{2}}}c^2+{\textstyle{\frac{1}{2}}}d^2)\,\delta_9
-
({\textstyle{\frac{1}{2}}}ad+{\textstyle{\frac{1}{2}}}bc)\,\delta_{10}
+
\\
&
\ \ \ \ \ \ \ \ \ \ \ \ \ \ \ \ \ \ \ \ \ \ \ \ \ \ \ \ \ \ \ \ \ \ \
\ \ \ \ \ \ \ \ \ \ \ \ \ \ \ \ \ \ \ \ \ \ \ \ \ \ \ \ \ \ \ \ \ \ \
+
({\textstyle{\frac{1}{2}}}c^2+{\textstyle{\frac{1}{2}}}d^2)\,\delta_{19}
+
3b^2+3a^2,
\endaligned
\]
\[
\footnotesize
\aligned
\alpha_{ti_1}
&
=
-({\textstyle{\frac{1}{4}}}a^2d+{\textstyle{\frac{1}{4}}}abc)\,\delta_1
+
(-{\textstyle{\frac{1}{4}}}abd+{\textstyle{\frac{1}{4}}}a^2c)\,\delta_2
+
({\textstyle{\frac{1}{24}}}d^3+{\textstyle{\frac{1}{24}}}c^2d)\,\delta_3
+
({\textstyle{\frac{1}{8}}}d^2b+{\textstyle{\frac{1}{8}}}bc^2)\,\delta_4
+
\\
&
\ \ \ \ \
+
({\textstyle{\frac{1}{8}}}c^3+{\textstyle{\frac{1}{8}}}cd^2)\,\delta_5
+
({\textstyle{\frac{1}{8}}}ac^2+{\textstyle{\frac{1}{8}}}ad^2)\,\delta_6
+
(-{\textstyle{\frac{1}{2}}}db^2+{\textstyle{\frac{1}{2}}}bca)\,\delta_7
-
({\textstyle{\frac{1}{4}}}dcb+{\textstyle{\frac{1}{4}}}ad^2)\,\delta_8
+
\\
&
\ \ \ \ \
+
({\textstyle{\frac{1}{4}}}bc^2+{\textstyle{\frac{1}{2}}}d^2b
-{\textstyle{\frac{1}{4}}}dca)\,\delta_9
-
({\textstyle{\frac{1}{2}}}bda+{\textstyle{\frac{1}{2}}}cb^2)\,\delta_{10}
+
({\textstyle{\frac{1}{4}}}d^2a+{\textstyle{\frac{1}{4}}}ac^2)\,\delta_{15}
+
\\
&
\ \ \ \ \
+
({\textstyle{\frac{1}{4}}}c^3+{\textstyle{\frac{1}{4}}}cd^2)\,\delta_{16}
+
({\textstyle{\frac{1}{4}}}c^2d+{\textstyle{\frac{1}{4}}}d^3)\,\delta_{17}
+
({\textstyle{\frac{1}{2}}}d^2b+{\textstyle{\frac{1}{2}}}bc^2)\,\delta_{19}
+
2a^2b+2b^3,
\endaligned
\]
\[
\footnotesize
\aligned
\alpha_{ti_2}
&
=
-({\textstyle{\frac{1}{4}}}bda+{\textstyle{\frac{1}{4}}}cb^2)\,\delta_1
+
(-{\textstyle{\frac{1}{4}}}db^2+{\textstyle{\frac{1}{4}}}bca)\,\delta_2
+
({\textstyle{\frac{1}{24}}}c^3+{\textstyle{\frac{1}{24}}}cd^2)\,\delta_3
+
(-{\textstyle{\frac{1}{8}}}ac^2-{\textstyle{\frac{1}{8}}}ad^2)\,\delta_4
+
\\
&
\ \ \ \ \
+
(-{\textstyle{\frac{1}{8}}}c^2d-{\textstyle{\frac{1}{8}}}d^3)\,\delta_5
+
({\textstyle{\frac{1}{8}}}d^2b+{\textstyle{\frac{1}{8}}}bc^2)\,\delta_6
+
(-{\textstyle{\frac{1}{2}}}a^2c+{\textstyle{\frac{1}{2}}}bda)\,\delta_7
-
({\textstyle{\frac{1}{4}}}dca+{\textstyle{\frac{1}{4}}}bc^2)\,\delta_8
+
\\
&
\ \ \ \ \
+
(-{\textstyle{\frac{1}{2}}}ac^2+{\textstyle{\frac{1}{4}}}dcb
-{\textstyle{\frac{1}{4}}}d^2a)\,\delta_9
+
({\textstyle{\frac{1}{2}}}a^2d+{\textstyle{\frac{1}{2}}}bca)\,\delta_{10}
+
({\textstyle{\frac{1}{4}}}bc^2+{\textstyle{\frac{1}{4}}}d^2b)\,\delta_{15}
-
\\
&
\ \ \ \ \
-
({\textstyle{\frac{1}{4}}}c^2d+{\textstyle{\frac{1}{4}}}d^3)\,\delta_{16}
+
({\textstyle{\frac{1}{4}}}c^3+{\textstyle{\frac{1}{4}}}cd^2)\,\delta_{17}
+
(-{\textstyle{\frac{1}{2}}}ac^2
-{\textstyle{\frac{1}{2}}}d^2a)\,\delta_{19}
-
2a^3-2ab^2,
\endaligned
\]
\[
\footnotesize
\aligned
\alpha_{tj}
&
=
-(acb+ade)\,\delta_1+(-edb+eca)\,\delta_2
+(cb^2a+ca^3-db^3-dba^2)\,\delta_7
+
\\
&
\ \ \ \ \
+
(-d^2ab+c^2ab-dcb^2+dca^2)\,\delta_8
+
(-2dcba+c^2a^2+d^2b^2)\,\delta_9
\\
&
\ \ \ \ \
+
(-da^3-db^2a-cb^3-cba^2)\,\delta_{10}
+(ed^2+ec^2)\,\delta_{15}
+
\\
&
\ \ \ \ \
+
(ad^3+c^3b+dac^2+cbd^2)\,\delta_{16}+(-cad^2+bd^3+dbc^2-c^3a)\,\delta_{17}
+
\\
&
\ \ \ \ \
+
(d^4+c^4+2d^2c^2)\,\delta_{18}+(b^2c^2+d^2b^2+c^2a^2+a^2d^2)\,\delta_{19}
+
6a^2b^2
-
\\
&
\ \ \ \ \ \ \ \ \ \ \ \ \ \ \ \ \ \ \ \ \ \ \ \ \ \ \ \ \ \ \ \ \ \ \ 
\ \ \ \ \ \ \ \ \ \ \ \ \ \ \ \ \ \ \ \ \ \ \ \ \ \ \ \ \ \ \ \ \ \ \ 
\ \ \ \ \ \ \ \ \ \ \ \ \ \ \ \ 
-
4e^2+3a^4+3b^4,
\endaligned
\]
\[
\footnotesize
\aligned
\alpha_{h_1h_1}
&
=
(d)\,\delta_{13}-(c)\,\delta_{14},
\\
\alpha_{h_1h_2}
&
=
(d)\,\delta_{11}-(c)\,\delta_{12},
\\
\alpha_{h_1d}
&
=
({\textstyle{\frac{1}{4}}}d)\,\delta_1
-
({\textstyle{\frac{1}{4}}}c)\,\delta_2
-
2b,
\\
\alpha_{h1r}
&
=
-({\textstyle{\frac{1}{2}}}c)\,\delta_7
+
({\textstyle{\frac{1}{2}}}d)\,\delta_{10}
-
6a,
\\
\alpha_{h_1i_1}
&
=
({\textstyle{\frac{1}{4}}}ad)\,\delta_1
-
({\textstyle{\frac{1}{4}}}ac)\,\delta_2
-
({\textstyle{\frac{1}{8}}}c^2+{\textstyle{\frac{1}{8}}}d^2)\,\delta_6
-
({\textstyle{\frac{1}{2}}}bc)\,\delta_7
+
({\textstyle{\frac{1}{4}}}d^2)\,\delta_8
+
({\textstyle{\frac{1}{4}}}cd)\,\delta_9
+
\\
&
\ \ \ \ \ \ \ \ \ \ \ \ \ \ \ \ \ \ \ \ \ \ \ \ \ \ \ \ \ \ \ \ \ \ \ 
\ \ \ \ \ \ \ \ \ \ \ \ \ \ \ \ \ \ \ \ \ \ \ \ \ \ \ \ \ \ \ \ \ \ \ 
\ \ \ \ \ \ \ \ \ \ \ \ \ \ \ \ 
+
({\textstyle{\frac{1}{2}}}bd)\,\delta_{10}
-
4ab-2e,
\endaligned
\]
\[
\footnotesize
\aligned
\alpha_{h_1i_2}
&
=
({\textstyle{\frac{1}{4}}}bd)\delta_1
-
({\textstyle{\frac{1}{4}}}bc)\,\delta_2
-
({\textstyle{\frac{1}{8}}}c^2+{\textstyle{\frac{1}{8}}}d^2)\,\delta_4
+
({\textstyle{\frac{1}{2}}}ac)\,\delta_7
+
({\textstyle{\frac{1}{4}}}cd)\,\delta_8
-
({\textstyle{\frac{1}{4}}}d^2)\,\delta_9
-
\\
&
\ \ \ \ \ \ \ \ \ \ \ \ \ \ \ \ \ \ \ \ \ \ \ \ \ \ \ \ \ \ \ \ \ \ \ 
\ \ \ \ \ \ \ \ \ \ \ \ \ \ \ \ \ \ \ \ \ \ \ \ \ \ \ \ \ \ \ \ \ \ \ 
\ \ \ \ \ \ \ \ \ \ \ \ \ \ \ \ 
-
({\textstyle{\frac{1}{2}}}ad)\,\delta_{10}
+
3a^2-b^2,
\endaligned
\]
\[
\footnotesize
\aligned
\alpha_{h_1j}
&
=
(de)\,\delta_1-(ce)\,\delta_2
-
({\textstyle{\frac{1}{6}}}c^3+{\textstyle{\frac{1}{6}}}cd^2)\,\delta_3
+
({\textstyle{\frac{1}{2}}}ac^2+{\textstyle{\frac{1}{2}}}d^2a)\,\delta_4
+
({\textstyle{\frac{1}{2}}}d^3+{\textstyle{\frac{1}{2}}}c^2d)\,\delta_5
-
\\
&
\ \ \ \ \
-
({\textstyle{\frac{1}{2}}}d^2b+{\textstyle{\frac{1}{2}}}bc^2)\,\delta_6
-
(a^2c+cb^2)\,\delta_7+(-dca+d^2b)\,\delta_8+(dcb+d^2a)\,\delta_9
+
\\
&
\ \ \ \ \
+(a^2d+db^2)\,\delta_{10}-8be-4a^3-4ab^2,
\endaligned
\]
\[
\footnotesize
\aligned
\alpha_{h_2h_1}
&
=
(c)\,\delta_{13}+(d)\,\delta_{14},
\\
\alpha_{h_2h_2}
&
=
(c)\,\delta_{11}+(d)\,\delta_{12},
\\
\alpha_{h_2d}
&
=
({\textstyle{\frac{1}{4}}}c)\,\delta_1
+
({\textstyle{\frac{1}{4}}}d)\,\delta_2
+
2a,
\\
\alpha_{h_2r}
&
=
({\textstyle{\frac{1}{2}}}d)\,\delta_7
+
({\textstyle{\frac{1}{2}}}c)\,\delta_{10}
-
6b,
\\
\alpha_{h_2i_1}
&
=
({\textstyle{\frac{1}{4}}}ac)\,\delta_1
+
({\textstyle{\frac{1}{4}}}ad)\,\delta_2
+
({\textstyle{\frac{1}{8}}}c^2+{\textstyle{\frac{1}{8}}}d^2)\,\delta_4
+
({\textstyle{\frac{1}{2}}}bd)\,\delta_7
+
({\textstyle{\frac{1}{4}}}cd)\,\delta_8
+
({\textstyle{\frac{1}{4}}}c^2)\,\delta_9
+
\\
&
\ \ \ \ \ \ \ \ \ \ \ \ \ \ \ \ \ \ \ \ \ \ \ \ \ \ \ \ \ \ \ \ \ \ \ 
\ \ \ \ \ \ \ \ \ \ \ \ \ \ \ \ \ \ \ \ \ \ \ \ \ \ \ \ \ \ \ \ \ \ \ 
\ \ \ \ \ \ \ \ \ \ \ \ \ \ \ \ 
({\textstyle{\frac{1}{2}}}bc)\,\delta_{10}
-
3b^2+a^2,
\endaligned
\]
\[
\footnotesize
\aligned
\alpha_{h_1i_2}
&
=
({\textstyle{\frac{1}{4}}}cb)\,\delta_1
+
({\textstyle{\frac{1}{4}}}bd)\,\delta_2
-
({\textstyle{\frac{1}{8}}}c^2+{\textstyle{\frac{1}{8}}}d^2)\,\delta_6
-
({\textstyle{\frac{1}{2}}}ad)\,\delta_7
+
({\textstyle{\frac{1}{4}}}c^2)\,\delta_8
-
({\textstyle{\frac{1}{4}}}dc)\,\delta_9
-
\\
&
\ \ \ \ \ \ \ \ \ \ \ \ \ \ \ \ \ \ \ \ \ \ \ \ \ \ \ \ \ \ \ \ \ \ \ 
\ \ \ \ \ \ \ \ \ \ \ \ \ \ \ \ \ \ \ \ \ \ \ \ \ \ \ \ \ \ \ \ \ \ \ 
\ \ \ \ \ \ \ \ \ \ \ \ \ \ \ \ 
-
({\textstyle{\frac{1}{2}}}ac)\,\delta_{10}
+
4ab-2e,
\endaligned
\]
\[
\footnotesize
\aligned
\alpha_{h_2j}
&
=
(ce)\,\delta_1+(ed)\,\delta_2
+
({\textstyle{\frac{1}{6}}}d^3+{\textstyle{\frac{1}{6}}}c^2d)\,\delta_3
+
({\textstyle{\frac{1}{2}}}bc^2+{\textstyle{\frac{1}{2}}}d^2b)\,\delta_4
+
({\textstyle{\frac{1}{2}}}c^3+{\textstyle{\frac{1}{2}}}cd^2)\,\delta_5
+
\\
&
\ \ \ \ \
+
({\textstyle{\frac{1}{2}}}d^2a+{\textstyle{\frac{1}{2}}}ac^2)\,\delta_6
+
(db^2+da^2)\,\delta_7+(-ac^2+dcb)\,\delta_8+(bc^2+dca)\,\delta_9
+
\\
&
\ \ \ \ \
+
(ca^2+cb^2)\,\delta_{10}-4a^2b+8ae-4b^3,
\endaligned
\]
where the coefficients $\delta_{k}$ are only with respect to the
horizontal coordinates $x,y,u$ and are independent of the fibre
variables.
\end{Lemma}

If we would assume that the functions $\alpha_{ {}_\bullet {}_\bullet}$
would
be independent of the horizontal coordinates, {\em i.e.} that
the functions $\delta_k$ would be 
constant, then condition {\bf (c1)} would
hold. However, a major problem would occur when we would try to
fulfill the normality condition
{\bf (c3)}. Hence, the dependence of the
$\alpha_{ {}_\bullet {}_\bullet}$'s with respect to horizontal
coordinates should better be determined by 
annihilating all curvatures $\kappa^{[
h]}\equiv 0$ of homogeneities $h = 0, 1, 2, 3$. We
will achieve this task in the next subsections
and it will require quite hard elimination computations.

\subsection{Graded differential structure}
\label{Graded-Diff}
Now, the choice of the seven coefficients $\alpha_{ tt}, \alpha_{
th_1}, \alpha_{ th_2}, \alpha_{ h_1h_1}, \alpha_{ h_1 h_2}, \alpha_{
h_2 h_1}, \alpha_{ h_2 h_2}$ is governed by the geometry\footnote{\,
The authors would like to thank Ben McLaughlin and Gerd Schmalz for
their helpful explanations in this regard.
} 
of the graded tangent bundle $T^c M \oplus (TM / T^cM)$. In fact, the
five coefficients $\alpha_{ tt}$ must be fiber coordinates with
respect to some fixed trivialization of $T^c M \oplus (TM / T^cM)$,
which is a principal bundle. Given our two fixed complex-tangent local
vector fields $H_1, H_2 \in \Gamma ( T^cM)$ spanning $T^c M$, the
$(\partial_x, \partial_y, \partial_u)$-part of a first lift $\widehat{
H}_1$ must take account, in terms of the coordinates $(a, b, c, d, e)$
of the principal bundle, of the non-uniqueness of the choice of a
first vector field in $\Gamma ( T^cM)$. Thus, the $(\partial_x,
\partial_y, \partial_u)$-part of $\widehat{ H}_1$ must be of the form
$c\, H_1 + d\, H_2$. Next, the $(\partial_x, \partial_y,
\partial_u)$-part of $\widehat{ H}_2$ must be equal to $J (c\, H_1 +
d\, H_2) = -d\, H_1 + c\, H_2$. With this, the coefficient $\alpha_{
tt}$ must be equal to the $(\partial_x, \partial_y, \partial_u)$-part
of $\big[ \widehat{ H}_1, \, \widehat{ H}_2 \big]$, whence $\alpha_{
tt} = c^2 + d^2$ necessarily. Finally, the choice of the $(\partial_x,
\partial_y, \partial_u)$-part of $\widehat{ T}$, as a section of the
quotient $TM / T^cM$, can still be made up to an arbitrary linear
combination $-a\, H_1 - b\, H_2$, whence $\alpha_{ th_1} = bd - ac$ and
$\alpha_{ th_2} = -ad - bc$. In summary, for geometric reasons, we
must have:
\[
\left\{
\aligned 
\alpha_{tt} 
& 
= 
c^2+d^2, \ \ \ \ \ 
\alpha_{th_1}
= 
bd-ac, \ \ \ \ \
\alpha_{th_2}
=
-ad-bc,
\\
\alpha_{h_1h_1} 
& 
= 
c, \ \ \ \ \ \ \ \ \ \ \ \ \ 
\alpha_{h_1h_2} = d,
\\
\alpha_{h_2h_1} 
& 
= 
-d, \ \ \ \ \ \ \ \ \ \ \alpha_{h_2h_2}=c,
\endaligned\right.
\]
which means, equivalently, that seven of the functions $\delta_{
{}_\bullet {}_\bullet}$ are already completely determined.
\[
\delta_{11}=\delta_{22}=1, \ \ \ \
\delta_{12}=\delta_{13}=\delta_{20}=\delta_{21}=0, \ \ \ \
\delta_{14}=-1.
\]

The remaining 15 undetermined coefficient functions $\delta_k$ will be
determined progressively (and uniquely) by subjecting them to the
conditions that the wanted Cartan connection be normal and regular. In
particular they should be determined such that all curvatures
$\kappa_{[h]}$ in the four homogeneities $h= 0, 1, 2, 3$ should be
zero. The next subsections are devoted to inspecting these conditions,
until we examine homogeneities $h = 4, 5$, for
which we shall take account of the
Bianchi-Tanaka identities, too.

\subsection{Brackets between horizontal vector fields}
For our access to Cartan curvatures, we need as a tool to compute the
Lie brackets between the three horizontal vector fields $\widehat{
H_1}$, $\widehat{ H}_2$, $\widehat{ T}$ whose coefficients $\alpha_{
tt}, \dots, \alpha_{ h2j}$, still unknown, are to be determined so as
to simplify curvatures. Thus three brackets $\big[ \widehat{ H}_1, \,
\widehat{ H}_2 \big]$, $\big[ \widehat{ H}_1, \, \widehat{ T} \big]$,
$\big[ \widehat{ H}_2, \, \widehat{ T} \big]$ must be considered, and
they all are of the general form:
\[
\small
\aligned
&
\big[
\alpha\,T+\beta\,H_1+\gamma\,H_2
+
\delta\,D+\rho\,R+\lambda\,I_1+\mu\,I_2+\nu\,J,
\ \
\\
&
\ \ \ \ \ \ \ \ \ \ \ \ \ \ \ \ \ \ \ \ \ \ \ \ \ \ \ \ \ \ \ \
\alpha'\,T+\beta'\,H_1+\gamma'\,H_2
+
\delta'\,D+\rho'\,R+\lambda'\,I_1+\mu'\,I_2+\nu'\,J
\big]
=
\\
&
=
\big[
\alpha\,T+\beta\,H_1+\gamma\,H_2,
\ \
\alpha'\,T+\beta'\,H_1+\gamma'\,H_2
\big]
+
\\
&
\ \ \ \ \
+
\big[
\alpha\,T+\beta\,H_1+\gamma\,H_2,
\ \
\delta'\,D+\rho'\,R+\lambda'\,I_1+\mu'\,I_2+\nu'\,J
\big]
+
\\
&
\ \ \ \ \
+
\big[
\delta\,D+\rho\,R+\lambda\,I_1+\mu\,I_2+\nu\,J,
\ \
\alpha'\,T+\beta'\,H_1+\gamma'\,H_2
\big]
+
\\
&
\ \ \ \ \
+
\big[
\delta\,D+\rho\,R+\lambda\,I_1+\mu\,I_2+\nu\,J,
\ \
\delta'\,D+\rho'\,R+\lambda'\,I_1+\mu'\,I_2+\nu'\,J
\big].
\endaligned
\]
Applying bilinearity, any obtained bracket of the form $\big[ \phi \,
X, \, \psi Y]$ then expands as: 
\[
\big[ \phi \,
X, \, \psi Y]
=
\phi\psi\,[X,\,Y]+\phi X(\psi)\,Y
- 
\psi Y(\phi)\,X,
\] 
and for the brackets in lines 2 and 3, all first
terms $\phi \psi \, [X, Y]$ vanish. After a reorganization using the
commutator rules between $H_1$, $H_2$, $T$ and those between $D$, $R$,
$I_1$, $I_2$, $J$, we obtain that this general bracket equals:
\[
\scriptsize
\aligned
&
\Big(
\alpha T(\alpha')+\beta H_1(\alpha')+\gamma H_2(\alpha')
-
\alpha'T(\alpha)-\beta'H_1(\alpha)-\gamma'H_2(\alpha)
-
\\
&
-
(\alpha\beta'-\alpha'\beta)\Phi_1
-
(\alpha\gamma'-\alpha'\gamma)\Phi_2
+
4(\beta\gamma'-\beta'\gamma)
+
\\
&
+
\delta D(\alpha')+\rho R(\alpha')+\lambda I_1(\alpha')
+
\mu I_2(\alpha')+\nu J(\alpha')
-
\delta'D(\alpha)-\rho'R(\alpha)-\lambda'I_1(\alpha)
-
\mu'I_2(\alpha)-\nu'J(\alpha)
\Big)\,T
+
\endaligned
\]
\[
\scriptsize
\aligned
&
+
\Big(
\alpha T(\beta')+\beta H_1(\beta')+\gamma H_2(\beta')
-
\alpha'T(\beta)-\beta'H_1(\beta)-\gamma'H_2(\beta)
+
\\
&
+
\delta D(\beta')+\rho R(\beta')+\lambda I_1(\beta')
+
\mu I_2(\beta')+\nu J(\beta')
-
\delta'D(\beta)-\rho'R(\beta)-\lambda'I_1(\beta)
-
\mu'I_2(\beta)-\nu'J(\beta)
\Big)\,H_1
+
\endaligned
\]
\[
\scriptsize
\aligned
&
+
\Big(
\alpha T(\gamma')+\beta H_1(\gamma')+\gamma H_2(\gamma')
-
\alpha'T(\gamma)-\beta'H_1(\gamma)-\gamma'H_2(\gamma)
+
\\
&
+
\delta D(\gamma')+\rho R(\gamma')+\lambda I_1(\gamma')
+
\mu I_2(\gamma')+\nu J(\gamma')
-
\delta'D(\gamma)-\rho'R(\gamma)-\lambda'I_1(\gamma)
-
\mu'I_2(\gamma)-\nu'J(\gamma)
\Big)\,H_2
+
\endaligned
\]
\[
\scriptsize
\aligned
&
+
\Big(
\alpha T(\delta')+\beta H_1(\delta')+\gamma H_2(\delta')
-
\alpha'T(\delta)-\beta'H_1(\delta)-\gamma'H_2(\delta)
+
\\
&
+
\delta D(\delta')+\rho R(\delta')+\lambda I_1(\delta')
+
\mu I_2(\delta')+\nu J(\delta')
-
\delta'D(\delta)-\rho'R(\delta)-\lambda'I_1(\delta)
-
\mu'I_2(\delta)-\nu'J(\delta)
\Big)\,D
+
\endaligned
\]
\[
\scriptsize
\aligned
&
+
\Big(
\alpha T(\rho')+\beta H_1(\rho')+\gamma H_2(\rho')
-
\alpha'T(\rho)-\beta'H_1(\rho)-\gamma'H_2(\rho)
+
\\
&
+
\delta D(\rho')+\rho R(\rho')+\lambda I_1(\rho')
+
\mu I_2(\rho')+\nu J(\rho')
-
\delta'D(\rho)-\rho'R(\rho)-\lambda'I_1(\rho)
-
\mu'I_2(\rho)-\nu'J(\rho)
\Big)\,R
+
\endaligned
\]
\[
\scriptsize
\aligned
&
+
\Big(
\alpha T(\lambda')+\beta H_1(\lambda')+\gamma H_2(\lambda')
-
\alpha'T(\lambda)-\beta'H_1(\lambda)-\gamma'H_2(\lambda)
+
\\
&
+
\delta D(\lambda')+\rho R(\lambda')+\lambda I_1(\lambda')
+
\mu I_2(\lambda')+\nu J(\lambda')
-
\delta'D(\lambda)-\rho'R(\lambda)-\lambda'I_1(\lambda)
-
\mu'I_2(\lambda)-\nu'J(\lambda)
+
\\
&
+
\delta\lambda'-\delta'\lambda
-
\rho\mu'+\rho'\mu
\Big)\,I_1
+
\endaligned
\]
\[
\scriptsize
\aligned
&
+
\Big(
\alpha T(\mu')+\beta H_1(\mu')+\gamma H_2(\mu')
-
\alpha'T(\mu)-\beta'H_1(\mu)-\gamma'H_2(\mu)
+
\\
&
+
\delta D(\mu')+\rho R(\mu')+\lambda I_1(\mu')
+
\mu I_2(\mu')+\nu J(\mu')
-
\delta'D(\mu)-\rho'R(\mu)-\lambda'I_1(\mu)
-
\mu'I_2(\mu)-\nu'J(\mu)
+
\\
&
+
\delta\mu'-\delta'\mu
-
\rho\lambda'+\rho'\lambda
\Big)\,I_2
+
\endaligned
\]
\[
\scriptsize
\aligned
&
+
\Big(
\alpha T(\nu')+\beta H_1(\nu')+\gamma H_2(\nu')
-
\alpha'T(\nu)-\beta'H_1(\nu)-\gamma'H_2(\nu)
+
\\
&
+
\delta D(\nu')+\rho R(\nu')+\lambda I_1(\nu')
+
\mu I_2(\nu')+\nu J(\nu')
-
\delta'D(\nu)-\rho'R(\nu)-\lambda'I_1(\nu)
-
\mu'I_2(\nu)-\nu'J(\nu)
+
\\
&
+
\delta\nu'-\delta'\nu
-
4\lambda\mu'-4\lambda'\mu
\Big)\,J.
\endaligned
\]
As will be useful later, in order to get for instance $\big[ \widehat{
H}_1, \, \widehat{ T}\big]$, it suffices to replace in this expression
$\alpha, \beta, \dots, \nu$ by $0, \alpha_{ h_1 h1}, \dots, \alpha_{
h_1j}$ and $\alpha', \beta', \dots, \nu'$ by $\alpha_{ tt}, \alpha_{
th_1}, \dots, \alpha_{ tj}$.

Now we are ready to start the computation of the curvature
components. To do this, recall at first that
the curvature function $\kappa$ as an
element of the space $\mathcal C^2(\frak g_-,\frak g)$, splits up in 
components of various
homogeneities. In the case of our Lie algebra:
\[
\frak g
=
\frak g_{-2}\oplus\frak g_{-1}
\oplus
\frak g_0\oplus\frak g_1\oplus\frak g_2,
\]
the minimal homogeneity occurs when we one considers the value of
$\kappa$ on $({\sf h}_1,{\sf h}_2)$ in $\frak g_{-2}=\frak
g_{-1-1+0}$. So the minimal homogeneity of $\kappa$ is zero. On the
other hand, the maximal homogeneity occurs when one considers the
value of $\kappa$ on $({\sf h}_i,{\sf t})$, $i=1,2$, in $\frak
g_2=\frak g_{-1-2+5}$. Hence the maximal homogeneity is five. Now let,
$\kappa^{p_{j_1}p_{j_2}}_{q_j}$ be the coefficient of ${\sf q}_j$ in
$\kappa({\sf p}_{j_1},{\sf p}_{j_2})$, for ${\sf p}_{j_1}\in\frak
g_{j_1},{\sf p}_{j_2}\in\frak g_{j_2}$ and ${\sf q}_j\in\frak g_j$,
where naturally $j_1,j_2 < 0$. Hence if $h = j - (j_1 + j_2)$, then
clearly the $h$-homogeneous component of the value of $\kappa({\sf
p}_{j_1},{\sf p}_{j_2})$ is:
\[
\kappa_{[h]}({\sf p}_{j_1},{\sf
p}_{j_2})=\sum_{h=j-(j_1+j_2)}\kappa^{p_{j_1}p_{j_2}}_{q_j} {\sf
q}_j.
\]
Any coefficient $\kappa^{p_{j_1}p_{j_2}}_{q_j}$ is called a {\em
curvature coefficient of homogeneity h}. 
From Proposition~\ref{representation-curvature-coefficients}, we
know that every occurring curvature
coefficient $\kappa^{p_{j_1}p_{j_2}}_{q_j}$
is equal to:
\[
\aligned
 \kappa^{p_{j_1}p_{j_2}}_{q_j}&=\widehat{Q}_j^\ast\big([\omega^{-1}{\sf
p}_{j_1},\omega^{-1}{\sf p}_{j_2}]-\omega^{-1}[{\sf p}_{j_1},{\sf p}_{j_2}]\big)\\
&=\widehat{Q}_j^\ast\big([\widehat{P}_{j_1},\widehat{P}_{j_2}]-\widehat{[{\sf
p}_{j_1},{\sf
p}_{j_2}]}_{\frak g}\big).
\endaligned
\]
Let us list the curvature coefficients corresponding to our sought
Cartan connection,
according to their homogeneities:
\[\footnotesize
\aligned &\fbox{\tiny 0} \ \ \
\kappa^{h_1h_2}_t={\widehat{T}}^\ast([\widehat{H}_1,\widehat{H}_2]-4\widehat{T})
\\
& \fbox{\tiny 1} \ \ \
\kappa^{h_1h_2}_{h_1}={\widehat{H}_1}^\ast([\widehat{H}_1,\widehat{H}_2]-4\widehat{T})
\\
& \fbox{\tiny 1} \ \ \
\kappa^{h_1h_2}_{h_2}={\widehat{H}_2}^\ast([\widehat{H}_1,\widehat{H}_2]-4\widehat{T})
\\
& \fbox{\tiny 1} \ \ \
\kappa^{h_1t}_{t}={\widehat{T}}^\ast([\widehat{H}_1,\widehat{T}])
\\
& \fbox{\tiny 1} \ \ \
\kappa^{h_2t}_{t}={\widehat{T}}^\ast([\widehat{H}_2,\widehat{T}])
\\
& \fbox{\tiny 2} \ \ \
\kappa^{h_1h_2}_{d}={\widehat{D}}^\ast([\widehat{H}_1,\widehat{H}_2]-4\widehat{T})
\\
& \fbox{\tiny 2} \ \ \
\kappa^{h_1h_2}_{r}={\widehat{R}}^\ast([\widehat{H}_1,\widehat{H}_2]-4\widehat{T})
\\
& \fbox{\tiny 2} \ \ \
\kappa^{h_1t}_{h_1}={\widehat{H}_1}^\ast([\widehat{H}_1,\widehat{T}])
\endaligned
 \ \ \ \ \ \
\aligned
& \fbox{\tiny 2} \ \ \
\kappa^{h_1t}_{h_2}={\widehat{H}_2}^\ast([\widehat{H}_1,\widehat{T}])
\\
& \fbox{\tiny 2} \ \ \
\kappa^{h_2t}_{h_1}={\widehat{H}_1}^\ast([\widehat{H}_2,\widehat{T}])
\\
& \fbox{\tiny 2} \ \ \
\kappa^{h_2t}_{h_2}={\widehat{H}_2}^\ast([\widehat{H}_2,\widehat{T}])
\\
& \fbox{\tiny 3} \ \ \
\kappa^{h_1h_2}_{i_1}={\widehat{I}_1}^\ast([\widehat{H}_1,\widehat{H}_2]-4\widehat{T})
\\
& \fbox{\tiny 3} \ \ \
\kappa^{h_1h_2}_{i_2}={\widehat{I}_2}^\ast([\widehat{H}_1,\widehat{H}_2]-4\widehat{T})
\\
& \fbox{\tiny 3} \ \ \
\kappa^{h_1t}_{d}={\widehat{D}}^\ast([\widehat{H}_1,\widehat{T}])
\\
& \fbox{\tiny 3} \ \ \
\kappa^{h_2t}_{d}={\widehat{D}}^\ast([\widehat{H}_2,\widehat{T}])
\\
& \fbox{\tiny 3} \ \ \
\kappa^{h_1t}_{r}={\widehat{R}}^\ast([\widehat{H}_1,\widehat{T}])
\endaligned
\ \ \ \ \ \
\aligned
& \fbox{\tiny 3} \ \ \
\kappa^{h_2t}_{r}={\widehat{R}}^\ast([\widehat{H}_2,\widehat{T}])
\\
& \fbox{\tiny 4} \ \ \
\kappa^{h_1h_2}_{j}={\widehat{J}}^\ast([\widehat{H}_1,\widehat{H}_2]-4\widehat{T})
\\
& \fbox{\tiny 4}\ \ \
\kappa^{h_1t}_{i_1}={\widehat{I}_1}^\ast([\widehat{H}_1,\widehat{T}])
\\
& \fbox{\tiny 4} \ \ \
\kappa^{h_1t}_{i_2}={\widehat{I}_2}^\ast([\widehat{H}_1,\widehat{T}])
\\
& \fbox{\tiny 4} \ \ \
\kappa^{h_2t}_{i_1}={\widehat{I}_1}^\ast([\widehat{H}_2,\widehat{T}])
\\
& \fbox{\tiny 4} \ \ \
\kappa^{h_2t}_{i_2}={\widehat{I}_2}^\ast([\widehat{H}_2,\widehat{T}])
\\
& \fbox{\tiny 5} \ \ \
\kappa^{h_1t}_{j}={\widehat{J}}^\ast([\widehat{H}_1,\widehat{T}])
\\
& \fbox{\tiny 5} \ \ \
\kappa^{h_2t}_{j}={\widehat{J}}^\ast([\widehat{H}_2,\widehat{T}])
\endaligned
\]
From now on, our aim will be to compute all of these 24 curvature
coefficients and to determine the functions 
$\alpha_{ {}_\bullet {}_\bullet}$ in conformity with 
the properties they are subjected to.

\subsection{Homogeneity 0}
\label{Hom0}
In this homogeneity, we encounter only one curvature coefficient:
\[
\kappa^{h_1h_2}_t 
= 
{\widehat{T}}^\ast([\widehat{H}_1,\widehat{H}_2]-4\widehat{T})
=
\frac{1}{\alpha_{tt}}
(4\alpha_{h_1h_1}\alpha_{h_2h_2}-4\alpha_{h_1h_2}\alpha_{h_2h_1})-4,
\]
which is in fact the $\sf t$-component of the value of the curvature
on $({\sf h}_1,{\sf h}_2)$ as a $\frak g$-valued bilinear
map. In order to satisfy
the regularity condition \textbf{(c4)}, this curvature component should
vanish. Hence we should have:
\[
\alpha_{tt}
=
\alpha_{h_1h_1}\alpha_{h_2h_2}-\alpha_{h_1h_2}\alpha_{h_2h_1}.
\]
But this equality is automatically satisfied, as one sees by coming
back to the expressions of $\alpha_{ {}_\bullet {}_\bullet}$
introduced after Lemma~\ref{fiber-type}. Therefore, the desired
condition $\kappa_{ [0] } = 0$ holds.

\subsection{Homogeneity 1}
\label{Hom1}
In this homogeneity, we have four curvature coefficients, while,
according to the regularity condition \textbf{(c4)}, all of them
should vanish. A detailed and latex-ed calculation gives:
\begin{eqnarray*}
\footnotesize\aligned
\kappa^{h_1h_2}_{h_1}
&
=
\widehat{H}^\ast_1([\widehat{H}_1,\widehat{H}_2]-4\widehat{T})
\\
&=\alpha_{h_1r}+\beta_{h_1h_1}\big(\alpha_{h_1h_1}
\zero{H_1(\alpha_{h_2h_1})}-\alpha_{h_2h_1}\zero{H_1(\alpha_{h_1h_1})}-
\alpha_{h_2h_2}\zero{H_2(\alpha_{h_1h_1})}-
\\
&-\alpha_{h_2d}\underbrace{\widehat{D}
(\alpha_{h_1h_1})}_{-\alpha_{h_1h_1}}-
\alpha_{h_2r}\underbrace{\widehat{R}
(\alpha_{h_1h_1})}_{-\alpha_{h_2h_1}}
-
\alpha_{h_2i_1}\zero{\widehat{I}_1(\alpha_{h_1h_1})}-
\alpha_{h_2i_2}\zero{\widehat{I}_2(\alpha_{h_1h_1})}-
\\
&
-
\alpha_{h_2j}\zero{\widehat{J}(\alpha_{h_1h_1})}+ 
\alpha_{h_1h_2}\zero{H_2(\alpha_{h_2h_1})}\big)+
\beta_{h_1h_2}\big(\alpha_{h_1h_1}
\zero{H_1(\alpha_{h_2h_2})}-\alpha_{h_2h_1}\zero{H_1(\alpha_{h_1h_2})}+
\\
&
+
\alpha_{h_1h_2}\zero{H_2(\alpha_{h_2h_2})} 
-\alpha_{h_2h_2}\zero{H_2(\alpha_{h_1h_2})}
-\alpha_{h_2d}
\underbrace{\widehat{D}(\alpha_{h_1h_2})}_{-\alpha_{h_1h_2}}
-
\alpha_{h_2r}
\underbrace{\widehat{R}(\alpha_{h_1h_2})}_{-\alpha_{h_2h_2}}-
\alpha_{h_2i_1}\zero{\widehat{I}_1
(\alpha_{h_1h_2})}
-
\\
&-
\alpha_{h_2i_2}\zero{\widehat{I}_2(\alpha_{h_1h_2})}-
\alpha_{h_2j}\zero{\widehat{J}
(\alpha_{h_1h_2})}\big)+
\beta_{h_1t}\big(4\alpha_{h_1h_1}
\alpha_{h_2h_2}-4\alpha_{h_1h_2}\alpha_{h_2h_1}\big).
\endaligned
\end{eqnarray*}
\begin{eqnarray*}
\footnotesize\aligned
\kappa^{h_1h_2}_{h_2}&=\widehat{H}^\ast_1
([\widehat{H}_1,\widehat{H}_2]-4\widehat{T})
\\
&=-\alpha_{h_1d}+\beta_{h_2h_1}
\big(\alpha_{h_1h_1}\zero{H_1(\alpha_{h_2h_1})}
-
\alpha_{h_2h_1}\zero{H_1(\alpha_{h_1h_1})}-
\alpha_{h_2h_2}\zero{H_2(\alpha_{h_1h_1})}-
\\
&-\alpha_{h_2d}\underbrace{\widehat{D}
(\alpha_{h_1h_1})}_{-\alpha_{h_1h_1}}-
\alpha_{h_2r}\underbrace{\widehat{R}
(\alpha_{h_1h_1})}_{-\alpha_{h_2h_1}}
-
\alpha_{h_2i_1}\zero{\widehat{I}_1(\alpha_{h_1h_1})}-
\alpha_{h_2i_2}\zero{\widehat{I}_2(\alpha_{h_1h_1})}-
\alpha_{h_2j}\zero{\widehat{J}(\alpha_{h_1h_1})}+
\\
&+ \alpha_{h_1h_2}\zero{H_2(\alpha_{h_2h_1})}\big)+
\beta_{h_2h_2}\big(\alpha_{h_1h_1}
\zero{H_1(\alpha_{h_2h_2})}-\alpha_{h_2h_1}\zero{H_1(\alpha_{h_1h_2})}
+
\alpha_{h_1h_2}\zero{H_2(\alpha_{h_2h_2})} -
\\
&-\alpha_{h_2h_2}\zero{H_2(\alpha_{h_1h_2})}-\alpha_{h_2d}
\underbrace{\widehat{D}(\alpha_{h_1h_2})}_
{-\alpha_{h_1h_2}}-\alpha_{h_2r}
\underbrace{\widehat{R}(\alpha_{h_1h_2})}_{-\alpha_{h_2h_2}}-
\alpha_{h_2i_1}\zero{\widehat{I}_1(\alpha_{h_1h_2})}
-
\alpha_{h_2i_2}\zero{\widehat{I}_2(\alpha_{h_1h_2})}-
\\
&-\alpha_{h_2j}\zero{\widehat{J}(\alpha_{h_1h_2})}\big)+
\beta_{h_2t}\big(4\alpha_{h_1h_1}\alpha_{h_2h_2}
-
4\alpha_{h_1h_2}\alpha_{h_2h_1}\big).
\endaligned
\end{eqnarray*}
\begin{eqnarray*}
\footnotesize\aligned \kappa^{h_1t}_t
&
=
\widehat{T}^\ast([\widehat{H}_1,\widehat{T}])
\\
&=-2\alpha_{h_1d}+\beta_{tt}
\big(4\alpha_{h_1h_1}\alpha_{th_2}
+
\alpha_{h_1h_1}\alpha_{tt}\Phi_1+
\alpha_{h_1h_1}\zero{H_1(\alpha_{tt})}-4\alpha_{h_1h_2}\alpha_{th_1}+
\\
&+\alpha_{h_1h_2}\alpha_{tt}\Phi_2
+
\alpha_{h_1h_2}\zero{H_2(\alpha_{tt})}\big).
\endaligned
\end{eqnarray*}
\begin{eqnarray*}
\footnotesize\aligned
\kappa^{h_2t}_t&=\widehat{T}^\ast
([\widehat{H}_2,\widehat{T}])\\
&=-2\alpha_{h_2d}+\beta_{tt}\big(4\alpha_{h_2h_1}\alpha_{th_2}+
\alpha_{h_2h_1}\alpha_{tt}
\Phi_1\alpha_{h_2h_1}\zero{H_1(\alpha_{tt})}-
4\alpha_{h_2h_2}\alpha_{th_1}+
\\
&+\alpha_{h_2h_2}\alpha_{tt}\Phi_2+
\alpha_{h_2h_2}\zero{H_2(\alpha_{tt})}\big).
\endaligned
\end{eqnarray*}
Here remind that in Lemma~\ref{fiber-type}, we saw the exact
expressions of the functions $\alpha_{tt}, \alpha_{th_1},
\alpha_{th_2}, \alpha_{h_1h_1}, \alpha_{h_1h_2}, \alpha_{h_2h_1},
\alpha_{h_2h_2}$ and visibly, they were independent of the horizontal
coordinates. Hence the value of the horizontal vector fields $H_1$,
$H_2$ and $T$ on these functions vanish, as is made visible 
by a specific underlining in the
above calculations. Moreover, we have replaced the the values of the
vector fields $\widehat{D}, \widehat{R}, \widehat{I}_1,
\widehat{I}_2$, $\widehat{J}$ on the concerned functions
$\alpha_{{}_\bullet {}_\bullet}$ by just using the 110 equations
stated before Lemma~\ref{fiber-type}. Furthermore, we notice here that
the functions $\beta_{{}_\bullet {}_\bullet}$ are the coefficients of
the dual basis introduced in Subsection~\ref{Dual-section}. After
simplifying carefully these four expressions, we get:
\[
\scriptsize
\aligned
k^{h_1h_2}_{h_1}
&
=
\big(-4\alpha_{th_1}c-4\alpha_{th_2}d
+
(\alpha_{h_1r}+\alpha_{h_2d})c^2+(\alpha_{h_1r}
+
\alpha_{h_2d})d^2\big)/\big(c^2+d^2\big)
\\
k^{h_1h_2}_{h_2}
&
=
-
\big(4\alpha_{th_2}c-4\alpha_{th_1}d
+
(-\alpha_{h_2r}+\alpha_{h_1d})c^2+(-\alpha_{h_2r}
+
\alpha_{h_1d})d^2\big)/\big(c^2+d^2\big)
\\
k^{h_1t}_t
&
=
\big(4\alpha_{th_2}c-2\alpha_{h_1d}c^2
+
\Phi_1c^3-4\alpha_{th_1}d-2\alpha_{h_1d}d^2
+
\Phi_1cd^2+\Phi_2d^3+\Phi_2c^2d\big)/\big(c^2+d^2\big)
\\
k^{h_2t}_t
&
=
-\big(4\alpha_{th_1}c
+
2\alpha_{h_2d}c^2-\Phi_2c^3+4\alpha_{th_2}d
+
2\alpha_{h_2d}d^2+\Phi_1d^3+\Phi_1c^2d-\Phi_2cd^2\big)/\big(c^2+d^2\big).
\endaligned
\]
Now as said, all these curvature coefficients should vanish. Looking
at the above expressions, we realize that one encounters exactly six
undetermined functions $\alpha_{h_1 d}$, $\alpha_{h_2 d}$,
$\alpha_{h_1 r}$, $\alpha_{h_2 r}, \alpha_{t h_1}$ and $\alpha_{t
h_2}$. But if we replace $\alpha_{t h_1} = bd - ac$ and $\alpha_{t
h_2} = -ad - bc$ by the values which were already ascribed in
Subsection~\ref{Graded-Diff}), we are left with exactly four
$\alpha_{{}_\bullet {}_\bullet}$. We therefore see that for these
four curvatures to vanish, it is necessary and sufficient that:
\[
\aligned
\alpha_{h_1d}&=-2b+{\textstyle{\frac{1}{2}}}\Phi_1c
+
{\textstyle{\frac{1}{2}}}\Phi_2d,
\ \
\ \
\alpha_{h_2d}=2a+{\textstyle{\frac{1}{2}}}\Phi_2c
-
{\textstyle{\frac{1}{2}}}\Phi_1d,
\\
\alpha_{h_1r}&=
-6a-{\textstyle{\frac{1}{2}}}\Phi_2c+{\textstyle{\frac{1}{2}}}\Phi_1d, 
\ \ \ \
\alpha_{h_2r}=-6b+{\textstyle{\frac{1}{2}}}\Phi_1c
+
{\textstyle{\frac{1}{2}}}\Phi_2d.
\endaligned
\]
Due to the existence of the functions $\Phi_1$ and $\Phi_2$ in the
last four expressions, one recognizes that the four functions
$\alpha_{ h_1, d}$, $\alpha_{ h_2, d}$, $\alpha_{ h_1, r}$, 
$\alpha_{ h_2, r}$ really depend on the horizontal
coordinates for $i=1,2$.

On the other hand, these four functions $\alpha_{h_1d}$,
$\alpha_{h_2d}$, $\alpha_{h_1r}$, $\alpha_{h_2r}$ should also be of
the form introduced in Lemma~\ref{fiber-type}, 
namely we should have:
\[
\footnotesize
\left[ \aligned &
\boxed{\alpha_{h_1d}}: \ \
-2b+{\textstyle{\frac{1}{2}}}\Phi_1c+{\textstyle{\frac{1}{2}}}\Phi_2d=-{\textstyle{\frac{1}{4}}}c\delta_2-2b+{\textstyle{\frac{1}{4}}}d\delta_1,
\\
 &\boxed{\alpha_{h_2d}}: \ \
2a+{\textstyle{\frac{1}{2}}}\Phi_2c-{\textstyle{\frac{1}{2}}}\Phi_1d=({\textstyle{{\textstyle{\frac{1}{4}}}}}c)\,\delta_1
+ ({\textstyle{{\textstyle{\frac{1}{4}}}}}d)\,\delta_2 + 2a,
\\
 &\boxed{\alpha_{h_1r}}: \ \
-6a-{\textstyle{\frac{1}{2}}}\Phi_2c+{\textstyle{\frac{1}{2}}}\Phi_1d=-({\textstyle{{\textstyle{\frac{1}{2}}}}}c)\,\delta_7
+ ({\textstyle{{\textstyle{\frac{1}{2}}}}}d)\,\delta_{10} - 6a
\\
& \boxed{\alpha_{h_2r}}: \ \
-6b+{\textstyle{\frac{1}{2}}}\Phi_1c+{\textstyle{\frac{1}{2}}}\Phi_2d=({\textstyle{{\textstyle{\frac{1}{2}}}}}d)\,\delta_7
+ ({\textstyle{{\textstyle{\frac{1}{2}}}}}c)\,\delta_{10} - 6b.
\endaligned\right.
\]
By plain identification, these four equations immediately determine
the values of four of the functions $\delta_{{}_\bullet}$ as follows:
\[
\delta_1=2\Phi_2, \ \ \ \ \
\delta_2=-2\Phi_1, \ \ \ \ \ 
\delta_7=\Phi_2, \ \ \ \ \
\delta_{10}=\Phi_1.
\]

\subsection{Homogeneity 2}
\label{Hom2}
In this homogeneity, we encounter the following six curvature
coefficients that should be annihilated in order to satisfy
the regularity condition
\textbf{(c4)}:
\begin{eqnarray*}
\footnotesize\aligned
\kappa^{h_1h_2}_d
&
=
\widehat{D}^\ast([\widehat{H}_1,\widehat{H}_2]-4\widehat{T})
\\
&=2\alpha_{h_1i_1}-\widehat{H}_2(a_{h_1d})
+
\alpha_{h_1h_2}H_2(ah_2d)+\alpha_{h_1h_1}H_1(\alpha_{h_2d})+\beta_{dh_1}\big(
\alpha_{h_1h_1}\zero{H_1(\alpha_{h_2h_1})}-
\\
&-\alpha_{h_2h_1}\zero{H_1(\alpha_{h_1h_1})}-\alpha_{h_2h_2}\zero{H_2(\alpha_{h_1h_1})}-
\alpha_{h_2d}\underbrace{\widehat{D}(\alpha_{h_1h_1})}_{-\alpha_{h_1h_1}}-\alpha_{h_2r}\underbrace{\widehat{R}(\alpha_{h_1h_1})}_{-\alpha_{h_2h_1}}-
\alpha_{h_2i_1}\zero{\widehat{I}_1(\alpha_{h_1h_1})}-
\\
&-
\alpha_{h_2i_2}\zero{\widehat{I}_2(\alpha_{h_1h_1})}-\alpha_{h_2j}\zero{\widehat{J}(\alpha_{h_1h_1})}+\alpha_{h_1h_2}\zero{H_2(\alpha_{h_2h_1})}\big)+
\beta_{dh_2}\big(\alpha_{h_1h_1}\zero{H_1(\alpha_{h_2h_2})}-
\\
&-\alpha_{h_2h_1}\zero{H_1(\alpha_{h_1h_2})}-\alpha_{h_1h_2}\zero{H_2(\alpha_{h_2h_2})}-\alpha_{h_2h_2}\zero{H_2(\alpha_{h_1h_2})}-
\alpha_{h_2d}\underbrace{\widehat{D}(\alpha_{h_1h_2})}_{-\alpha_{h_1h_2}}-
\alpha_{h_2r}\underbrace{\widehat{R}(\alpha_{h_1h_2})}_{-\alpha_{h_2h_2}}-
\\
&-\alpha_{h_2i_1}\zero{\widehat{I}_1(\alpha_{h_1h_2})}-
\alpha_{h_2i_2}\zero{\widehat{I}_2(\alpha_{h_1h_2})}-\alpha_{h_2j}\zero{\widehat{J}(\alpha_{h_1h_2})}\big)+
\beta_{dt}\big(4\alpha_{h_1h_1}\alpha_{h_2h_2}-4\alpha_{h_1h_2}\alpha_{h_2h_1}\big),
\endaligned
\end{eqnarray*}
\begin{eqnarray*}
\footnotesize\aligned
\kappa^{h_1h_2}_r&=\widehat{R}^\ast([\widehat{H}_1,\widehat{H}_2]-4\widehat{T})
\\
&=-6\alpha_{h_1i_2}-\widehat{H}_2(\alpha_{h_1r})+\alpha_{h_1h_2}H_2(\alpha_{h_2r})+\alpha_{h_1h_1}H_1(\alpha_{h_2r})+\beta_{rh_1}\big(
\alpha_{h_1h_1}\zero{H_1(\alpha_{h_2h_1})}-
\\
&-\alpha_{h_2h_1}\zero{H_1(\alpha_{h_1h_1})}-\alpha_{h_2h_2}\zero{H_2(\alpha_{h_1h_1})}-
\alpha_{h_2d}\underbrace{\widehat{D}(\alpha_{h_1h_1})}_{-\alpha_{h_1h_1}}-\alpha_{h_2r}\underbrace{\widehat{R}(\alpha_{h_1h_1})}_{-\alpha_{h_2h_1}}-
\alpha_{h_2i_1}\zero{\widehat{I}_1(\alpha_{h_1h_1})}-
\\
&-\alpha_{h_2i_2}\zero{\widehat{I}_2(\alpha_{h_1h_1})}-
\alpha_{h_2j}\zero{\widehat{J}(\alpha_{h_1h_1})}+\alpha_{h_1h_2}\zero{H_2(\alpha_{h_2h_1})}\big)+
\beta_{rh_2}\big(\alpha_{h_1h_1}\zero{H_1(\alpha_{h_2h_2})}-
\\
&-\alpha_{h_2h_1}\zero{H_1(\alpha_{h_1h_2})}-\alpha_{h_1h_2}\zero{H_2(\alpha_{h_2h_2})}-\alpha_{h_2h_2}\zero{H_2(\alpha_{h_1h_2})}-
\alpha_{h_2d}\underbrace{\widehat{D}(\alpha_{h_1h_2})}_{-\alpha_{h_1h_2}}-
\alpha_{h_2r}\underbrace{\widehat{R}(\alpha_{h_1h_2})}_{-\alpha_{h_2h_2}}-
\\
&-\alpha_{h_2i_1}\zero{\widehat{I}_1(\alpha_{h_1h_2})}-
\alpha_{h_2i_2}\zero{\widehat{I}_2(\alpha_{h_1h_2})}-\alpha_{h_2j}\zero{\widehat{J}(\alpha_{h_1h_2})}\big)+
\beta_{rt}\big(4\alpha_{h_1h_1}\alpha_{h_2h_2}-4\alpha_{h_1h_2}\alpha_{h_2h_1}\big),
\endaligned
\end{eqnarray*}
\begin{eqnarray*}
\footnotesize\aligned
\kappa^{h_1t}_{h_1}&=\widehat{H}^\ast_1([\widehat{H}_1,\widehat{T}])
\\
&=-\alpha_{h_1i_1}+\beta_{h_1h_1}\big(\alpha_{h_1h_1}\zero{H_1(\alpha_{th_1})}-\alpha_{th_1}\zero{H_1(\alpha_{h_1h_1})}-
\alpha_{th_2}\zero{H_2(\alpha_{h_1h_1})}-\alpha_{tt}\zero{T(\alpha_{h_1h_1})}-
\\
&-\alpha_{td}\underbrace{\widehat{D}(\alpha_{h_1h_1})}_{-\alpha_{h_1h_1}}-
\alpha_{ti_1}\zero{\widehat{I}_1(\alpha_{h_1h_1})}-\alpha_{ti_2}\zero{\widehat{I}_2(\alpha_{h_1h_1})}-
\alpha_{tr}\underbrace{\widehat{R}(\alpha_{h_1h_1})}_{-\alpha_{h_2h_1}}-
\alpha_{tj}\zero{\widehat{J}(\alpha_{h_1h_1})}\big)+
\\
&+
\beta_{h_1h_2}\big(\alpha_{h_1h_1}\zero{H_1(\alpha_{th_2})}-\alpha_{th_1}\zero{H_1(\alpha_{h_1h_2})}
+\alpha_{h_1h_2}\zero{H_2(\alpha_{th_2})}-\alpha_{th_2}\zero{H_2(\alpha_{h_1h_2})}-
\\
&-\alpha_{tt}\zero{T(\alpha_{h_1h_2})}-\alpha_{td}\underbrace{\widehat{D}(\alpha_{h_1h_2})}_{-\alpha_{h_1h_2}}-
\alpha_{tr}\underbrace{\widehat{R}(\alpha_{h_1h_2})}_{-\alpha_{h_2h_2}}-\alpha_{ti_1}\zero{\widehat{I}_1(\alpha_{h_1h_2})}-
\alpha_{ti_2}\zero{\widehat{I}_2(\alpha_{h_1h_2})}-
\\
&-\alpha_{tj}\zero{\widehat{J}(\alpha_{h_1h_2})}\big)+
\beta_{h_1t}\big(4\alpha_{h_1h_1}\alpha_{th_2}+\alpha_{h_1h_1}\alpha_{tt}\Phi_1-4\alpha_{h_1h_2}\alpha_{th_1}+\alpha_{h_1h_2}\alpha_{tt}\Phi_2\big),
\endaligned
\end{eqnarray*}
\begin{eqnarray*}
\footnotesize\aligned
\kappa^{h_1t}_{h_2}&=\widehat{H}^\ast_2([\widehat{H}_1,\widehat{T}])
\\
&=-\alpha_{h_1i_2}+\beta_{h_2h_1}\big(\alpha_{h_1h_1}\zero{H_1(\alpha_{th_1})}-\alpha_{th_1}\zero{H_1(\alpha_{h_1h_1})}-
\alpha_{th_2}\zero{H_2(\alpha_{h_1h_1})}-\alpha_{tt}\zero{T(\alpha_{h_1h_1})}-
\\
&-\alpha_{td}\underbrace{\widehat{D}(\alpha_{h_1h_1})}_{-\alpha_{h_1h_1}}-
\alpha_{ti_1}\zero{\widehat{I}_1(\alpha_{h_1h_1})}-\alpha_{ti_2}\zero{\widehat{I}_2(\alpha_{h_1h_1})}-
\alpha_{tr}\underbrace{\widehat{R}(\alpha_{h_1h_1})}_{-\alpha_{h_2h_1}}-
\alpha_{tj}\zero{\widehat{J}(\alpha_{h_1h_1})}\big)+
\\
&+
\beta_{h_2h_2}\big(\alpha_{h_1h_1}\zero{H_1(\alpha_{th_2})}-\alpha_{th_1}\zero{H_1(\alpha_{h_1h_2})}
+\alpha_{h_1h_2}\zero{H_2(\alpha_{th_2})}-\alpha_{th_2}\zero{H_2(\alpha_{h_1h_2})}-
\\
&-\alpha_{tt}\zero{T(\alpha_{h_1h_2})}-\alpha_{td}\underbrace{\widehat{D}(\alpha_{h_1h_2})}_{-\alpha_{h_1h_2}}-
\alpha_{tr}\underbrace{\widehat{R}(\alpha_{h_1h_2})}_{-\alpha_{h_2h_2}}-\alpha_{ti_1}\zero{\widehat{I}_1(\alpha_{h_1h_2})}-
\alpha_{ti_2}\zero{\widehat{I}_2(\alpha_{h_1h_2})}-
\\
&-\alpha_{tj}\zero{\widehat{J}(\alpha_{h_1h_2})}\big)+
\beta_{h_2t}\big(4\alpha_{h_1h_1}\alpha_{th_2}+\alpha_{h_1h_1}\alpha_{tt}\Phi_1-4\alpha_{h_1h_2}\alpha_{th_1}+\alpha_{h_1h_2}\alpha_{tt}\Phi_2\big),
\endaligned
\end{eqnarray*}
\begin{eqnarray*}
\footnotesize\aligned
\kappa^{h_2t}_{h_1}&=\widehat{H}^\ast_1([\widehat{H}_2,\widehat{T}])
\\
&=-\alpha_{h_2i_1}+\beta_{h_1h_1}\big(\alpha_{h_2h_1}\zero{H_1(\alpha_{th_1})}-\alpha_{th_1}\zero{H_1(\alpha_{h_2h_1})}-
\alpha_{th_2}\zero{H_2(\alpha_{h_2h_1})}-\alpha_{tt}\zero{T(\alpha_{h_2h_1})}-
\\
&-\alpha_{td}\underbrace{\widehat{D}(\alpha_{h_2h_1})}_{-\alpha_{h_2h_1}}-
\alpha_{tr}\underbrace{\widehat{R}(\alpha_{h_2h_1})}_{\alpha_{h_1h_1}}-\alpha_{ti_1}\zero{\widehat{I}_1(\alpha_{h_2h_1})}-
\alpha_{ti_2}\zero{\widehat{I}_2(\alpha_{h_2h_1})}-
\alpha_{tj}\zero{\widehat{J}(\alpha_{h_2h_1})}+
\\
&+\alpha_{h_2h_1}\zero{H_2(\alpha_{th_1})}\big)+
\beta_{h_1h_2}\big(\alpha_{h_2h_1}\zero{H_1(\alpha_{th_2})}-\alpha_{th_1}\zero{H_1(\alpha_{h_2h_2})}
+\alpha_{h_2h_2}\zero{H_2(\alpha_{th_2})}-
\\
&-\alpha_{th_2}\zero{H_2(\alpha_{h_2h_2})}-
\alpha_{tt}\zero{T(\alpha_{h_2h_2})}-\alpha_{td}\underbrace{\widehat{D}(\alpha_{h_2h_2})}_{-\alpha_{h_2h_2}}-
\alpha_{tr}\underbrace{\widehat{R}(\alpha_{h_2h_2})}_{\alpha_{h_1h_2}}-\alpha_{ti_1}\zero{\widehat{I}_1(\alpha_{h_2h_2})}-
\\
&-\alpha_{ti_2}\zero{\widehat{I}_2(\alpha_{h_2h_2})}-\alpha_{tj}\zero{\widehat{J}(\alpha_{h_2h_2})}\big)+
\beta_{h_1t}\big(4\alpha_{h_2h_1}\alpha_{th_2}+\alpha_{h_2h_1}\alpha_{tt}\Phi_1-4\alpha_{h_2h_2}\alpha_{th_1}+\alpha_{h_2h_2}\alpha_{tt}\Phi_2\big),
\endaligned
\end{eqnarray*}
\begin{eqnarray*}
\footnotesize\aligned
\kappa^{h_2t}_{h_2}&=\widehat{H}^\ast_2([\widehat{H}_2,\widehat{T}])\\
&=-\alpha_{h_2i_2}+\beta_{h_2h_1}\big(\alpha_{h_2h_1}\zero{H_1(\alpha_{th_1})}-\alpha_{th_1}\zero{H_1(\alpha_{h_2h_1})}-
\alpha_{th_2}\zero{H_2(\alpha_{h_2h_1})}-\alpha_{tt}\zero{T(\alpha_{h_2h_1})}-
\\
&-\alpha_{td}\underbrace{\widehat{D}(\alpha_{h_2h_1})}_{-\alpha_{h_2h_1}}-
\alpha_{tr}\underbrace{\widehat{R}(\alpha_{h_2h_1})}_{\alpha_{h_1h_1}}-\alpha_{ti_1}\zero{\widehat{I}_1(\alpha_{h_2h_1})}-
\alpha_{ti_2}\zero{\widehat{I}_2(\alpha_{h_2h_1})}-
\alpha_{tj}\zero{\widehat{J}(\alpha_{h_2h_1})}+
\\
&+\alpha_{h_2h_1}\zero{H_2(\alpha_{th_1})}\big)+
\beta_{h_2h_2}\big(\alpha_{h_2h_1}\zero{H_1(\alpha_{th_2})}-\alpha_{th_1}\zero{H_1(\alpha_{h_2h_2})}
+\alpha_{h_2h_2}\zero{H_2(\alpha_{th_2})}-
\\
&-\alpha_{th_2}\zero{H_2(\alpha_{h_2h_2})}-
\alpha_{tt}\zero{T(\alpha_{h_2h_2})}-\alpha_{td}\underbrace{\widehat{D}(\alpha_{h_2h_2})}_{-\alpha_{h_2h_2}}-
\alpha_{tr}\underbrace{\widehat{R}(\alpha_{h_2h_2})}_{\alpha_{h_1h_2}}-\alpha_{ti_1}\zero{\widehat{I}_1(\alpha_{h_2h_2})}-
\\
&-\alpha_{ti_2}\zero{\widehat{I}_2(\alpha_{h_2h_2})}-\alpha_{tj}\zero{\widehat{J}(\alpha_{h_2h_2})}\big)+
\beta_{h_2t}\big(4\alpha_{h_2h_1}\alpha_{th_2}+\alpha_{h_2h_1}\alpha_{tt}\Phi_1-4\alpha_{h_2h_2}\alpha_{th_1}+\alpha_{h_2h_2}\alpha_{tt}\Phi_2\big),
\endaligned
\end{eqnarray*}
Before simplifying the above expressions, we notice that the
expressions of the four functions $\alpha_{h_1d}$, $\alpha_{h_2d}$,
$\alpha_{h_1r}$, $\alpha_{h_2r}$ obtained in the previous subsection
enable us to identify how the three horizontal vector fields $H_1$,
$H_2$, $T$ act on them as first-order differential operators. Because
these vector fields only differentiate with respect to $(x, y, u)$ and
not with respect to $(a, b, c, d, e)$, we have, for any
for $Y=H_1,H_2,T$:
\[
\footnotesize\aligned
 Y(a_{h_1d})&={\textstyle{\frac{1}{2}}}Y(\Phi_1)c+{\textstyle{\frac{1}{2}}}Y(\Phi_2)d-2b,
 \ \ \ \ \ \ \
 Y(a_{h_2d}) =
{\textstyle{\frac{1}{2}}}(\Phi_2)c-{\textstyle{\frac{1}{2}}}Y(\Phi_1)d+2a,
\\
Y(a_{h_1r})&=
-{\textstyle{\frac{1}{2}}}Y(\Phi_2)c+{\textstyle{\frac{1}{2}}}Y(\Phi_1)d-6a, \ \ \
Y(a_{h_2r})={\textstyle{\frac{1}{2}}}Y(\Phi_1)c+{\textstyle{\frac{1}{2}}}Y(\Phi_2)d-6b.
\endaligned
\]
Moreover, we can find the expressions of
$\widehat{H}_2(\alpha_{h_1d})$ and $\widehat{H}_2(\alpha_{h_1r})$
which appear above in the expressions of the two curvature
coefficients $\kappa^{h_1h_2}_d$, $\kappa^{h_1h_2}_r$, coming back to
the definition of the vector field $\widehat{H}_2$ in
Subsection~\ref{Dual-section}:
\[
\footnotesize\aligned
\widehat{H}_2(\alpha_{h_1d})&=\alpha_{h_2h_1}H_1(\alpha_{h_1d})+\alpha_{h_2h_2}H_2(\alpha_{h_1d})+
\alpha_{h_2d}\underbrace{\widehat{D}(\alpha_{h_1d})}_{-\alpha_{h_1d}}+
\alpha_{h_2r}\underbrace{\widehat{R}(\alpha_{h_1d})}_{-\alpha_{h_2d}}+
\alpha_{h_2i_2}\underbrace{\widehat{I}_2(\alpha_{h_1d})}_{-2}+
\\
&+\alpha_{h_2i_1}\zero{\widehat{I}_1(\alpha_{h_1d})}+\alpha_{h_2j}\zero{\widehat{J}(\alpha_{h_1d})}
\\
&-2\alpha_{h_2i_2}+({\textstyle{\frac{1}{2}}}H_2(\Phi_1)-{\textstyle{\frac{1}{2}}}\Phi_1\Phi_2)c^2+
(-{\textstyle{\frac{1}{2}}}H_1(\Phi_2)+{\textstyle{\frac{1}{2}}}\Phi_1\Phi_2)d^2+16ab-
2\Phi_2ad+
\\
&+4\Phi_2bc-2\Phi_1ac-4\Phi_1bd+({\textstyle{\frac{1}{2}}}\Phi_1^2+{\textstyle{\frac{1}{2}}}H_2(\Phi_2)-
{\textstyle{\frac{1}{2}}}\Phi_2^2-{\textstyle{\frac{1}{2}}}H_1(\Phi_1))cd
\endaligned
\]
\[
\footnotesize\aligned
\widehat{H}_2(\alpha_{h_1r})&=\alpha_{h_2h_1}H_1(\alpha_{h_1r})+\alpha_{h_2h_2}H_2(\alpha_{h_1r})+
\alpha_{h_2d}\underbrace{\widehat{D}(\alpha_{h_1r})}_{-\alpha_{h_1r}}+
\alpha_{h_2r}\underbrace{\widehat{R}(\alpha_{h_1r})}_{-\alpha_{h_2r}}+\alpha_{h_2i_1}\underbrace{\widehat{I}_1(\alpha_{h_1r})}_{-6}+
\\
&+\alpha_{h_2i_2}\zero{\widehat{I}_2(\alpha_{h_1r})}+\alpha_{h_2j}\zero{\widehat{J}(\alpha_{h_1r})}
\\
&=- 6\alpha_{h_2i_1}+(-{\textstyle{\frac{1}{2}}}H_2(\Phi_2)-
{\textstyle{\frac{1}{4}}}\Phi_1^2+{\textstyle{\frac{1}{4}}}\Phi_2^2)c^2+
({\textstyle{\frac{1}{4}}}\Phi_1^2-{\textstyle{\frac{1}{4}}}\Phi_1^2-{\textstyle{\frac{1}{2}}}H_1(\Phi_1))d^2+(-\Phi_1\Phi_2+
\\
&+{\textstyle{\frac{1}{2}}}H_1(\Phi_2)+{\textstyle{\frac{1}{2}}}H_2(\Phi_1))cd+
12a^2-36b^2+4\Phi_2ac-4\Phi_1ad+6\Phi_1bc+6\Phi_2bd.
\endaligned
\]
Thanks to these preparations, we are 
now in a position to simplify the 
six curvature coefficients of homogeneity $2$, and
careful calculations give at the end: 
\[
\footnotesize\aligned
\kappa^{h_1h_2}_d&=2\alpha_{h_2i_2}-4\alpha_{td}+2\alpha_{h_1i_1}+{\textstyle{\frac{1}{2}}}\big(-H_2(\Phi_1)+
H_1(\Phi_2)\big)c^2+ {\textstyle{\frac{1}{2}}}\big(-H_2(\Phi_1)+H_1(\Phi_2)\big)d^2+
\\
&+2\big(\Phi_1bd-\Phi_2bc-\Phi_1ac-\Phi_2ad\big),
\\
\kappa^{h_1h_2}_r&=-6\alpha_{h_1i_2}+6\alpha_{h_2i_1}-4\alpha_{tr}+
{\textstyle{\frac{1}{2}}}\big(H_2(\Phi_2)+H_1(\Phi_1)\big)c^2+
{\textstyle{\frac{1}{2}}}\big(H_2(\Phi_2)+(H_1)\Phi_1\big)d^2+
\\
&+24a^2+24b^2+2\Phi_2ac-2\Phi_1bc-2\Phi_1ad-2\Phi_2bd,
\\
\kappa^{h_1t}_{h_1}&=\alpha_{td}-\alpha_{h_1i_1}+\Phi_1ac+\Phi_2ad-4ab, 
\\
\kappa^{h_1t}_{h_2}&=\alpha_{tr}-\alpha_{h_1i_2}+\Phi_1bc+\Phi_2bd-4b^2,
\\
\kappa^{h_2t}_{h_1}&=-\alpha_{tr}-\alpha_{h_2i_1}+\Phi_2ac-\Phi_1ad+4a^2,
\\
\kappa^{h_2t}_{h_2}&=\alpha_{td}-\alpha_{h_2i_2}+\Phi_2bc-\Phi_1bd+4ab.
\\
\endaligned
\]
Inspecting these six equations, we see that there are exactly six
undetermined functions $\alpha_{td}$, $\alpha_{tr}$,
$\alpha_{h_1i_1}$, $\alpha_{h_1i_2}$ ,$\alpha_{h_2i_1}$,
$\alpha_{h_2i_2}$. Thus, we might to fully annihilate the curvature
of homogeneity $2$ by solving this system of six equations in six
unknowns. But unfortunately, the solution set of this system, as it is
written here, happens to be empty!

Temporarily, let us omit the first equation involving the curvature
coefficient $\kappa^{h_1h_2}_d$ and let us solve the remaining system
of five equations with six unknowns. The obtained solution set of this
system reveals that vanishing of this new system is independent of the
function $\alpha_{h_2i_2}$. Consequently, $\alpha_{h_2i_2}$ is free
and we then determine the remaining five functions $\alpha_{td}$,
$\alpha_{tr}$, $\alpha_{h_1i_1}$, $\alpha_{h_1i_2}$, $\alpha_{h_2i_1}$
according to this subsystem of five equations.

On the other hand, reminding Lemma~\ref{fiber-type}, 
$\alpha_{ h_2 i_2}$ should be of the form:
\[
\footnotesize\aligned
\alpha_{h_2i_2}
=
(-{\textstyle{\frac{1}{8}}}c^2
-
{\textstyle{\frac{1}{8}}}d^2)\delta_6
+
{\textstyle{\frac{1}{4}}}\delta_8c^2
-
{\textstyle{\frac{1}{4}}}\delta_9dc
-
{\textstyle{\frac{1}{2}}}\Phi_1ac
+
4ab-{\textstyle{\frac{1}{2}}}\Phi_2ad+
{\textstyle{\frac{1}{2}}}\Phi_2bc
-
{\textstyle{\frac{1}{2}}}\Phi_1bd-2e.
\endaligned
\] 
and furthermore, the obtained values of the five functions in question
should identify to:
\[\footnotesize\aligned
\alpha_{td}&=\alpha_{h_2i_2
}
-\Phi_2bc+\Phi_1bd-4ab
\\
&
=
-{\textstyle{\frac{1}{8}}}(c^2+d^2)\delta_6
+
{\textstyle{\frac{1}{4}}}\delta_8c^2-
{\textstyle{\frac{1}{4}}}\delta_9cd
+
{\textstyle{\frac{1}{2}}}\Phi_1bd
-
{\textstyle{\frac{1}{2}}}\Phi_2bc
-
{\textstyle{\frac{1}{2}}}\Phi_1ac-{\textstyle{\frac{1}{2}}}\Phi_2ad-2e,
\\
\alpha_{tr}&=3(a^2+b^2)+{\textstyle{\frac{1}{32}}}\big(H_1(\Phi_1)+H_2(\Phi_2)\big)c^2+
{\textstyle{\frac{1}{32}}}\big(H_1(\Phi_1)+H_2(\Phi_2)\big)d^2+
\\
&+ {\textstyle{\frac{1}{2}}}\big(\Phi_2ac-\Phi_1bc-\Phi_2bd-\Phi_1ad\big),
\\
\alpha_{h_1i_1}&=\alpha_{h_2i_2}+\Phi_1ac-\Phi_2bc+\Phi_1bd-8ab+\Phi_2ad
\\
&
=
(-{\textstyle{\frac{1}{8}}}c^2
-
{\textstyle{\frac{1}{8}}}d^2)\delta_6
+
{\textstyle{\frac{1}{4}}}\delta_8c^2
-
{\textstyle{\frac{1}{4}}}\delta_9cd
+
{\textstyle{\frac{1}{2}}}\Phi_1bd
-
{\textstyle{\frac{1}{2}}}\Phi_2bc
+
{\textstyle{\frac{1}{2}}}\Phi_1ac
+
{\textstyle{\frac{1}{2}}}\Phi_2ad-4ab-2e,
\\
\alpha_{h_1i_2}&=3a^2-b^2
+
{\textstyle{\frac{1}{32}}}\big(H_1(\Phi_1)
+
H_2(\Phi_2)\big)c^2+
{\textstyle{\frac{1}{32}}}\big(H_1(\Phi_1)
+
H_2(\Phi_2)\big)d^2+
\\
&+{\textstyle{\frac{1}{2}}}\big(\Phi_2ac+\Phi_2bd+ \Phi_1bc-\Phi_1ad\big),
\\
\alpha_{h_2i_1}
&
=
a^2-3b^2+\big(-{\textstyle{\frac{1}{32}}}H_1(\Phi_1)
-
{\textstyle{\frac{1}{32}}}H_2(\Phi_2)\big)c^2+
\big(-{\textstyle{\frac{1}{32}}}H_1(\Phi_1)
-
{\textstyle{\frac{1}{32}}}H_2(\Phi_2)\big)d^2
+
\\
&
+
{\textstyle{\frac{1}{2}}}\Phi_2bd+{\textstyle{\frac{1}{2}}}\Phi_1bc
+
{\textstyle{\frac{1}{2}}}\Phi_2ac-{\textstyle{\frac{1}{2}}}\Phi_1ad,
\endaligned
\]
after the necessary simplifications. These five equations guarantee
the vanishing of the five curvature coefficients $\kappa^{ h_1
h_2}_r,\kappa^{ h_1t}_{ h_1}, \kappa^{ h_1t}_{ h_2}$, $\kappa^{ h_2t
}_{ h_1}$, but still, it remains to also annihilate the first, left
aside, curvature coefficient $\kappa^{ h_1 h_2}_d$. But replacing the
above expressions in the expression of this curvature coefficient, we
get:
\[
\footnotesize\aligned
\kappa^{h_1h_2}_d
=
\big(-{\textstyle{\frac{1}{2}}}H_2(\Phi_1)
+
{\textstyle{\frac{1}{2}}}H_1(\Phi_2)\big)c^2
+
\big(-{\textstyle{\frac{1}{2}}}H_2(\Phi_1)
+
{\textstyle{\frac{1}{2}}}H_1(\Phi_2)\big)d^2.
\endaligned
\]
Fortunately, this last expression vanishes thanks to the fact that the
two functions $H_1(\Phi_2) = H_2(\Phi_1)$ are equal, which was already
seen in Lemma~\ref{H1-Phi2-H2-Phi1}.

Now, we are sure that the above determination of the functions
$\alpha_{td}$, $\alpha_{tr}$, $\alpha_{h_1i_1}$, $\alpha_{h_1i_2}$,
$\alpha_{h_2i_1}$, $\alpha_{h_2i_2}$ annihilates
all the curvature components
of homogeneity two, which is what was announced in the regularity
condition \textbf{(c4)}. Lastly, we also have to take care of the
condition \textbf{(c1)}. Similarly as in homogeneity one,
reminding Lemma~\ref{fiber-type}, an identification gives:
\[
\footnotesize
 \left[ \aligned
&\boxed{\alpha_{td}}: \ \
(-{\textstyle{\frac{1}{8}}}c^2-{\textstyle{\frac{1}{8}}}d^2)\delta_6+{\textstyle{\frac{1}{4}}}\delta_8c^2-
{\textstyle{\frac{1}{4}}}\delta_9cd+{\textstyle{\frac{1}{2}}}\Phi_1bd-{\textstyle{\frac{1}{2}}}\Phi_2bc-
{\textstyle{\frac{1}{2}}}\Phi_1ac-{\textstyle{\frac{1}{2}}}\Phi_2ad-2e
\\
& \ \ \ \ \ \ \ \ \ \
=2(-{\textstyle{{\textstyle{\frac{1}{4}}}}}bc-{\textstyle{{\textstyle{\frac{1}{4}}}}}ad)\,\Phi_2
-2
({\textstyle{{\textstyle{\frac{1}{4}}}}}ac-{\textstyle{{\textstyle{\frac{1}{4}}}}}bd)\,\Phi_1
+
({\textstyle{{\textstyle{\frac{1}{4}}}}}c^2+{\textstyle{{\textstyle{\frac{1}{4}}}}}d^2)\,\delta_{15}
- 2e,
\\
 &\boxed{\alpha_{tr}}: \ \
3a^2+3b^2+\big({\textstyle{\frac{1}{32}}}H_1(\Phi_1)+{\textstyle{\frac{1}{32}}}H_2(\Phi_2)\big)c^2+
 \big({\textstyle{\frac{1}{32}}}H_1(\Phi_1)+{\textstyle{\frac{1}{32}}}H_2(\Phi_2)\big)d^2+
{\textstyle{\frac{1}{2}}}\Phi_2ac-{\textstyle{\frac{1}{2}}}\Phi_1bc-
\\
& \ \ \ \ \ \ \ \ \ \
-{\textstyle{\frac{1}{2}}}\Phi_2bd-{\textstyle{\frac{1}{2}}}\Phi_1ad
\\
& \ \ \ \ \ \ \ \ \ \ 
=
({\textstyle{{\textstyle{\frac{1}{4}}}}}c^2+{\textstyle{{\textstyle{\frac{1}{4}}}}}d^2)\,\delta_4
+
({\textstyle{{\textstyle{\frac{1}{2}}}}}ac-{\textstyle{{\textstyle{\frac{1}{2}}}}}bd)\,\Phi_2
+
({\textstyle{{\textstyle{\frac{1}{2}}}}}c^2+{\textstyle{{\textstyle{\frac{1}{2}}}}}d^2)\,\delta_9
-
({\textstyle{{\textstyle{\frac{1}{2}}}}}ad+{\textstyle{{\textstyle{\frac{1}{2}}}}}bc)\,\Phi_1
+
({\textstyle{{\textstyle{\frac{1}{2}}}}}c^2+{\textstyle{{\textstyle{\frac{1}{2}}}}}d^2)\,\delta_{19}
+ 3b^2+3a^2,
\\
&\boxed{\alpha_{h_1i_1}}: \ \
(-{\textstyle{\frac{1}{8}}}c^2-{\textstyle{\frac{1}{8}}}d^2)\delta_6+
{\textstyle{\frac{1}{4}}}\delta_8c^2-{\textstyle{\frac{1}{4}}}\delta_9cd+{\textstyle{\frac{1}{2}}}\Phi_1bd-{\textstyle{\frac{1}{2}}}\Phi_2bc+
{\textstyle{\frac{1}{2}}}\Phi_1ac+{\textstyle{\frac{1}{2}}}\Phi_2ad-4ab-2e
\\
& \ \ \ \ \ \ \ \ \ \ \ \ \ \ \ \
=2({\textstyle{{\textstyle{\frac{1}{4}}}}}ad)\,\Phi_2 +
2({\textstyle{{\textstyle{\frac{1}{4}}}}}ac)\,\Phi_1 -
({\textstyle{{\textstyle{\frac{1}{8}}}}}c^2+{\textstyle{{\textstyle{\frac{1}{8}}}}}d^2)\,\delta_6
-
({\textstyle{{\textstyle{\frac{1}{2}}}}}bc)\,\Phi_2 +
({\textstyle{{\textstyle{\frac{1}{4}}}}}d^2)\,\delta_8 +
({\textstyle{{\textstyle{\frac{1}{4}}}}}cd)\,\delta_9 +
({\textstyle{{\textstyle{\frac{1}{2}}}}}bd)\,\Phi_1 - 4ab-2e,
\\
& \boxed{\alpha_{h_1i_2}}: \ \
3a^2-b^2+\big({\textstyle{\frac{1}{32}}}H_1(\Phi_1)+{\textstyle{\frac{1}{32}}}H_2(\Phi_2)\big)c^2+\big({\textstyle{\frac{1}{32}}}H_1(\Phi_1)+
{\textstyle{\frac{1}{32}}}H_2(\Phi_2)\big)d^2+{\textstyle{\frac{1}{2}}}\Phi_2ac+{\textstyle{\frac{1}{2}}}\Phi_2bd+
\\
& \ \ \ \ \ \ \ \ \ \ \ \ \ \ \ \ \
{\textstyle{\frac{1}{2}}}\Phi_1bc-{\textstyle{\frac{1}{2}}}\Phi_1ad
\\
&\ \ \ \ \ \ \ \ \ \ =2({\textstyle{{\textstyle{\frac{1}{4}}}}}bd)\Phi_2 +
2({\textstyle{{\textstyle{\frac{1}{4}}}}}bc)\,\Phi_1 -
({\textstyle{{\textstyle{\frac{1}{8}}}}}c^2+{\textstyle{{\textstyle{\frac{1}{8}}}}}d^2)\,\delta_4
+
({\textstyle{{\textstyle{\frac{1}{2}}}}}ac)\,\Phi_2 +
({\textstyle{{\textstyle{\frac{1}{4}}}}}cd)\,\delta_8 -
({\textstyle{{\textstyle{\frac{1}{4}}}}}d^2)\,\delta_9 -
({\textstyle{{\textstyle{\frac{1}{2}}}}}ad)\,\Phi_1 + 3a^2-b^2,
\\
& \boxed{\alpha_{h_2i_1}}: \ \
a^2-3b^2+\big(-{\textstyle{\frac{1}{32}}}H_1(\Phi_1)-{\textstyle{\frac{1}{32}}}H_2(\Phi_2)\big)c^2+\big(-{\textstyle{\frac{1}{32}}}H_1(\Phi_1)-
{\textstyle{\frac{1}{32}}}H_2(\Phi_2)\big)d^2+{\textstyle{\frac{1}{2}}}\Phi_2bd+{\textstyle{\frac{1}{2}}}\Phi_1bc+
\\
& \ \ \ \ \ \ \ \ \ \ \ \ \ \ \ \
{\textstyle{\frac{1}{2}}}\Phi_2ac-{\textstyle{\frac{1}{2}}}\Phi_1ad
\\
& \ \ \ \ \ \ \ \ \ \ = 2({\textstyle{{\textstyle{\frac{1}{4}}}}}ac)\,\Phi_2 -2
({\textstyle{{\textstyle{\frac{1}{4}}}}}ad)\,\Phi_1 +
({\textstyle{{\textstyle{\frac{1}{8}}}}}c^2+{\textstyle{{\textstyle{\frac{1}{8}}}}}d^2)\,\delta_4
+
({\textstyle{{\textstyle{\frac{1}{2}}}}}bd)\,\Phi_2 +
({\textstyle{{\textstyle{\frac{1}{4}}}}}cd)\,\delta_8 +
({\textstyle{{\textstyle{\frac{1}{4}}}}}c^2)\,\delta_9 +
({\textstyle{{\textstyle{\frac{1}{2}}}}}bc)\,\Phi_1 - 3b^2+a^2.
\endaligned\right.
\]
The right hand sides of these equations are the expressions of the
mentioned functions as in Lemma~\ref{fiber-type} after the possible
simplifications. Inspection of this system shows that it will be
satisfied whenever one has:
\[
\footnotesize
\left[ \aligned
&\boxed{\alpha_{td}}: \ \ \
-{\textstyle{\frac{1}{8}}}\delta_6+{\textstyle{\frac{1}{4}}}\delta_8={\textstyle{\frac{1}{4}}}\delta_{15},
\ \ \ -{\textstyle{\frac{1}{8}}}\delta_6={\textstyle{\frac{1}{4}}}\delta_{15}, \ \ \
{\textstyle{\frac{1}{4}}}\delta_9=0,
\\
&\boxed{\alpha_{tr}}: \ \ \
{\textstyle{\frac{1}{32}}}H_1(\Phi_1)+{\textstyle{\frac{1}{32}}}H_2(\Phi_2)={\textstyle{\frac{1}{4}}}\delta_4+{\textstyle{\frac{1}{2}}}\delta_9+{\textstyle{\frac{1}{2}}}\delta_{19},
\\
&\boxed{\alpha_{h_1i_1}}: \ \
-{\textstyle{\frac{1}{8}}}\delta_6+{\textstyle{\frac{1}{4}}}\delta_8=-{\textstyle{\frac{1}{8}}}\delta_6,
\ \ \
-{\textstyle{\frac{1}{8}}}\delta_6={\textstyle{\frac{1}{8}}}\delta_6+{\textstyle{\frac{1}{4}}}\delta_8,
\ \ \ -{\textstyle{\frac{1}{4}}}\delta_9={\textstyle{\frac{1}{4}}}\delta_9,
\\
&\boxed{\alpha_{h_1i_2}}: \ \
{\textstyle{\frac{1}{32}}}H_1(\Phi_1)+{\textstyle{\frac{1}{32}}}H_2(\Phi_2)=-{\textstyle{\frac{1}{8}}}\delta_4,
\\
&\boxed{\alpha_{h_2i_1}}: \ \
-{\textstyle{\frac{1}{32}}}H_1(\Phi_1)-{\textstyle{\frac{1}{32}}}H_2(\Phi_2)={\textstyle{\frac{1}{8}}}\delta_4+{\textstyle{\frac{1}{4}}}\delta_9,
\ \ \
-{\textstyle{\frac{1}{32}}}H_1(\Phi_1)-{\textstyle{\frac{1}{32}}}H_2(\Phi_2)={\textstyle{\frac{1}{8}}}\delta_4,
\ \ \ {\textstyle{\frac{1}{4}}}\delta_8=0.
\endaligned\right.
\]
One immediately checks that this system has the following solution set
which guarantees that our computations in this homogeneity are in
direction of satisfying both of the conditions \textbf{(c1)} and
\textbf{(c4)}:
\[
\delta_6=\delta_8=\delta_9=\delta_{15}=0, \ \ \
\delta_4=-{\textstyle{\frac{1}{4}}}H_1(\Phi_1)-{\textstyle{\frac{1}{4}}}H_2(\Phi_2),
\ \ \
\delta_{19}=\frac{3}{16}H_1(\Phi_1)+\frac{3}{16}H_2(\Phi_2).
\]
 Finally we obtain:
\[
\footnotesize\aligned
\alpha_{td}&={\textstyle{\frac{1}{2}}}(bd-ac)\Phi_1-{\textstyle{\frac{1}{2}}}\Phi_2(bc+ad)\Phi_2-2e,
\\
\alpha_{tr}&={\textstyle{\frac{1}{32}}}\big(H_1(\Phi_1)H_2(\Phi_2)\big)c^2+
 {\textstyle{\frac{1}{32}}}\big(H_1(\Phi_1)+H_2(\Phi_2)\big)d^2-{\textstyle{\frac{1}{2}}}(ad+bc)\Phi_1+{\textstyle{\frac{1}{2}}}(ac-bd)\Phi_2
 +3a^2+3b^2,
\\
\alpha_{h_1i_1}&={\textstyle{\frac{1}{2}}}(bd+ac)\Phi_1-{\textstyle{\frac{1}{2}}}(bc-ad)\Phi_2-4ab-2e,
\\
\alpha_{h_1i_2}&={\textstyle{\frac{1}{32}}}\big(H_1(\Phi_1)+H_2(\Phi_2)\big)c^2+{\textstyle{\frac{1}{32}}}\big(H_1(\Phi_1)+
H_2(\Phi_2)\big)d^2+
{\textstyle{\frac{1}{2}}}(bc-ad)\Phi_1+{\textstyle{\frac{1}{2}}}(ac+bd)\Phi_2+3a^2-b^2,
\\
\alpha_{h_2i_1}&=-{\textstyle{\frac{1}{32}}}H_1\big((\Phi_1)+H_2(\Phi_2)\big)c^2-{\textstyle{\frac{1}{32}}}\big(H_1(\Phi_1)+
H_2(\Phi_2)\big)d^2+{\textstyle{\frac{1}{2}}}(bc-ad)\Phi_1+{\textstyle{\frac{1}{2}}}(ac+bd)\Phi_2+
a^2-3b^2,
\\
\alpha_{h_2i_2}&=-{\textstyle{\frac{1}{2}}}(ac+bd)\Phi_1-{\textstyle{\frac{1}{2}}}(ad-bc)\Phi_2+4ab-2e.
\endaligned
\]

\subsection{Homogeneity 3}
\label{Hom3}
In this homogeneity, we have exactly six curvature coefficients:
\begin{eqnarray*}
\footnotesize\aligned
\kappa^{h_1h_2}_{i_1}&=\widehat{I}^\ast_1([\widehat{H}_1,\widehat{H}_2]-4\widehat{T})\\
&=-\widehat{H}_2(\alpha_{h_1i_1})+\alpha_{h_1h_2}H_2(\alpha_{h_2i_1})+\alpha_{h_1h_1}H_1(\alpha_{h_2i_1})+
\beta_{i_1h_1}\big( \alpha_{h_1h_1}\zero{H_1(\alpha_{h_2h_1})}-
\\
&-\alpha_{h_2h_1}\zero{H_1(\alpha_{h_1h_1})}-\alpha_{h_2h_2}\zero{H_2(\alpha_{h_1h_1})}-
\alpha_{h_2d}\underbrace{\widehat{D}(\alpha_{h_1h_1})}_{-\alpha_{h_1h_1}}-\alpha_{h_2r}\underbrace{\widehat{R}(\alpha_{h_1h_1})}_{-\alpha_{h_2h_1}}
-\alpha_{h_2i_1}\zero{\widehat{I}_1(\alpha_{h_1h_1})}-
\\
&-\alpha_{h_2i_2}\zero{\widehat{I}_2(\alpha_{h_1h_1})}-\alpha_{h_2j}\zero{\widehat{J}(\alpha_{h_1h_1})}+\alpha_{h_1h_2}\zero{H_2(\alpha_{h_2h_1})}\big)
+\beta_{i_1h_2}\big(\alpha_{h_1h_1}\zero{H_1(\alpha_{h_2h_2})}-
\\
&-\alpha_{h_2h_1}\zero{H_1(\alpha_{h_1h_2})}
-\alpha_{h_1h_2}\zero{H_2(\alpha_{h_2h_2})}-\alpha_{h_2h_2}\zero{H_2(\alpha_{h_1h_2})}-
\alpha_{h_2d}\underbrace{\widehat{D}(\alpha_{h_1h_2})}_{-\alpha_{h_1h_2}}-
\alpha_{h_2r}\underbrace{\widehat{R}(\alpha_{h_1h_2})}_{-\alpha_{h_2h_2}}-
\\
&-\alpha_{h_2i_1}\zero{\widehat{I}_1(\alpha_{h_1h_2})}-\alpha_{h_2i_2}\zero{\widehat{I}_2(\alpha_{h_1h_2})}-\zero{\alpha_{h_2j}\widehat{J}(\alpha_{h_1h_2})}\big)+
\beta_{i_1t}\big(4\alpha_{h_1h_1}\alpha_{h_2h_2}-4\alpha_{h_1h_2}\alpha_{h_2h_1}\big),
\endaligned
\end{eqnarray*}
\begin{eqnarray*}
\footnotesize\aligned
\kappa^{h_1h_2}_{i_2}&=\widehat{I}^\ast_2([\widehat{H}_1,\widehat{H}_2]-4\widehat{T})\\
&-\alpha_{h_1j}-
\widehat{H}_2(\alpha_{h_1i_2})+\alpha_{h_1h_2}H_2(\alpha_{h_2i_2})+\alpha_{h_1h_1}H_1(\alpha_{h_2i_2})+\beta_{i_2h_1}\big(
\alpha_{h_1h_1}\zero{H_1(\alpha_{h_2h_1})}-
\\
&-\alpha_{h_2h_1}\zero{H_1(\alpha_{h_1h_1})}-\alpha_{h_2h_2}\zero{H_2(\alpha_{h_1h_1})}-
\alpha_{h_2d}\underbrace{\widehat{D}(\alpha_{h_1h_1})}_{-\alpha_{h_1h_1}}-\alpha_{h_2r}\underbrace{\widehat{R}(\alpha_{h_1h_1})}_{-\alpha_{h_2h_1}}
-\alpha_{h_2i_1}\zero{\widehat{I}_1(\alpha_{h_1h_1})}-
\\
&-\alpha_{h_2i_2}\zero{\widehat{I}_2(\alpha_{h_1h_1})}-\alpha_{h_2j}\zero{\widehat{J}(\alpha_{h_1h_1})}+
\alpha_{h_1h_2}\zero{H_2(\alpha_{h_2h_1})}\big)
+\beta_{i_2h_2}\big(\alpha_{h_1h_1}\zero{H_1(\alpha_{h_2h_2})}-
\\
&-\alpha_{h_2h_1}\zero{H_1(\alpha_{h_1h_2})}-\alpha_{h_1h_2}\zero{H_2(\alpha_{h_2h_2})}-\alpha_{h_2h_2}\zero{H_2(\alpha_{h_1h_2})}-
\alpha_{h_2d}\underbrace{\widehat{D}(\alpha_{h_1h_2})}_{-\alpha_{h_1h_2}}-
\alpha_{h_2r}\underbrace{\widehat{R}(\alpha_{h_1h_2})}_{-\alpha_{h_2h_2}}-
\\
&-\alpha_{h_2i_1}\zero{\widehat{I}_1(\alpha_{h_1h_2})}-\alpha_{h_2i_2}\zero{\widehat{I}_2(\alpha_{h_1h_2})}-
\zero{\alpha_{h_2j}\widehat{J}(\alpha_{h_1h_2})}\big)+
\beta_{i_2t}\big(4\alpha_{h_1h_1}\alpha_{h_2h_2}-4\alpha_{h_1h_2}\alpha_{h_2h_1}\big),
\endaligned
\end{eqnarray*}
\begin{eqnarray*}
\footnotesize\aligned
\kappa^{h_1t}_d&=\widehat{D}^\ast([\widehat{H}_1,\widehat{T}])=-\alpha_{h_1j}-\widehat{T}(\alpha_{h_1d})+
\alpha_{h_1h_2}H_2(\alpha_{td})+\alpha_{h_1h_1}H_1(\alpha_{td})+\beta_{dh_1}\big(\alpha_{h_1h_1}\zero{H_1(\alpha_{th_1})}-
\\
&-\alpha_{th_1}\zero{H_1(\alpha_{h_1h_1})}-
\alpha_{th_2}\zero{H_2(\alpha_{h_1h_1})}-\alpha_{tt}\zero{T(\alpha_{h_1h_1})}-\alpha_{ti_1}\zero{\widehat{I}_1(\alpha_{h_1h_1})}-
\alpha_{ti_2}\zero{\widehat{I}_2(\alpha_{h_1h_1})}-
\\
&-
\alpha_{td}\underbrace{\widehat{D}(\alpha_{h_1h_1})}_{-\alpha_{h_1h_1}}-\alpha_{tr}\zero{\widehat{R}(\alpha_{h_1h_1})}-
\alpha_{tj}\zero{\widehat{J}(\alpha_{h_1h_1})}\big)-\beta_{dh_2}\big(\alpha_{td}\underbrace{\widehat{D}(\alpha_{h_1h_2})}_{-\alpha_{h_1h_2}}+
\alpha_{tr}\underbrace{\widehat{R}(\alpha_{h_1h_2})}_{-\alpha_{h_2h_2}}+\alpha_{h_1h_1}\zero{H_1(\alpha_{th_2})}-
\\
&+ \alpha_{th_1}\zero{H_1(\alpha_{h_1h_2})}
-\alpha_{h_1h_2}\zero{H_2(\alpha_{th_2})}+\alpha_{th_2}\zero{H_2(\alpha_{h_1h_2})}+
\alpha_{tt}\zero{T(\alpha_{h_1h_2})}+\alpha_{ti_1}\zero{\widehat{I}_1(\alpha_{h_1h_2})}+
\\
&+
\alpha_{ti_2}\zero{\widehat{I}_2(\alpha_{h_1h_2})}+\alpha_{tj}\zero{\widehat{J}(\alpha_{h_1h_2})}\big)+
\beta_{dt}\big(4\alpha_{h_1h_1}\alpha_{th_2}+
\alpha_{h_1h_1}\alpha_{tt}\Phi_1-4\alpha_{h_1h_2}\alpha_{th_1}+\alpha_{h_1h_2}\alpha_{tt}\Phi_2\big),
\endaligned
\end{eqnarray*}
\begin{eqnarray*}
\footnotesize\aligned
\kappa^{h_2t}_d&=\widehat{D}^\ast([\widehat{H}_2,\widehat{T}])=-\alpha_{h_2j}-\widehat{T}(\alpha_{h_2d})+
\alpha_{h_2h_2}H_2(\alpha_{td})+\alpha_{h_2h_1}H_1(\alpha_{td})+\beta_{dh_1}\big(\alpha_{h_2h_1}\zero{H_1(\alpha_{th_1})}-
\\
&-\alpha_{th_1}\zero{H_1(\alpha_{h_2h_1})}-
\alpha_{th_2}\zero{H_2(\alpha_{h_2h_1})}-\alpha_{tt}\zero{T(\alpha_{h_2h_1})}-
\alpha_{ti_2}\zero{\widehat{I}_2(\alpha_{h_2h_1})}-\alpha_{ti_1}\zero{\widehat{I}_1(\alpha_{h_2h_1})}-
\\
&-\alpha_{td}\underbrace{\widehat{D}(\alpha_{h_2h_1})}_{-\alpha_{h_2h_1}}-\alpha_{tr}\underbrace{\widehat{R}(\alpha_{h_2h_1})}_{\alpha_{h_1h_1}}-
\alpha_{tj}\zero{\widehat{J}(\alpha_{h_2h_1})}\big)-\beta_{dh_2}\big(
\alpha_{td}\underbrace{\widehat{D}(\alpha_{h_2h_2})}_{-\alpha_{h_2h_2}}+
\alpha_{tr}\underbrace{\widehat{R}(\alpha_{h_2h_2})}_{\alpha_{h_1h_2}}-\alpha_{h_2h_1}\zero{H_1(\alpha_{th_2})}
+
\\
&+\alpha_{th_1}\zero{H_1(\alpha_{h_2h_2})}
-\alpha_{h_2h_2}\zero{H_2(\alpha_{th_2})}+\alpha_{th_2}\zero{H_2(\alpha_{h_2h_2})}+
\alpha_{tt}\zero{T(\alpha_{h_2h_2})}+\alpha_{ti_1}\zero{\widehat{I}_1(\alpha_{h_2h_2})}+
\\
&+
\alpha_{ti_2}\zero{\widehat{I}_2(\alpha_{h_2h_2})}+\alpha_{tj}\zero{\widehat{J}(\alpha_{h_2h_2})}\big)+
\beta_{dt}\big(4\alpha_{h_2h_1}\alpha_{th_2}+
\alpha_{h_2h_1}\alpha_{tt}\Phi_1-4\alpha_{h_2h_2}\alpha_{th_1}+\alpha_{h_2h_2}\alpha_{tt}\Phi_2\big),
\endaligned
\end{eqnarray*}
\begin{eqnarray*}
\footnotesize\aligned
\kappa^{h_1t}_r&=\widehat{R}^\ast([\widehat{H}_1,\widehat{T}])=-\widehat{T}(\alpha_{h_1r})+
\alpha_{h_1h_2}H_2(\alpha_{tr})+\alpha_{h_1h_1}H_1(\alpha_{tr})+\beta_{rh_1}\big(\alpha_{h_1h_1}\zero{H_1(\alpha_{th_1})}-
\\
&-\alpha_{th_1}\zero{H_1(\alpha_{h_1h_1})}-
\alpha_{th_2}\zero{H_2(\alpha_{h_1h_1})}-\alpha_{tt}\zero{T(\alpha_{h_1h_1})}-\alpha_{ti_1}\zero{\widehat{I}_1(\alpha_{h_1h_1})}-
\alpha_{ti_2}\zero{\widehat{I}_2(\alpha_{h_1h_1})}-
\\
&-\alpha_{td}\underbrace{\widehat{D}(\alpha_{h_1h_1})}_{-\alpha_{h_1h_1}}-
\alpha_{tr}\underbrace{\widehat{R}(\alpha_{h_1h_1})}_{-\alpha_{h_2h_1}}-
\alpha_{tj}\zero{\widehat{J}(\alpha_{h_1h_1})}\big)+\beta_{rh_2}\big(\alpha_{h_1h_1}\zero{H_1(\alpha_{th_2})}-
\alpha_{td}\underbrace{\widehat{D}(\alpha_{h_1h_2})}_{-\alpha_{h_1h_2}}-
\alpha_{tr}\underbrace{\widehat{R}(\alpha_{h_1h_2})}_{-\alpha_{h_2h_2}}-
\\
&- \alpha_{th_1}\zero{H_1(\alpha_{h_1h_2})}
+\alpha_{h_1h_2}\zero{H_2(\alpha_{th_2})}-\alpha_{th_2}\zero{H_2(\alpha_{h_1h_2})}-
\alpha_{tt}\zero{T(\alpha_{h_1h_2})}-\alpha_{ti_1}\zero{\widehat{I}_1(\alpha_{h_1h_2})}-
\\
&-
\alpha_{ti_2}\zero{\widehat{I}_2(\alpha_{h_1h_2})}-\alpha_{tj}\zero{\widehat{J}(\alpha_{h_1h_2})}\big)+
\beta_{rt}\big(4\alpha_{h_1h_1}\alpha_{th_2}+
\alpha_{h_1h_1}\alpha_{tt}\Phi_1-4\alpha_{h_1h_2}\alpha_{th_1}+\alpha_{h_1h_2}\alpha_{tt}\Phi_2\big),
\endaligned
\end{eqnarray*}
\begin{eqnarray*}
\footnotesize\aligned
\kappa^{h_2t}_r&=\widehat{R}^\ast([\widehat{H}_2,\widehat{T}])=-\widehat{T}(\alpha_{h_2r})+
\alpha_{h_2h_2}H_2(\alpha_{tr})+\alpha_{h_2h_1}H_1(\alpha_{tr})+\beta_{rh_1}\big(\alpha_{h_2h_1}\zero{H_1(\alpha_{th_1})}-
\\
&-\alpha_{th_1}\zero{H_1(\alpha_{h_2h_1})}-
\alpha_{th_2}\zero{H_2(\alpha_{h_2h_1})}-\alpha_{tt}\zero{T(\alpha_{h_2h_1})}-\alpha_{ti_1}\zero{\widehat{I}_1(\alpha_{h_2h_1})}-
\alpha_{ti_2}\zero{\widehat{I}_2(\alpha_{h_2h_1})}-
\\
&-\alpha_{td}\underbrace{\widehat{D}(\alpha_{h_2h_1})}_{-\alpha_{h_2h_1}}-
\alpha_{tr}\underbrace{\widehat{R}(\alpha_{h_2h_1})}_{\alpha_{h_1h_1}}-
\alpha_{tj}\zero{\widehat{J}(\alpha_{h_2h_1})}\big)+\beta_{rh_2}\big(\alpha_{h_2h_1}\zero{H_1(\alpha_{th_2})}-
\alpha_{td}\underbrace{\widehat{D}(\alpha_{h_2h_2})}_{-\alpha_{h_2h_2}}
-\alpha_{tr}\underbrace{\widehat{R}(\alpha_{h_2h_2})}_{\alpha_{h_1h_2}}-
\\
&- \alpha_{th_1}\zero{H_1(\alpha_{h_2h_2})}
+\alpha_{h_2h_2}\zero{H_2(\alpha_{th_2})}-\alpha_{th_2}\zero{H_2(\alpha_{h_2h_2})}-
\alpha_{tt}\zero{T(\alpha_{h_2h_2})}-\alpha_{ti_1}\zero{\widehat{I}_1(\alpha_{h_2h_2})}-
\\
&-
\alpha_{ti_2}\zero{\widehat{I}_2(\alpha_{h_2h_2})}-\alpha_{tj}\zero{\widehat{J}(\alpha_{h_2h_2})}\big)+
\beta_{rt}\big(4\alpha_{h_2h_1}\alpha_{th_2}+
\alpha_{h_2h_1}\alpha_{tt}\Phi_1-4\alpha_{h_2h_2}\alpha_{th_1}+\alpha_{h_2h_2}\alpha_{tt}\Phi_2\big),
\endaligned
\end{eqnarray*}
According to the expressions of the functions $\alpha_{{}_\bullet
{}_\bullet}$, obtained in homogeneity two and according to the
properties of the horizontal vector field $Y=H_1,H_2,T$ we have:
\[
\footnotesize\aligned
Y(\alpha_{td})&=-\frac{1}{2}(bc+ad)\,Y(\Phi_2) -\frac{1}{2} (ac-bd)\,Y(\Phi_1),
\\
Y(\alpha_{tr})&={\textstyle{\frac{1}{32}}}\big[Y(H_1(\Phi_1))+Y(H_2(\Phi_2))\big]c^2+
 {\textstyle{\frac{1}{32}}}\big[Y(H_1(\Phi_1))+Y(H_2(\Phi_2))\big]d^2+
{\textstyle{\frac{1}{2}}}Y(\Phi_2)ac-
\\
&-{\textstyle{\frac{1}{2}}}Y(\Phi_1)bc-{\textstyle{\frac{1}{2}}}Y(\Phi_2)bd-{\textstyle{\frac{1}{2}}}Y(\Phi_1)ad,
\\
Y(\alpha_{h_1i_1})&={\textstyle{\frac{1}{2}}}(bd+ac)Y(\Phi_1)+{\textstyle{\frac{1}{2}}}(ad-bc)Y(\Phi_2),
\\
Y(\alpha_{h_1i_2})&={\textstyle{\frac{1}{32}}}\big[Y(H_1(\Phi_1))+Y(H_2(\Phi_2))\big]c^2+{\textstyle{\frac{1}{32}}}\big[Y(H_1(\Phi_1))+
Y(H_2(\Phi_2))\big]d^2+{\textstyle{\frac{1}{2}}}Y(\Phi_2)ac+
\\
&+{\textstyle{\frac{1}{2}}}Y(\Phi_2)bd+
{\textstyle{\frac{1}{2}}}Y(\Phi_1)bc-{\textstyle{\frac{1}{2}}}Y(\Phi_1)ad,
\\
Y(\alpha_{h_2i_1})&=-{\textstyle{\frac{1}{32}}}(Y(H_1(\Phi_1))+Y(H_2(\Phi_2)))c^2-{\textstyle{\frac{1}{32}}}(Y(H_1(\Phi_1))+
Y(H_2(\Phi_2)))d^2+{\textstyle{\frac{1}{2}}}Y(\Phi_2)bd+
\\
&+
{\textstyle{\frac{1}{2}}}Y(\Phi_1)bc+{\textstyle{\frac{1}{2}}}Y(\Phi_2)ac-{\textstyle{\frac{1}{2}}}Y(\Phi_1)ad,
\\
Y(\alpha_{h_2i_2})&=-{\textstyle{\frac{1}{2}}}Y(\Phi_1)ac-{\textstyle{\frac{1}{2}}}Y(\Phi_2)ad+
{\textstyle{\frac{1}{2}}}Y(\Phi_2)bc-{\textstyle{\frac{1}{2}}}Y(\Phi_1)bd.
\endaligned
\]
Furthermore, we can simplify the following functions
which appear just above by using the 110 equations introduced
before Lemma~\ref{fiber-type}:
\[
\footnotesize\aligned
\widehat{H}_2(\alpha_{h_1i_1})&=\alpha_{h_2h_1}H_1(\alpha_{h_1i_1})+\alpha_{h_2h_2}H_2(\alpha_{h_1i_1})+
\alpha_{h_2d}\underbrace{\widehat{D}(\alpha_{h_1i_1})}_{-2\alpha_{h_1i_1}}+
\alpha_{h_2r}\underbrace{\widehat{R}(\alpha_{h_1i_1})}_{-\alpha_{h_1i_2}-\alpha_{h_2i_1}}+
\\
&+\alpha_{h_2i_1}\underbrace{\widehat{I}_1(\alpha_{h_1i_1})}_{\alpha_{h_1d}}+
\alpha_{h_2i_2}\underbrace{\widehat{I}_2(\alpha_{h_1i_1})}_{\alpha_{h_1r}}+
\alpha_{h_2j}\underbrace{\widehat{J}(\alpha_{h_1i_1})}_{-1}
\endaligned
\]
\[
\footnotesize\aligned
\widehat{H}_2(\alpha_{h_1i_2})&=\alpha_{h_2h_1}H_1(\alpha_{h_1i_2})+\alpha_{h_2h_2}H_2(\alpha_{h_1i_2})+
\alpha_{h_2d}\underbrace{\widehat{D}(\alpha_{h_1i_2})}_{-2\alpha_{h_1i_2}}+
\alpha_{h_2r}\underbrace{\widehat{R}(\alpha_{h_1i_2})}_{\alpha_{h_1i_1}-\alpha_{h_2i_2}}+
\\
&+\alpha_{h_2i_1}\underbrace{\widehat{I}_1(\alpha_{h_1i_2})}_{-\alpha_{h_1r}}+
\alpha_{h_2i_2}\underbrace{\widehat{I}_2(\alpha_{h_1i_2})}_{\alpha_{h_1d}}+
\alpha_{h_2j}\zero{\widehat{J}(\alpha_{h_1i_2})},
\\
\widehat{T}(\alpha_{h_1d})&=\alpha_{tt}T(\alpha_{h_1d})+\alpha_{th_1}H_1(\alpha_{h_1d})+\alpha_{th_2}H_2(\alpha_{h_1d})+
\alpha_{td}\underbrace{\widehat{D}(\alpha_{h_1d})}_{-\alpha_{h_1d}}+
\alpha_{tr}\underbrace{\widehat{R}(\alpha_{h_1d})}_{-\alpha_{h_2d}}+
\alpha_{ti_2}\underbrace{\widehat{I}_2(\alpha_{h_1d})}_{-2}+
\\
&+\alpha_{ti_1}\zero{\widehat{I}_1(\alpha_{h_1d})}+
\alpha_{tj}\zero{\widehat{J}(\alpha_{h_1d})},
\endaligned
\]
\[
\footnotesize\aligned
\widehat{T}(\alpha_{h_2d})&=\alpha_{tt}T(\alpha_{h_2d})+\alpha_{th_1}H_1(\alpha_{h_2d})+\alpha_{th_2}H_2(\alpha_{h_2d})+
\alpha_{td}\underbrace{\widehat{D}(\alpha_{h_2d})}_{-\alpha_{h_2d}}+
\alpha_{tr}\underbrace{\widehat{R}(\alpha_{h_2d})}_{\alpha_{h_1d}}+\alpha_{ti_1}\underbrace{\widehat{I}_1(\alpha_{h_2d})}_{2}+
\\
&+\alpha_{ti_2}\zero{\widehat{I}_2(\alpha_{h_2d})}+\alpha_{tj}\zero{\widehat{J}(\alpha_{h_2d})},
\\
\widehat{T}(\alpha_{h_1r})&=\alpha_{tt}T(\alpha_{h_1r})+\alpha_{th_1}H_1(\alpha_{h_1r})+\alpha_{th_2}H_2(\alpha_{h_1r})+
\alpha_{td}\underbrace{\widehat{D}(\alpha_{h_1r})}_{-\alpha_{h_1r}}+
\alpha_{tr}\underbrace{\widehat{R}(\alpha_{h_1r})}_{-\alpha_{h_2r}}+\alpha_{ti_1}\underbrace{\widehat{I}_1(\alpha_{h_1r})}_{-6}+
\\
&+\alpha_{ti_2}\zero{\widehat{I}_2(\alpha_{h_1r})}+\alpha_{tj}\zero{\widehat{J}(\alpha_{h_1r})},
\\
\widehat{T}(\alpha_{h_2r})&=\alpha_{tt}T(\alpha_{h_2r})+\alpha_{th_1}H_1(\alpha_{h_2r})+\alpha_{th_2}H_2(\alpha_{h_2r})+
\alpha_{td}\underbrace{\widehat{D}(\alpha_{h_2r})}_{-\alpha_{h_2r}}+
\alpha_{tr}\underbrace{\widehat{R}(\alpha_{h_2r})}_{\alpha_{h_1r}}+\alpha_{ti_2}\underbrace{\widehat{I}_2(\alpha_{h_2r})}_{-6}+
\\
&+\alpha_{ti_1}\zero{\widehat{I}_1(\alpha_{h_2r})}+\alpha_{tj}\zero{\widehat{J}(\alpha_{h_2r})}.
\endaligned
\]
Before replacing the above expressions in the curvature coefficients
in order to simplify, we should be aware of the following fact, which
helps us to substitute the vector field $T$ in terms of the two
basic sections $H_1$ and $H_2$ of $T^cM$:
\[\footnotesize\aligned
T(\Phi_i)=4[H_1,H_2](\Phi_i)=4\big[H_1(H_2((\Phi_i)))-H_2(H_1((\Phi_i)))\big], \ \ \
\ \ \ \ i=1,2.
\endaligned
\]
This is natural and this helps us to get more simplified
expressions. Thus, replacing the above expressions in the curvature
coefficients of homogeneity $3$ and carefully simplifying, we get:
\[
\footnotesize\aligned
\kappa^{h_1h_2}_{i_1}&=\alpha_{h_2j}-4\alpha_{ti_1}+{\textstyle{\frac{1}{32}}}\big[H_1(\Phi_1)\Phi_1-H_1(H_2(\Phi_2))+
H_2(\Phi_2)\Phi_1-H_1(H_1(\Phi_1))\big]c^3+
\\
&+{\textstyle{\frac{1}{32}}}\big[H_2(\Phi_2)\Phi_2+H_1(\Phi_1)\Phi_2-H_2(H_2(\Phi_2))-
H_2(H_1(\Phi_1))\big]d^3+12a^2b-2\Phi_2ce+12b^3-
\\
&-3\Phi_1a^2c-3\Phi_1b^2c+2\Phi_1de-3\Phi_2a^2d-3\Phi_2b^2d-
8ae+{\textstyle{\frac{1}{32}}}\big[H_2(\Phi_2)\Phi_2+
\\
&+H_1(\Phi_1)\Phi_2-H_2(H_2(\Phi_2))-H_2(H_1(\Phi_1))\big]c^2d+
{\textstyle{\frac{1}{32}}}\big[H_1(\Phi_1)\Phi_1-
\\
&-H_1(H_2(\Phi_2))+H_2(\Phi_2)\Phi_1-H_1(H_1(\Phi_1))\big]cd^2+{\textstyle{\frac{3}{8}}}\big[H_1\Phi_1+H_2(\Phi_2)\big]bc^2+
{\textstyle{\frac{3}{8}}}\big[H_1(\Phi_1)+H_2(\Phi_2)\big]bd^2,
\endaligned
\]
\[
\footnotesize\aligned
\kappa^{h_1h_2}_{i_2}&=-\alpha_{h_1j}-4\alpha_{ti_2}+{\textstyle{\frac{1}{32}}}\big[H_2(\Phi_2)\Phi_2+
H_1(\Phi_1)\Phi_2-H_2(H_2(\Phi_2))-H_2(H_1(\Phi_1))\big]c^3-
\\
&-{\textstyle{\frac{3}{8}}}\big[H_2(\Phi_2)+H_1(\Phi_1)\big]ad^2+{\textstyle{\frac{1}{32}}}\big[H_1(H_1(\Phi_1))-H_2(\Phi_2)\Phi_1-
H_1(\Phi_1)\Phi_1+H_1(H_2(\Phi_2))\big]d^3+
\\
&+2\Phi_1ce-8be-3\Phi_2a^2c-3\Phi_2b^2c+2\Phi_2de+
3\Phi_1a^2d+3\Phi_1b^2d-12a^3-12ab^2-
\\
&-{\textstyle{\frac{3}{8}}}\big[H_2(\Phi_2)+H_1(\Phi_1)\big]ac^2+
{\textstyle{\frac{1}{32}}}\big[H_1(H_1(\Phi_1))-H_2(\Phi_2)\Phi_1-H_1(\Phi_1)\Phi_1+H_1(H_2)\Phi_2\big]c^2d+
\\
&+{\textstyle{\frac{1}{32}}}\big[H_2(\Phi_2)\Phi_2+H_1(\Phi_1)\Phi_2-H_2(H_2(\Phi_2))-H_2(H_1(\Phi_1))\big]cd^2,
\endaligned
\]
\[
\footnotesize\aligned
\kappa^{h_1t}_d&=2\alpha_{ti_2}-\alpha_{h_1j}-{\textstyle{\frac{1}{8}}}\big[H_1(H_1(\Phi_2))-H_2(H_1(\Phi_1))\big]c^3+
{\textstyle{\frac{1}{8}}}\big[H_2(H_1(\Phi_2))-H_1(H_2(\Phi_2))\big]d^3+
\\
&+2\Phi_1ce-
8be+2\Phi_2de+{\textstyle{\frac{1}{8}}}\big[H_2(H_1(\Phi_2))-H_1(H_2(\Phi_2))\big]c^2d-
{\textstyle{\frac{1}{8}}}\big[H_1(H_1(\Phi_2))-H_2(H_1(\Phi_1))\big]cd^2,
\endaligned
\]
\[
\footnotesize\aligned
\kappa^{h_2t}_d&=-2\alpha_{ti_1}-\alpha_{h_2j}+{\textstyle{\frac{1}{8}}}\big[H_2(H_1(\Phi_2))-H_1(H_2(\Phi_2))\big]c^3+
{\textstyle{\frac{1}{8}}}\big[H_1(H_1(\Phi_2))-H_2(H_1(\Phi_1))\big]d^3+
\\
&+2\Phi_2ce-2\Phi_1de+8ae+
{\textstyle{\frac{1}{8}}}\big[H_1(H_1(\Phi_2))-H_2(H_1(\Phi_1))\big]c^2d+{\textstyle{\frac{1}{8}}}\big[H_2(H_1(\Phi_2))-H_1(H_2(\Phi_2))\big]cd^2,
\endaligned
\]
\[
\footnotesize\aligned
\kappa^{h_1t}_r&=6\alpha_{ti_1}+{\textstyle{\frac{1}{32}}}\big[-H_1(\Phi_1)\Phi_1-H_2(\Phi_2)\Phi_1-{\textstyle{4}}H_2(H_1(\Phi_2))+
{\textstyle{5}}H_1(H_2(\Phi_2))+ {\textstyle{\frac{1}{16}}}H_1(H_1(\Phi_1))\big]c^3-
\\
&-
{\textstyle{\frac{3}{8}}}\big[H_2(\Phi_2)+H_1(\Phi_1)\big]bd^2-{\textstyle{\frac{3}{8}}}\big[H_2(\Phi_2)+H_1(\Phi_1)\big]bc^2+
{\textstyle{\frac{1}{32}}}\big[-H_2(\Phi_2)\Phi_2+ 5H_2(H_1(\Phi_1))-
\\
&-
H_1(\Phi_1)\Phi_2-4H_1(H_1(\Phi_2))+H_2(H_2(\Phi_2))\big]c^2d+{\textstyle{\frac{1}{32}}}\big[-H_2(\Phi_2)\Phi_2+{\textstyle{5}}H_2(H_1(\Phi_1))-
\\
&- H_1(\Phi_1)\Phi_2-
{\textstyle{4}}H_1(H_1(\Phi_2))+H_2(H_2(\Phi_2))\big]d^3-12a^2b-12b^3+3\Phi_1a^2c+3\Phi_1cb^2+3\Phi_2a^2d+3\Phi_2b^2d+
\\
&+ {\textstyle{\frac{1}{32}}}\big[-H_1(\Phi_1)\Phi_1-H_2(\Phi_2)\Phi_1-
{\textstyle{4}}H_2(H_1(\Phi_2))+{\textstyle{5}}H_1(H_2(\Phi_2))+H_1(H_1(\Phi_1))\big]d^2c
\endaligned
\]
\[
\footnotesize\aligned
\kappa^{h_2t}_r&=6\alpha_{ti_2}-{\textstyle{\frac{1}{32}}}\big[H_2(\Phi_2)\Phi_2-{\textstyle{5}}H_2(H_1(\Phi_1))+H_1(\Phi_1)\Phi_2+
{\textstyle{4}}H_1(H_1(\Phi_2))-H_2(H_2(\Phi_2))\big]c^3+
\\
&+
{\textstyle{\frac{1}{32}}}\big[-{\textstyle{5}}H_1(H_2(\Phi_2))+{\textstyle{4}}H_2(H_1(\Phi_2))+H_2(\Phi_2)\Phi_1-H_1(H_1(\Phi_1))+
H_1(\Phi_1)\Phi_1\big]d^3+ 3ca^2\Phi_2+
\\
&+3cb^2\Phi_2-
3a^2d\Phi_1-3b^2d\Phi_1+12a^3+12ab^2+{\textstyle{\frac{3}{8}}}\big[H_2\Phi_2+H_1(\Phi_1)big]c^2a+
{\textstyle{\frac{1}{32}}}\big[-{\textstyle{5}}H_1(H_2(\Phi_2))+
\\
&+{\textstyle{4}}H_2(H_1(\Phi_2))+H_2(\Phi_2)\Phi_1-H_1(H_1(\Phi_1))+H_1(\Phi_1)\Phi_1\big]c^2d+
{\textstyle{\frac{1}{32}}}\big[-H_2(\Phi_2)\Phi_2+{\textstyle{5}}H_2(H_1(\Phi_1))-
\\
&-H_1(\Phi_1)\Phi_2-{\textstyle{4}}H_1(H_1(\Phi_2))+H_2(H_2(\Phi_2))\big]cd^2+
{\textstyle{\frac{3}{8}}}\big[H_2(\Phi_2)+H_1(\Phi_1)\big]ad^2.
\endaligned
\]
We therefore see here exactly four undetermined functions
$\alpha_{ti_1}$, $\alpha_{ti_2}$, $\alpha_{h_1j}$, $\alpha_{h_2j}$
within these six expressions. Although the number (six) of equations
is greater than the number (four) of unknowns, we can 
annihilate all
the six curvature coefficients by
making the following appropriate
determinations:
\[
\footnotesize\aligned \alpha_{ti_1}&={\textstyle{\frac{1}{192}}}\big[-H_2(H_2(\Phi_2))+
H_1(\Phi_1)\Phi_2-{\textstyle{5}}H_2(H_1(\Phi_1))+H_2(\Phi_2)\Phi_2+{\textstyle{4}}H_1(H_1(\Phi_2))\big]d^3+
\\
&+ {\textstyle{\frac{1}{192}}}\big[{\textstyle{4}}H_2(H_1(\Phi_2))+
H_2(\Phi_2)\Phi_1-H_1(H_1(\Phi_1))-{\textstyle{5}}H_1(H_2(\Phi_2))+H_1(\Phi_1)\Phi_1\big]c^3+
\\
&+ {\textstyle{\frac{1}{192}}}\big[{\textstyle{4}}H_2(H_1(\Phi_2))+
H_2(\Phi_2)\Phi_1-H_1(H_1(\Phi_1))-{\textstyle{5}}H_1(H_2(\Phi_2))+H_1(\Phi_1)\Phi_1\big]cd^2+
\\
&+ {\textstyle{\frac{1}{16}}}\big[H_2(\Phi_2)+
H_1(\Phi_1)\big]bc^2+{\textstyle{\frac{1}{192}}}\big[-H_2(H_2(\Phi_2))+H_1(\Phi_1)\Phi_2-{\textstyle{5}}H_2(H_1(\Phi_1))+
\\
&+ H_2(\Phi_2)\Phi_2+{\textstyle{4}}H_1(H_1(\Phi_2))\big]c^2d+
{\textstyle{\frac{1}{16}}}\big[H_2(\Phi_2)+H_1\Phi_1\big]bd^2+
\\
&+
{\textstyle{\frac{1}{2}}}\big[-\Phi_1a^2c+{\textstyle{4}}b^3-\Phi_1b^2c+{\textstyle{4}}ba^2-\Phi_2b^2d-\Phi_2a^2d\big],
\endaligned
\]
\[\footnotesize\aligned
\alpha_{ti_2}&={\textstyle{\frac{1}{192}}}\big[-H_2(H_2(\Phi_2))+H_1(\Phi_1)\Phi_2-{\textstyle{5}}H_2(H_1(\Phi_1))+H_2(\Phi_2)\Phi_2+
{\textstyle{4}}H_1(H_1(\Phi_2))\big]c^3-
\\
&-{\textstyle{\frac{1}{16}}}\big[H_1(\Phi_1)+H_2(\Phi_2)\big]ac^2-{\textstyle{\frac{1}{16}}}\big[H_1(\Phi_1)+H_2(\Phi_2)\big]ad^2+
{\textstyle{\frac{1}{192}}}\big[-H_2(H_2(\Phi_2))+H_1(\Phi_1)\Phi_2-
\\
&-{\textstyle{5}}H_2(H_1(\Phi_1))+H_2(\Phi_2)\Phi_2+{\textstyle{4}}H_1(H_1(\Phi_2))\big]cd^2+
{\textstyle{\frac{1}{192}}}\big[-{\textstyle{4}}H_2(H_1(\Phi_2))-H_1(\Phi_1)\Phi_1-
\\
&-H_2(\Phi_2)\Phi_1+{\textstyle{5}}H_1(H_2(\Phi_2))+
H_1(H_1(\Phi_1))\big]c^2d+{\textstyle{\frac{1}{192}}}\big[-{\textstyle{4}}H_2(H_1(\Phi_2))-H_1(\Phi_1)\Phi_1-
H_2(\Phi_2)\Phi_1+
\\
&+{\textstyle{5}}H_1(H_2(\Phi_2))+H_1(H_1(\Phi_1))\big]d^3-
{\textstyle{\frac{1}{2}}}\big[\Phi_2a^2c+\Phi_2b^2c-\Phi_1b^2d-
\Phi_1a^2d-{\textstyle{4}}ab^2+{\textstyle{4}}a^3\big],
\endaligned
\]
\[
\footnotesize\aligned
\alpha_{h_1j}&={\textstyle{\frac{1}{96}}}\big[-H_2(H_2(\Phi_2)+H_2(\Phi_2)\Phi_2+H_1(\Phi_1)\Phi_2+
{\textstyle{7}}H_2(H_1(\Phi_1))-{\textstyle{8}}H_1(H_1(\Phi_2))\big]c^3-
\\
&+{\textstyle{\frac{1}{96}}}\big[-H_1(\Phi_1)\Phi_1+
{\textstyle{8}}H_2(H_1(\Phi_2))-{\textstyle{7}}H_1(H_2(\Phi_2))-H_2(\Phi_2)\Phi_1+
H_1(H_1(\Phi_1))\big]c^2d+
\\
&+{\textstyle{\frac{1}{96}}}\big[-H_2(H_2(\Phi_2))+H_2(\Phi_2)\Phi_2+H_1(\Phi_1)\Phi_2+
{\textstyle{\frac{7}{16}}}H_2(H_1(\Phi_1))- {\textstyle{8}}H_1(H_1(\Phi_2))\big]cd^2+
\\
&+{\textstyle{\frac{1}{96}}}\big[-H_1(\Phi_1)\Phi_1+{\textstyle{8}}H_2(H_1(\Phi_2))-{\textstyle{7}}H_1(H_2(\Phi_2))-H_2(\Phi_2)\Phi_1+
H_1(H_1(\Phi_1))\big]d^3-
\\
&-{\textstyle{\frac{1}{8}}}\big[H_2(\Phi_2)+H_1(\Phi_1)\big]ac^2-{\textstyle{\frac{1}{8}}}\big[\frac{1}{8}H_2(\Phi_2)+H_1(\Phi_1)\big]ad^2
-\Phi_2a^2c-\Phi_2b^2c+2\Phi_1ce-
\\
&-8be+2\Phi_2de+\Phi_1b^2d+ \Phi_1a^2d-4ab^2-4a^3,
\endaligned
\]
\[
\footnotesize\aligned
\alpha_{h_2j}&={\textstyle{\frac{1}{96}}}\big[-H_1(\Phi_1)\Phi_1+{\textstyle{8}}H_2(H_1(\Phi_2))-{\textstyle{7}}H_1(H_2(\Phi_2))-H_2(\Phi_2)\Phi_1+
H_1(H_1(\Phi_1))\big]c^3+
\\
&-{\textstyle{\frac{1}{8}}}\big[H_2(\Phi_2)+
H_1(\Phi_1)\big]bd^2+{\textstyle{\frac{1}{96}}}\big[-H_2(\Phi_2)\Phi_2-
H_1(\Phi_1)\Phi_2+H_2(H_2(\Phi_2))+8H_1(H_1(\Phi_2))-
\\
&-7H_2(H_1(\Phi_1))\big]d^3-
{\textstyle{\frac{1}{8}}}\big[H_2(\Phi_2)+H_1(\Phi_1)\big]bc^2+{\textstyle{\frac{1}{96}}}\big[-H_1(\Phi_1)\Phi_1+8H_2(H_1(\Phi_2))-
\\
&-
7H_1(H_2(\Phi_2))-H_2(\Phi_2)\Phi_1+H_1(H_1(\Phi_1))\big]cd^2+{\textstyle{\frac{1}{96}}}\big[-H_2(\Phi_2)\Phi_2-
H_1(\Phi_1)\Phi_2+H_2(H_2(\Phi_2))+
\\
&+8H_1(H_1(\Phi_2))-7H_2(H_1(\Phi_1))\big]c^2d+
\Phi_1a^2c-2\Phi_1de+\Phi_2b^2d+\Phi_2a^2d+8ae{\textstyle{4}}b^3+\Phi_1b^2c-
\\
&-4a^2b+2\Phi_2ce,
\endaligned
\]
Lastly, reminding Lemma~\ref{fiber-type}, 
we must have (after possible simplifications) by 
identification:
\[
\footnotesize
 \left[ \aligned
\alpha_{ti_1}&=-2({\textstyle{\frac{1}{4}}}a^2d+{\textstyle{\frac{1}{4}}}abc)\,\Phi_2
-2
(-{\textstyle{\frac{1}{4}}}abd+{\textstyle{\frac{1}{4}}}a^2c)\,\Phi_1 +
({\textstyle{\frac{1}{24}}}d^3+{\textstyle{\frac{1}{24}}}c^2d)\,\delta_3 -
\\
& \ \ \ \ \ -\frac{1}{4}
({\textstyle{\frac{1}{8}}}d^2b+{\textstyle{\frac{1}{8}}}bc^2)\,(H_1(\Phi_1)+H_2(\Phi_2))
+
({\textstyle{\frac{1}{8}}}c^3+{\textstyle{\frac{1}{8}}}cd^2)\,\delta_5+
(-{\textstyle{\frac{1}{2}}}db^2+{\textstyle{\frac{1}{2}}}bca)\,\Phi_2 -
({\textstyle{\frac{1}{2}}}bda+{\textstyle{\frac{1}{2}}}cb^2)\,\Phi_1 +
\\
& \ \ \ \ \ +
({\textstyle{\frac{1}{4}}}c^3+{\textstyle{\frac{1}{4}}}cd^2)\,\delta_{16} +
({\textstyle{\frac{1}{4}}}c^2d+{\textstyle{\frac{1}{4}}}d^3)\,\delta_{17} +\frac{3}{16}
({\textstyle{\frac{1}{2}}}d^2b+{\textstyle{\frac{1}{2}}}bc^2)\,(H_1(\Phi_1)+H_2(\Phi_2))
+
2a^2b+2b^3,
\\
\alpha_{ti_2}&=
-2({\textstyle{\frac{1}{4}}}bda+{\textstyle{\frac{1}{4}}}cb^2)\,\Phi_2 -2
(-{\textstyle{\frac{1}{4}}}db^2+{\textstyle{\frac{1}{4}}}bca)\,\Phi_1 +
({\textstyle{\frac{1}{24}}}c^3+{\textstyle{\frac{1}{24}}}cd^2)\,\delta_3 -\frac{1}{4}
(-{\textstyle{\frac{1}{8}}}ac^2-{\textstyle{\frac{1}{8}}}ad^2)\,(H_1(\Phi_1)
+
\\
& \ \ \ \ \ +H_2(\Phi_2))
+ (-{\textstyle{\frac{1}{8}}}c^2d-{\textstyle{\frac{1}{8}}}d^3)\,\delta_5 +
(-{\textstyle{\frac{1}{2}}}a^2c+{\textstyle{\frac{1}{2}}}bda)\,\Phi_2 +
({\textstyle{\frac{1}{2}}}a^2d+{\textstyle{\frac{1}{2}}}bca)\,\Phi_1 -
({\textstyle{\frac{1}{4}}}c^2d+{\textstyle{\frac{1}{4}}}d^3)\,\delta_{16} +
\\
& \ \ \ \ \ + ({\textstyle{\frac{1}{4}}}c^3+{\textstyle{\frac{1}{4}}}cd^2)\,\delta_{17}
+\frac{3}{16} (-{\textstyle{\frac{1}{2}}}ac^2
-{\textstyle{\frac{1}{2}}}d^2a)\,(H_1(\Phi_1)+H_2(\Phi_2)) - 2a^3-2ab^2,
\\
\alpha_{h_1j}&= 2(de)\,\Phi_2+2(ce)\,\Phi_1 -
({\textstyle{\frac{1}{6}}}c^3+{\textstyle{\frac{1}{6}}}cd^2)\,\delta_3 -\frac{1}{4}
({\textstyle{\frac{1}{2}}}ac^2+{\textstyle{\frac{1}{2}}}d^2a)\,(H_1(\Phi_1)+H_2(\Phi_2))
+
\\
& \ \ \ \ \ +({\textstyle{\frac{1}{2}}}d^3+{\textstyle{\frac{1}{2}}}c^2d)\,\delta_5 -
(a^2c+cb^2)\,\Phi_2+(a^2d+db^2)\,\Phi_1-8be-4a^3-4ab^2,
\\
\alpha_{h_2j}&= 2(ce)\,\Phi_2-2(ed)\,\Phi_1 +
({\textstyle{\frac{1}{6}}}d^3+{\textstyle{\frac{1}{6}}}c^2d)\,\delta_3 -\frac{1}{4}
({\textstyle{\frac{1}{2}}}bc^2+{\textstyle{\frac{1}{2}}}d^2b)\,(H_1(\Phi_1)+H_2(\Phi_2))
+
\\
& \ \ \ \ \ + ({\textstyle{\frac{1}{2}}}c^3+{\textstyle{\frac{1}{2}}}cd^2)\,\delta_5
+(db^2+da^2)\,\Phi_2+(ca^2+cb^2)\,\Phi_1-4a^2b+8ae-4b^3.
\endaligned\right.
\]
These equations will be satisfied if and only if the functions
$\delta_{{}_\bullet}$ are determined as follows:
\[
\footnotesize\aligned
\delta_3&={\textstyle{\frac{1}{2}}}H_1(H_1(\Phi_2))-{\textstyle{\frac{1}{16}}}\Phi_2H_2(\Phi_2)-{\textstyle{\frac{7}{16}}}H_2(H_1(\Phi_1))+
{\textstyle{\frac{1}{16}}}H_2(H_2(\Phi_2))-{\textstyle{\frac{1}{16}}}\Phi_2H_1(\Phi_1),
\\
\delta_5&=-{\textstyle{\frac{1}{48}}}\Phi_1H_2(\Phi_2)-{\textstyle{\frac{7}{48}}}H_1(H_2(\Phi_2))+{\textstyle{\frac{1}{48}}}H_1(H_1(\Phi_1))-
{\textstyle{\frac{1}{48}}}\Phi_1H_1(\Phi_1)+{\textstyle{\frac{1}{6}}}H_2(H_1(\Phi_2)),
\\
\delta_{16}&={\textstyle{\frac{1}{32}}}\Phi_1H_1(\Phi_1)-{\textstyle{\frac{1}{32}}}H_1(H_1(\Phi_1))+
{\textstyle{\frac{1}{32}}}\Phi_1H_2(\Phi_2)-{\textstyle{\frac{1}{32}}}H_1(H_2(\Phi_2)),
\\
\delta_{17}&={\textstyle{\frac{1}{32}}}\Phi_2H_1(\Phi_1)-{\textstyle{\frac{1}{32}}}H_2(H_2(\Phi_2))+
{\textstyle{\frac{1}{32}}}\Phi_2H_2(\Phi_2)-{\textstyle{\frac{1}{32}}}H_2(H_1(\Phi_1)).
\endaligned
\]

\subsection{Homogeneity 4}
\label{Hom4}
In this homogeneity, we see exactly five curvature coefficients:
\begin{eqnarray*}
\footnotesize\aligned
\kappa^{h_1h_2}_j&=\widehat{J}^\ast([\widehat{H}_1,\widehat{H}_2]-4\widehat{T})-\widehat{H}_2(\alpha_{h_1j})+\alpha_{h_1h_2}H_2(\alpha_{h_2j})+\alpha_{h_1h_1}H_1(\alpha_{h_2j})+
\beta_{jh_1}\big(
\alpha_{h_1h_1}\zero{H_1(\alpha_{h_2h_1})}-
\\
&-\alpha_{h_2h_1}\zero{H_1(\alpha_{h_1h_1})}-\alpha_{h_2h_2}\zero{H_2(\alpha_{h_1h_1})}-
\alpha_{h_2i_1}\zero{\widehat{I}_1(\alpha_{h_1h_1})}-
\alpha_{h_2i_2}\zero{\widehat{I}_2(\alpha_{h_1h_1})}-\alpha_{h_2j}\zero{\widehat{J}(\alpha_{h_1h_1})}+
\\
&+\alpha_{h_1h_2}\zero{H_2(\alpha_{h_2h_1})}-
\alpha_{h_2d}\underbrace{\widehat{D}(\alpha_{h_1h_1})}_{-\alpha_{h_1h_1}}-
\alpha_{h_2r}\underbrace{\widehat{R}(\alpha_{h_1h_1})}_{-\alpha_{h_2h_1}}\big)+
\beta_{jh_2}\big(-
\alpha_{h_2d}\underbrace{\widehat{D}(\alpha_{h_1h_2})}_{-\alpha_{h_1h_2}}-
\alpha_{h_2r}\underbrace{\widehat{R}(\alpha_{h_1h_2})}_{-\alpha_{h_2h_2}}+
\\
&+\alpha_{h_1h_1}\zero{H_1(\alpha_{h_2h_2})}-\alpha_{h_2h_1}\zero{H_1(\alpha_{h_1h_2})}-
\alpha_{h_1h_2}\zero{H_2(\alpha_{h_2h_2})}-\alpha_{h_2h_2}\zero{H_2(\alpha_{h_1h_2})}-
\\
&-\alpha_{h_2i_1}\zero{\widehat{I}_1(\alpha_{h_1h_2})}-
\alpha_{h_2i_2}\zero{\widehat{I}_2(\alpha_{h_1h_2})}-\alpha_{h_2j}\zero{\widehat{J}(\alpha_{h_1h_2})}\big)+
\beta_{jt}\big(4\alpha_{h_1h_1}\alpha_{h_2h_2}-4\alpha_{h_1h_2}\alpha_{h_2h_1}\big),
\endaligned
\end{eqnarray*}
\begin{eqnarray*}
\footnotesize\aligned
\kappa^{h_1t}_{i_1}&=\widehat{I}^\ast_1([\widehat{H}_1,\widehat{T}])-\widehat{T}(\alpha_{h_1i_1})+\alpha_{h_1h_2}H_2(\alpha_{ti_1})+\alpha_{h_1h_1}H_1(\alpha_{ti_1})+
\beta_{i_1h_1}\big(\alpha_{h_1h_1}\zero{H_1(\alpha_{th_1})}-
\\
&-\alpha_{th_1}\zero{H_1(\alpha_{h_1h_1})}-
\alpha_{th_2}\zero{H_2(\alpha_{h_1h_1})}-\alpha_{tt}\zero{T(\alpha_{h_1h_1})}-
\alpha_{ti_1}\zero{\widehat{I}_1(\alpha_{h_1h_1})}-\alpha_{ti_2}\zero{\widehat{I}_2(\alpha_{h_1h_1})}-
\\
&-
\alpha_{td}\underbrace{\widehat{D}(\alpha_{h_1h_1})}_{-\alpha_{h_1h_1}}-\alpha_{tr}\underbrace{\widehat{R}(\alpha_{h_1h_1})}_{-\alpha_{h_2h_1}}
\alpha_{tj}\zero{\widehat{J}(\alpha_{h_1h_1})}\big)+\beta_{i_1h_2}(-
\alpha_{td}\underbrace{\widehat{D}(\alpha_{h_1h_2})}_{-\alpha_{h_1h_2}}-
\alpha_{tr}\underbrace{\widehat{R}(\alpha_{h_1h_2})}_{-\alpha_{h_2h_2}}+
\\
&+\alpha_{h_1h_1}\zero{H_1(\alpha_{th_2})}-\alpha_{th_1}\zero{H_1(\alpha_{h_1h_2})}
+\alpha_{h_1h_2}\zero{H_2(\alpha_{th_2})}-\alpha_{th_2}\zero{H_2(\alpha_{h_1h_2})}-
\alpha_{tt}\zero{T(\alpha_{h_1h_2})}-
\\
&-\alpha_{ti_1}\zero{\widehat{I}_1(\alpha_{h_1h_2})}-
\alpha_{ti_2}\zero{\widehat{I}_2(\alpha_{h_1h_2})}-\alpha_{tj}\zero{\widehat{J}(\alpha_{h_1h_2})}\big)+
\beta_{i_1t}\big(4\alpha_{h_1h_1}\alpha_{th_2}+\alpha_{h_1h_1}\alpha_{tt}\Phi_1-
\\
&- 4\alpha_{h_1h_2}\alpha_{th_1}+\alpha_{h_1h_2}\alpha_{tt}\Phi_2\big),
\endaligned
\end{eqnarray*}
\begin{eqnarray*}
\footnotesize\aligned
\kappa^{h_1t}_{i_2}&=\widehat{I}^\ast_2([\widehat{H}_1,\widehat{T}])
=-\widehat{T}(\alpha_{h_1i_2})+\alpha_{h_1h_2}H_2(\alpha_{ti_2})+\alpha_{h_1h_1}H_1(\alpha_{ti_2})+
\beta_{i_2h_1}\big(\alpha_{h_1h_1}\zero{H_1(\alpha_{th_1})}-
\\
&-\alpha_{th_1}\zero{H_1(\alpha_{h_1h_1})}-
\alpha_{th_2}\zero{H_2(\alpha_{h_1h_1})}-\alpha_{tt}\zero{T(\alpha_{h_1h_1})}-
\alpha_{ti_1}\zero{\widehat{I}_1(\alpha_{h_1h_1})}-\alpha_{ti_2}\zero{\widehat{I}_2(\alpha_{h_1h_1})}-
\\
&-
\alpha_{td}\underbrace{\widehat{D}(\alpha_{h_1h_1})}_{-\alpha_{h_1h_1}}-\alpha_{tr}\underbrace{\widehat{R}(\alpha_{h_1h_1})}_{-\alpha_{h_2h_1}}
\alpha_{tj}\zero{\widehat{J}(\alpha_{h_1h_1})}\big)+\beta_{i_2h_2}(-
\alpha_{td}\underbrace{\widehat{D}(\alpha_{h_1h_2})}_{-\alpha_{h_1h_2}}-
\alpha_{tr}\underbrace{\widehat{R}(\alpha_{h_1h_2})}_{-\alpha_{h_2h_2}}+
\\
&+\alpha_{h_1h_1}\zero{H_1(\alpha_{th_2})}-\alpha_{th_1}\zero{H_1(\alpha_{h_1h_2})}
+\alpha_{h_1h_2}\zero{H_2(\alpha_{th_2})}-\alpha_{th_2}\zero{H_2(\alpha_{h_1h_2})}-
\alpha_{tt}\zero{T(\alpha_{h_1h_2})}-
\\
&-\alpha_{ti_1}\zero{\widehat{I}_1(\alpha_{h_1h_2})}-
\alpha_{ti_2}\zero{\widehat{I}_2(\alpha_{h_1h_2})}-\alpha_{tj}\zero{\widehat{J}(\alpha_{h_1h_2})}\big)+
\beta_{i_2t}\big(4\alpha_{h_1h_1}\alpha_{th_2}+\alpha_{h_1h_1}\alpha_{tt}\Phi_1-
\\
&- 4\alpha_{h_1h_2}\alpha_{th_1}+\alpha_{h_1h_2}\alpha_{tt}\Phi_2\big),
\endaligned
\end{eqnarray*}
\begin{eqnarray*}
\footnotesize\aligned
 \kappa^{h_2t}_{i_1} &=\widehat{I}^\ast_1([\widehat{H}_2,\widehat{T}])
-\widehat{T}(\alpha_{h_2i_1})+\alpha_{h_2h_2}H_2(\alpha_{ti_1})+\alpha_{h_2h_1}H_1(\alpha_{ti_1})+
\beta_{i_1h_1}\big(\alpha_{h_2h_1}\zero{H_1(\alpha_{th_1})}-
\\
&-\alpha_{th_1}\zero{H_1(\alpha_{h_2h_1})}-
\alpha_{th_2}\zero{H_2(\alpha_{h_2h_1})}-\alpha_{tt}\zero{T(\alpha_{h_2h_1})}-\alpha_{ti_1}\zero{\widehat{I}_1(\alpha_{h_2h_1})}-
\alpha_{ti_2}\zero{\widehat{I}_2(\alpha_{h_2h_1})}-
\\
&-\alpha_{td}\underbrace{\widehat{D}(\alpha_{h_2h_1})}_{-\alpha_{h_2h_1}}-
\alpha_{tr}\underbrace{\widehat{R}(\alpha_{h_2h_1})}_{\alpha_{h_1h_1}}-
\alpha_{tj}\zero{\widehat{J}(\alpha_{h_2h_1})}+
\alpha_{h_2h_1}\zero{H_2(\alpha_{th_1})}\big)+\beta_{i_1h_2}(\alpha_{h_2h_1}\zero{H_1(\alpha_{th_2})}-
\\
&-\alpha_{th_1}\zero{H_1(\alpha_{h_2h_2})}
+\alpha_{h_2h_2}\zero{H_2(\alpha_{th_2})}-\alpha_{th_2}\zero{H_2(\alpha_{h_2h_2})}-\alpha_{tt}\zero{T(\alpha_{h_2h_2})}-
\\
 &-
\alpha_{td}\underbrace{\widehat{D}(\alpha_{h_2h_2})}_{-\alpha_{h_2h_2}}-
\alpha_{tr}\underbrace{\widehat{R}(\alpha_{h_2h_2})}_{\alpha_{h_1h_2}}-\alpha_{ti_1}\zero{\widehat{I}_1(\alpha_{h_2h_2})}-
\alpha_{ti_2}\zero{\widehat{I}_2(\alpha_{h_2h_2})}-\alpha_{tj}\zero{\widehat{J}(\alpha_{h_2h_2})}\big)+
\\
&+ \beta_{i_1t}\big(4\alpha_{h_2h_1}\alpha_{th_2}+\alpha_{h_2h_1}\alpha_{tt}\Phi_1-
4\alpha_{h_2h_2}\alpha_{th_1}+\alpha_{h_2h_2}\alpha_{tt}\Phi_2\big),
\endaligned
\end{eqnarray*}
\begin{eqnarray*}
\footnotesize\aligned
\kappa^{h_2t}_{i_2}
&
=
\widehat{I}^\ast_2([\widehat{H}_2,\widehat{T}])=-\widehat{T}(\alpha_{h_2i_2})+\alpha_{h_2h_2}H_2(\alpha_{ti_2})+\alpha_{h_2h_1}H_1(\alpha_{ti_2})+
\beta_{i_2h_1}\big(\alpha_{h_2h_1}\zero{H_1(\alpha_{th_1})}-
\\
&-\alpha_{th_1}\zero{H_1(\alpha_{h_2h_1})}-
\alpha_{th_2}\zero{H_2(\alpha_{h_2h_1})}-\alpha_{tt}\zero{T(\alpha_{h_2h_1})}-\alpha_{ti_1}\zero{\widehat{I}_1(\alpha_{h_2h_1})}-
\alpha_{ti_2}\zero{\widehat{I}_2(\alpha_{h_2h_1})}-
\\
&-\alpha_{td}\underbrace{\widehat{D}(\alpha_{h_2h_1})}_{-\alpha_{h_2h_1}}-
\alpha_{tr}\underbrace{\widehat{R}(\alpha_{h_2h_1})}_{\alpha_{h_1h_1}}-
\alpha_{tj}\zero{\widehat{J}(\alpha_{h_2h_1})}+
\alpha_{h_2h_1}\zero{H_2(\alpha_{th_1})}\big)+\beta_{i_2h_2}(\alpha_{h_2h_1}\zero{H_1(\alpha_{th_2})}-
\\
&-\alpha_{th_1}\zero{H_1(\alpha_{h_2h_2})}
+\alpha_{h_2h_2}\zero{H_2(\alpha_{th_2})}-\alpha_{th_2}\zero{H_2(\alpha_{h_2h_2})}-\alpha_{tt}\zero{T(\alpha_{h_2h_2})}-
\\
 &-
\alpha_{td}\underbrace{\widehat{D}(\alpha_{h_2h_2})}_{-\alpha_{h_2h_2}}-
\alpha_{tr}\underbrace{\widehat{R}(\alpha_{h_2h_2})}_{\alpha_{h_1h_2}}-\alpha_{ti_1}\zero{\widehat{I}_1(\alpha_{h_2h_2})}-
\alpha_{ti_2}\zero{\widehat{I}_2(\alpha_{h_2h_2})}-\alpha_{tj}\zero{\widehat{J}(\alpha_{h_2h_2})}\big)+
\\
&+ \beta_{i_2t}\big(4\alpha_{h_2h_1}\alpha_{th_2}+\alpha_{h_2h_1}\alpha_{tt}\Phi_1-
4\alpha_{h_2h_2}\alpha_{th_1}+\alpha_{h_2h_2}\alpha_{tt}\Phi_2\big),
\endaligned
\end{eqnarray*}
As for the preceding homogeneities, it is possible to simplify further
these curvature coefficients by applying the horizontal vector fields
$Y=H_1,H_2,T$ on the expressions of $\alpha_{h_1j}$, $\alpha_{h_2j}$,
$\alpha_{ti_1}$, $\alpha_{ti_2}$ already determined. For this, we
employ the following equality from elementary differential geometry:
\[
Y(fg)
=
fY(g)+gY(f), 
\ \ \ \ \ 
{\scriptstyle{f,g\in C^\infty (M)}}.
\]
Moreover, we can compute the values of the vector fields
$\widehat{H}_1$, $\widehat{H}_2$, $\widehat{T}$ on the functions
$\alpha_{{}_\bullet {}_\bullet}$ that are visible in the above curvature
coefficients. After that, we will be ready to simplify the expressions of
the curvature coefficients of homogeneity four. The obtained
expressions are a bit too long for them to appear here, 
{\em see} instead~\cite{AMSMaple}. On the other hand, inspection
of the set of functions $\alpha_{{}_\bullet {}_\bullet}$ shows that
the only undetermined one is $\alpha_{tj}$. This function appears
in the expressions of $\beta_{jt}$ ({\em see}
subsection~\ref{Dual-section}), visible in $\kappa^{h_1h_2}_j$ and
then we may annihilate at least this curvature
coefficient by choosing:
\[
\footnotesize\aligned
&
\alpha_{tj}=3a^4+3b^4-4e^2-\Phi_1a^2bc+ca\Phi_2b^2-\Phi_1ab^2d-
\Phi_2a^2bd-2\Phi_2bce-2\Phi_1ace-2\Phi_2ade+2\Phi_1bde-
\\
&-\Phi_1a^3d+\Phi_2a^3c-\Phi_1b^3c-\Phi_2b^3d+6a^2b^2+\big[{\textstyle{\frac{3}{16}}}H_1(\Phi_1)+{\textstyle{\frac{3}{16}}}H_2(\Phi_2)\big]b^2d^2+
\\
&+\big[-{\textstyle{\frac{11}{1536}}}H_2(\Phi_2)H_1(\Phi_1)-
{\textstyle{\frac{1}{192}}}H_1(H_1(\Phi_1))\Phi_1-{\textstyle{\frac{11}{3072}}}H_2({\Phi_2}^2)+{\textstyle{\frac{1}{384}}}{\Phi_2}^2H_2(\Phi_2)-
{\textstyle{\frac{11}{3072}}}H_1({\Phi_1^2})+
\\
&+{\textstyle{\frac{1}{384}}}\Phi_1^2H_1(\Phi_1)+{\textstyle{\frac{1}{48}}}H_1(H_2(H_1(\Phi_2)))+
{\textstyle{\frac{1}{384}}}H_2(H_2(H_2(\Phi_2)))+{\textstyle{\frac{1}{384}}}H_1(H_1(H_1(\Phi_1)))+
{\textstyle{\frac{1}{384}}}{\Phi_2}^2H_1(\Phi_1)-
\endaligned
\]
\[
\footnotesize\aligned
&-{\textstyle{\frac{1}{192}}}H_2(H_2(\Phi_2))\Phi_2+
{\textstyle{\frac{1}{48}}}H_2(H_1(H_1(\Phi_2)))+{\textstyle{\frac{1}{64}}}H_2(H_1(\Phi_1))\Phi_2-{\textstyle{\frac{1}{48}}}\Phi_1H_2(H_1(\Phi_2))+
{\textstyle{\frac{1}{384}}}{\Phi_1}^2H_2(\Phi_2)-
\\
&-{\textstyle{\frac{7}{384}}}H_2(H_2(H_1(\Phi_1)))+{\textstyle{\frac{1}{64}}}H_1(H_2(\Phi_2))\Phi_1-
{\textstyle{\frac{7}{384}}}H_1(H_1(H_2(\Phi_2)))-{\textstyle{\frac{1}{48}}}\Phi_2H_1(H_1(\Phi_2))\big]d^4+
\\
&+\big[-{\textstyle{\frac{11}{768}}}H_2(\Phi_2)H_1(\Phi_1)-{\textstyle{\frac{7}{192}}}H_2(H_2(H_1(\Phi_1)))+
{\textstyle{\frac{1}{192}}}H_2(H_2(H_2(\Phi_2)))+
{\textstyle{\frac{1}{192}}}H_1(H_1(H_1(\Phi_1)))+
\\
&+{\textstyle{\frac{1}{24}}}H_1(H_2(H_1(\Phi_2)))-
{\textstyle{\frac{1}{96}}}H_2(H_2(\Phi_2))\Phi_2+{\textstyle{\frac{1}{32}}}H_1(H_2(\Phi_2))\Phi_1+
{\textstyle{\frac{1}{192}}}\Phi_2^2H_1(\Phi_1)-{\textstyle{\frac{7}{192}}}H_1(H_1(H_2(\Phi_2)))+
\endaligned
\]
\[
\footnotesize\aligned
&+ {\textstyle{\frac{1}{192}}}\Phi_2^2H_2(\Phi_2)-
{\textstyle{\frac{11}{1536}}}H_1(\Phi_1^2)-{\textstyle{\frac{1}{24}}}\Phi_2H_1(H_1(\Phi_2))-
{\textstyle{\frac{11}{1536}}}H_2(\Phi_2^2)+{\textstyle{\frac{1}{32}}}H_2(H_1(\Phi_1))\Phi_2-{\textstyle{\frac{1}{96}}}H_1(H_1(\Phi_1))\Phi_1+
\\
&+
{\textstyle{\frac{1}{192}}}\Phi_1^2H_2(\Phi_2)+{\textstyle{\frac{1}{192}}}\Phi_1^2H_1(\Phi_1)-
{\textstyle{\frac{1}{24}}}\Phi_1H_2(H_1(\Phi_2))+{\textstyle{\frac{1}{24}}}H_2(H_1(H_1(\Phi_2)))\big]c^2d^2+
\big[-{\textstyle{\frac{1}{32}}}H_1(H_1(\Phi_1))+
\\
&+{\textstyle{\frac{1}{32}}}H_2(\Phi_2)\Phi_1-{\textstyle{\frac{1}{32}}}H_1(H_2(\Phi_2))+
{\textstyle{\frac{1}{32}}}H_1(\Phi_1)\Phi_1\big]bcd^2+
\big[{\textstyle{\frac{1}{32}}}H_2(H_1(\Phi_1))+{\textstyle{\frac{1}{32}}}H_2(H_2(\Phi_2))-
\\
&-{\textstyle{\frac{1}{32}}}H_2(\Phi_2)\Phi_2-
{\textstyle{\frac{1}{32}}}H_1(\Phi_1)\Phi_2\big]acd^2+
\big[-{\textstyle{\frac{1}{32}}}H_1(H_1(\Phi_1))+{\textstyle{\frac{1}{32}}}H_2(\Phi_2)\Phi_1-{\textstyle{\frac{1}{32}}}H_1(H_2(\Phi_2))+
\endaligned
\]
\[
\footnotesize\aligned
& +{\textstyle{\frac{1}{32}}}H_1(\Phi_1)\Phi_1\big]ad^3+
\big[{\textstyle{\frac{1}{32}}}H_2(H_1(\Phi_1))+{\textstyle{\frac{1}{32}}}H_2(H_2(\Phi_2))-{\textstyle{\frac{1}{32}}}H_2(\Phi_2)\Phi_2-
{\textstyle{\frac{1}{32}}}H_1(\Phi_1)\Phi_2\big]ac^3+
\\
&+ {\textstyle{\frac{3}{16}}}\big[H_1(\Phi_1)+H_2(\Phi_2)\big]a^2d^2+
{\textstyle{\frac{1}{32}}}\big[H_2(\Phi_2)\Phi_2-H_2(H_1(\Phi_1))-H_2(H_2(\Phi_2))+H_1(\Phi_1)\Phi_2\big]bd^3+
\\
&+\big[-{\textstyle{\frac{1}{32}}}H_1(H_1(\Phi_1))+
{\textstyle{\frac{1}{32}}}H_2(\Phi_2)\Phi_1-
{\textstyle{\frac{1}{32}}}H_1(H_2(\Phi_2))+{\textstyle{\frac{1}{32}}}H_1(\Phi_1)\Phi_1\big]bc^3+
\\
&+{\textstyle{\frac{3}{16}}}\big[H_1(\Phi_1)+H_2(\Phi_2)\big]a^2c^2+
{\textstyle{\frac{3}{16}}}\big[H_1(\Phi_1)+H_2(\Phi_2)\big]b^2c^2+{\textstyle{\frac{1}{32}}}\big[H_2(\Phi_2)\Phi_2-
H_2(H_1(\Phi_1))-
\endaligned
\]
\[
\footnotesize\aligned
&-
H_2(H_2(\Phi_2))+H_1(\Phi_1)\Phi_2\big]dbc^2+{\textstyle{\frac{1}{32}}}\big[-H_1(H_1(\Phi_1))+
H_2(\Phi_2)\Phi_1- H_1(H_2(\Phi_2))+H_1(\Phi_1)\Phi_1\big]ac^2d+
\\
&+\big[-{\textstyle{\frac{11}{1536}}}H_2(\Phi_2)H_1(\Phi_1)-
{\textstyle{\frac{1}{192}}}H_1(H_1(\Phi_1))\Phi_1-{\textstyle{\frac{11}{3072}}}H_2(\Phi_2^2)+
{\textstyle{\frac{1}{384}}}\Phi_2^2H_2(\Phi_2)-{\textstyle{\frac{11}{3072}}}H_1(\Phi_1^2)+
\\
&+
{\textstyle{\frac{1}{384}}}\Phi_1^2H_1(\Phi_1)+{\textstyle{\frac{1}{48}}}H_1(H_2(H_1(\Phi_2)))+{\textstyle{\frac{1}{384}}}H_2(H_2(H_2(\Phi_2)))+
{\textstyle{\frac{1}{384}}}H_1(H_1(H_1(\Phi_1)))+{\textstyle{\frac{1}{384}}}\Phi_2^2H_1(\Phi_1)-
\\
&-
{\textstyle{\frac{1}{192}}}H_2(H_2(\Phi_2))\Phi_2+{\textstyle{\frac{1}{48}}}H_2(H_1(H_1(\Phi_2)))+{\textstyle{\frac{1}{64}}}H_2(H_1(\Phi_1))\Phi_2-
{\textstyle{\frac{1}{48}}}\Phi_1H_2(H_1(\Phi_2))+{\textstyle{\frac{1}{384}}}\Phi_1^2H_2(\Phi_2)-
\\
&-{\textstyle{\frac{7}{384}}}H_2(H_2(H_1(\Phi_1)))+{\textstyle{\frac{1}{64}}}H_1(H_2(\Phi_2))\Phi_1-{\textstyle{\frac{7}{384}}}H_1(H_1(H_2(\Phi_2)))-
{\textstyle{\frac{1}{48}}}\Phi_2H_1(H_1(\Phi_2))\big]c^4.
\endaligned
\]
However, this choice does not annihilate the remaining four curvature
coefficients. By a careful examination (either by hand or with the
help of a computer), we realize that the remaining four curvature
coefficients have the following expressions:
\[
\footnotesize\aligned \kappa^{h_1t}_{i_1}&=\mathbf{\Delta_1}d^4+
(\mathbf{\Delta_2}-\mathbf{\Delta_1})c^4+\mathbf{\Delta_2}c^2d^2+\mathbf{\Delta_3}c^3d+\mathbf{\Delta_3}cd^3,
\\
 \kappa^{h_1t}_{i_2}&=\mathbf{\Delta_4}d^4+
(\mathbf{\Delta_3}+\mathbf{\Delta_4})c^4+(\mathbf{\Delta_3}+2\mathbf{\Delta_4})c^2d^2+
(2\mathbf{\Delta_1}-\mathbf{\Delta_2})c^3d+(2\mathbf{\Delta_1}-\mathbf{\Delta_2})cd^3,
\endaligned
\]
\[
\footnotesize\aligned
\kappa^{h_2t}_{i_1}&=-\kappa^{h_1t}_{i_2}-\widehat{R}(\kappa^{h_1t}_{i_1}),
\ \ \ \ \ \ \ \ \ \ \ \ \ \ \ \ \ \ \ \ \ \ \ \ \ \ \ \ \ \ \ \ \
\ \ \ \ \ \ \ \ \ \ \ \ \ \ \ \ \ \ \ \ \ \ \ \ \ \ \ \ \ \ \ \ \
\ \ \ \ \ \ \ \ \ \ \ \ \ \ \ \ \ \ \ \ \ \ \ \ \ \ \ \ \ \ \ \
\ \ \ \ \ \ \ \ \ \ \ \
\\
\kappa^{h_2t}_{i_2}&=\kappa^{h_1t}_{i_1}-\widehat{R}(\kappa^{h_1t}_{i_2}),
\endaligned
\]
where:
\[
\footnotesize\aligned
 \mathbf{\Delta_1}&:=\textstyle{\frac{1}{384}}\big[-20\Phi_2H_1(H_1(\Phi_2))-(H_1(\Phi_1))^2-
2\Phi_2^2H_1(\Phi_1)+8H_1(H_2(H_1(\Phi_2)))+ 2\Phi_1^2H_1(\Phi_1)-
\\
&-7H_1(H_1(H_2(\Phi_2)))-
4\Phi_1H_2(H_1(\Phi_2))+H_1(H_1(H_1(\Phi_1)))+\Phi_1H_1(H_2(\Phi_2))+
\\
&+23\Phi_2H_2(H_1(\Phi_1))+(H_2(\Phi_2))^2
-3\Phi_1H_1(H_1(\Phi_1))+3\Phi_2H_2(H_2(\Phi_2))- 2\Phi_2^2H_2(\Phi_2)-
\\
&-17H_2(H_2(H_1(\Phi_1)))+2\Phi_1^2H_2(\Phi_2)+16H_2(H_1(H_1(\Phi_2)))-H_2(H_2(H_2(\Phi_2)))\big],
 \endaligned
\]
\[
\footnotesize\aligned
\mathbf{\Delta_2}&:=\textstyle{\frac{1}{384}}\big[24H_1(H_2(H_1(\Phi_2)))-
24\Phi_1H_2(H_1(\Phi_2))+24\Phi_1H_1(H_2(\Phi_2))+ 24H_2(H_1(H_1(\Phi_2)))-
\\
&-24 H_1(H_1(H_2(\Phi_2)))-24\Phi_2H_1(H_1(\Phi_2))+24\Phi_2H_2(H_1(\Phi_1))-
24H_2(H_2(H_1(\Phi_1)))\big], \ \ \ \ \ \ \ \ \ \ \ \ \ \ \ \ \ \ \ \
\endaligned
\]
\[
\footnotesize\aligned
\mathbf{\Delta_3}&:={\textstyle{\frac{1}{384}}}\big[-2H_2(H_1(H_1(\Phi_1)))+8H_1(H_1(H_1(\Phi_2)))-2\Phi_1\Phi_2H_1(\Phi_1)-
8\Phi_1\Phi_2H_2(\Phi_2)-
\\
& -2H_1(H_2(H_2(\Phi_2)))- 10H_2(H_1(H_2(\Phi_2)))-16\Phi_1H_1(H_1(\Phi_2))+
8H_1(\Phi_2)H_2(\Phi_2)+
\\
&+6\Phi_1H_2(H_2(\Phi_2))+
8H_2(H_2(H_1(\Phi_2)))+22\Phi_2H_1(H_2(\Phi_2))-16\Phi_2H_2(H_1(\Phi_2))+
\\
& +22\Phi_1H_2(H_1(\Phi_1))- 10H_1(H_2(H_1(\Phi_1)))+4H_1(\Phi_1)H_1(\Phi_2)+
6\Phi_2H_1(H_1(\Phi_1))\big], \ \ \ \ \ \ \ \ \ \ \ \ \ \ \ \ \ \ \ \ \ \ \ \ \ \ \
\ \ \ \ \ \ \
\ \
\endaligned
\]
\[
\footnotesize\aligned
\mathbf{\Delta_4}&:={\textstyle{\frac{1}{384}}}\big[4\Phi_1H_1(H_1(\Phi_2))-2H_1(\Phi_2)H_2(\Phi_2)-
2H_1(\Phi_1)H_1(\Phi_2)+ 13H_2(H_1(H_2(\Phi_2)))-
\\
&-3H_1(H_2(H_2(\Phi_2)))-3\Phi_2H_1(H_1(\Phi_1))-15\Phi_2H_1(H_2(\Phi_2))+
4\Phi_1\Phi_2H_1(\Phi_1)-
\\
&-8H_2(H_2(H_1(\Phi_2)))-3H_1(H_2(H_1(\Phi_1)))+12\Phi_2H_2(H_1(\Phi_2))-
3\Phi_1H_2(H_2(\Phi_2))-
\\
&-7\Phi_1H_2(H_1(\Phi_1))+ 4\Phi_1\Phi_2H_2(\Phi_2)+5H_2(H_1(H_1(\Phi_1)))\big].
\endaligned
\]

\begin{Lemma}
One in fact has, identically as functions of $(x, y, u)$:
\[
\boxed{
0
\equiv
\Delta_2}
\ \ \ \ \ \ \ \ \ \
\text{\rm and}
\ \ \ \ \ \ \ \ \ \
\boxed{0
\equiv
\Delta_3+2\,\Delta_4}\,.
\]
\end{Lemma}

\proof
These two nontrivial relations were already prepared 
in advance, {\em cf.} the Corollary~\ref{0-Delta2-Delta3-2-Delta4}.
\endproof

Furthermore, by taking account of the relations listed in
Proposition~\ref{relations-Phi-H-I-V} and of $H_1 ( \Phi_2) = H_2 (
\Phi_1)$, one sees that the expressions of the two remaining functions
$\mathbf{\Delta_1}$ and $\mathbf{\Delta_4}$ of $(x, y, u)$ can be
given better, completely symmetric forms, as is stated by the
following summarizing proposition.

\begin{Proposition}
The four remaining curvature coefficients of homogeneity
$h = 4$ express explicitly as follows:
\[
\footnotesize
\aligned
\kappa_{i_1}^{h_1t}
&
=
-\,\mathbf{\Delta_1}\,c^4
-
2\,\mathbf{\Delta_4}\,c^3d
-
2\,\mathbf{\Delta_4}\,cd^3
+
\mathbf{\Delta_1}\,d^4,
\\
\kappa_{i_2}^{h_1t}
&
=
-\,\mathbf{\Delta_4}\,c^4
+
2\,\mathbf{\Delta_1}\,c^3d
+
2\,\mathbf{\Delta_1}\,cd^3
+
\mathbf{\Delta_4}\,d^4,
\\
\kappa_{i_1}^{h_2t}
&
=
\kappa_{i_2}^{h_1t},
\\
\kappa_{i_2}^{h_2t}
&
=
-\,\kappa_{i_1}^{h_1t},
\endaligned
\]
where the two functions $\mathbf{ \Delta_1}$ and $\mathbf{ \Delta_4}$
of the three horizontal variables $(x, y, u)$ have the following
explicit expressions:
\[
\boxed{
\footnotesize\aligned
\mathbf{\Delta_1}
&
=
{\textstyle{\frac{1}{384}}}
\Big[
H_1(H_1(H_1(\Phi_1)))
-
H_2(H_2(H_2(\Phi_2)))
+
11\,H_1(H_2(H_1(\Phi_2)))
-
11\,H_2(H_1(H_2(\Phi_1)))
+
\\
&
\ \ \ \ \ \ \ \ \ \ \ \ \
+6\,\Phi_2\,H_2(H_1(\Phi_1))
-
6\,\Phi_1\,H_1(H_2(\Phi_2))
-
3\,\Phi_2\,H_1(H_1(\Phi_2))
+
3\,\Phi_1\,H_2(H_2(\Phi_1))
-
\\
&
\ \ \ \ \ \ \ \ \ \ \ \ \
-\,3\,\Phi_1\,H_1(H_1(\Phi_1))
+
3\,\Phi_2\,H_2(H_2(\Phi_2))
-
2\,\Phi_1\,H_1(\Phi_1)
+
2\,\Phi_2\,H_2(\Phi_2)
-
\\
&
\ \ \ \ \ \ \ \ \ \ \ \ \
-\,2\,(\Phi_2)^2\,H_1(\Phi_1)
+
2\,(\Phi_1)^2\,H_2(\Phi_2)
-
2\,(\Phi_2)^2\,H_2(\Phi_2)
+
2\,(\Phi_1)^2\,H_1(\Phi_1)
\Big],
\\
\mathbf{\Delta_4}
&
=
{\textstyle{\frac{1}{384}}}
\Big[
-\,3\,H_2(H_1(H_2(\Phi_2)))
-
3\,H_1(H_2(H_1(\Phi_1)))
+
5\,H_1(H_2(H_2(\Phi_2)))
+
5\,H_2(H_1(H_1(\Phi_1)))
+
\\
&
\ \ \ \ \ \ \ \ \ \ \ \ \
+4\,\Phi_1\,H_1(H_1(\Phi_2))
+
4\,\Phi_2\,H_2(H_1(\Phi_2))
-
3\,\Phi_2\,H_1(H_1(\Phi_1))
-
3\,\Phi_1\,H_2(H_2(\Phi_2))
-
\\
&
\ \ \ \ \ \ \ \ \ \ \ \ \
-\,7\,\Phi_2\,H_1(H_2(\Phi_2))
-
7\,\Phi_1\,H_2(H_1(\Phi_1))
-
2\,H_1(\Phi_1)\,H_1(\Phi_2)
-
2\,H_2(\Phi_2)\,H_2(\Phi_1)
+
\\
&
\ \ \ \ \ \ \ \ \ \ \ \ \
+4\,\Phi_1\Phi_2\,H_1(\Phi_1)
+
4\,\Phi_1\Phi_2\,H_2(\Phi_2)
\Big].
\endaligned}
\]
\end{Proposition}

\proof
As said, one uses the relations listed in
Proposition~\ref{relations-Phi-H-I-V} until formal expressions
show up symmetries.
\endproof

In the next subsection, we will establish that the obtained
Cartan connection is actually normal. Up to know, all the functions
$\alpha_{{}_\bullet {}_\bullet}$ are determined, 
but still there is last
function of type $\delta_{{}_\bullet}$, namely
$\delta_{18}$ which is
yet undetermined. To determine it 
it suffices to equate the above expression of $\alpha_{ tj}$ to the
corresponding one in Lemma~\ref{fiber-type},
after making the possible simplification:
\[
\footnotesize
 \left[ \aligned
\alpha_{tj}&=(c^4+2c^2d^2+d^4)\delta_{18}+6a^2b^2-
2\Phi_2bce-2\Phi_1ace-\Phi_1a^2bc+\Phi_2ab^2c-\Phi_2a^2bd- 2\Phi_2ade-
\\
&-\Phi_1ab^2d+2\Phi_1bde+3a^4+3b^4-4e^2-\Phi_1b^3c+
\Phi_2a^3c-\Phi_1a^3d-\Phi_2b^3d+{\textstyle{\frac{3}{16}}}\big[H_1(\Phi_1)+H_2(\Phi_2)\big]a^2c^2+
\\
&+ {\textstyle{\frac{1}{32}}}\big[-H_1(\Phi_1)\Phi_2+H_2(H_1(\Phi_1))+
H_2(H_2(\Phi_2))-H_2(\Phi_2)\Phi_2\big]ac^3+
{\textstyle{\frac{3}{16}}}\big[H_1(\Phi_1)+H_2(\Phi_2)\big]b^2c^2+
\\
&+ {\textstyle{\frac{1}{32}}}\big[H_1(\Phi_1)\Phi_1-H_1(H_1(\Phi_1))+
H_2(\Phi_2)\Phi_1-H_1(H_2(\Phi_2))\big]bc^3+
{\textstyle{\frac{3}{16}}}\big[H_1(\Phi_1)+H_2(\Phi_2)\big]b^2d^2+
\\
&+ {\textstyle{\frac{3}{16}}}\big[H_1(\Phi_1)+H_2(\Phi_2)\big]a^2d^2+
{\textstyle{\frac{1}{32}}}\big[H_1(\Phi_1)\Phi_2-H_2(H_2(\Phi_2))+
H_2(\Phi_2)\Phi_2-H_2(H_1(\Phi_1))\big]bd^3+
\\
&+ {\textstyle{\frac{1}{32}}}\big[H_1(\Phi_1)\Phi_1-H_1(H_1(\Phi_1))+
H_2(\Phi_2)\Phi_1-H_1(H_2(\Phi_2))\big]ad^3+
{\textstyle{\frac{1}{32}}}\big[H_1(\Phi_1)\Phi_2-H_2(H_2(\Phi_2))+
\\
&+H_2(\Phi_2)\Phi_2-
H_2(H_1(\Phi_1))\big]bc^2d+{\textstyle{\frac{1}{32}}}\big[H_1(\Phi_1)\Phi_1-
H_1(H_1(\Phi_1))+H_2(\Phi_2)\Phi_1-H_1(H_2(\Phi_2))\big]ac^2d+
\\
&+{\textstyle{\frac{1}{32}}}\big[-H_1(\Phi_1)\Phi_2+ H_2(H_1(\Phi_1))+H_2(H_2(\Phi_2))-
H_2(\Phi_2)\Phi_2\big]acd^2+{\textstyle{\frac{1}{32}}}\big[H_1(\Phi_1)\Phi_1-
\\
&- H_1(H_1(\Phi_1))+H_2(\Phi_2)\Phi_1-H_1(H_2(\Phi_2))\big]bcd^2,
\endaligned\right.
\]
By identification, 
this equation holds when one makes the following assignment:
\[
\footnotesize\aligned
\delta_{18}&={\textstyle{\frac{1}{64}}}\Phi_2H_2(H_1(\Phi_1))-{\textstyle{\frac{11}{1536}}}H_2(\Phi_2)H_1(\Phi_1)-
{\textstyle{\frac{1}{192}}}\Phi_1H_1(H_1(\Phi_1))-
{\textstyle{\frac{11}{3072}}}H_1(\Phi_1^2)+{\textstyle{\frac{1}{48}}}H_1(H_2(H_1(\Phi_2)))-
\\
&-
{\textstyle{\frac{7}{384}}}H_1(H_1(H_2(\Phi_2)))+{\textstyle{\frac{1}{384}}}\Phi_2^2H_2(\Phi_2)-
{\textstyle{\frac{1}{48}}}\Phi_2H_1(H_1(\Phi_2))-{\textstyle{\frac{1}{192}}}\Phi_2H_2(H_2(\Phi_2))+
{\textstyle{\frac{1}{64}}}\Phi_1H_1(H_2(\Phi_2))+
\\
&+{\textstyle{\frac{1}{384}}}\Phi_1^2H_1(\Phi_1)-
{\textstyle{\frac{11}{3072}}}H_2(\Phi_2^2)+{\textstyle{\frac{1}{384}}}\Phi_2^2H_1(\Phi_1)-{\textstyle{\frac{1}{48}}}\Phi_1H_2(H_1(\Phi_2))+
{\textstyle{\frac{1}{48}}}H_2(H_1(H_1(\Phi_2)))+
\\
&+{\textstyle{\frac{1}{384}}}\Phi_1^2H_2(\Phi_2)+{\textstyle{\frac{1}{384}}}H_2(H_2(H_2(\Phi_2)))+
{\textstyle{\frac{1}{384}}}H_1(H_1(H_1(\Phi_1)))-{\textstyle{\frac{7}{384}}}H_2(H_2(H_1(\Phi_1))).
\endaligned
\]

\subsection{Homogeneity 5}
\label{Hom5}
In this homogeneity, it is possible to express curvature coefficients
in terms of the curvature coefficients of lower homogeneities. More
precisely, consider the restricted
$h$-homogeneous differential operators:
\[
\partial_{[h]}
\colon\ \ \
{\mathcal
C}_{[h]}^2(\frak g_-,\frak g)
\rightarrow{\mathcal C}_{[h]}^3(\frak g_-,\frak g).
\]
We will use the graded Bianchi-Tanaka identities
of Proposition~\ref{Bianchi-Tanaka-graded} to
identify $\partial \kappa^{(5)}$ in terms of the lower components of
$\kappa$.

In homogeneity $5$, 
we encounter two curvature coefficients $\kappa^{h_1t}_j$ and
$\kappa^{h_2t}_j$ and according to the notations introduced in 
Section~\ref{Cohomology-section}, we have:

\begin{equation}
\label{k5} \aligned 
\kappa_{[5]}
=
\kappa^{h_1t}_j\,{\sf h}_1^\ast\wedge {\sf
t}^\ast\otimes {\sf j} +
\kappa^{h_2t}_j\,{\sf h}_2^\ast\wedge {\sf t}^\ast\otimes {\sf j}.
\endaligned
\end{equation}
Applying the differential gives:
\begin{equation}
\label{partial-def}
\footnotesize\aligned
\partial
\kappa_{[5]}({\sf h}_1,{\sf h}_2,{\sf t})
&
=
\big[{\sf
h}_1,\underbrace{\kappa_{[5]}({\sf h}_2,
{\sf t})}_{\kappa^{h_2t}_j{\sf j}}\big]-[{\sf h}_2,
\underbrace{\kappa_{[5]}({\sf h}_1,{\sf
t})}_{\kappa^{h_1t}_j{\sf j}}]
+
[{\sf t},\zero{\kappa_{[5]}({\sf h}_1,{\sf h}_2)}]-
\\
&-[\kappa_{[5]}(\underbrace{[{\sf h}_1,{\sf h}_2]}_{\sf t},{\sf
t})]+[\kappa_{[5]}(\zero{[{\sf
h}_1,{\sf t}]},{\sf h_2})]
-
[\kappa_{[5]}(\zero{[{\sf h}_2,{\sf t}]},{\sf h}_1)]
\\
&=\kappa^{h_2t}_j\,{\sf i}_1-\kappa^{h_1t}_j\,{\sf i}_2
+0+0+0+0.
\endaligned
\end{equation}

On the other hand, the graded Bianchi-Tanaka
identities of Proposition~\ref{Bianchi-Tanaka-graded} assert that:
\[
\footnotesize\aligned
\partial 
\kappa_{[5]}({\sf h}_1,{\sf h}_2,{\sf t})
&
=
\underbrace{-\,\sum_{j=1}^4\,
\kappa_{[5-j]}
\Big(
{\sf proj}_{\mathfrak{g}_-}\big(
\kappa_{[j]}
({\sf h}_1,{\sf h}_2)\big),
{\sf t}\Big)
-
\big(\widehat{T}\zero{\kappa_{[3]}}\big)
({\sf h}_1,{\sf h}_2)
}_{-\,P_1}
-
\\
& 
\underbrace{
-\,\sum_{j=1}^4\,
\kappa_{[5-j]}
\Big(
{\sf proj}_{\mathfrak{g}_-}\big(
\kappa_{[j]}
({\sf t},{\sf h}_1)\big),
{\sf h}_2)
\Big)
-
\big(\widehat{H}_2\kappa_{[4]}\big)
({\sf t},{\sf h}_1)
}_{-\,P_2}
-
\\
&
\underbrace{
-\,\sum_{j=1}^4\,
\kappa_{[5-j]}
\Big(
{\sf proj}_{\mathfrak{g}_-}\big(
\kappa_{[j]}
({\sf h}_2,{\sf t})\big),
{\sf h}_1)
\Big)
-
\big(\widehat{H}_1\kappa_{[4]}\big)
({\sf h}_2,{\sf t})
}_{-\,P_3}.
\endaligned
\]
Let us compute for example the last term $P_3$,
taking account of the 
vanishing of the curvature components $\kappa_{[1]}$,
$\kappa_{[2]}$ and
$\kappa_{[3]}$. At first we compute 
its $\sum_{ j=1}^4\,$ part:
\[
\footnotesize\aligned
\boxed{{\scriptstyle{j=1}}} 
\ \ \ \ \ \ \
\kappa_{[4]}
\Big(
\zero{
{\sf proj}_{\mathfrak{g}_-}\big(
\kappa_{[1]}({\sf h}_2,{\sf t})\big)},\,
{\sf h}_1
\Big)
&
=
0,
\\
\boxed{{\scriptstyle{j=2}}} 
\ \ \ \ \ \ \
\kappa_{[3]}
\Big(
\zero{
{\sf proj}_{\mathfrak{g}_-}\big(
\kappa_{[2]}({\sf h}_2,{\sf t})\big)},\,
{\sf h}_1
\Big)
&
=
0,
\\
\boxed{{\scriptstyle{j=3}}} 
\ \ \ \ \ \ \
\kappa_{[2]}
\Big(
\zero{
{\sf proj}_{\mathfrak{g}_-}\big(
\kappa_{[3]}({\sf h}_2,{\sf t})\big)},\,
{\sf h}_1
\Big)
&
=
0,
\\
\boxed{{\scriptstyle{j=4}}} \ \ \ \ \ \ \ 
\zero{\kappa_{[1]}}
\Big(
{\sf proj}_{\mathfrak{g}_-}
\big(\kappa_{[4]}({\sf h}_2,{\sf t})\big),
{\sf h}_1
\Big)
&
=
0.
\endaligned
\]
So the $\sum_{ j=1}^4\,$ part of 
$P_3$ is zero. For the remaining part, reminding 
that:
\begin{equation}
\label{k4}
\aligned
\kappa_{[4]}
=
\kappa^{h_1t}_{i_1}\,
{\sf h}_1^\ast\wedge{\sf t}^\ast\otimes {\sf i}_1
+ 
\kappa^{h_2t}_{i_1}\,{\sf h}_2^\ast\wedge
{\sf t}^\ast\otimes {\sf i}_1+\kappa^{h_1t}_{i_2}\,
{\sf h}_1^\ast\wedge {\sf
t}^\ast\otimes {\sf i}_2
+ 
\kappa^{h_2t}_{i_2}\,{\sf h}_2^\ast\wedge {\sf t}^\ast\otimes {\sf i}_2,
\endaligned
\end{equation}
it is clear that:
\[
\kappa_{[4]}({\sf h}_2,{\sf t})
=
\kappa^{h_2t}_{i_1}\,{\sf i}_1
+
\kappa^{h_2t}_{i_2}\,{\sf i}_2,
\]
whence:
\[
\big(\widehat{H}_1\kappa_{[4]}\big)
({\sf h}_2,{\sf t})
=
\widehat{H}_1\big(\kappa^{h_2t}_{i_1}\big)\,
{\sf i}_1
+
\widehat{H}_1\big(\kappa^{h_2t}_{i_2}\big)\,
{\sf i}_2,
\]
which is the expression of $P_3$. Similar computations provide:
\[
\aligned 
P_1
=
0,
\ \ \ \ \ \
P_2
=
-{\widehat H}_2\big(\kappa^{h_1t}_{i_1}\big)\,
{\sf i}_1
-
{\widehat H}_2\big(\kappa^{h_1t}_{i_2}\big)\,{\sf i}_2,
\endaligned
\]
and consequently:
\begin{equation}
\label{Bianchi}
\aligned
\partial
\kappa_{[5]}
({\sf h}_1,{\sf h}_2,{\sf t})
&
=
-\,\big(P_1+P_2+P_3\big)
\\
&
=
\big(\widehat{H}_2\big(\kappa^{h_1t}_{i_1}\big)
-
\widehat{H}_1\big(\kappa^{h_2t}_{i_1}\big)\big)\,
{\sf i}_1
+
\big(\widehat{H}_2\big(\kappa^{h_1t}_{i_2}\big)
-
\widehat{H}_1\big(\kappa^{h_2t}_{i_2}\big)\big)\,
{\sf i}_2.
\endaligned
\end{equation}
Now comparison of~\thetag{ \ref{partial-def}} 
and~\thetag{ \ref{Bianchi}} implies
that:
\[
\aligned
\kappa^{h_1t}_j
&
=
\widehat{H}_1\big(\kappa^{h_2t}_{i_2}\big)
-
\widehat{H}_2\big(\kappa^{h_1t}_{i_2}\big),
\\
\kappa^{h_2t}_j
&
=
-\widehat{H}_1\big(\kappa^{h_2t}_{i_1}\big)
+
\widehat{H}_2\big(\kappa^{h_1t}_{i_1}\big).
\endaligned
\]

\subsection{Conclusion}
\label{conclusion}
A review of the results obtained so far
shows that the only non-zero curvature
coefficients are:
\[
\boxed{{\scriptstyle{\text{\rm Hom}\,4}}} 
\ \ \ \ 
\kappa^{h_1t}_{i_1}, \ \ \
\kappa^{h_1t}_{i_2}, \ \ \
\kappa^{h_2t}_{i_1}, \ \ \ 
\kappa^{h_2t}_{i_2};
\]
\[
\boxed{{\scriptstyle{\text{\rm Hom}\,5}}} 
\ \ \ \ 
\kappa^{h_1t}_j, \ \ \ \kappa^{h_2t}_j. \
\ \ \ \ \ \ \ \
\ \ \ \ \ \ \ \ \ \ \ \ \
\]
All these curvature coefficients can be expressed as the
combinations of $\kappa^{h_1t}_{i_1}$ and $\kappa^{h_1t}_{i_2}$ and
the values of the constant vector fields on them. These two curvature
coefficients are called {\em essential curvatures}. A Cartan geometry
is homogeneous if and only if all of its essential curvatures
vanish (\cite{ Sharpe, EMS}). Hence a consequence
of our results is the following

\begin{Theorem}
The Cartan geometry associated to any strongly pseudoconvex
deformation $M^3 \subset \C^2$ the Heisenberg sphere $\mathbb{ H}^3
\subset \C^2$ 
having curvature function equal to:
\begin{equation}
\label{kappa} 
\aligned
\kappa
&
=
\kappa_{[4]}+\kappa_{[5]}
\\
&
=
\kappa^{h_1t}_{i_1}\,{\sf h}_1^\ast\wedge{\sf t}^\ast\otimes{\sf
i}_1+\kappa^{h_1t}_{i_2}\,{\sf
h}_1^\ast\wedge{\sf t}^\ast\otimes {\sf i}_2
+
\kappa^{h_2t}_{i_1}\,{\sf
h}_2^\ast\wedge{\sf
t}^\ast\otimes {\sf i}_1+
\\
&
+
\kappa^{h_2t}_{i_2}\,{\sf h}_2^\ast\wedge{\sf t}^\ast\otimes {\sf
i}_2+\kappa^{h_1t}_j\,{\sf
h}_1^\ast\wedge{\sf t}^\ast\otimes{\sf j}
+
\kappa^{h_2t}_j\,{\sf h}_2^\ast\wedge{\sf
t}^\ast\otimes
{\sf j},
\endaligned
\end{equation}
with:
\[
\aligned
\kappa_{i_1}^{h_2t}
&
=
\kappa_{i_2}^{h_1t},
\\
\kappa_{i_2}^{h_2t}
&
=
-\,\kappa_{i_1}^{h_1t},
\\
\kappa^{h_1t}_j
&
=
\widehat{H}_1\big(\kappa^{h_2t}_{i_2}\big)
-
\widehat{H}_2\big(\kappa^{h_1t}_{i_2}\big),
\\
\kappa^{h_2t}_j
&
=
-\widehat{H}_1\big(\kappa^{h_2t}_{i_1}\big)
+
\widehat{H}_2\big(\kappa^{h_1t}_{i_1}\big),
\endaligned
\]
is locally homogeneous if and only if its two essential
curvatures:
\[
\footnotesize
\aligned
\kappa_{i_1}^{h_1t}
&
=
-\,\mathbf{\Delta_1}\,c^4
-
2\,\mathbf{\Delta_4}\,c^3d
-
2\,\mathbf{\Delta_4}\,cd^3
+
\mathbf{\Delta_1}\,d^4,
\\
\kappa_{i_2}^{h_1t}
&
=
-\,\mathbf{\Delta_4}\,c^4
+
2\,\mathbf{\Delta_1}\,c^3d
+
2\,\mathbf{\Delta_1}\,cd^3
+
\mathbf{\Delta_4}\,d^4
\endaligned
\]
vanish identically; equivalently, 
the following two
{\em explicit} functions $\mathbf{ \Delta_1}$ and $\mathbf{ \Delta_4}$
of only the three horizontal variables $(x, y, u)$:
\[
\footnotesize\aligned
\mathbf{\Delta_1}
&
=
{\textstyle{\frac{1}{384}}}
\Big[
H_1(H_1(H_1(\Phi_1)))
-
H_2(H_2(H_2(\Phi_2)))
+
11\,H_1(H_2(H_1(\Phi_2)))
-
11\,H_2(H_1(H_2(\Phi_1)))
+
\\
&
\ \ \ \ \ \ \ \ \ \ \ \ \
+6\,\Phi_2\,H_2(H_1(\Phi_1))
-
6\,\Phi_1\,H_1(H_2(\Phi_2))
-
3\,\Phi_2\,H_1(H_1(\Phi_2))
+
3\,\Phi_1\,H_2(H_2(\Phi_1))
-
\\
&
\ \ \ \ \ \ \ \ \ \ \ \ \
-\,3\,\Phi_1\,H_1(H_1(\Phi_1))
+
3\,\Phi_2\,H_2(H_2(\Phi_2))
-
2\,\Phi_1\,H_1(\Phi_1)
+
2\,\Phi_2\,H_2(\Phi_2)
-
\\
&
\ \ \ \ \ \ \ \ \ \ \ \ \
-\,2\,(\Phi_2)^2\,H_1(\Phi_1)
+
2\,(\Phi_1)^2\,H_2(\Phi_2)
-
2\,(\Phi_2)^2\,H_2(\Phi_2)
+
2\,(\Phi_1)^2\,H_1(\Phi_1)
\Big],
\\
\mathbf{\Delta_4}
&
=
{\textstyle{\frac{1}{384}}}
\Big[
-\,3\,H_2(H_1(H_2(\Phi_2)))
-
3\,H_1(H_2(H_1(\Phi_1)))
+
5\,H_1(H_2(H_2(\Phi_2)))
+
5\,H_2(H_1(H_1(\Phi_1)))
+
\\
&
\ \ \ \ \ \ \ \ \ \ \ \ \
+4\,\Phi_1\,H_1(H_1(\Phi_2))
+
4\,\Phi_2\,H_2(H_1(\Phi_2))
-
3\,\Phi_2\,H_1(H_1(\Phi_1))
-
3\,\Phi_1\,H_2(H_2(\Phi_2))
-
\\
&
\ \ \ \ \ \ \ \ \ \ \ \ \
-\,7\,\Phi_2\,H_1(H_2(\Phi_2))
-
7\,\Phi_1\,H_2(H_1(\Phi_1))
-
2\,H_1(\Phi_1)\,H_1(\Phi_2)
-
2\,H_2(\Phi_2)\,H_2(\Phi_1)
+
\\
&
\ \ \ \ \ \ \ \ \ \ \ \ \
+4\,\Phi_1\Phi_2\,H_1(\Phi_1)
+
4\,\Phi_1\Phi_2\,H_2(\Phi_2)
\Big]
\endaligned
\]
vanish identically.
\qed
\end{Theorem}

Up to now, we have achieved 
all the necessary computations. We saw that
the regularity condition~\textbf{(c4)} is satisfied automatically and
condition~\textbf{(c1)} holds by applying Lemma~\ref{fiber-type} in
the computations of various homogeneities. The only remaining task is
to verify that both
the isomorphism condition~\textbf{(c2)} and the normality
condition~\textbf{(c3)} hold. 
We inspect these two final conditions via the
following propositions, respectively:

\begin{Proposition}
For any element $p = (a,b,c,d,e,x,y,u)$ of $\mathcal G$, 
the $\mathfrak{ g}$-valued map $\omega_p
\colon T_p{ \mathcal P} \longrightarrow\frak g$ is an isomorphism.
\end{Proposition}

\proof
According to
the expressions of $\widehat{T},\ldots,\widehat{J}$ as 
components of $\omega_p^{-1}$, the matrix corresponding to 
$\omega_p^{-1}$ is:
\[
\left(%
\begin{array}{cccccccc}
\alpha_{tt} & \alpha_{th_1} & \alpha_{th_2} & \alpha_{td} & \alpha_{tr} &
\alpha_{ti_1} & \alpha_{ti_2} & \alpha_{tj} \\
0 & \alpha_{h_1h_1} & \alpha_{h_1h_2} & \alpha_{h_1d} & \alpha_{h_1r} &
\alpha_{h_1i_1} & \alpha_{h_1i_2} & \alpha_{h_1j} \\
0 & \alpha_{h_2h_1} & \alpha_{h_2h_2} & \alpha_{h_2d} & \alpha_{h_2r} &
\alpha_{h_2i_1} & \alpha_{h_2i_2} & \alpha_{h_2j} \\
0 & 0 & 0 & 1 & 0 & 0 & 0 & 0 \\
0 & 0 & 0 & 0 & 1 & 0 & 0 & 0 \\
0 & 0 & 0 & 0 & 0 & 1 & 0 & 0 \\
0 & 0 & 0 & 0 & 0 & 0 & 1 & 0 \\
0 & 0 & 0 & 0 & 0 & 0 & 0 & 1 \\
\end{array}%
\right),
\]
and it is has determinant:
\[
\alpha_{tt}\big(\alpha_{h_1h_1}
\alpha_{h_2h_2}-\alpha_{h_2h_1}\alpha_{h_1h_2}\big)=(c^2+d^2)^2,
\]
which is nonzero by assumption.
\endproof

\begin{Proposition}
The Cartan connection constructed in the
preceding paragraphs in a completely effective
way is normal.
\end{Proposition}

\proof
According to ~\thetag{\ref{kappa}}, the ${\sf t}$-, ${\sf h}_1$-,
${\sf h}_2$-, ${\sf d}$- and ${\sf r}$-components of the Cartan
curvature $\kappa$ vanish together. Vanishing of its ${\sf t},{\sf
h}_1$ and ${\sf h}_2$-components means that this curvature is torsion
free. Moreover, the ${\sf d}$- and ${\sf r}$-components of $\kappa$
constitute its $\frak g_0$-component and consequently $\kappa_{[ 0]}
\equiv 0$ by construction. Therefore the Cartan connection is normal
according to Definition 1.6.7 page 128 of \cite{Slovak}.
\endproof

\section{General Formulas for the Second Cohomology 
of Graded Lie Algebras}
\label{second-cohomology-graded}

\HEAD{General Formulas for the Second Cohomology of Graded Lie Algebras}{
Mansour Aghasi, Joël Merker, and Masoud Sabzevari}

Throughout cohomology considerations, the ground field $\K$ will be a
commutative field of characteristic zero, while in most expected
applications to exterior differential systems, $\K$ will be either
$\Q$, $\R$ or $\C$.

\subsection{Arbitrary abstract Lie algebra}
Let $\mathfrak{ g}$ be a Lie algebra over $\K$ of dimension $r
\geqslant 2$ containing a Lie subalgebra $\mathfrak{ g}_- \subset
\mathfrak{ g}$ of dimension $n$ with $1 \leqslant n \leqslant r - 1$,
so that $[ \mathfrak{ g}_-, \mathfrak{ g}_- ]_{\mathfrak{g}} \subset
\mathfrak{ g}_-$. Let ${\sf x}_1, \dots, {\sf x}_n$ be an arbitrary
but fixed basis of $\mathfrak{ g}_-$ which is completed by means of
vectors ${\sf x}_{ n+1}, \dots, {\sf x}_r$ to produce a basis of
$\mathfrak{ g}$. To any such pair of bases are associated the
so-called {\sl structure constants} $c_{ k_1,k_2}^s \in \K$ encoding
the Lie bracket:
\begin{equation}
\label{general-structure-g}
\aligned
&
\big[{\sf x}_{k_1},\,{\sf x}_{k_2}\big]_{\mathfrak{g}}
=
\sum_{s=1}^r\,c_{k_1,k_2}^s\,{\sf x}_s
\\
&
\ \ \ \ \ \
{\scriptstyle{(k_1,\,k_2\,=\,1\,\cdots\,n,\,n\,+\,1,\dots,r)}},
\endaligned
\end{equation}
and because $[ \mathfrak{ g}_-, \mathfrak{ g}_- ]_{\mathfrak{g}}
\subset \mathfrak{ g}_-$, we must naturally have $c_{ k_1,k_2}^s = 0$
for $s = n+1, \dots, r$ whenever $1 \leqslant k_1, k_2 \leqslant n =
\dim \mathfrak{ g}_-$. Furthermore, one may adopt the convention that
$c_{ k_1, k_2}^s = 0$ whenever one does not have $1 \leqslant k_1, k_2
\leqslant r$ and $1 \leqslant s \leqslant r$. Of course, the
skew-symmetry and the Jacobi identity:
\[
\aligned
0
&
=
[{\sf x}_{k_1},{\sf x}_{k_2}]_{\mathfrak{g}}
+
[{\sf x}_{k_1},{\sf x}_{k_2}]_{\mathfrak{g}}
\\
0
&
=
\big[
[{\sf x}_{k_1},{\sf x}_{k_2}]_{\mathfrak{g}},
{\sf x}_{k_3}\big]_{\mathfrak{g}}
+
\big[
[{\sf x}_{k_3},{\sf x}_{k_1}]_{\mathfrak{g}},
{\sf x}_{k_2}\big]_{\mathfrak{g}}
+
\big[
[{\sf x}_{k_2},{\sf x}_{k_3}]_{\mathfrak{g}},
{\sf x}_{k_1}\big]_{\mathfrak{g}}
\endaligned
\]
read at the level of structure constants as:
\[
\aligned
0
&
=
c_{k_1,k_2}^s
+
c_{k_2,k_1}^s
\\
0
&
=
\sum_{s=1}^r\,
\big(
c_{k_1,k_2}^s\,c_{s,k_3}^l
+
c_{k_3,k_1}^s\,c_{s,k_2}^l
+
c_{k_2,k_3}^s\,c_{s,k_1}^l
\big)
\\
&
\ \ \ \ \ \ \ \ \ \ \ \ \ \ \ \ \ \ \ \ \ \ \
{\scriptstyle{(k_1,\,\,k_2,\,\,k_3,\,\,l\,=\,1\,\cdots\,r)}}.
\endaligned
\]

\subsection{Exterior algebra}
Given any integer $\ell \geqslant 1$, consider the $\ell$-th exterior power
$\Lambda^\ell \mathfrak{ g}_-$. Whenever $\ell \leqslant n$, it is a
nonzero vector space generated over $\K$ by the basis consisting of
the $\binom{ n}{ \ell} = \frac{ n!}{ \ell!\,(n-\ell)!}$ linearly independent
$\ell$-fold exterior products:
\[
\big(
{\sf x}_{j_1}\wedge {\sf x}_{j_2}
\wedge\cdots\wedge
{\sf x}_{j_\ell}
\big)_{1\leqslant j_1<j_2<\cdots<j_\ell\leqslant n}, 
\]
while $\Lambda^\ell \mathfrak{ g}_- = \{ 0\}$ for all $\ell \geqslant
n+1$. Next, let $\mathfrak{ g}^* = {\rm Lin} ( \mathfrak{ g}, \K)$
denote the dual of the Lie algebra $\mathfrak{ g}$, viewed as a plain
vector space (it has no natural Lie bracket structure). If 
we introduce
the basis: 
\[
{\sf x}_1^*,\dots,{\sf x}_n^*,{\sf x}_{ n+1}^*,\dots,{\sf x}_r^*
\] 
of $\mathfrak{ g}^*$ which is dual to the previously fixed basis ${\sf
x}_1, \dots, {\sf x}_n, {\sf x}_{ n+1}, \dots, {\sf x}_r$ of
$\mathfrak{ g}$, then by definition, for any $i, j = 1, \dots, n, n+1,
\dots, r$, we have:
\[
{\sf x}_i^*({\sf x}_j)
=
\delta_j^i
:=
\left\{
\aligned
&
1\ \ 
\text{\rm if}\ \
i=j,
\\
&
0\ \ 
\text{\rm otherwise}.
\endaligned\right. 
\]
For any $\ell \geqslant 1$, let us define 
({\em cf.}~\cite{Goze}, Chap.~3) 
the space $\mathcal{ C}^\ell ( \mathfrak{ g}_-,
\mathfrak{ g})$ of {\sl $\ell$-cochains} as the space of linear maps from
$\Lambda^k \mathfrak{ g}_-$ to $\mathfrak{ g}$, that is to say:
\[
\mathcal{C}^\ell(\mathfrak{g}_-,\mathfrak{g})
\overset{\text{\rm def}}{=}
{\rm Lin}\big(\Lambda^\ell\mathfrak{g}_-,\,\mathfrak{g}\big). 
\]
Thanks to the canonical identifications: 
\[
\aligned
{\rm Lin}\big(\Lambda^\ell\mathfrak{g}_-,\,\mathfrak{g}\big)
&
\simeq
\big(\Lambda^\ell\mathfrak{g}_-\big)^*\otimes\mathfrak{g}
\\
&
\simeq
\Lambda^\ell\mathfrak{g}_-^*
\otimes
\mathfrak{g},
\endaligned 
\]
valid for any $\ell$, an arbitrary $\ell$-cochain writes, in term of bases
for $\mathfrak{ g}_-^*$ and for $\mathfrak{ g}$, under the general
form:
\[
\Phi
=
\sum_{1\leqslant i_1<i_2<\cdots<i_\ell\leqslant n}\,
\sum_{l=1}^r\,\phi_{i_1,i_2,\dots,i_\ell}^k\,
\big({\sf x}_{i_1}^*\wedge{\sf x}_{i_2}^*
\wedge\cdots\wedge
{\sf x}_{i_\ell}^*\big)
\otimes
{\sf x}_k,
\] 
where the $\phi_{ i_1, i_2, \dots, i_\ell}^k$ are coefficients in the
ground field $\K$.

On the other hand, without any reference to bases, we recall from
basic algebra that a $\ell$-cochain $\Phi \in \mathcal{ C}^\ell ( \mathfrak{
g}_-, \mathfrak{ g})$ may be seen either as being a $\mathfrak{
g}$-valued linear map acting on exterior $\ell$-vectors ${\sf z}_1 \wedge
{\sf z}_2 \wedge \cdots \wedge {\sf z}_\ell$ with the ${\sf z}_i$
belonging to $\mathfrak{ g}_-$, or equivalently as being a {\em
multilinear} map from the $\ell$-fold product $\mathfrak{ g}_- \times
\mathfrak{ g}_- \times \cdots \times \mathfrak{ g}_-$ to $\mathfrak{
g}$ which has the property that:
\[
\Phi\big({\sf z}_{\sigma(1)},{\sf z}_{\sigma(2)},
\dots,
{\sf z}_{\sigma(\ell)}\big)
=
(-1)^{{\rm sgn}(\sigma)}\,
\Phi\big({\sf z}_1,{\sf z}_2,
\dots,
{\sf z}_\ell\big)
\]
for every permutation $\sigma$ of $\{ 1, 2, \dots, \ell\}$. 
This last property is easily seen to be equivalent to the property
that the value of $\Phi$ vanishes as soon as two at least of
its arguments coincide.

\subsection{Differential complex and cohomology} From $\ell$-cochains 
to $(\ell+1)$-cochains, there is a canonical {\sl boundary operator}: 
\[
\partial^\ell
\colon\ \ \ 
\mathcal{C}^\ell\big(\mathfrak{g}_-,\,\mathfrak{g}\big)
\longrightarrow
\mathcal{C}^{\ell+1}\big(\mathfrak{g}_-,\,\mathfrak{g}\big).
\]
which, to any $\Phi \in \mathcal{ C}^\ell \big( \mathfrak{ g}_-,
\mathfrak{ g}\big)$ associates a $(\ell+1)$-cochain $\partial^\ell 
\Phi$ whose
action on any collection of $\ell+1$ vectors ${\sf z}_0, {\sf z}_1,
\dots, {\sf z}_\ell$ of $\mathfrak{ g}_-$ is defined by the specific
formula:
\begin{equation}
\label{cochain-definition}
\small
\aligned
(\partial^\ell\Phi)
\big({\sf z}_0,{\sf z}_1,\dots,{\sf z}_\ell\big)
&
:=
\sum_{i=0}^\ell\,(-1)^i
\big[{\sf z}_i,\,\Phi({\sf z}_0,
\dots,
\widehat{\sf z}_i,\dots,{\sf z}_\ell)\big]_{\mathfrak{g}}
+
\\
&
\ \ \ \ \
+
\sum_{0\leqslant i<j\leqslant \ell}\,
(-1)^{i+j}\,
\Phi\big([{\sf z}_i,{\sf z}_j]_{\mathfrak{g}},
{\sf z}_0,\dots,\widehat{\sf z}_i,
\dots,
\widehat{\sf z}_j,\dots,{\sf z}_n\big),
\endaligned
\end{equation}
where, as usual, $\widehat{\sf z}_l$ means removal of the term ${\sf
z}_l$; first, second and third cohomology spaces ({\em see} below)
associated to this specific differential $\partial$ occur naturally as
providing one-dimensional extensions of Lie algebras ($H^1$), as
parametrizing infinitesimal deformations of Lie algebras ($H^2$), as
obstruction to their deformations ($H^3$), or as the algebraic
skeleton of curvatures of Cartan connections on principal bundles
($H^2$ and Bianchi-type identities). Let us check that $\partial^\ell
\Phi$ really belongs to $\mathcal{ C}^{ \ell+1} (\mathfrak{ g}_-,
\mathfrak{ g})$.

The so-defined action of $\partial^\ell \Phi$ is clearly linear with
respect to each argument ${\sf z}_i$, for $i= 0, 1, \dots, \ell$.
Furthermore, from the assumption that $\Phi$ vanishes when two of its
arguments coincide, one immediately infers that $(\partial^\ell \Phi)
({\sf z}_0, {\sf z}_1, \dots, {\sf z}_\ell) = 0$ vanishes as soon as
${\sf z}_{ i_1} = {\sf z}_{ i_2}$ for at least two distinct indices
$i_1 \neq i_2$, whence by a standard elementary reasoning, we have the
skew-symmetry:
\[
(\partial^\ell\Phi)\big({\sf z}_{\sigma(0)},{\sf z}_{\sigma(1)},\dots,
{\sf z}_{\sigma(\ell)}\big)
=
(-1)^{{\rm sgn}(\sigma)}\,
(\partial^\ell\Phi)
\big({\sf z}_0,{\sf z}_1,\dots,{\sf z}_\ell\big), 
\]
for every permutation $\sigma$ of $\{ 0, 1, \dots, \ell\}$. Consequently,
$(\partial^\ell \Phi)$ effectively identifies to a certain linear map:
\[
\partial^\ell\Phi\colon\ \ \
\Lambda^{\ell+1}\mathfrak{g}_-
\longrightarrow
\mathfrak{g}
\]
from the $( \ell+1)$-th exterior product of $\mathfrak{ g}_-$ into
$\mathfrak{ g}$, namely it truly is a $(\ell+1)$-cochain. Precisely, for
any element ${\sf z}_0 \wedge {\sf z}_1 \wedge \cdots \wedge {\sf z}_\ell
\in \Lambda^{ \ell+1} \mathfrak{ g}_-$ with ${\sf z}_i \in \mathfrak{
g}_-$, one simply sets:
\[
(\partial^\ell\Phi)\big({\sf z}_0\wedge {\sf z}_1
\wedge\cdots\wedge 
{\sf z}_\ell\big)
:=
(\partial^\ell\Phi)\big({\sf z}_0,{\sf z}_1,\dots,{\sf z}_\ell).
\]

On the other hand, it is usually left as an exercise to verify that
from any level $\ell$ to the level $\ell+2$, one has
$\partial^{\ell+1} ( \partial^\ell \Phi) \equiv 0$, so that the datum:
\begin{equation}
\label{complex-partial}
0
\overset{\partial^0}{\longrightarrow}
\mathcal{C}^1
\overset{\partial^1}{\longrightarrow}
\mathcal{C}^2
\overset{\partial^2}{\longrightarrow}
\cdots
\overset{\partial^{n-2}}{\longrightarrow}
\mathcal{C}^{n-1}
\overset{\partial^{n-1}}{\longrightarrow}
\mathcal{C}^n
\overset{\partial^n}{\longrightarrow}
0
\end{equation}
forms what is called a {\sl complex}, namely the composition 
$\partial^{ \ell + 1}
\circ \partial^\ell = 0$ from any $\mathcal{ C}^\ell$ to $\mathcal{ C}^{
\ell+2}$ always vanishes. Equivalently, one has:
\[
{\rm im}\big(\partial^{\ell-1}
\colon
\mathcal{C}^{\ell-1}\to\mathcal{C}^\ell\big)
\subset
{\rm ker}\big(\partial^\ell
\colon
\mathcal{C}^\ell\to\mathcal{C}^{\ell+1}\big),
\]
and the classical terminology is to call:
\[
\mathcal{Z}^\ell(\mathfrak{g}_-,\mathfrak{g}\big)
:=
{\rm ker}\big(\partial^\ell
\colon
\mathcal{C}^\ell\to\mathcal{C}^{\ell+1}\big)
\]
the {\sl space of cocycles} of order $\ell$, and also to call:
\[
\mathcal{B}^\ell(\mathfrak{g}_-,\mathfrak{g}\big)
:=
{\rm im}\big(
\partial^{\ell-1}
\colon\mathcal{C}^{\ell-1}\to\mathcal{C}^\ell\big)
\] 
the {\sl space of coboundaries} of order $\ell$, which is thus always
a vector subspace of $\mathcal{ Z}^\ell( \mathfrak{ g}_-, \mathfrak{
g} \big)$.

\begin{Definition}
The quotient space:
\[
H^\ell\big(\mathfrak{g}_-,\mathfrak{g}\big)
:=
\frac{\mathcal{Z}^\ell(\mathfrak{g}_-,\mathfrak{g}\big)}{
\mathcal{B}^\ell(\mathfrak{g}_-,\mathfrak{g}\big)}
\]
is called the $\ell$-th cohomology space of $\mathfrak{ g}_-$
in $\mathfrak{ g}$.
\end{Definition}

For applications to either deformations of Lie algebras or to the
explicit constructions of Cartan connections, we will mainly be
interested in computing the second cohomology:
\[
H^2\big(\mathfrak{g}_-,\mathfrak{g}\big)
=
\frac{\mathcal{Z}^2(\mathfrak{g}_-,\mathfrak{g}\big)}{
\mathcal{B}^2(\mathfrak{g}_-,\mathfrak{g}\big)},
\]
which is a plain finite-dimensional vector space over $\K$, the
complexity of which will depend on the geometric situation under
study.

\subsection{Basis for $2$-cochains}
As we already saw above, a $2$-cochain writes, in term of bases, under
the general form (we shall from now on regularly omit the parentheses
in $({\sf x}_{ i_1}^* \wedge {\sf x}_{ i_2}^*) \otimes {\sf x}_k$):
\[
\Phi
=
\sum_{1\leqslant i_1<i_2\leqslant n}\,
\sum_{k=1}^r\,
\phi_{i_1,i_2}^k\,
{\sf x}_{i_1}^*\wedge {\sf x}_{i_2}^*\otimes {\sf x}_k,
\]
where the $\phi_{ i_1, i_2}^k \in \K$ are arbitrary 
coefficients, namely it is a linear combination of
the $\frac{ n ( n-1)}{ 2} \, r$ basic elements:
\[
{\sf x}_{i_1}^*\wedge {\sf x}_{i_2}^*\otimes {\sf x}_k
\ \ \ \ \ \ \ \ \ \ \ \ \
{\scriptstyle{(1\,\leqslant\,i_1\,<\,i_2\,\leqslant\,n\,;\,\,
k\,=\,1,\,\dots,\,n,\,n+1,\,\dots,\,r)}},
\]
which visibly form a basis of $\Lambda^2 \mathfrak{ g}_-^*
\otimes \mathfrak{ g}$. We remind that if $E$ is a finite-dimensional
$\K$-vector space and if $\omega^*$, $\pi^*$ are one-forms
belonging to its dual $E^* = 
{\rm Lin} ( E, \K)$, then the two-form $\omega^* \wedge \pi^*$
acts on pairs $(e, f) \in E^2$ by definition as:
\[
\omega^*\wedge\pi^*(e,f)
\overset{def}{=}
\omega^*(e)\,\pi^*(f)-\omega^*(f)\,\pi^*(e).
\] 
In particular, for any $i_1$, $i_2$ with 
$i_1 < i_2$ and for any $j_1$, $j_2$ without restriction, we have: 
\begin{equation}
\label{i-1-i-2-j-1-j-2-bis}
\aligned
{\sf x}_{i_1}^*\wedge{\sf x}_{i_2}^*
({\sf x}_{j_1},{\sf x}_{j_2})
&
=
{\sf x}_{i_1}^*({\sf x}_{j_1})\,
{\sf x}_{i_2}^*({\sf x}_{j_2})
-
{\sf x}_{i_1}^*({\sf x}_{j_2})\,
{\sf x}_{i_2}^*({\sf x}_{j_1})
\\
&
=
\delta_{j_1}^{i_1}\,\delta_{j_2}^{i_2}
-
\delta_{j_2}^{i_1}\,\delta_{j_1}^{i_2}.
\endaligned
\end{equation}
However, we observe {\em passim} that the second product of Kronecker
deltas necessarily vanishes whenever $j_1 < j_2$, a natural
restriction we will sometimes make, though not always.

\subsection{Boundary of a basic $2$-cochain}
According to the definition~\thetag{ \ref{cochain-definition}}, for
any triple of indices $j_1, j_2, j_3$ with $1 \leqslant j_1 < j_2 <
j_3 \leqslant n$, we have:
\[
\small
\aligned
(\partial^2\Phi)\big({\sf x}_{j_1}\wedge 
{\sf x}_{j_2}\wedge {\sf x}_{j_3}\big)
&
=
(\partial^2\Phi)\big({\sf x}_{j_1},{\sf x}_{j_2},{\sf x}_{j_3}\big)
\\
&
=
\big[{\sf x}_{j_1},\,\Phi({\sf x}_{j_2},{\sf x}_{j_3})\big]_{\mathfrak{g}}
-
\big[{\sf x}_{j_2},\,\Phi({\sf x}_{j_1},{\sf x}_{j_3})\big]_{\mathfrak{g}}
+
\big[{\sf x}_{j_3},\,\Phi({\sf x}_{j_1},{\sf x}_{j_2})\big]_{\mathfrak{g}}
-
\\
&
\ \ \ \ \ 
-
\Phi\big([{\sf x}_{j_1},\,{\sf x}_{j_2}]_{\mathfrak{g}},\,{\sf x}_{j_3}\big)
+
\Phi\big([{\sf x}_{j_1},\,{\sf x}_{j_3}]_{\mathfrak{g}},\,{\sf x}_{j_2}\big)
-
\Phi\big([{\sf x}_{j_2},\,{\sf x}_{j_3}]_{\mathfrak{g}},\,{\sf x}_{j_1}\big).
\endaligned
\] 
Let us hence apply this formula to any basic $2$-cochain $\Phi_k^{
i_1, i_2} = {\sf x}_{ i_1}^* \wedge {\sf x}_{ i_2}^* \otimes {\sf
x}_k$ and perform a few natural computational transformations,
using the Lie algebra structure~\thetag{ \ref{general-structure-g}} and 
applying formulas~\thetag{ \ref{i-1-i-2-j-1-j-2-bis}}:
\[
\scriptsize
\aligned
&
\partial^2({\sf x}_{i_1}^*\wedge{\sf x}_{i_2}^*\otimes{\sf x}_k)
\big({\sf x}_{j_1}\wedge {\sf x}_{j_2}\wedge {\sf x}_{j_3}\big)
=
\\
&
=
\big[{\sf x}_{j_1},\,
({\sf x}_{i_1}^*\wedge{\sf x}_{i_2}^*\otimes{\sf x}_k)
({\sf x}_{j_2},{\sf x}_{j_3})\big]_{\mathfrak{g}}
-
\big[{\sf x}_{j_2},\,
({\sf x}_{i_1}^*\wedge{\sf x}_{i_2}^*\otimes{\sf x}_k)
({\sf x}_{j_1},{\sf x}_{j_3})\big]_{\mathfrak{g}}
+
\big[{\sf x}_{j_3},\,
({\sf x}_{i_1}^*\wedge{\sf x}_{i_2}^*\otimes{\sf x}_k)
({\sf x}_{j_1},{\sf x}_{j_2})\big]_{\mathfrak{g}}
-
\\
&
\ \ \ \ \ 
-
({\sf x}_{i_1}^*\wedge{\sf x}_{i_2}^*\otimes{\sf x}_k)
\big([{\sf x}_{j_1},\,{\sf x}_{j_2}]_{\mathfrak{g}},\,{\sf x}_{j_3}\big)
+
({\sf x}_{i_1}^*\wedge{\sf x}_{i_2}^*\otimes{\sf x}_k)
\big([{\sf x}_{j_1},\,{\sf x}_{j_3}]_{\mathfrak{g}},\,{\sf x}_{j_2}\big)
-
({\sf x}_{i_1}^*\wedge{\sf x}_{i_2}^*\otimes{\sf x}_k)
\big([{\sf x}_{j_2},\,{\sf x}_{j_3}]_{\mathfrak{g}},\,{\sf x}_{j_1}\big)
\\
&
=
\big[{\sf x}_{j_1},\,
\delta_{j_2}^{i_1}\delta_{j_3}^{i_2}\,{\sf x}_k\big]_{\mathfrak{g}}
-
\big[{\sf x}_{j_2},\,
\delta_{j_1}^{i_1}\delta_{j_3}^{i_2}\,{\sf x}_k\big]_{\mathfrak{g}}
+
\big[{\sf x}_{j_3},\,
\delta_{j_1}^{i_1}\delta_{j_2}^{i_2}\,{\sf x}_k\big]_{\mathfrak{g}}
-
\\
&
\ \ \ \ \
-
{\sf x}_{i_1}^*\wedge{\sf x}_{i_2}^*\otimes{\sf x}_k
\bigg(
\sum_{l=1}^r\,c_{j_1,j_2}^l\,{\sf x}_l,\,\,{\sf x}_{j_3}
\bigg)
+
{\sf x}_{i_1}^*\wedge{\sf x}_{i_2}^*\otimes{\sf x}_k
\bigg(
\sum_{l=1}^r\,c_{j_1,j_3}^l\,{\sf x}_l,\,\,{\sf x}_{j_2}
\bigg)
-
{\sf x}_{i_1}^*\wedge{\sf x}_{i_2}^*\otimes{\sf x}_k
\bigg(
\sum_{l=1}^r\,c_{j_2,j_3}^l\,{\sf x}_l,\,\,{\sf x}_{j_1}
\bigg)
\\
&
=
\delta_{j_2}^{i_1}\delta_{j_3}^{i_2}\,
[{\sf x}_{j_1},{\sf x}_k]_{\mathfrak{g}}
-
\delta_{j_1}^{i_1}\delta_{j_3}^{i_2}\,
[{\sf x}_{j_2},{\sf x}_k]_{\mathfrak{g}}
+
\delta_{j_1}^{i_1}\delta_{j_2}^{i_2}\,
[{\sf x}_{j_3},{\sf x}_k]_{\mathfrak{g}}
-
\\
&
\ \ \ \ \ \ 
-
\bigg(
\sum_{l=1}^r\,c_{j_1,j_2}^l\,
\big(
\delta_l^{i_1}\delta_{j_3}^{i_2}-\delta_l^{i_2}\delta_{j_3}^{i_1}
\big)
+
\sum_{l=1}^r\,c_{j_1,j_3}^l\,
\big(
\delta_l^{i_1}\delta_{j_2}^{i_2}-\delta_l^{i_2}\delta_{j_2}^{i_1}
\big)
-
\sum_{l=1}^r\,c_{j_2,j_3}^l\,
\big(
\delta_l^{i_1}\delta_{j_1}^{i_2}-\delta_l^{i_2}\delta_{j_1}^{i_1}
\big)
\bigg)\,{\sf x}_k
\\
&
=
\sum_{l=1}^r\,
\Big(
c_{j_1,k}^l\,\delta_{j_2}^{i_1}\delta_{j_3}^{i_2}
-
c_{j_2,k}^l\,\delta_{j_1}^{i_1}\delta_{j_3}^{i_2}
+
c_{j_3,k}^l\,\delta_{j_1}^{i_1}\delta_{j_2}^{i_2}
\Big)\,{\sf x}_l
+
\\
&
\ \ \ \ \
+
\big(
-c_{j_1,j_2}^{i_1}\,\delta_{j_3}^{i_2}
+
c_{j_1,j_2}^{i_2}\,\delta_{j_3}^{i_1}
+
c_{j_1,j_3}^{i_1}\,\delta_{j_2}^{i_2}
-
c_{j_1,j_3}^{i_2}\,\delta_{j_2}^{i_1}
-
c_{j_2,j_3}^{i_1}\,\delta_{j_1}^{i_2}
+
c_{j_2,j_3}^{i_2}\,\delta_{j_1}^{i_1}
\big)\,{\sf x}_k.
\endaligned
\]
At this point, in order to reach a neat formula, let us replace ${\sf
x}_k$ in the last line by $\sum_{ l=1}^r\, \delta_k^l\, {\sf x}_l$ and
reorganize:
\[
\footnotesize
\aligned
(\partial^2{\sf x}_{i_1}^*\wedge{\sf x}_{i_2}^*\otimes{\sf x}_k)
&
\big({\sf x}_{j_1}\wedge {\sf x}_{j_2}\wedge {\sf x}_{j_3}\big)
=
\sum_{l=1}^r\,
\Big(
c_{j_1,k}^l\,\delta_{j_2}^{i_1}\delta_{j_3}^{i_2}
-
c_{j_2,k}^l\,\delta_{j_1}^{i_1}\delta_{j_3}^{i_2}
+
c_{j_3,k}^l\,\delta_{j_1}^{i_1}\delta_{j_2}^{i_2}
+
\\
&
\ \ \ \ \ \ \ \ \ \ \ \
+
\delta_k^l
\big(
-c_{j_1,j_2}^{i_1}\,\delta_{j_3}^{i_2}
+
c_{j_1,j_2}^{i_2}\,\delta_{j_3}^{i_1}
+
c_{j_1,j_3}^{i_1}\,\delta_{j_2}^{i_2}
-
c_{j_1,j_3}^{i_2}\,\delta_{j_2}^{i_1}
-
c_{j_2,j_3}^{i_1}\,\delta_{j_1}^{i_2}
+
c_{j_2,j_3}^{i_2}\,\delta_{j_1}^{i_1}
\big)
\Big)\,{\sf x}_l.
\endaligned
\]
Since we may {\em a priori} represent without arguments
this $3$-cochain $\partial^2
({\sf x}_{i_1}^* \wedge{\sf x}_{i_2}^* \otimes{\sf x}_k)$ as:
\[
\partial^2({\sf x}_{i_1}^*\wedge{\sf x}_{i_2}^*\otimes{\sf x}_k)
=
\sum_{1\leqslant j_1<j_2<j_3\leqslant n}\,
\sum_{l=1}^r\,
{\sf coefficient}_{k;j_1,j_2,j_3}^{i_1,i_2;l}
\cdot
{\sf x}_{j_1}^*\wedge {\sf x}_{j_2}^*\wedge {\sf x}_{j_3}^*
\otimes{\sf x}_l, 
\]
the preceding computations precisely give us that: 
\[
\footnotesize
\aligned
{\sf coefficient}_{k;j_1,j_2,j_3}^{i_1,i_2;l}
&
=
c_{j_1,k}^l\,\delta_{j_2}^{i_1}\delta_{j_3}^{i_2}
-
c_{j_2,k}^l\,\delta_{j_1}^{i_1}\delta_{j_3}^{i_2}
+
c_{j_3,k}^l\,\delta_{j_1}^{i_1}\delta_{j_2}^{i_2}
-
\\
&
\ \ \ \ \
-
c_{j_1,j_2}^{i_1}\,\delta_k^l\delta_{j_3}^{i_2}
+
c_{j_1,j_2}^{i_2}\,\delta_k^l\delta_{j_3}^{i_1}
+
c_{j_1,j_3}^{i_1}\,\delta_k^l\delta_{j_2}^{i_2}
-
c_{j_1,j_3}^{i_2}\,\delta_k^l\delta_{j_2}^{i_1}
-
c_{j_2,j_3}^{i_1}\,\delta_k^l\delta_{j_1}^{i_2}
+
c_{j_2,j_3}^{i_2}\,\delta_k^l\delta_{j_1}^{i_1}.
\endaligned
\]
As a result, we obtain explicitly that: 
\[
\footnotesize
\aligned
\partial^2\big({\sf x}_{i_1}^*\wedge{\sf x}_{i_2}^*
\otimes{\sf x}_k\big)
&
=
\sum_{1\leqslant j_1<j_2<j_3\leqslant n}\,
\sum_{l=1}^r\,
{\sf x}_{j_1}^*\wedge {\sf x}_{j_2}^*\wedge {\sf x}_{j_3}^*
\otimes {\sf x}_l
\Big(
c_{j_1,k}^l\,\delta_{j_2}^{i_1}\delta_{j_3}^{i_2}
-
c_{j_2,k}^l\,\delta_{j_1}^{i_1}\delta_{j_3}^{i_2}
+
c_{j_3,k}^l\,\delta_{j_1}^{i_1}\delta_{j_2}^{i_2}
-
\\
&
\ \ \ \ \
-
\delta_k^l
\big[
c_{j_1,j_2}^{i_1}\,\delta_{j_3}^{i_2}
+
c_{j_1,j_2}^{i_2}\,\delta_{j_3}^{i_1}
+
c_{j_1,j_3}^{i_1}\,\delta_{j_2}^{i_2}
-
c_{j_1,j_3}^{i_2}\,\delta_{j_2}^{i_1}
-
c_{j_2,j_3}^{i_1}\,\delta_{j_1}^{i_2}
+
c_{j_2,j_3}^{i_2}\,\delta_{j_1}^{i_1}
\big]
\Big)
\\
&
\ \ \ \ \ \ \ \ \ \ \ \ \ \ \ \ \ \ 
{\scriptstyle{(1\,\leqslant\,i_1\,<\,i_2\,\leqslant\,n\,;\,\,
k\,=\,1,\,\dots,\,n,\,n+1,\,\dots\,r)}}.
\endaligned
\]

\subsection{Boundary of a general $2$-cochain}
Next, with:
\[
\Phi 
= 
\sum_{1\leqslant i_1<i_2\leqslant n}\, 
\sum_{k=1}^r\, 
\phi_{i_1,i_2}^k\, 
{\sf x}_{i_1}^*\wedge{\sf x}_{ i_2}^*
\otimes {\sf x}_k,
\] 
we may now compute $\partial^2 \Phi$ by linearity: 
\[
\scriptsize
\aligned
\partial^2\Phi
&
=
\sum_{1\leqslant i_1<i_2\leqslant n}\,
\sum_{k=1}^r\,\phi_{i_1,i_2}^k\,
\partial^2\big(
{\sf x}_{i_1}^*\wedge{\sf x}_{i_2}^*\otimes{\sf x}_k
\big)
\\
&
=
\sum_{1\leqslant i_1<i_2\leqslant n}\,
\sum_{k=1}^r\,
\phi_{i_1,i_2}^k\,
\sum_{1\leqslant j_1<j_2<j_3\leqslant n}\,
\sum_{l=1}^r\,
{\sf x}_{j_1}^*\wedge {\sf x}_{j_2}^*\wedge{\sf x}_{j_3}^*
\otimes {\sf x}_l
\Big(
c_{j_1,k}^l\,\delta_{j_2}^{i_1}\delta_{j_3}^{i_2}
-
c_{j_2,k}^l\,\delta_{j_1}^{i_1}\delta_{j_3}^{i_2}
+
c_{j_3,k}^l\,\delta_{j_1}^{i_1}\delta_{j_2}^{i_2}
-
\\
&
\ \ \ \ \
-
\delta_k^l
\big[
c_{j_1,j_2}^{i_1}\,\delta_{j_3}^{i_2}
+
c_{j_1,j_2}^{i_2}\,\delta_{j_3}^{i_1}
+
c_{j_1,j_3}^{i_1}\,\delta_{j_2}^{i_2}
-
c_{j_1,j_3}^{i_2}\,\delta_{j_2}^{i_1}
-
c_{j_2,j_3}^{i_1}\,\delta_{j_1}^{i_2}
+
c_{j_2,j_3}^{i_2}\,\delta_{j_1}^{i_1}
\big]
\Big)
\\
&
=
\sum_{1\leqslant j_1<j_2<j_3\leqslant n}\,
\sum_{l=1}^r\,
{\sf x}_{j_1}^*\wedge {\sf x}_{j_2}^*\wedge{\sf x}_{j_3}^*
\otimes {\sf x}_l
\bigg(
\sum_{k=1}^r\,
\Big(
c_{j_1,k}^l\,\phi_{j_2,j_3}^k
-
c_{j_2,k}^l\,\phi_{j_1,j_3}^k
+
c_{j_3,k}^l\,\phi_{j_1,j_2}
\Big)
+
\\
&
\ \ \ \ \
+
\sum_{1\leqslant i_1<i_2\leqslant n}\,
\phi_{i_1,i_2}^l
\big(
-c_{j_1,j_2}^{i_1}\,\delta_{j_3}^{i_2}
+
c_{j_1,j_2}^{i_2}\,\delta_{j_3}^{i_1}
\big)
+
\sum_{1\leqslant i_1<i_2\leqslant n}\,
\phi_{i_1,i_2}^l
\big(
c_{j_1,j_3}^{i_1}\,\delta_{j_2}^{i_2}
-
c_{j_1,j_3}^{i_2}\,\delta_{j_2}^{i_1}
\big)
+
\\
&
\ \ \ \ \ 
+
\sum_{1\leqslant i_1<i_2\leqslant n}\,
\phi_{i_1,i_2}^l
\big(
-c_{j_2,j_3}^{i_1}\,\delta_{j_1}^{i_2}
+
c_{j_2,j_3}^{i_2}\,\delta_{j_1}^{i_1}
\big).
\endaligned
\]
At this point, we must finish the computation of the three sums
appearing in the last two lines. In fact, any general triangle-like 
sum of the form:
\[
\sum_{1\leqslant i_1<i_2\leqslant n}\,
\mu_{i_1,i_2}\,
\big(
-\nu^{i_1}\,\delta_{j_3}^{i_2}
+
\nu^{i_2}\,\delta_{j_3}^{i_1}
\big),
\]
where the $\mu_{ \cdot, \cdot}$, $\nu^{\cdot}$ are indexed
numbers, has the property that its general term within 
parentheses is zero unless $i_2 = j_3$ or
$i_1 = j_3$, whence it decomposes symbolically just as
two simple sums:
\[
\sum_{1\leqslant i_1<i_2\leqslant n}\,
=
\sum_{i_1=1}^{j_3-1}\,
\bigg\vert_{i_2=j_3}
+
\sum_{i_2=j_3+1}^n\,
\bigg\vert_{i_1=j_3},
\]
so that the sum in question expands as:
\[
-\sum_{i_1=1}^{j_3-1}\,
\mu_{i_1,j_3}\,\nu^{i_1}
+
\sum_{i_2=j_3+1}^n\,\mu_{j_3,i_2}\,\nu^{i_2}.
\]
Applying this formula, we may finish the computation of the
three sums mentioned above and they are equal to:
\[
\footnotesize
\aligned
&
-\sum_{i_1=1}^{j_3-1}\,
c_{j_1,j_2}^{i_1}\,\phi_{j_1,j_3}^l
+
\sum_{i_2=j_3+1}^n\,c_{j_1,j_2}^{i_2}\,\phi_{j_3,i_2}^l
+
\sum_{i_1=1}^{j_2-1}\,c_{j_1,j_3}^{i_1}\,\phi_{i_1,j_2}^l
-
\sum_{i_2=j_2+1}^n\,c_{j_1,j_3}^{i_2}\,\phi_{j_2,i_2}^l
-
\\
&
-
\sum_{i_1=1}^{j_1-1}\,
c_{j_2,j_3}^{i_1}\,\phi_{i_1,j_1}^l
+
\sum_{i_2=j_1+1}^n\,
c_{j_2,j_3}^{i_2}\,\phi_{j_1,i_2}.
\endaligned
\]
As a result, we may explicitly characterize the condition that $\Phi$
be a {\sl $2$-cocycle}, stating the initial hypothesis for
self-contentness reasons.

\begin{Proposition}
\label{characterization-2-cocycles}
Let $\mathfrak{ g}$ be an $r$-dimensional Lie algebra over $\K$, let
$\mathfrak{ g}_- \subset \mathfrak{ g}$ be a proper $n$-dimensional
Lie subalgebra with $2 \leqslant n \leqslant r - 1$, and let ${\sf
x}_1, \dots, {\sf x}_n, {\sf x}_{ n+1}, \dots, {\sf x}_r$ be a basis
of $\mathfrak{ g}$, its first $n$ terms ${\sf x}_1, \dots, {\sf x}_n$
simultaneously constituting a basis of $\mathfrak{ g}_-$ so that the
${\sf x}_{ i_1}^* \wedge {\sf x}_{ i_2}^* \otimes {\sf x}_k$ with $1
\leqslant i_1 < i_2 \leqslant n$ and $k = 1, \dots, n, n+1, \dots, r$
make a naturally associated basis for $\Lambda^2 \mathfrak{ g}_-^*
\otimes \mathfrak{ g}$. Then a general $2$-cochain in $\Lambda^2
\mathfrak{ g}_-^* \otimes \mathfrak{ g}$:
\[
\Phi 
= 
\sum_{1\leqslant i_1<i_2\leqslant n}\, 
\sum_{k=1}^r\,
\phi_{i_1, i_2}^k\,
{\sf x}_{ i_1}^*\wedge{\sf x}_{ i_2}^* 
\otimes{\sf x}_k
\]
having arbitrary undetermined coefficients $\phi_{ i_1, i_2}^k \in \K$
is a cocycle, namely satisfies $\partial^2 \Phi = 0$, if and only if the
following $r \, \binom{ n}{ 3}$ linear equations hold:
\[
\footnotesize
\aligned
0
&
=
\sum_{k=1}^r\,
\Big(
c_{j_1,k}^l\,\phi_{j_2,j_3}^k
-
c_{j_2,k}^l\,\phi_{j_1,j_3}^k
+
c_{j_3,k}^l\,\phi_{j_1,j_2}^k
\Big)
-
\\
&
\ \ \ \ \
-\sum_{i_1=1}^{j_3-1}\,
c_{j_1,j_2}^{i_1}\,\phi_{i_1,j_3}^l
+
\sum_{i_2=j_3+1}^n\,c_{j_1,j_2}^{i_2}\,\phi_{j_3,i_2}^l
+
\sum_{i_1=1}^{j_2-1}\,c_{j_1,j_3}^{i_1}\,\phi_{i_1,j_2}^l
-
\\
&
\ \ \ \ \
-
\sum_{i_2=j_2+1}^n\,
c_{j_1,j_3}^{i_2}\,\phi_{j_2,i_2}^l
-
\sum_{i_1=1}^{j_1-1}\,
c_{j_2,j_3}^{i_1}\,\phi_{i_1,j_1}^l
+
\sum_{i_2=j_1+1}^n\,
c_{j_2,j_3}^{i_2}\,\phi_{j_1,i_2}^l
\ \ \
\\
&
\ \ \ \ \ \ \ \ \ \ \ \ \ \ \ \ \ \ \ \ \ \ \ \ \ \ \ \ \ 
{\scriptstyle{(1\,\leqslant\,j_1\,<\,j_2\,<\,j_3\,\leqslant\,n\,;\,\,
l\,=\,1,\,\dots\,n,\,n+1,\,\dots,\,r)}},
\endaligned
\]
where the $c_{ j,k}^s$ are the structure constants:
$[ {\sf x}_j, \, {\sf x}_k ]_{\mathfrak{g}} = 
\sum_{ s=1}^r\, c_{ j,k}^s\, {\sf x}_s$. 
\end{Proposition}

\subsection{Basis for $1$-cochains}
Now, we want to characterize, in terms of the structure constants of
$\mathfrak{ g}$, the condition that the $2$-cochain $\Phi$ identifies
to the differential $\partial^1 \Psi$ of a $1$-cochain $\Psi \in
\mathcal{ C}^1 ( \mathfrak{ g}_-, \, \mathfrak{ g})$. Let us
therefore write down such a general $1$-cochain: 
\[
\Psi
=
\sum_{1\leqslant i\leqslant n}\,
\sum_{k=1}^r\,\psi_k^i\,\,{\sf x}_i^*\otimes {\sf x}_k,
\]
as being the linear combination, with arbitrary coefficients
$\psi_k^i \in \K$, of the $n \cdot n$ basic $1$-cochains:
\[
{\sf x}_i^*\otimes {\sf x}_k
\ \ \ \ \ \ \ \ \ \ \ \ \ 
{\scriptstyle{(i,\,\,j\,=\,1\,\cdots\,n)}}
\]
which clearly form a basis of $\mathcal{ C}^1 ( \mathfrak{ g}_-,
\mathfrak{ g}) = \mathfrak{ g}_-^* \otimes \mathfrak{ g}$ over
$\K$. 

\subsection{Boundary of a basic $1$-cochain}
Applying the definition~\thetag{ \ref{cochain-definition}}, the
differential $\partial^1$ acts as follows on such a general $1$-cochain:
\[
(\partial^1\Psi)\big({\sf x}_{j_1},{\sf x}_{j_2}\big)
=
\big[{\sf x}_{j_1},\,\Psi({\sf x}_{j_2})\big]_{\mathfrak{g}}
-
\big[{\sf x}_{j_2},\,\Psi({\sf x}_{j_1})\big]_{\mathfrak{g}}
-
\Psi\big([{\sf x}_{j_1},\,{\sf x}_{j_2}]_{\mathfrak{g}}\big).
\]
Applying this formula to the basic forms, we may compute for any two
indices $j_1, j_2$ with $1 \leqslant j_1 < j_2 \leqslant n$:
\[
\footnotesize
\aligned
\big(\partial^1({\sf x}_i^*\otimes{\sf x}_k)\big)
\big({\sf x}_{j_1},\,{\sf x}_{j_2}\big)
&
=
\big[{\sf x}_{j_1},\,\delta_{j_2}^i\,{\sf x}_k\big]_{\mathfrak{g}}
-
\big[{\sf x}_{j_2},\,\delta_{j_1}^i\,{\sf x}_k\big]_{\mathfrak{g}}
-
{\sf x}_i^*\otimes {\sf x}_k
\Big(
\sum_{1\leqslant s\leqslant n}\,
c_{j_1,j_2}^s\,{\sf x}_s
\Big)
\\
&
=
\sum_{s=1}^r\,{\sf x}_s
\Big(
-\,c_{k,j_1}^s\delta_{j_2}^i
+
c_{k,j_2}^s\delta_{j_1}^i
-
c_{j_1,j_2}^i\delta_k^s
\Big).
\endaligned
\]
This means that we have obtained the following representation of
the differentials of all basic $1$-cochains:
\[
\aligned
\partial^1\big({\sf x}_i^*\otimes {\sf x}_k\big)
=
\sum_{1\leqslant j_1<j_2\leqslant n}\,
\sum_{s=1}^r\,
&
{\sf x}_{j_1}^*\wedge {\sf x}_{j_2}^*\otimes {\sf x}_s
\Big(
-\,c_{k,j_1}^s\delta_{j_2}^i
+
c_{k,j_2}^s\delta_{j_1}^i
-
c_{j_1,j_2}^i\delta_k^s
\Big)
\\
&
{\scriptstyle{(1\,\leqslant\,i\,\leqslant\,n\,;\,\,
k\,=\,1,\,\dots,\,n,\,n+1,\,\dots,\,r).}}
\endaligned
\]

\subsection{Boundary of a general $1$-cochain}
Thanks to these formulas, we may then 
compute $\partial^1 \Psi$ by linearity:
\[
\footnotesize
\aligned
\partial^1\Psi
&
=
\sum_{1\leqslant i\leqslant n}\,
\sum_{k=1}^r\,\psi_k^i\,
\partial^1\big({\sf x}_i^*\otimes{\sf x}_k\big)
\\
&
=
\sum_{1\leqslant i\leqslant n}\,
\sum_{k=1}^r\,
\psi_k^i\,
\sum_{1\leqslant j_1<j_2\leqslant n}\,
\sum_{s=1}^r\,
{\sf x}_{j_1}^*\wedge {\sf x}_{j_2}^*\otimes {\sf x}_s
\Big(
-\,c_{k,j_1}\delta_{j_2}^i
+
c_{k,j_2}^s\delta_{j_1}^i
-
c_{j_1,j_2}^i\delta_k^s
\Big)
\\
&
=
\sum_{1\leqslant j_1<j_2\leqslant n}\,
\sum_{s=1}^r\,
{\sf x}_{j_1}^*\wedge {\sf x}_{j_2}^*\otimes {\sf x}_s
\bigg(
\sum_{1\leqslant i\leqslant n}\,
\sum_{k=1}^r\,
-\,\psi_k^ic_{k,j_1}^s\delta_{j_2}^i
+
\psi_k^ic_{k,j_2}^s\delta_{j_1}^{i_1}
-
\psi_k^ic_{j_1,j_2}^i\delta_k^s
\bigg)
\\
&
=
\sum_{1\leqslant j_1<j_2\leqslant n}\,
\sum_{s=1}^r\,
{\sf x}_{j_1}^*\wedge {\sf x}_{j_2}^*\otimes {\sf x}_s
\bigg(
-\,\sum_{k=1}^r\,
\psi_k^{j_2}c_{k,j_1}^s
+
\sum_{k=1}^r\,
\psi_k^{j_1}c_{k,j_2}^s
-
\sum_{1\leqslant i\leqslant n}\,
\psi_s^ic_{j_1,j_2}^i
\bigg).
\endaligned
\] 
As a result, we may explicitly characterize the condition 
that $\Phi$ be a {\sl $2$-coboundary}.

\begin{Proposition}
Under the assumptions of Proposition~\ref{characterization-2-cocycles},
a $2$-cochain:
\[
\Phi
=
\sum_{1\leqslant i_1<i_2\leqslant n}\,
\sum_{k=1}^r\,
\phi_{i_1,i_2}^k\,
{\sf x}_{i_1}^*\wedge{\sf x}_{i_2}^*
\otimes{\sf x}_k
\]
is the boundary $\Phi = \partial^1 \Psi$ of a $1$-cochain:
\[
\Psi
=
\sum_{1\leqslant i\leqslant n}\,
\sum_{k=1}^r\,\psi_k^i\,\,
{\sf x}_i^*\otimes{\sf x}_k
\]
if and only if all its coefficients $\phi_{ l_1, l_2}^k$
are uniquely determined as the following linear 
combinations of the $\psi_\cdot^\cdot$:
\[
\aligned
\phi_{j_1,j_2}^s
=
-\,\sum_{k=1}^r\,
&
\psi_k^{j_2}c_{k,j_1}^s
+
\sum_{k=1}^r\,
\psi_k^{j_1}c_{k,j_2}^s
-
\sum_{1\leqslant i\leqslant n}\,
\psi_s^ic_{j_1,j_2}^i
\\
&
{\scriptstyle{(1\,\leqslant\,j_1\,<\,j_2\,\leqslant\,n\,;\,\,
s\,=\,1,\,\dots\,n,\,n+1,\,\dots,\,r).}}
\endaligned
\]
\end{Proposition}

\subsection{Combinatorial assumptions for a general grading}
\label{form-of-a-graded-Lie-algebra}
Now, we want to show that the complexity of cohomological computations
splits when the linear system of equations considered above may be
decomposed as a direct sum of blocks of linear systems in smaller
dimensions. The typical and quite general case where such a splitting
is available holds when $\mathfrak{ g}$ is endowed with the
supplementary structure of a {\sl grading} in the sense that
$\mathfrak{ g}$, viewed as a vector space, writes down as a direct
sum:
\[
\aligned
\mathfrak{g}
&
=
\mathfrak{g}_{-a}
\oplus\cdots\oplus
\mathfrak{g}_{-1}
\oplus
\mathfrak{g}_0
\oplus
\mathfrak{g}_1
\oplus\cdots\oplus
\mathfrak{g}_b
\\
&
=
\bigoplus_{-a\leqslant k\leqslant b}\,
\mathfrak{g}_k
\endaligned
\]
of nonzero vector subspaces $\mathfrak{ g}_k$, where 
$a \geqslant 1$ and $b \geqslant 0$ are certain integers,
when one assumes that:
\[
\big[\mathfrak{g}_{k_1},\,\mathfrak{g}_{k_2}\big]_{\mathfrak{g}}
\subset \mathfrak{ g}_{k_1+k_2},
\] 
for all $k_1, k_2 \in \Z$, after prolonging trivially $\mathfrak{
g}_k := \{ 0\}$ for either $k \leqslant - a - 1$ or $k \geqslant b
+1$. In this setting, one naturally considers:
\[
\mathfrak{g}_-
:=
\mathfrak{g}_{-a}
\oplus\cdots\oplus
\mathfrak{g}_{-1}
\]
for computing the second cohomology $H^2 ( \mathfrak{ g}_-, \mathfrak{
g})$ in the sense of the preceding sections. As before, we shall
denote:
\[
r
:=
\dim_\K\,\mathfrak{g}
\ \ \ \ \ \ \ \ \ \
\text{\rm and}
\ \ \ \ \ \ \ \ \ \
n
:=
\dim_\K\,\mathfrak{g}_-.
\]
Working abstractly and in full generality, we shall not assume that
$\mathfrak{ g}_-$ is generated by $\mathfrak{ g}_{ -1}$ in the sense
that $\mathfrak{ g}_{ - i - 1} = \big[ \mathfrak{ g}_{ -1}, \,
\mathfrak{ g}_{ -i} \big]_{\mathfrak{g}}$ for all $i = 1, \dots, \mu
-1$, an assumption which, however, comes naturally in Tanaka's theory
of graded differential systems.

Now, we need more notation. For any $k$ with $-a \leqslant k \leqslant
b$, each $\mathfrak{ g}_k$ has a certain positive dimension, call it:
\[
d_{(k)}
:=
\dim_\K\,\mathfrak{g}_k,
\]
so that one naturally has:
\[
\aligned
r
&
=
d_{(\!-a)}+\cdots+d_{(\!-1)}+d_{(0)}+d_{(1)}+\cdots+d_{(b)},
\\
n
&
=
d_{(\!-a)}+\cdots+d_{(\!-1)}.
\endaligned
\]
Let us introduce, for each $k$ with $-a \leqslant k \leqslant b$, an
arbitrary fixed basis: 
\[
\big({\sf x}_k^{i_k} 
\big)^{1\leqslant i_k\leqslant d_{(k)}}
\]
of the $\K$-vector subspace $\mathfrak{ g}_k$ of $\mathfrak{
g}$. The lower index $k$ refers to the graded part $\mathfrak{ g}_k$
to which all the ${\sf x}_k^{ i_k}$ belong, for $i_k = 1, \dots, d_{
(k)}$. Accordingly, because $\big[ \mathfrak{ g}_{ k_1}, \,
\mathfrak{ g}_{ k_2} \big]_{\mathfrak{g}} \subset \mathfrak{ g}_{
k_1+k_2}$, the structure constants of $\mathfrak{ g}$ are of a
certain specific form such that:
\[
\big[{\sf x}_{k_1}^{i_1},\,{\sf x}_{k_2}^{i_2}\big]_{\mathfrak{g}}
=
\sum_{i'=1}^{d_{(k_1+k_2)}}\,
c_{k_1,k_2}^{i_1,i_2,i'}\,
{\sf x}_{k_1+k_2}^{i'}.
\]
Furthermore, one may adopt the convention that $c_{ k_1, k_2}^{ i_1,
i_2, i'} = 0$ whenever one does not simultaneously have $-a \leqslant
k_1, k_2, k_1 + k_2 \leqslant b$, $1 \leqslant i_1 \leqslant d_{
(k_1)}$, $1 \leqslant i_2 \leqslant d_{ (k_2)}$ and $1
\leqslant i' \leqslant d_{ (k_1 + k_2)}$. 

\subsection{Splitting of cochains, of cocycles, 
of coboundaries and of cohomologies according to homogeneity}
Each vector space $\mathcal{ C}^\ell (\mathfrak{ g}_-, \mathfrak{ g})$
naturally splits into a direct sum of so-called {\sl homogeneous
cochains} as follows: an $\ell$-cochain $\Phi\in\mathcal{C}^\ell (
\mathfrak{ g}_-, \mathfrak{ g})$ is said to be {\sl of homogeneity} a
certain integer $h \in \Z$ whenever for any $\ell$ vectors:
\[
{\sf z}_{i_1}\in\mathfrak{g}_{k_1}, 
\ldots\ldots, 
{\sf z}_{i_\ell}\in\mathfrak{g}_{k_\ell}
\]
belonging to certain arbitrary but determined $\mathfrak{ g}$-components, 
its value: 
\[
\Phi({\sf z}_{i_1},\dots,{\sf z}_{i_\ell})
\in
\mathfrak{g}_{i_1+\cdots+i_\ell+h}
\]
belongs to the $(i_1 + \cdots + i_\ell + h)$-th component of
$\mathfrak{ g}$. In fact, one easily convinces oneself that any 
$\ell$-cochain $\Phi \in \mathcal{ C}
( \mathfrak{ g}_-, \mathfrak{ g})$ 
splits as a direct sum of $\ell$-cochains of
fixed homogeneity:
\[
\Phi
=
\cdots+
\Phi^{[h-1]}
+
\Phi^{[h]}
+
\Phi^{[h+1]}
+\cdots,
\]
where we denote the completely $h$-homogeneous component of
$\Phi$ just by $\Phi^{ [h]}$. In other words:
\[
\mathcal{C}^\ell(\mathfrak{g}_-,\mathfrak{g})
=
\bigoplus_{h\in\Z}\,
\mathcal{C}_{[h]}^\ell(\mathfrak{g}_-,\mathfrak{g}),
\]
where of course the spaces $\mathcal{ C}_{ [h]}^\ell ( \mathfrak{
g}_-, \mathfrak{ g})$ reduce to $\{ 0\}$ for all large $\vert h
\vert$. Furthermore, applying the definition~\thetag{
\ref{cochain-definition}}, one verifies the important fact that
$\partial^\ell$ respects homogeneity for all $\ell = 0, 1, \dots, n$,
that is to say, for any $h \in \Z$, one has $\partial^\ell ( \mathcal{
C}_{ [h]}^\ell ) \subset \mathcal{ C}_{ [h]}^{ \ell+1}$, whence the
complex~\thetag{ \ref{complex-partial}} splits up as a direct sum of
complexes:
\[
0
\overset{\partial_{[h]}^0}{\longrightarrow}
\mathcal{C}^1
\overset{\partial_{[h]}^1}{\longrightarrow}
\mathcal{C}^2
\overset{\partial_{[h]}^2}{\longrightarrow}
\cdots
\overset{\partial_{[h]}^{n-2}}{\longrightarrow}
\mathcal{C}^{n-1}
\overset{\partial_{[h]}^{n-1}}{\longrightarrow}
\mathcal{C}^n
\overset{\partial_{[h]}^n}{\longrightarrow}
0
\]
indexed by $h \in \Z$, where $\partial_{ [h]}^\ell$ naturally denotes the
restriction:
\[
\partial_{[h]}^\ell
:=
\partial^\ell\big\vert_{\mathcal{C}_{[h]}^\ell}
\colon
\mathcal{C}_{[h]}^\ell
\longrightarrow
\mathcal{C}_{[h]}^{\ell+1}.
\]
Consequently, one may introduce the spaces of
{\sl $h$-homogeneous cocycles} of order $\ell$: 
\[
\mathcal{Z}_{[h]}^\ell(\mathfrak{g}_-,\mathfrak{g}\big)
:=
{\rm ker}\big(\partial_{[h]}^\ell
\colon
\mathcal{C}_{[h]}^\ell\to\mathcal{C}_{[h]}^{\ell+1}\big),
\]
together with the spaces of 
{\sl $h$-homogeneous coboundaries} of order $\ell$:
\[
\mathcal{B}_{[h]}^\ell(\mathfrak{g}_-,\mathfrak{g}\big)
:=
{\rm im}\big(\partial_{[h]}^{\ell-1}
\colon
\mathcal{C}_{[h]}^{\ell-1}\to\mathcal{C}_{[h]}^\ell\big).
\]

\begin{Definition}
The quotient space:
\[
H_{[h]}^\ell\big(\mathfrak{g}_-,\mathfrak{g}\big)
:=
\frac{\mathcal{Z}_{[h]}^\ell(\mathfrak{g}_-,\mathfrak{g}\big)}{
\mathcal{B}_{[h]}^\ell(\mathfrak{g}_-,\mathfrak{g}\big)}
\]
is called the $h$-homogeneous $\ell$-th cohomology space of
$\mathfrak{ g}_-$ in $\mathfrak{ g}$.
\end{Definition}

In the sequel, we will mainly be interested in showing how to compute
$h$-homogeneous second cohomologies:
\[
H_{[h]}^2\big(\mathfrak{g}_-,\mathfrak{g}\big)
=
\frac{\mathcal{Z}_{[h]}^2(\mathfrak{g}_-,\mathfrak{g}\big)}{
\mathcal{B}_{[h]}^2(\mathfrak{g}_-,\mathfrak{g}\big)},
\]
so that the task of computing the full cohomology spaces:
\[
H^2(\mathfrak{g}_-,\mathfrak{g})
=
\bigoplus_{h\in\Z}\,
H_{[h]}^2(\mathfrak{g}_-,\mathfrak{g})
\]
requires to deal with vector (sub)spaces of smaller dimensions.

\subsection{Collecting $1$-cochains and $2$-cochains
according to constant homogeneity}
In terms of the bases ${\sf x}_l^{ j_l}$ for $\mathfrak{ g}_-$ and
${\sf x}_k^{ i_k}$ for $\mathfrak{ g}$, the collection ${\sf
x}_l^{j_l*} \otimes {\sf x}_k^i$ where $-a \leqslant l \leqslant -1$,
$j_l = 1, \dots, d_{ (l)}$ and where $- a \leqslant k \leqslant b$,
$i_k = 1, \dots, d_{ (k)}$ makes an obvious basis over $\K$ of the
space $1$-cochains. Similarly, the ${\sf x}_{ l_1}^{ j_{ l_1}*} \wedge
{\sf x}_{ l_2}^{ j_{ l_2}*} \otimes {\sf x}_k^{ i_k}$, where either
$-a \leqslant l_1 < l_2 \leqslant -1$ or $l_1 = l_2$ but $1 \leqslant
j_{ l_1} < j_{ l_2} \leqslant d_{ l_1} = d_{ l_2}$, makes too a basis
over $\K$ of the space of $2$-cochains. Thus, if we mind the fact
that any double sum $\sum_{ k = -a}^b \, \sum_{ i_k = 1}^{ d_{(k)}}$
may also be written without mentioning that the second index $i = i_k$
depends upon the first index $k$ (provided the order of summation is
not permuted), it follows that a general $1$-cochain $\Psi \in
\mathcal{ C}^1 (\mathfrak{ g}_-, \mathfrak{ g})$ and a general $2$
cochain $\Psi \in \mathcal{ C}^2 ( \mathfrak{ g}_-, \mathfrak{ g})$
write, in terms of these natural bases, as the following linear
combinations:
\[
\aligned
\Psi
&
=
\sum_{-a\leqslant l\leqslant-1}\,
\sum_{j=1}^{d_{(l)}}\,
\sum_{k=-a}^b\,\sum_{i=1}^{d_{(k)}}\,
\psi_{l,j}^{k,i}\,\,
{\sf x}_l^{j*}\otimes{\sf x}_k^i
\\
\Phi
&
=
\sum_{-a\leqslant l_1<l_2\leqslant -1}\,
\sum_{j_1=1}^{d_{(l_1)}}\,\sum_{j_2=1}^{d_{(l_2)}}\,
\sum_{k=-a}^b\,\sum_{i=1}^{d_{(k)}}\,
\phi_{l_1,j_1,l_2,j_2}^{k,i}\,\,
{\sf x}_{l_1}^{j_1*}\wedge{\sf x}_{l_2}^{j_2*}
\otimes
{\sf x}_k^i
+
\\
&
\ \ \ \ \
+
\sum_{l=-a}^{-1}\,\sum_{1\leqslant j'<j''\leqslant d_{(l)}}\,
\sum_{k=-a}^b\,\sum_{i=1}^{d_{(k)}}\,
\phi_{l,j',j''}^{k,i}\,
{\sf x}_l^{j'*}\wedge{\sf x}_l^{j''*}\otimes{\sf x}_k^i,
\endaligned
\]
where the $\psi_\cdot^\cdot$ and the $\phi_\cdot^\cdot$ are arbitrary
constants in $\K$. However, we must at first improve such a
preliminary representation.

The homogeneity of any basic $1$-cochain ${\sf x}_l^{j*} \otimes {\sf
x}_k^i$ and of any basic $2$-cochain ${\sf x}_{l_1}^{j_1*} \wedge {\sf
x}_{ l_2}^{ j_2*} \otimes {\sf x}_k^i$ is clearly given by:
\[
\aligned
\text{\rm homogeneity}
\big({\sf x}_l^{j*}\otimes{\sf x}_k^i\big)
&
=
-l+k,
\\
\text{\rm homogeneity}
\big({\sf x}_{l_1}^{j_1*}\wedge{\sf x}_{l_2}^{j_2*}
\otimes{\sf x}_k^i\big)
&
=
-l_1-l_2+k,
\endaligned
\]
so that minimal values and maximal values of homogeneities are equal to:
\[
\aligned
\text{\rm $1$-cochains $\Psi$}\,\,
&
\colon\,\,
-a+1\leqslant\text{\rm homogeneity}\leqslant a+b
\\
\text{\rm $2$-cochains $\Phi$}\,\,
&
\colon\,\,
-a+2\leqslant\text{\rm homogeneity}\leqslant 2a+b.
\endaligned
\]
Thus, if we denote by the letter $h$ the homogeneity of a cochain, in
order to split our two cochains:
\[
\Psi
=
\sum_{h=-a+1}^{a+b}\,\Psi^{[h]}
\ \ \ \ \
\text{\rm and}
\ \ \ \ \
\Phi
=
\sum_{h=-a+2}^{2a+b}\,\Phi^{[h]}
\]
as cochains having constant homogeneity $h$ (for any $h \in \Z$), 
and if we introduce the following two sets of integers:
\[
\small
\aligned
\Delta_1^{\![h]}
&:=
\big\{
(l,j)\in\Z\times\N
\colon
-a\leqslant l\leqslant -1,\,\,
j=1,\dots,d_{(l)},\,\,
-a\leqslant l+h\leqslant b
\big\},
\\
\Delta_2^{\![h]}
&:=
\big\{
(l_1,j_1,l_2,j_2)\in\Z\times\N\times\Z\times\N
\colon
\\
&
\ \ \ \ \ \ \ \ \ \ \ \ \ \ 
\aligned
&
-a\leqslant l_1\leqslant -1,\,\,
j_1=1,\dots,d_{(l_1)},
\\
&
-a\leqslant l_2\leqslant -1,\,\,
j_2=1,\dots,d_{(l_2)},
\endaligned\,\,
-a\leqslant l_1+l_2+h\leqslant b
\big\},
\endaligned
\]
we can expand $\Psi^{ [ h]}$ for $-a + 1 \leqslant h \leqslant a+b$
and $\Phi^{ [ h]}$ for $-a + 2 \leqslant h \leqslant 2a + b$ as:
\[
\small
\aligned
\Psi^{[h]}
&
=
\sum_{(l,j)\in\Delta_1^{\![h]}}\,
\sum_{i=1}^{d_{(l+h)}}\,
\psi_{l,j}^{l+h,i}\,\,
{\sf x}_l^{j*}\otimes{\sf x}_{l+h}^i
\ \ \ \ \ \ \ \ \ \ \ \ \ {\scriptstyle{(-\,a\,+\,1\,\leqslant\,
h\,\leqslant\,a\,+\,b)}}
\\
\Phi^{[h]}
&
=
\sum_{(l_1,j_1,l_2,j_2)\in\Delta_2^{\![h]}
\atop(l_1,j_1)<_{\sf lex}(l_2,j_2)}\,
\sum_{k=1}^{d_{(l_1+l_2+h)}}\,\,
\phi_{l_1,j_1,l_2,j_2}^{l_1+l_2+h,k}\,
{\sf x}_{l_1}^{j_1*}\wedge{\sf x}_{l_2}^{j_2*}
\otimes{\sf x}_{l_1+l_2+h}^k
\ \ \ \ \ \ \ \ \ \ \ \ \ {\scriptstyle{(-\,a\,+\,2\,\leqslant\,
h\,\leqslant\,2\,a\,+\,b)}},
\endaligned
\]
where $(l_1, j_1) <_{ \sf lex} (l_2, j_2)$ means either $l_1 < l_2$ or
$l_1 = l_2$ but $j_1 < j_2$.

\subsection{Boundary of $h$-homogeneous $1$-cochains}
Thus, fix a homogeneity $h$ for a $1$-cochain with $-a + 1 \leqslant h
\leqslant a + b$, let $(l, j) \in \Delta_1^{\! [h]}$, let also $(l_1,
j_1, l_2, j_2) \in \Delta_2^{ \! [h]}$ and compute the boundary of an
arbitrary basic $h$-homogeneous $1$-cochain:
\[
\footnotesize
\aligned
\partial_{[h]}^1({\sf x}_l^{j*}\otimes{\sf x}_{l+h}^i)
\big({\sf x}_{l_1}^{j_1},{\sf x}_{l_2}^{j_2}\big)
&
=
\big[
{\sf x}_{l_1}^{j_1},\,
({\sf x}_l^{j*}\otimes{\sf x}_{l+h}^i)({\sf x}_{l_2}^{j_2})
\big]_{\mathfrak{g}}
-
\big[
{\sf x}_{l_2}^{j_2},\,
({\sf x}_l^{j*}\otimes{\sf x}_{l+h}^i)({\sf x}_{l_1}^{j_1})
\big]_{\mathfrak{g}}
-
({\sf x}_l^{j*}\otimes{\sf x}_{l+h}^i)
\big(
\big[
{\sf x}_{l_1}^{j_1},{\sf x}_{l_2}^{j_2}
\big]_{\mathfrak{g}}
\big)
\\
&
=
[{\sf x}_{l_1}^{j_1},\,\delta_{l_2}^l\delta_{j_2}^j{\sf x}_{l+h}^i]_{\mathfrak{g}}
-
[{\sf x}_{l_2}^{j_2},\,\delta_{l_1}^l\delta_{j_1}^j{\sf x}_{l+h}^i]_{\mathfrak{g}}
-
\big({\sf x}_l^{j*}\otimes{\sf x}_{l+h}^i\big)
\Big(
\sum_{k=1}^{d_{(l_1+l_2)}}\,
c_{l_1,l_2}^{j_1,j_2,k}\,{\sf x}_{l_1+l_2}^k
\Big)
\\
&
=
\delta_{l_2}^l\delta_{j_2}^j\!
\sum_{k=1}^{d_{l_1+l+h}}\,
c_{l_1,l+h}^{j_1,i,k}\,{\sf x}_{l_1+l+h}^k
-
\delta_{l_1}^l\delta_{j_1}^j\!
\sum_{k=1}^{d_{(l_2+l+h)}}\,
c_{l_2,l+h}^{j_2,i,k}\,{\sf x}_{l_2+l+h}^k
-
\big(\delta_{l_1+l_2}^lc_{l_1,l_2}^{j_1,j_2,j}\big)\,{\sf x}_{l+h}^i
\\
&
=
\sum_{k=1}^{d_{(l_1+l_2+h)}}\,
\big(
\delta_{l_2}^l\delta_{j_2}^jc_{l_1,l+h}^{j_1,i,k}
-
\delta_{l_1}^l\delta_{j_1}^jc_{l_2,l+h}^{j_2,i,k}
-
\delta_{l_1+l_2}^l\delta_i^kc_{l_1,l_2}^{j_1,j_2,j}
\big)\,{\sf x}_{l_1+l_2+h}^k.
\endaligned
\]
In other words, we have obtained the following representation
for the differential of any basic $h$-homogeneous
$1$-cochain:
\[
\footnotesize
\aligned
\partial_{[h]}^1\big({\sf x}_l^{j*}\otimes{\sf x}_{l+h}^i\big)
&
=
\sum_{(l_1,j_1,l_2,j_2)\in\Delta_2^{\![h]}
\atop
(l_1,j_1)<_{\sf lex}(l_2,j_2)}\,
\sum_{k=1}^{d_{(l_1+l_2+h)}}\,
\big(
\delta_{l_2}^l\delta_{j_2}^jc_{l_1,l+h}^{j_1,i,k}
-
\delta_{l_1}^l\delta_{j_1}^jc_{l_2,l+h}^{j_2,i,k}
-
\delta_{l_1+l_2}^l\delta_i^kc_{l_1,l_2}^{j_1,j_2,j}
\big)\cdot
\\
&
\ \ \ \ \ \ \ \ \ \ \ \ \ \ \ \ \ \ \ \ \ \ \ \ \ \ \ \ \ \ 
\ \ \ \ \ \ \ \ \ \ \ \ \ \ \ \ \ 
\cdot
{\sf x}_{l_1}^{j_1*}\wedge{\sf x}_{l_2}^{j_2*}
\otimes
{\sf x}_{l_1+l_2+h}^k
\ \ \ \ \ \ \
{\scriptstyle{\big((j,\,l)\,\in\,\Delta_1^{\![h]},\,\,
i\,=\,1\,\cdots\,d_{(l+h)}\big)}}.
\endaligned
\]
Now by linearity, we can compute the boundary of a general
$h$-homogeneous $1$-cochain:
\[
\footnotesize
\aligned
\partial_{[h]}^1\Psi^{[h]}
&
=
\sum_{(l,j)\in\Delta_1^{\![h]}}\,
\sum_{i=1}^{d_{(l+h)}}\,
\psi_{l,j}^{l+h,i}\,
\partial_{[h]}^1({\sf x}_l^{j*}\otimes{\sf x}_{l+h}^i)
\\
&
=
\sum_{(l_1,j_1,l_2,j_2)\in\Delta_2^{\![h]}}\,
\sum_{k=1}^{d_{(l_1+l_2+h)}}\,
\bigg(
\sum_{(l,j)\in\Delta_1^{\![h]}}\,
\sum_{i=1}^{d_{(l+h)}}\,
\\
&
\ \ \ \ \ \ \ \ \ \ \ \ \ \ \ \ \ \ \ \ \ \ 
\Big(
\delta_{l_2}^l\delta_{j_2}^jc_{l_1,l+h}^{j_1,i,k}
-
\delta_{l_1}^l\delta_{j_1}^jc_{l_2,l+h}^{j_2,i,k}
-
\delta_{l_1+l_2}^l\delta_k^ic_{l_1,l_2}^{j_1,j_2,j}
\Big)
\psi_{l,j}^{l+h,i}
\bigg)\,
{\sf x}_{l_1}^{j_1*}\wedge{\sf x}_{l_2}^{j_2*}
\otimes{\sf x}_{l_1+l_2+h}^k.
\endaligned
\]

\begin{Proposition}
Under the above assumptions, an arbitrary $h$-homogeneous 
$2$-cochain:
\[
\aligned
\Phi^{[h]}
=
\sum_{(l_1,j_1,l_2,j_2)\in\Delta_2^{\![h]}
\atop(l_1,j_1)<_{\sf lex}(l_2,j_2)}\,
\sum_{k=1}^{d_{(l_1+l_2+h)}}\,\,
\phi_{l_1,j_1,l_2,j_2}^{l_1+l_2+h,k}\,
{\sf x}_{l_1}^{j_1*}\wedge{\sf x}_{l_2}^{j_2*}
\otimes{\sf x}_{l_1+l_2+h}^k
\endaligned
\]
with $-a + 2 \leqslant h \leqslant 2a + b$ is the boundary $\Phi^{
[h]} = \partial_{[h]}^1 \Psi^{ [ h]}$ of a $h$-homogeneous $1$-cochain:
\[
\aligned
\Psi^{[h]}
=
\sum_{(l,j)\in\Delta_1^{\![h]}}\,
\sum_{i=1}^{d_{(l+h)}}\,
\psi_{l,j}^{l+h,i}\,\,
{\sf x}_l^{j*}\otimes{\sf x}_{l+h}^i
\endaligned
\]
if and only if its homogeneous degree $h$ satisfies in fact $-a + 1
\leqslant h \leqslant a + b$ and if all its coefficients $\phi_{ (l_1,
j_1, l_2, j_2)}^{ l_1 + l_2 + h, k}$ are uniquely determined as the
following linear combinations of the $\psi_\cdot^\cdot$:
\[
\aligned
\phi_{l_1,j_1,l_2,j_2}^{l_1+l_2+h,k}
&
=
\sum_{(l,j)\in\Delta_1^{\![h]}}\,
\sum_{i=1}^{d_{(l+h)}}\,
\Big(
\delta_{l_2}^l\delta_{j_2}^jc_{l_1,l+h}^{j_1,i,k}
-
\delta_{l_1}^l\delta_{j_1}^jc_{l_2,l+h}^{j_2,i,k}
-
\delta_{l_1+l_2}^l\delta_k^ic_{l_1,l_2}^{j_1,j_2,j}
\Big)
\psi_{l,j}^{l+h,i}.
\endaligned
\]
\end{Proposition}

\subsection{Boundary of $h$-homogeneous $2$-cochains}
Next, for any $h$ with $-a + 2 \leqslant h \leqslant 2 a + b$, for any
$(l_1, j_1, l_2, j_2) \in \Delta_2^{ [h]}$, for any $k = 1, \dots, d_{
(l_1 + l_2 + h)}$ and for any $(l_1', j_1', l_2', j_2', l_3', j_3')
\in \Delta_3^{[h]}$, belonging to the set:
\[
\footnotesize
\aligned
\Delta_3^{[h]}
&
:=
\big\{
(l_1',j_1',l_2',j_2',l_3',j_3')
\in
\Z\times\N\times\Z\times\N\times\Z\times\N
\colon\,\,
\\
&
\ \ \ \ \ \ \ \
\aligned
&
-a\leqslant l_1'\leqslant -1,\,\,
j_1'=1,\dots,d_{(l_1')},
\\
&
-a\leqslant l_2'\leqslant -1,\,\,
j_2'=1,\dots,d_{(l_2')},
\\
&
-a\leqslant l_3'\leqslant -1,\,\,
j_3'=1,\dots,d_{(l_3')},\,\,\,\,
-a\leqslant l_1'+l_2'+l_3'+h\leqslant b
\big\},
\endaligned
\endaligned
\]
by applying the definitional formula~\thetag{
\ref{cochain-definition}}, we obtain the value of
the boundary of a basic $2$-cochain:
\[
\footnotesize
\aligned
&
\partial_{[h]}^2\big({\sf x}_{l_1}^{j_1*}\wedge{\sf x}_{l_2}^{j_2*}
\otimes{\sf x}_{l_1+l_2+h}^k\big)
\big(
{\sf x}_{l_1'}^{j_1'},{\sf x}_{l_2'}^{j_2'},{\sf x}_{l_3'}^{j_3'}
\big)
=
\ \ \ \ \,
\big[{\sf x}_{l_1'}^{j_1'},\,
({\sf x}_{l_1}^{j_1*}\wedge{\sf x}_{l_2}^{j_2*}
\otimes{\sf x}_{l_1+l_2+h}^k)
({\sf x}_{l_2'}^{j_2'},{\sf x}_{l_3'}^{j_3'})\big]_{\mathfrak{g}}
-
\\
&
\ \ \ \ \
-
\big[{\sf x}_{l_2'}^{j_2'},\,
({\sf x}_{l_1}^{j_1*}\wedge{\sf x}_{l_2}^{j_2*}
\otimes{\sf x}_{l_1+l_2+h}^k)
({\sf x}_{l_1'}^{j_1'},{\sf x}_{l_3'}^{j_3'})\big]_{\mathfrak{g}}
+
\big[{\sf x}_{l_3'}^{j_3'},\,
({\sf x}_{l_1}^{j_1*}\wedge{\sf x}_{l_2}^{j_2*}
\otimes{\sf x}_{l_1+l_2+h}^k)
({\sf x}_{l_1'}^{j_1'},{\sf x}_{l_2'}^{j_2'})\big]_{\mathfrak{g}}
-
\\
&
\ \ \ \ \
-
({\sf x}_{l_1}^{j_1*}\wedge{\sf x}_{l_2}^{j_2*}
\otimes{\sf x}_{l_1+l_2+h}^k)
\big(
\big[{\sf x}_{l_1'}^{j_1'},{\sf x}_{l_2'}^{j_2'}\big]_{\mathfrak{g}},\,
{\sf x}_{l_3'}^{j_3'}
\big)
+
({\sf x}_{l_1}^{j_1*}\wedge{\sf x}_{l_2}^{j_2*}
\otimes{\sf x}_{l_1+l_2+h}^k)
\big(
\big[{\sf x}_{l_1'}^{j_1'},{\sf x}_{l_3'}^{j_3'}\big]_{\mathfrak{g}},\,
{\sf x}_{l_2'}^{j_2'}
\big)
-
\\
&
\ \ \ \ \
-
({\sf x}_{l_1}^{j_1*}\wedge{\sf x}_{l_2}^{j_2*}
\otimes{\sf x}_{l_1+l_2+h}^k)
\big(
\big[{\sf x}_{l_2'}^{j_2'},{\sf x}_{l_3'}^{j_3'}\big]_{\mathfrak{g}},\,
{\sf x}_{l_1'}^{j_1'}
\big).
\endaligned
\]
Let us focus attention on terms at the first and fourth lines, since
other terms are obtained by obvious permutations of triples $\big(
(l_1', j_1'), (l_2', j_2'), (l_3', j_3') \big)$. The term in the first
line continues as:
\[
\footnotesize
\aligned
\big[{\sf x}_{l_1'}^{j_1'},\,
({\sf x}_{l_1}^{j_1*}\wedge{\sf x}_{l_2}^{j_2*}
\otimes{\sf x}_{l_1+l_2+h}^k)
({\sf x}_{l_2'}^{j_2'},{\sf x}_{l_3'}^{j_3'})\big]_{\mathfrak{g}}
&
=
\big[{\sf x}_{l_1'}^{j_1'},\,
\delta_{l_2'}^{l_1}\delta_{j_2'}^{j_1}
\delta_{l_3'}^{l_2}\delta_{j_3'}^{j_2}\,
{\sf x}_{l_1+l_2+h}\big]_{\mathfrak{g}}
\\
&
=
\delta_{l_2'}^{l_1}\delta_{j_2'}^{j_1}
\delta_{l_3'}^{l_2}\delta_{j_3'}^{j_2}\,
\sum_{s=1}^{d_{(l_1'+l_1+l_2+h)}}\,
c_{l_1',l_1+l_2+h}^{j_1',k,s}\,
{\sf x}_{l_1+l_2+l_1'+h}^s
\\
&
=
\sum_{s=1}^{d_{(l_1'+l_2'+l_3'+h)}}\,
\Big(
\delta_{l_2'}^{l_1}\delta_{j_2'}^{j_1}
\delta_{l_3'}^{l_2}\delta_{j_3'}^{j_2}\,
c_{l_1',l_1+l_2+h}^{j_1',k,s}
\Big)\,
{\sf x}_{l_1'+l_2'+l_3'+h}^s.
\endaligned
\]
The term in the fourth line continues as:
\[
\footnotesize
\aligned
-
({\sf x}_{l_1}^{j_1*}\wedge{\sf x}_{l_2}^{j_2*}
&
\otimes{\sf x}_{l_1+l_2+h}^k)
\big(
\big[{\sf x}_{l_1'}^{j_1'},{\sf x}_{l_2'}^{j_2'}\big]_{\mathfrak{g}},\,
{\sf x}_{l_3'}^{j_3'}
\big)
=
\\
&
=
-
({\sf x}_{l_1}^{j_1*}\wedge{\sf x}_{l_2}^{j_2*}
\otimes{\sf x}_{l_1+l_2+h}^k)
\Big(
\sum_{s=1}^{d_{(l_1'+l_2')}}\,
c_{l_1',l_2'}^{j_1',j_2',s}\,
{\sf x}_{l_1'+l_2'}^s,\,\,
{\sf x}_{l_3'}^{j_3'}
\Big)
\\
&
=
-
\bigg(
\sum_{s=1}^{d_{(l_1'+l_2')}}\,
c_{l_1',l_2'}^{j_1',j_2',s}\,
\big(
\delta_{l_1'+l_2'}^{l_1}\delta_s^{j_1}
\delta_{l_3'}^{l_2}\delta_{j_3'}^{j_2}
-
\delta_{l_3'}^{l_1}\delta_{j_3'}^{j_1}
\delta_{l_1'+l_2'}^{l_2}\delta_s^{j_2}
\big)
\bigg)\,{\sf x}_{l_1+l_2+h}^k
\\
&
=
\sum_{s=1}^{d_{(l_1'+l_2'+l_3'+h)}}\,
\Big(
\delta_k^s\,
\big[
-\delta_{l_1'+l_2'}^{l_1}\delta_{l_3'}^{l_2}
\delta_{j_3'}^{j_2}\,c_{l_1',l_2'}^{j_1',j_2',j_1}
+
\delta_{l_3'}^{l_1}\delta_{j_3'}^{j_1}
\delta_{l_1'+l_2'}^{l_2}\,c_{l_1',l_2'}^{j_1',j_2',j_2}
\big]
\Big)\,{\sf x}_{l_1'+l_2'+l_3'+h}^s.
\endaligned
\]
Finally, bringing back the two pairs of terms left above which are
obtained by simple permutations, we obtain the value of the boundary
of a basic $2$-cochain:
\[
\footnotesize
\aligned
&
\partial_{[h]}^2\big({\sf x}_{l_1}^{j_1*}\wedge{\sf x}_{l_2}^{j_2*}
\otimes{\sf x}_{l_1+l_2+h}^k\big)
\big(
{\sf x}_{l_1'}^{j_1'},{\sf x}_{l_2'}^{j_2'},{\sf x}_{l_3'}^{j_3'}
\big)
=
\\
&
=
\sum_{s=1}^{d_{(l_1'+l_2'+l_3'+h)}}\,
\Big(
\delta_{l_2'}^{l_1}\delta_{j_2'}^{j_1}
\delta_{l_3'}^{l_2}\delta_{j_3'}^{j_2}\,
c_{l_1',l_1+l_2+h}^{j_1',k,s}
-
\delta_{l_1'}^{l_1}\delta_{j_1'}^{j_1}
\delta_{l_3'}^{l_2}\delta_{j_3'}^{j_2}\,
c_{l_2',l_1+l_2+h}^{j_2',k,s}
+
\delta_{l_1'}^{l_1}\delta_{j_1'}^{j_1}
\delta_{l_2'}^{l_2}\delta_{j_2'}^{j_2}\,
c_{l_3',l_1+l_2+h}^{j_3',k,s}
+
\\
&
\ \ \ \ \ \ \ \ \ \ \ \ \ \ \ \ \ \ \ \ \ \ \ \ \ \ \ 
+
\delta_k^s\,
\big[
-\delta_{l_1'+l_2'}^{l_1}\delta_{l_3'}^{l_2}\delta_{j_3'}^{j_2}\,
c_{l_1',l_2'}^{j_1',j_2',j_1}
+
\delta_{l_3'}^{l_1}\delta_{j_3'}^{j_1}\delta_{l_1'+l_2'}^{l_2}\,
c_{l_1',l_2'}^{j_1',j_2',j_2}
+
\delta_{l_1'+l_3'}^{l_1}\delta_{l_2'}^{l_2}\delta_{j_2'}^{j_2}\,
c_{l_1',l_3'}^{j_1',j_3',j_1}
-
\\
&
\ \ \ \ \ \ \ \ \ \ \ \ \ \ \ \ \ \ \ \ \ \ \ \ \ \ \ 
-
\delta_{l_2'}^{l_1}\delta_{j_2'}^{j_1}\delta_{l_1'+l_3'}^{l_2}\,
c_{l_1',l_3'}^{j_1',j_3',j_2}
-
\delta_{l_2'+l_3'}^{l_1}\delta_{l_1'}^{l_2}\delta_{j_1'}^{j_2}\,
c_{l_2',l_3'}^{j_2',j_3',j_1}
+
\delta_{l_1'}^{l_1}\delta_{j_1'}^{j_1}\delta_{l_2'+l_3'}^{l_2}\,
c_{l_2',l_3'}^{j_2',j_3',j_2}
\big]
\Big)\,{\sf x}_{l_1'+l_2'+l_3'+h}^s.
\endaligned
\]
Equivalently:
\[
\footnotesize
\aligned
&
\partial_{[h]}^2\big({\sf x}_{l_1}^{j_1*}\wedge{\sf x}_{l_2}^{j_2*}
\otimes{\sf x}_{l_1+l_2+h}^k\big)
=
\sum_{(l_1',j_1',l_2',j_2',l_3',j_3')\in\Delta_3^{[h]}}\,
\sum_{s=1}^{d_{(l_1'+l_2'+l_3'+h)}}\,
\\
&
\ \ \ \ \ 
\Big(
\delta_{l_2'}^{l_1}\delta_{j_2'}^{j_1}
\delta_{l_3'}^{l_2}\delta_{j_3'}^{j_2}\,
c_{l_1',l_1+l_2+h}^{j_1',k,s}
-
\delta_{l_1'}^{l_1}\delta_{j_1'}^{j_1}
\delta_{l_3'}^{l_2}\delta_{j_3'}^{j_2}\,
c_{l_2',l_1+l_2+h}^{j_2',k,s}
+
\delta_{l_1'}^{l_1}\delta_{j_1'}^{j_1}
\delta_{l_2'}^{l_2}\delta_{j_2'}^{j_2}\,
c_{l_3',l_1+l_2+h}^{j_3',k,s}
+
\\
&
\ \ \ \ \ \ \
+
\delta_k^s\,
\big[
-\delta_{l_1'+l_2'}^{l_1}\delta_{l_3'}^{l_2}\delta_{j_3'}^{j_2}\,
c_{l_1',l_2'}^{j_1',j_2',j_1}
+
\delta_{l_3'}^{l_1}\delta_{j_3'}^{j_1}\delta_{l_1'+l_2'}^{l_2}\,
c_{l_1',l_2'}^{j_1',j_2',j_2}
+
\delta_{l_1'+l_3'}^{l_1}\delta_{l_2'}^{l_2}\delta_{j_2'}^{j_2}\,
c_{l_1',l_3'}^{j_1',j_3',j_1}
-
\\
&
\ \ \ \ \ \ \ 
-
\delta_{l_2'}^{l_1}\delta_{j_2'}^{j_1}\delta_{l_1'+l_3'}^{l_2}\,
c_{l_1',l_3'}^{j_1',j_3',j_2}
-
\delta_{l_2'+l_3'}^{l_1}\delta_{l_1'}^{l_2}\delta_{j_1'}^{j_2}\,
c_{l_2',l_3'}^{j_2',j_3',j_1}
+
\delta_{l_1'}^{l_1}\delta_{j_1'}^{j_1}\delta_{l_2'+l_3'}^{l_2}\,
c_{l_2',l_3'}^{j_2',j_3',j_2}
\big]
\Big)\,
{\sf x}_{l_1'}^{j_1'*}\wedge{\sf x}_{l_2'}^{j_2'*}
\wedge{\sf x}_{l_3'}^{j_3'*}
\otimes{\sf x}_{l_1'+l_2'+l_3'+h}^s
\endaligned
\]
Now, by linearity, we deduce the boundary $\partial_{[h]}^2
\Phi^{ [h]}$ of a general $h$-homogeneous $2$-cochain $\Phi^{[h]}$,
and this yields the following statement.

\begin{Proposition}
Under the above assumptions, the boundary of a general
$h$-homogeneous $2$-cochain:
\[
\Phi^{[h]}
=
\sum_{(l_1,j_1,l_2,j_2)\in\Delta_2^{[h]}
\atop
(l_1,j_1)<_{\sf lex}(l_2,j_2)}\,
\sum_{k=1}^{d_{(l_1+l_2+h)}}\,
\phi_{l_1,j_1,l_2,j_2}^{l_1+l_2+h,k}\,\,
{\sf x}_{l_1}^{j_1*}\wedge{\sf x}_{l_2}^{j_2*}
\otimes{\sf x}_{l_1+l_2+h}^k,
\]
where $h$ satisfies $-a + 2 \leqslant h \leqslant 2a + b$ which has
arbitrary coefficients $\phi_{ l_1, j_1, l_2, j_2}^{ l_1 + l_2 + h, k}
\in \K$, is a cocycle, namely satisfies $0 = \partial_{ [h]}^2 \Phi^{
[h]}$, if and only if all the following linear equations hold:
\[
\footnotesize
\aligned
0
&
=
\sum_{(l_1,j_1,l_2,j_2)\in\Delta_2^{[h]}
\atop
(l_1,j_1)<_{\sf lex}(l_2,j_2)}\,
\sum_{k=1}^{d_{(l_1+l_2+h)}}\,
\Big(
\delta_{l_2'}^{l_1}\delta_{j_2'}^{j_1}
\delta_{l_3'}^{l_2}\delta_{j_3'}^{j_2}\,
c_{l_1',l_1+l_2+h}^{j_1',k,s}
-
\delta_{l_1'}^{l_1}\delta_{j_1'}^{j_1}
\delta_{l_3'}^{l_2}\delta_{j_3'}^{j_2}\,
c_{l_2',l_1+l_2+h}^{j_2',k,s}
+
\delta_{l_1'}^{l_1}\delta_{j_1'}^{j_1}
\delta_{l_2'}^{l_2}\delta_{j_2'}^{j_2}\,
c_{l_3',l_1+l_2+h}^{j_3',k,s}
+
\\
&
\ \ \ \ \ \ \
+
\delta_k^s\,
\big[
-\delta_{l_1'+l_2'}^{l_1}\delta_{l_3'}^{l_2}\delta_{j_3'}^{j_2}\,
c_{l_1',l_2'}^{j_1',j_2',j_1}
+
\delta_{l_3'}^{l_1}\delta_{j_3'}^{j_1}\delta_{l_1'+l_2'}^{l_2}\,
c_{l_1',l_2'}^{j_1',j_2',j_2}
+
\delta_{l_1'+l_3'}^{l_1}\delta_{l_2'}^{l_2}\delta_{j_2'}^{j_2}\,
c_{l_1',l_3'}^{j_1',j_3',j_1}
-
\\
&
\ \ \ \ \ \ \ 
-
\delta_{l_2'}^{l_1}\delta_{j_2'}^{j_1}\delta_{l_1'+l_3'}^{l_2}\,
c_{l_1',l_3'}^{j_1',j_3',j_2}
-
\delta_{l_2'+l_3'}^{l_1}\delta_{l_1'}^{l_2}\delta_{j_1'}^{j_2}\,
c_{l_2',l_3'}^{j_2',j_3',j_1}
+
\delta_{l_1'}^{l_1}\delta_{j_1'}^{j_1}\delta_{l_2'+l_3'}^{l_2}\,
c_{l_2',l_3'}^{j_2',j_3',j_2}
\big]
\Big)\,
\phi_{l_1,j_1,l_2,j_2}^{l_1+l_2+h,k}.
\endaligned
\]
\end{Proposition}

\noindent
Further considerations accompanied with an algorithm using Gr\"obner
bases may be found in~\cite{AAMS}.

\end{document}